\documentclass[11pt]{ip-journal}

 \usepackage{amssymb}
 
 \usepackage{amsmath,amsfonts}
 \usepackage{amscd}
 \usepackage[mathscr]{eucal}

 \usepackage[T3,T1]{fontenc}
 \DeclareSymbolFont{tipa}{T3}{cmr}{m}{n}
 \DeclareMathAccent{\invbreve}{\mathalpha}{tipa}{16}
 
\newenvironment{frcseries}{\fontfamily{frc}\selectfont}{}
 \newcommand{\textfrc}[1]{{\frcseries#1}}

\newcommand{\mathfrc}[1]{\text{\textfrc{\upshape\if h#1\scriptsize
      h\else\if b#1\scriptsize b\else\small#1\fi\fi}}}

 \usepackage{xy}
 \xyoption{all}
 \CompileMatrices
\usepackage{graphicx}
\usepackage{calc}

\newlength{\wcwidth}
\newlength{\wcheight}

\numberwithin{equation}{section}
 \renewcommand{\theequation}{{\thesection.\arabic{equation}}}
 
\theoremstyle{plain}
\newtheorem*{thm*}{Theorem}
\newtheorem{thm}{Theorem}[section]\newtheorem{prop}[thm]{Proposition}
\newtheorem{lemma}[thm]{Lemma}\newtheorem{cor}[thm]{Corollary}
\newtheorem{conj}[thm]{Conjecture}

\theoremstyle{definition}
\newtheorem{defn}[thm]{Definition}

\theoremstyle{remark}
\newtheorem{rem}[thm]{Remark}
\newtheorem{remarks}[thm]{Remarks}
\newtheorem{assumption}[thm]{Assumption}
\newtheorem{facts}[thm]{Facts}

 \newcommand{\bth}{\begin{*thm}}
 \newcommand{\ethm}{\end{*thm}}
 \newcommand{\bco}{\begin{*cor}}
 \newcommand{\eco}{\end{*cor}}
 \newcommand{\bcj}{\begin{*conj}}
 \newcommand{\ecj}{\end{*conj}}
 \newcommand{\bpr}{\begin{*prop}}
 \newcommand{\epr}{\end{*prop}}
 \newcommand{\bprs}{\begin{**prop}}
 \newcommand{\eprs}{\end{**prop}}
 \newcommand{\ble}{\begin{*lemma}}
 \newcommand{\ele}{\end{*lemma}}
 \newcommand{\bles}{\begin{**lemma}}
 \newcommand{\eles}{\end{**lemma}}
 \newcommand{\bsl}{\begin{*sublemma}}
 \newcommand{\esl}{\end{*sublemma}}
 \newcommand{\bre}{\begin{*rem}}
 \newcommand{\ere}{\end{*rem}}
 \newcommand{\bres}{\begin{**rem}}
 \newcommand{\eres}{\end{**rem}}
 \newcommand{\bnt}{\begin{*nota}}
 \newcommand{\ent}{\end{*nota}}
 \newcommand{\bnts}{\begin{**nota}}
 \newcommand{\ents}{\end{**nota}}
 \newcommand{\bde}{\begin{*defn}}
 \newcommand{\ede}{\end{*defn}}
 \newcommand{\bdes}{\begin{**defn}}
 \newcommand{\edes}{\end{**defn}}

\DeclareMathAlphabet\eusm{U}{eus}{m}{n}
\def\makebb#1{\expandafter\def\csname bb#1\endcsname{{\mathbb{#1}}}\ignorespaces}
\def\makerm#1{\expandafter\def\csname rm#1\endcsname{{\rm #1}}\ignorespaces}
\def\makebf#1{\expandafter\def\csname bf#1\endcsname{{\bf #1}}\ignorespaces}
\def\makegr#1{\expandafter\def\csname gr#1\endcsname{{\mathfrak{#1}}}\ignorespaces}
\def\makescr#1{\expandafter\def\csname scr#1\endcsname{{\mathscr{#1}}}\ignorespaces}
\def\makecal#1{\expandafter\def\csname cal#1\endcsname{{\cal #1}}\ignorespaces}
\def\makeudl#1{\expandafter\def\csname udl#1\endcsname{{\underline{#1}}}\ignorespaces}
\def\doLetters#1{%
  #1A #1B #1C #1D #1E #1F #1G #1H #1I #1J #1K #1L #1M
  #1N #1O #1P #1Q #1R #1S #1T #1U #1V #1W #1X #1Y #1Z}
\def\doletters#1{%
  #1a #1b #1c #1d #1e #1f #1g #1h #1i #1j #1k #1l #1m
  #1n #1o #1p #1q #1r #1s #1t #1u #1v #1w #1x #1y #1z}
\doLetters\makebb\doLetters\makecal\doLetters\makerm\doLetters\makebf\doLetters\makescr
\doletters\makerm \doletters\makebf \doLetters\makegr
\doletters\makegr \doLetters\makeudl \doletters\makeudl

 \newcommand{\op}{\operatorname}

\def\co{\colon\thinspace}
 \def\pf{\noindent{\it Proof.}\enspace}

 \def\Conn{\mathop{\mathrm{Conn}}\nolimits}
 \def\Hom{\mathop{\mathrm{Hom}\, }\nolimits}
 \def\im{\mathop{\mathrm{Im}\, }\nolimits}

\newcommand{\ul}{\textup{l}}

\def\zzz{\textsl{z}}
\def\pp{\textsl{p}}

\def\fh{\op{\mathfrc{h}}}
\def\fc{\op{\mathfrc{c}}}
\def\fd{\op{\mathfrc{d}}}
\def\fm{\op{\mathfrc{m}}}
\def\fs{\op{\mathfrc{s}}}

\def\hata{\op{\hat{a}}}
\def\ra{\mathrm{a}}

\def\hY{\hat{Y}}

 \def\ker{\mathop{\mathrm{Ker}\,}\nolimits}

 \def\dist{\mathop{\mathrm{dist}\, }\nolimits}

\def\hatl{\hat{\mathop{\mathrm{l}}}}

 \def\Spin{\mathop{\mathrm{Spin}}\nolimits}

 \newcommand{\epf}{\hfill$\Box$\medbreak}
 \newcommand{\td}{\tilde}
\newcommand{\ts}{\tilde{\sigma}}
\newcommand{\slp}{\mbox{\mbox{\(\partial\mkern-7.5mu\mbox{/}\mkern 2mu\)}}}
\newcommand{\ud}{\underline}
\newcommand{\ov}{\overline}

\def\ss{\textsc{s}}

\def\sO{\textsc{o}}

\def\sW{\textsc{w}}

\def\smE{\textsc{e}}

\def\BTitem#1\ETitem{\begin{equation}\hss\left\{\mkern-20mu\parbox{0.9\hsize}{%
\begin{itemize}#1\end{itemize}}\right.\end{equation}}
\def\BTTitem#1\ETTitem{\begin{equation*}\tag{\Ttag}\hss\left\{\mkern-20mu\parbox{0.9\hsize}{%
\begin{itemize}#1\end{itemize}}\right.\end{equation*}}

 \usepackage{pdfsync,hyperref}

\begin{document}

 \title[From SW to Gr: MCE  and Singular Symplectic Forms]{From Seiberg-Witten to Gromov: \\MCE  and Singular Symplectic Forms}

 \author{Yi-Jen Lee} 
 \address{Institute of Mathematical Sciences \\the Chinese University
   of Hong Kong \\ Shatin, N.T. Hong Kong\\ }
 \email{yjlee@math.cuhk.edu.hk}



 \begin{abstract}
Motivated by various possible generalizations of Taubes's \(SW=Gr\)
theorem \cite{T} to Floer-theoretic
setting,
we prove certain variants of Taubes's convergence
theorem in \cite{T} (the first part of his proof of \(SW=Gr\)). In
place of the closed symplectic 4-manifold considered in \cite{T},
this article considers non-compact manifolds with cylindrical ends, 
equipped with a self-dual harmonic 2-form with non-degenerate zeroes.  
This extends and simplifies some central technical ingredients of the
author's prior work in \cite{LT} and \cite{KLT5}.
 Other expected applications include: 
extending the \(HM=PFH\) theorem in \cite{LT} and the \(HM=HF\)
theorem in  \cite{KLT1}-\cite{KLT5} to TQFTs on both sides \cite{L1};
definitions of large-perturbation Seiberg-Witten analogs of Heegaard
Floer theory's link Floer homologies and link cobordism invariants.
\end{abstract}

 \maketitle
\tableofcontents

\section{Introduction}

In a series of ground breaking papers \cite{T}, Taubes proved the equivalence of
the Seiberg-Witten invariant of closed symplectic 4-manifolds and a version
of Gromov invariant. The first article of the series, also technically
the most important of the four, \(SW\Rightarrow Gr\), proves that by perturbing the Seiberg-Witten equations by a large multiple of the
symplectic form, the Seiberg-Witten solutions converge in a technical
sense to a union of
connected pseudo-holomorphic curves with weights.  
As the first step to the ultimate goal of generalizing the equivalence
theorem to smooth 4-manifolds with \(b^{2+}>0\), Taubes proved a
 generalization of this convergence theorem in \cite{Ts}. In the
language adopted in this paper, this is rephrased as follows.

Let \(X\) be a closed \(\Spin^c\) 4-manifold, and \(\bbS_X=\bbS=\bbS^+\oplus \bbS^-\) be its
associated spinor bundle. The Seiberg-Witten equations on \(X\) take
the following form:
\begin{equation}\label{eq:SW}
\begin{cases}\frac{1}{2}F_{A}^+-\rho^{-1}(\Psi\Psi^*)_0+\frac{i}{4}\mu^+=0,&\\ 
\slp_A^+\Psi=0, &\end{cases}
\end{equation}
where \(A\) is a connection on \(\det \bbS^+\), \(\Psi\) is a section
on \(\bbS^+\), and \(\rho\), \(\slp_A^+\co \bbS^+\to \bbS^-\) denote respectively the
Clifford multiplication and the Dirac operator. (We use the definition
given in Sections I.1.1-I.1.3 of \cite{KM}). In the last term of the first
equation, \(\mu\) denotes a 2-form,
and \(\mu^+\) its self-dual part. This is often regarded as a
perturbation.  
\begin{thm}
{\rm \cite{Ts}}\label{thm:T1}
Let \(X\) be a closed  \(\Spin^c\) 4-manifold with \(b^{2+}>0\)
as above, and let 
\(\omega\) be a self-dual harmonic 2-form on \(X\) vanishing
transversely along an embedded 1-submanifold. Take the perturbation
two form \(\mu^+\) in Equation (\ref{eq:SW}) above to be 
\begin{equation}\label{eq:mu-pert}
\mu^+_r=\frac{1}{2}r\omega+w^+_r, 
\end{equation}
where \(r>1\), and \(w^+_r\) is a smooth self-dual 2-form with
\(\|w_r^+\|_{C^2}\) is bounded by an \(r\)-independent constant. 
Let \(\{(A_r,
\Psi_r)\}_r\) be a sequence of corresponding solutions, with \(r\in
\{r_1, r_2, \ldots, r_n, \ldots\}\), \(r_n\to \infty\) as \(n\to \infty\).
Then there exists a t-curve \(\mathbf{C}\) on \(X-\omega^{-1}(0)\) with
respect to the almost complex structure determined by \(\omega\) and
the metric on \(X\), and a subsequence
of the gauge equivalence classes \(\{[(A_r, \Psi_r)]\}_r\) which t-converges to \(\mathbf{C}\).
\end{thm}
The terms  `t-curves' and `t-convergence' above are respectively
catch-all phrases we coined for the kind of pseudo-holomorphic subvarieties
and the technical notion of convergence in Taubes's
result\footnote{Admittedly, this is a poor choice of terminology. The
  ``t-convergence'' is not ``real'' convergence in the usual sense:
  The sequence \(\{[(A_r, \Psi_r)]\}_r\) does not converge in any
  reasonable topological space.}.  They will be defined in
\textsection \ref{sec:t-curve} and \textsection \ref{sec:t-conv}.

\begin{rem}
  Only the \(w_r^+=0\) case was discussed in
\cite{Ts}. However, the arguments in \cite{Ts} work to establish the
slightly more general statement in the preceding
theorem.
\end{rem}

In the special case when \(X=S^1\times Y\), where \(Y\) is a closed
\(\Spin^c\) 3-manifold with \(b_1>0\), one may choose a harmonic
Morse-Novikov 1-form \(\theta \) on \(Y\), and take \(\omega=2(ds \wedge
\theta)^+\) on \(X\). Here, \(s\) denotes the affine parameter of \(S^1=\bbR/\bbZ\). With
respect to the product metric, \(\omega\) is self-dual harmonic, and
the \(S^1\)-invariant Seiberg-Witten solutions on \(X\) correspond to
the Seiberg-Witten solutions on the 3-manifold \(Y\). Let
\(\check{\theta}\) denote the vector dual to \(\theta\). The
\(S^1\)-invariant t-curves are products of \(S^1\) with certain sets of
(weighted) trajectories of \(\check{\theta}\) on \(Y\). We call the
latter  `t-orbits', with the precise definition given in \textsection \ref{sec:t-orbit}. Thus, as a
corollary to the preceding theorem, an analogous statement for closed
3-manifolds also holds. (A more precise statement is given in
Proposition \ref{prop:t-conv3d} below). For the purpose of this article however, we
need a slightly stronger notion of convergence, which we call ``strong
t-convergence''. This will be defined in \textsection \ref{sec:strong-t}.

The 3-dimensional Seiberg-Witten
equation, arising from a dimensional reduction of the four dimensional
version, has the following general form: 
\begin{equation}
\label{eq:SW3}
\grF_{\mu}(B, \Phi):=\begin{cases}\frac{1}{2}*F_{B}+\rho^{-1}(\Phi\Phi^*)_0+\frac{i}{4}*\mu=0,&\\ 
\slp_B \Phi=0, &\end{cases}
\end{equation}
where \(B\) 
is a unitary connection on \(\det\mathbb{S}\), \(\mathbb{S}=\bbS_Y\)
being the spinor-bundle associated to the given \(\Spin^c\) structure,
and \(\Phi\) is a section of \(\mathbb{S}\). \(\rho\) denotes the
Clifford action, and \(\slp_B\) denotes the Dirac operator associated
to \(B\). \((\Phi\Phi^*)_0\) denotes the traceless part of
\(\Phi\Phi^*\). \(\mu\) is a closed 2-form on \(Y\). Two pairs \((B,
\Phi)\), \((B', \Phi')\) are said to be {\em gauge equivalent} if
there exists a \(u\in C^\infty(Y, S^1)\), such that \((B',
\Phi')=(B-2u^{-1}du, u\cdot \Phi)\). Recall also that a closed 1-form \(\nu\)
on a 3-manifold \(Y\) is said to be {\em Morse-Novikov} if \(\nu\) has
isolated, nondegenerate zeroes.

\begin{thm}\label{thm:strong-t}
Let \((Y, \grs)\) be a closed oriented \(\Spin^c\) 3-manifold with \(b_1>0\), and
\(\theta\) be a harmonic Morse-Novikov 1-form on \(Y\). 
Let \(\nu=*\theta\) and take \(\mu=\mu_r= r\nu+w_r\), where
\(w_r\) is a closed 2-form on \(Y\), and \(r>1\).  Suppose that
\(\|w_r\|_{C^2}\) is bounded by an \(r\)-independent constant. Let \(\{(B_r,
\Phi_r)\}_r\) be a sequence of solutions to the \(\mu =\mu _r\)'s version of (\ref{eq:SW3}), with \(r\in
\{r_1, r_2, \ldots, r_n, \ldots\}\), \(r_n\to \infty\) as \(n\to \infty\).
Then there exists a t-orbit \(\pmb{\gamma}\) on \(Y-\theta^{-1}(0)\) with
respect to the flow of \(\check{\theta}\), the vector field dual to \(\theta\), and a subsequence
of the gauge equivalence classes \(\{[(B_r, \Phi_r)]\}_r\) which
strongly t-converges to \(\pmb{\gamma}\).

Furthermore, \begin{equation}\label{z-theta}
\int_{\tilde{\gamma}}\theta=\frac{1}{2}\big(c_1(\grs)\cdot[\theta]-\zeta_\theta\big)\geq
0, 
\end{equation}
where
\(\zeta_\theta\) is a constant depending only on \(\theta\) (and not
the spin structure \(\grs\)).
\end{thm}
 An explicit formula for the constant
\(\zeta_{\theta}\) in (\ref{z-theta}) is given in Proposition
\ref{prop:t-conv3d}. 

See Section
\ref{sec:strong-t} for the proof of this theorem and terminologies in
its statement. The preceding theorem should be viewed as the first part of a (direct) proof of
a 3-dimensional variant of \(SW=Gr\). In fact, the aforementioned
variant has a precise formulation as Conjecture 1.9 in \cite{HL}. A
slightly weaker version of the latter conjecture is proved
(indirectly) by combining the results in \cite{HL, MT}. A outline of
proof for Conjecture 1.9 in \cite{HL} is given in Section 1.5 of \cite{G}. 

An explicit expression for the constant
\(\zeta_{\theta}\) in (\ref{z-theta}) will be given in
Section \ref{sec:t-conv} below. 
Suffices to say for now that when 
\(\theta=d\tilde{f}\) for a \(S^1\)-valued Morse
function \(\tilde{f}\), \(\zeta_\theta\) is the maximal Euler characteristic of the
regular fibers of \(\tilde{f}\co Y\to S^1\).

\subsection{Stating the main theorems}\label{sec:1.1}

The goal of this article is to prove an analog of  Taubes's theorem (Theorem
\ref{thm:T1} above)  for 4-manifolds with cylindrical 
ends, equipped with a certain type of harmonic self-dual 2-form \(\omega\).
\begin{defn}\label{def:mce}
A {\em manifold with cylindrical ends} (or ``MCE'' for short) is a
connected oriented
manifold \(X\) with a complete metric \(g\), so that there is a
compact manifold-with-boundary \(X_c\subset X\), and an isometry
\(\gri\) 
from the closure of \(X-X_c\) to \(\coprod_{i\in \grY} \bbR^{\geq 0}\times Y_i\), where \(Y_i\)
are closed connected manifolds and each half cylinder \(\bbR^+\times Y_i\)
is endowed with the product metric. The index set,  \(\grY=\grY_X\),  is a non-empty
set of finitely many elements. In the rest of this article, we shall
sometimes implicitly identify subspaces in the closure of \(X-X_c\)
and their image in \(\coprod_{i\in \grY} \bbR^{\geq 0}\times Y_i\)
under \(\gri\). 
\end{defn}
We call the preimage (under \(\gri\)) of the half cylinder \(\bbR^{\geq 0}\times Y_i\) above an {\em end} of
\(X\), or the {\em \(Y_i\)-end} when we wish to be specific. It will
be frequently denoted by \(\hat{Y}_i\). Let 
\[\mathfrc{s}_i\co \bbR^{\geq 0}\times Y_i\simeq \hat{Y}_i \subset X\, \to \bbR^{\geq 0}
\]
denote the projection to the first factor. Denote
\(\mathfrc{s}_i^{-1}[l, \infty)\subset \hat{Y}_i\) by
\(\hat{Y}_{i,l}\). Slightly abusing terminology, we sometimes also
call such \(\hat{Y}_{i, l}\) an end of \(X\). We call each \(Y_i\) an {\em
  ending 3-manifold} of \(X\). Unless specified otherwise, the ending
3-manifolds \(Y_i\) are oriented such that \(\gri \) 
is orientation preserving. 

The analog of the self-dual
harmonic 2-form \(\omega\) has to satisfy appropriate asymptotic
conditions. Correspondingly, the Seiberg-Witten solutions \((A_r,
\Psi_r)\) need also satisfy certain asymptotic conditions. 

Regular Morse-Novikov 1-forms are generic in an appropriate sense.

\begin{defn}\label{def:reg_form}
We call a Morse-Novikov closed 1-form \(\theta \) on a closed oriented
3-manifold {\em regular} if:  
\begin{itemize}
\item all
  finite-length flow lines \(\gamma\) of
  the vector field dual to \(\theta\)  are
  regular. Namely, either  \(\gamma\) is a flow line from an index 2 zero
   \(p\in \theta^{-1}(0)\) to an index 1 zero \(q\in \theta^{-1}(0)\), where the
  descending manifold from \(p\) and ascending manifold of \(q\)
  intersect transversely at \(\gamma\), or \(\gamma\) is a closed orbit with
  nondegenerate linearized Poincar\'e return map;
\item for any fixed \(l>0\), there are finitely many  such flow lines
  \(\gamma\) with \(\int_\gamma\theta\leq l\).
\end{itemize}
\end{defn}
\begin{defn}\label{def:adm}
A smooth 2-form \(\nu\) on a MCE \(X\) is said to be {\em admissible} if it is a
nontrivial harmonic 2-form satisfying:
\begin{itemize}
\item  \(\nu\) is asymptotic to a harmonic 2-form \(v_i\) on \(Y_i\)
  as \(\mathfrc{s}_i\to \infty\) on each end \(\hat{Y_i}\). 
By this we mean the following: \(\tau_{-L}\nu\) converges to
\(\pi_2^*\nu_i\) on \([-1, 1]\times Y_i\) in the
\(C^\infty\)-topology, where \(\tau_{-L}\) denotes translation by
\(-L\in \bbR\) in the \(\bbR\)-factor of \(\hat{Y}_i\simeq
\bbR^{\geq 0}\times Y_i\), and \(\pi_2\) denotes projection to the
\(Y_i\)-factor of \(\hat{Y}_i\simeq\bbR^{\geq 0}\times Y_i\).
The
pair \((Y_i, \nu_i)\) is called the {\em ending pair} associated to \(\hat{Y}_i\).
The harmonic form \(\theta_i:=*_3\nu_i\) is either Morse-Novikov and regular, in which case the
  corresponding end \(\hat{Y}_i\) is called a {\em Morse end}, or
  \(\nu_i=0\), in which case the corresponding end is called a {\em
    vanishing end}. Let \(\grY_m\subset \grY\), and
  \(\grY_v:=\grY-\grY_m\subset \grY\) denote respectively the set of
  Morse ends and the set of vanishing ends. 

\item There are  functions \({\bf l}_v, \bfl_v^+ \co \grY_v\to
  \bbR^+\),
  \(\bfl_v(j)=:\ul_j\),  \(\bfl^+_v(j)=:\ul^+_j>\ul_j+10\), so that
  the restriction of \(\nu\) on
  \[
X'':=X-\bigcup_{j\in \grY_v}\hat{Y}_{j,\ul^+_j}^\circ :\quad
\]
vanishes  transversely along a 1-submanifold \(\nu ^{-1}(0)\) (possibly with
  boundary).
(In the above, \(\hat{Y}_{j,l}^\circ=\mathfrc{s}_j^{-1}(l, \infty)\) denotes the
interior of \(\hat{Y}_{j,l}\).)
  Let \(X':=X-\bigcup_{j\in \grY_v}\hat{Y}_{j,\ul_j}^\circ\). 
  Without loss of generality, we take
  \(\ul_j^+=\ul_j+10\) for the rest of this article. 
  We shall also encounter  intermediate domains lying between \(X'\) and
\(X''\). These are defined as follows: Given \(a\in \bbR\), let 
\[
X^{'a}:=X-\bigcup_{j\in\grY_v}\hat{Y}_{\ul_j+a}^\circ.
\]
In particular,  \(X'\subset
X^{'a}\subset X''\) when \(0\leq a\leq 10\); \(X'=X^{'0}\) and 
\(X''=X^{'10}\).  The value of \(a\) in \(X^{'a}\) is left
unspecified in some statements in this article. In such case, the
statement holds for all value of \(a\) with  \(0\leq  a\leq 10\), and
the precise value of \(a\) does not matter.

\item \(X\) has at least one Morse end.
\end{itemize}
\end{defn}
A \(\Spin^c\) MCE \((X,\grs)\) together with an admissible form \(\nu\) on
it, denoted \((X, \nu)\), is called an {\em admissible pair}.

In the case when
\(X=\mathbb{R}\times Y\) and \((X, \nu)\) is invariant under the
\(\mathbb{R}\)-action, the admissible pair is said to be {\em
  cylindrical}. In this case \(X\) has two ends, \([1, \infty)\times
Y\) and \((-\infty, -1]\times Y\simeq \mathbb{R}^{\geq 0}\times
(-Y)\). We call the former the \(+\infty\)-end, and the latter the
\(-\infty\)-end. 

All pairs of MCE's and 2-forms \(\nu\), \((X, \nu)\) are assumed to be
admissible in this article. 
\begin{rem}\label{rem:cyl-end}
In Definition \ref{def:mce}'s terminology, the
ending 3-manifolds of a cylindrical \(X=\bbR\times Y\) are \(Y\) and
\(-Y\). 
\end{rem}

Note that because the ending manifolds \(Y_i\) are compact, \(\nu\) decays
exponentially to \(\nu_i\) in the 
\(C^k\)-norm. (Cf. e.g. (\ref{eq:xi-exp})  for a more precise
statement.) Meanwhile, the second bullet in the
Definition \ref{def:adm} can be met with a generic choice of \((X, \nu)\)
and an appropriately chosen \({\bf l}_v\). Fix such a choice. The
condition also implies that over \(X''\) 
the restriction of \(|\nabla\nu|\) to \(\nu^{-1}(0)\) is bounded below
by a positive number. 

 The precise definition of  {\em admissible} Seiberg-Witten solutions
 associated to an admissible pair \((X, \nu)\) and their \(Y_i\)-end
 limits are given in \textsection \ref{sec:adm-SW}
 below. It suffices to say for now that this is basically a finite
 energy condition that prescribes the asymptotic behavior of the
 solutions. For a Morse
end, the asymptotic conditions are specified by certain strongly t-convergent
Seiberg-Witten solutions associated to the ending pair. 

On ``vanishing ends'',
the asymptotic conditions are of the kind appearing in typical
Seiberg-Witten-Floer theories, as described in \cite{KM}.  

To ensure
transversality, an additional nonlocal perturbation term,
\(\hat{\grp}\), is added to the left hand side of the Seiberg-Witten
equation in \cite{KM}. For this purpose, $\hat{\grp}$ can be set
to be zero on the Morse ends, 
and may be taken to be arbitrarily small on the vanishing ends. 
Let $\bbS=\bbS_X=\bbS^+\oplus \bbS^-$ denote the spinor bundle associated to
the $\Spin^c$ structure on $X$. Given a \(\Spin^c\)-connection on
\(\bbS\), let $A$ denote the associated unitary connection on
$\det\bbS^+$, and let $\slp_A^+\co \bbS^+\to \bbS^-$ denote the
associated Dirac operator. Let $\Psi\in \Gamma (\bbS^+)$. 
The general Seiberg-Witten equation on \(X\) takes the form:  
\begin{equation}
\mathfrak{S}_{\mu, \hat{\grp}}\, (A,
\Psi)
:=\Big(\frac{F_{A}^+}{2}-\rho^{-1}(\Psi\Psi^*)_0+i\frac{\mu^+}{4},
\slp_A^+\Psi\Big)+\hat{\grp}\,( A, \Psi),\label{eq:SW4}
\end{equation}
where  \(\mu^+\) denotes
  the self-dual part of a 2-form \(\mu\). (Cf. \textsection \ref{sec:convention}).
Note that \(\mathfrak{S}_{\mu, \hat{\grp}}\, (A,
\Psi)\) is invariant under {\em gauge transformations} \((A, \Psi)\mapsto
(A-2u^{-1}d u, u\cdot\Psi)\), where \(u\) maps \(X\) to
\(S^1\). The gauge equivalence class of \((A, \Psi)\) represented by
\((A, \Psi)\) is denoted by \([(A, \Psi)]\).

\begin{thm}[{\bf local convergence}]\label{thm:l-conv}
Let \(X\) be a 4-dimensional \(\Spin^c\) MCE and \(\nu\) be an
admissible 2-form on \(X\).  
Given \(r\geq 1\), let \(w_r, \hat{\grp}_r\) be
respectively a closed 2-form and a nonlocal perturbation
of the type described in \cite{KM}.  Suppose they satisfy Assumption \ref{assume} below. 
Let \(\mu _r:=r\nu +w_r\), and write
\[ \text{
\(\omega=2\nu^+\); \(\mu_r^+=\frac{r}{2}\omega +w_r^+\). }
\]
Fix a
sequence of positive real numbers \(\Gamma:=\{r_n\}_n\), \(r_n\to\infty\) as
\(n\to \infty\),  and a corresponding sequence of admissible solutions
\((A_{r}, \Psi _{r})\) to the Seiberg-Witten equation
\(\mathfrak{S}_{\mu_{r}, \hat{\grp}_r}(A_r, \Psi_r)=0\) for
\(r\in \{r_n\}_n\). Suppose that the sequence \(\{(A_{r}, \Psi
_{r})\}_{r\in \Gamma }\) satisfies: 
\begin{itemize}
\item [(0)] For each \(i\in \grY\) and \(r\), The \(Y_i\)-end limit of
  \((A_r, \Psi _r)\), denoted \((B_{i,r}, \Phi _{i,r})\), is 
  regular in the sense that the linearization of \(\grF_{\mu _{i,r}}\)
  from (\ref{eq:SW3}) is surjective. Here, \(\mu _{i,r}\) denotes the
  \(Y_i\)-end limit of \(\mu _r\) in the sense of
  Section \ref{sec:convention} (6).  (Such a limit exists by the assumptions
  on \(\nu \) and \(w_r\).)
\item [(1)]  If \(\hat{Y}_i\) is a vanishing-end of \(X\), then the \(Y_i\)-end limit
  of \((A_r, \Psi_r)\), denoted \((B_i, \Phi _i)=(B_{i,r}, \Phi _{i,r})\), is independent of \(r\). If \(\hat{Y}_i\) is a
  Morse-end of \(X\), then the \(Y_i\)-end limit of \((A_r, \Psi_r)\) 
  strongly t-converges to a t-orbit \(\pmb{\gamma}_i\) as $r\to
  \infty$. The \(Y_i\)-end limits of \((A_r, \Psi _r)\) are
  regular in the sense that 
\item[(2)] \(\exists \, N\in \bbR^+\) such that \(\forall r>N\), the
  relative homotopy classes of  \([(A_r, \Psi_r)]\) are identified
  via the canonical isomorphisms in Lemma \ref{rem:rel_class}.
Denote this class by $\mathfrc{h}$. 
\end{itemize}
Then:
\begin{itemize}
 \item[(a)] There is a
 t-curve \(\mathbf{C}\) and an unbounded subsequence \(\Gamma_0\) of 
 \(\Gamma =\{r_n\}_n\subset \bbR^+\), such that \(\{(A_r, \Psi_r)\}_{r\in \Gamma_0}\)
 t-converges to \(\mathbf{C}\) over \(X'\). 
\item[(b)] Furthermore,
 \(\mathbf{C}\) has finite $\omega$-energy, which is bounded above by a positive constant determined by 
$\mathfrc{h}$ (which depends on  \(\pmb{\gamma}_i\)), the metric and
the \(\Spin^c\)-structure of \(X\), \(\nu \)
and the constants \(\varsigma_w, \zzz_\grp\) in Assumption \ref{assume}.
\end{itemize}
\end{thm}
\begin{rem}
(a) See \textsection \ref{sec:SW-class} for the definition of 
relative homotopy classes of Seiberg-Witten solutions. The notion of
  $\omega$-energy of t-curves, as well as the $\nu$-energy
  mentioned in the next theorem, are introduced in \textsection
  \ref{sec:energy} below. They are analogs of the notion of energies in symplectic field theory in \textsection 5.3 of
\cite{BEHWZ}.

(b) By Theorem \ref{thm:strong-t}, 
the sequence \(\{(A_{r}, \Psi
_{r})\}_{r\in \{r_n\}_n}\) in the statement of the preceding theorem
exists only when  \(c_1(\grs_i)\cdot[*\nu_i]\geq -\zeta _{*\nu _i}\) for
all Morse ends \(\hat{Y}_i\). 
\end{rem}
For Floer-theoretic purposes, a stronger version of convergence result
is often desired. In heuristic terms, a gauge equivalence class
of Seiberg-Witten solutions on a 4-dimensional cobordism is regarded as
a (generalized) flow line. 
This type of  results say that a
sequence of flow lines with the same end points and relative homotopy class 
``weakly converge'' to a ``broken trajectory'' that is  a
concatenation of flow lines. (See Theorem 5.1.1 in \cite{KM} for an
example.) The following theorem is intended
to play a similar role. The ``broken trajectories'' here in our
setting are the so-called ``chains of t-curves''. 

\begin{thm}[{\bf global convergence}]\label{thm:g-conv}
Adopt the assumptions and notations in Theorem \ref{thm:l-conv}; in
particular,  \(\{\pmb{\gamma}_i\}_{i\in \grY_m}\), \(\mathfrc{h}\) are as in
Conditions (1) and (2) in Theorem \ref{thm:l-conv} respectively. Let
\(\mathfrak{s}=\grs_X\) denote the \(\Spin^c\) structure of \(X\).  
Then
\begin{itemize}
\item[(a)] There exists a chain of t-curves \(\mathfrak{C}\) in \(X'\)
with \(\Spin^c\) structure \(\mathfrak{s}\),
\(Y_i\)-end limits \(\pmb{\gamma}_i\) for \(i\in \grY_m\), 
together with a
subsequence \(\Gamma'=\{r'_n\}_n\) of \(\Gamma\), \(r'_n\to \infty\) as
\(n\to\infty\), such that over \(X'\), the corresponding subsequence of Seiberg-Witten
solutions \(\{(A_r, \Psi_r)\}_{r\in\Gamma'}\) weakly t-converges to
\(\grC\) in the sense of Definition \ref{def:w-t-conv}. 

\item[(b)] 
\(\grC\) has finite \(\omega\)-energy, which has 
an upper bound determined by the same factors listed in Item (b) in the
statement of Theorem \ref{thm:l-conv}. 

\item[(c)] Suppose in addition that: 
\begin{equation}\label{b_1=0}
b^1(Y_i)=0 \quad \forall i\in \grY_v, 
\end{equation}
and let  \(\grh\) be the isomorphism  introduced in Lemma
  \ref{lem:htpy}. (This is a map from  the set of relative
  homotopy classes of Seiberg-Witten gauge equivalence classes to the
  set of relative homology classes of chains of t-curves.)  Then
  \(\grC\subset X'\) is of relative homology class \(\grh\,
  (\mathfrc{h})\).  

\item[(d)] Assuming that \(X\) has no vanishing ends.  Then the
\(\nu\)-energy of \(\grC\) is determined by \(\nu \), the
\(\Spin^c\)-structure, \(\{\pmb{\gamma }_i\}_i\) and \(\fh\) via the explicit formula
(\ref{eq:F_vgrC}). 
\end{itemize}
\end{thm}


See Definition \ref{def:chain} below for the definitions of chains of t-curves, as
well as their \(\Spin^c\) structures and relative homology
classes. 
\begin{rem}
As with typical Floer theories, constraints on the relative homotopy class (Condition (2) of
Theorem \ref{thm:l-conv} in our setting; also assumed in Theorem
\ref{thm:g-conv}) may be replaced by the more general condition of
a uniform upper bound on a certain ``topological energy''. Compare
e.g. \cite{KM}'s Theorem 24.6.2 and Proposition 24.6.4. The
corresponding ``energy-bound'' condition in our setting takes the form
of:  
\begin{equation}\label{eq:CSD-est0}
r^{-1}\scrE_{top}^{\mu_r, \hat{\grp}_r}(X)(A_r, \Psi_r)\leq \smE  \quad \text{for an \(r\)-independent constant \(\smE\),}
\end{equation}
where \(\scrE_{top}^{\mu, \hat{\grp}}\) is a natural generalization of
\cite{KM}'s ``perturbed topological energy''
\(\scrE_{top}^{\hat{\grp}}\). 
(Cf. (\ref{assume:EtopX-ubdd}) below.) Fixing the relative homotopy class \(\mathfrc{h}\) as in Condition (2) of
Theorem \ref{thm:l-conv} leads to a bound of the preceding form (with 
\(\smE_p\) depending on \(\mathfrc{h}\)); see Lemma \ref{lem:E_topX}
below. As explained in \textsection \ref{sec:energy-comp}, such a
bound on the (Seiberg-Witten) topological energy is virtually 
equivalent to a bound on the corresponding symplectic version of
energy, which is  required in Gromov-type compactness theorems. 
A special case of such a generalization is used in the proof of \(HM=PFH\) in \cite{LT}.
\end{rem}

\subsection{Motivations and applications}

This article is intended to provide some essential analytic foundations of
various speculated relations between Seiberg-Witten-Floer homologies
and \(Gr\)-type Floer homologies (such as Hutchings-Taubes's \(ECH\)/\(PFH\)
\cite{Def:ECH}). More specifically, it began as a part of the author's
program of proving the equivalence of Seiberg-Witten-Floer homology (\(\mathring{HM}\))
and Heegaard Floer homology (\(HF^\circ\)) \cite{L}; and the earliest incarnation of
the present paper dates back to 2006. Since then, several major progress
has been made 
on the relation between the two types of Floer theories. In
particular, due to certain technical difficulties, the program in
\cite{L} has been since been modified to take advantage of Taubes's
monumental work on the equivalence of \(HM\) and \(ECH\) in 2008
\cite{T:ech}. See \cite{KLT1}, \cite{KLT2}, \cite{KLT3},
\cite{KLT4}, \cite{KLT5}. This renders the results stated here not immediately
relevant to our original goal (which partially accounts for its delay
in appearance). However, the aforementioned relation comes in many
forms, (some even lacking precise formulation), and is still very far
from being 
completely understood. This article is expected to serve as technical
basis for possible results in this general direction. Here we briefly describe
some of them, as well as their background. 

\subsubsection*{\it (a) \(PFH\) and TQFT}
The definition of purported \(Gr\) analogs of \(HM\) is itself a difficult
subject. Certain variants, called Periodic Floer homology (\(PFH\)) and
Embedded Floer homology (\(ECH\)), are proposed by Hutchings, and shown to
be well-defined by Hutchings and Taubes \cite{Def:ECH}. The latter (\(ECH\)) is
associated to \(\Spin^c\) 3-manifolds equipped with a contact form. As
mentioned before, it is shown to be equivalent to the (unperturbed
version) of \(HM\) by Taubes. In fact, some of the expected properties of
\(ECH\) are difficult to establish directly, and this relation with
the better-developed \(HM\) are used to justify them, resulting in
important consequences.
For example, the independence of \(ECH\) on
the contact form and the definition of cobordism
maps in \(ECH\) both hinge on their counterparts in \(HM\), and these are
in turn the basis for the proofs of 3-dimensional Weinstein conjecture
\cite{T:W}, Arnolds's conjecture \cite{Arnold}, and Hutchings's
\(ECH\) capacities \cite{H}. 

The closely related \(PFH\) pertains to a special case of the setting studied in this
paper. These are defined for 3-manifolds that are mapping tori, and
they are shown to be equivalent to (a perturbed version) of \(HM\) in
the author's joint work with Taubes \cite{LT}. Given Taubes's prior
work on \(HM=ECH\), the only new ingredient of this proof 
consists of  a special case of this article's Theorem
\ref{thm:l-conv}, namely the case when \((X, \nu)=(\bbR\times Y,
\pi_2^**_3\theta)\) is cylindrical, and the harmonic 1-form \(\theta\) is nowhere
vanishing. One naturally expects that this equivalence extends to an
equivalence of the TQFT on both sides, and the theorems of this
article in that case when \(X\) is a symplectic cobordism would form
the main part of the proof. Unfortunately, like \(ECH\), cobordism
maps for \(PFH\) is not at present (directly) defined. In \(ECH\), it is defined
indirectly through its Seiberg-Witten counterparts by Taubes's
\(HM=ECH\)).  A similar indirect definition of PFH
cobordism maps via the main theorems of this article is very likely
possible. Alternatively, for heuristic reasons one expects the
purported TQFT for \(PFH \) to be equivalent to Usher's TQFT for what
he calls  the ``FCOB category'' \cite{U}. based on
Donaldson's result on the equivalence of (closed) symplectic
4-manifolds and Lefschetz pencils, Usher introduces in \cite{U} a
version of TQFT on what is called the ``FCOB category'' in
\cite{U}. The latter \(TQFT\) is rigorously defined; and therefore one
has a precisely formulated conjecture:
\begin{conj}\cite{U}\label{conj:U}
Usher's TQFT is equivalent to the restriction of
(perturbed versions of)
Seiberg-Witten \(TQFT\), restricted to the subcategory of \(FCOB\).
\end{conj}
The main results of this article are expected to play a key role in 
a proof of this conjecture. 

The aforementioned heuristic reasoning also leads to other variants of
this conjecture: In the context of closed symplectic manifolds, The
TQFTs of Usher's and of \(PFH\)'s are respectively the counterparts of 
Donaldson-Smith's invariant \cite{DS} for symplectic
Lefschetz pencils, and Taubes's Gromov invariant for symplectic
4-manifolds.  The equivalence of the two has been established by Usher
\cite{U1}.  Hence by Taubes's \(SW_4=Gr\) \cite{T}, the
Donaldson-Smith invariant is equivalent to the Seiberg-Witten
invariant for closed 4-manifolds. This is the counterpart of
Conjecture \ref{conj:U} above. 

As \(PFH\) concerns the setting of Theorem
\ref{thm:strong-t} in the special case when the harmonic \(1\)-form has no
zeroes, in full generality Theorems
\ref{thm:strong-t} and \ref{thm:l-conv}, \ref{thm:g-conv} suggest that there
should be a generalized version of the \(PFH\) TQFT, which is defined for
a larger category containing \(FCOB\) as a subcategory, whose objects
consist of closed 3-manifolds with \(b_1>0\) (equipped with a
nondegenerate harmonic 1-form \(\theta\)), and whose morphisms are 4-dimensional
cobordisms with \(b^+>0\), equipped with nondegenerate harmonic
self-dual 2-forms that appear as \(\omega\) in the setup of Theorems
\ref{thm:l-conv}, \ref{thm:g-conv}. The closed 4-manifold counterpart of
this larger category consists of pairs \((X, \omega)\) described in
the statement of  Taubes's Theorem \ref{thm:T1}. Such manifolds are
called  ``near symplectic manifolds'' by
Donaldson and some other authors. Donaldson's description of symplectic manifolds
as Lefschetz pencils has been generalized to near-symplectic
4-manifolds \cite{ADK}. This endows every near-symplectic manifold
with a ``singular Lefschetz pencil''. Correspondingly, Perutz defined
a ``Lagrangian-matching invariant'' for singular Lefschetz fibrations
\cite{P}, which can be regarded as a generalization of
Donaldson-Smith's invariant for symplectic Lefschetz invariants. A
TQFT version of this ``Lagrangian-matching invariant'' is briefly
outlined in Section 1 of \cite{P}, and currently being developed by
Lekili and Perutz \cite{LP}. This is supposed to be equivalent to the
putative ``generalized \(PFH\) TQFT mentioned above, and Conjecture
\ref{conj:U} has a corresponding generalization in terms of these
invariants.

A more modest project is to understand the relation between \(HM\)
and \(PFH\) cobordism maps in special cases. Hutchings communicated to us that \(PFH\) cobordism maps
can be well-defined for some simple cobordisms equipped with Lefschetz
fibration structure. G.H. Chen is working on generalizing these and
establishing the equivalence of \(HM\) and \(PFH\) cobordism maps in
these special cases. 
Also according Hutchings, one complication of directly
define \(ECH\)/\(PFH\)-cobordism maps is that, unlike the case of Floer
homologies, such maps would inevitably involve counting
(nontrivially-) broken trajectories (analogs of our ``chains of 
t-curves'' consisting of more than one t-curves), even after
generic perturbation. We hope that by
exploring some explicit examples of Theorem \ref{thm:g-conv} one may
gain some insight on how to correctly count holomorphic curves in
cobordisms for the definition of \(ECH\)/\(PFH\)-cobordism maps.

\subsubsection*{\it (b) Seiberg-Witten invariants and Heegaard Floer
  invariants.}
One common essential ingredient of both the proof of \(HM=HF\) in
\cite{KLT1}-\cite{KLT5} and the original program \cite{L} is a positivity result, exemplified by 
Propositions 3.4 and 3.7 in \cite{KLT5}. (A special case of this
result is also a major ingredient of \cite{KLT4}).  Their counterparts
in \cite{L} are Claim (1) of its Section 6.4, and the last part of its
sketch of proof for Theorem 8.2. These type of positive results serve
the purpose of showing that the differential and cobordism maps
between Floer complexes preserve the filtration, ensuring that the
``filtered'' Floer homologies and cobordism maps among them are
well-defined. 

The following proposition is a by-product of the proof of Theorem
\ref{thm:g-conv}, given in Sections \ref{sec:g-conv:a}, \ref{sec:g-conv:b}
below. The positivity results needed in \cite{L} are special cases of
this proposition.
\begin{prop}\label{cor:F-positive}
Let \(\{(A_r, \Psi_r)\}_{r\in \Gamma'}\) and \(\grC
\) be respectively the sequence of
Seiberg-Witten solutions and the chain of t-curves from the conclusion of Theorem \ref{thm:g-conv}. 
Suppose that \(P\subset \mathring{X}'\) is
a 
pseudo-holomorphic submanifold with the following properties: 
\BTitem\label{def:P}
\item \(P\) is disjoint from \(\nu ^{-1}(0)\), and its distance to the
  latter is bounded below by a positive number. 
  \ETitem
   \def\Ttag{\ref{def:P}}
   \BTTitem
   \item  On each Morse end \(\hY_i\), \(P\) is asymptotic to a union of
  flow lines of the vector field dual to \(*\nu _i\), denoted \(\pp_i\), 
  in this sense of Definition \ref{def:asymp-cur}. (\(\pp_i\) may be
  empty.) 
  For every \(i\in \grY_m\), \(\pp_i\) is disjoint from all the rest
  orbits of \(\grC\). 
   \item When \((X, \nu )\) is non-cylindrical: Write \(\grC=\{\bfC_0,
  \{\grC_i\}_{i\in \grY_m}\}\). Then \(P\) intersects \(C_0\) transversely in a discrete set in the interior of \(P\), and
  for each \(\grC_i\) with at least one component, 
  \(\bbR\times \pp_i\) intersects each component of \(\grC_i\)
  transversely in a discrete subset  in the interior of
  \(\bbR\times \pp_i\). When \((X, \nu )\) is cylindrical: Let
  \(\pp_\pm\) denote the limit of \(P\) in the \(\pm\infty\)-end. Then
  every component of \(\grC\) intersects all three of \(P\),
  \(\bbR\times \pp_\pm\)  transversely in discrete points in the interior; and for every
  rest orbit \(\pmb{\gamma }\) of \(\grC\), \(P\) intersects
  \(\bbR\times \pmb{\gamma }\) transversely in discrete points in the
  interior. 
   \ETTitem
Then
\(\lim_{r\to \infty}\int_P\frac{i}{2\pi}F_{A_r^E}\) exists and is a
non-negative integer. (See \textsection \ref{defn:t-convergence} for the definition of
\(A^E\)). Moreover, if \(\grY_v=\emptyset\), this limit is determined
by the relative homotopy class \(\mathfrc{h}\) (which depends
implicitly on \(\nu \) and \(\{\pmb{\gamma }_i\}_{i\in \grY}\)) and
the (relative) homology class of \(P\) relative to \(\{\pp_i\}_{i\in
  \grY}\). In the above, two pseudo-holomorphic subvarieties \(P, P'\)
are said to be homologous relative to \(\{\pp_i\}_{i\in\grY}\) if they
both satisfy (\ref{def:P}) with the same \(\{\pp_i\}_{i\in\grY}\), and
\(P-P'\) bounds a chain with closed support (cf. the footnote in
Section \ref{sec:gr-spin}).
\end{prop}
See Definition \ref{def:chain} for terminology and notation (e.g. rest
orbits, components) in the preceding statement. 

In fact, with this proposition in place, \cite{KLT5} can be
significantly shortened. This will be explained in more detail in the
next subsection. Its generality also renders the definition of
filtered Seiberg-Witten counterparts of some of the extensions of
\(HF^\circ\) immediate. For example, the knot
Floer homology \(HFK^\circ(Y, K)\) in \cite{OS:knot} has a \(HMT^\circ\)
counterpart \(HMT^\circ(Y, K)\) by a straightforward translation of Ozsvath-Szabo's 
construction to the large-perturbation Seiberg-Witten theory according
to the prescription of Section 6 of \cite{L}, by introduction to a second
filtration to \(CMT^\circ\) associated to \(\ud{\gamma}_w\), the
1-cycle in \(\ud{Y}\) corresponding to the second base point \(w\) in
\cite{OS:knot}. Modifying \(HMT^\circ(Y, K)\) in the same way as
\cite{KLT1}-\cite{KLT5}, one also has an extension of the latter's
\(H^\circ\) to knots. Theorem 1.1
of \cite{KLT5} states that \(H^\circ\) is isomorphic to \(HMT^\circ\)
tensored with copies of \(H_*(S^1;\bbZ)\). The aforementioned
extension of \(H^\circ\) is related to \(HMT^\circ(Y, K)\) in the same
manner.

\subsection{Relation to existing literature}\label{sec:literature}
The main theorems of this article may also be viewed as an extension
of  a counterpart of
\cite{Arnold2}'s Proposition 5.2, which partially generalizes the \(SW\Rightarrow Gr\)
part of Taubes's \(HM=ECH\) proof in \cite{T:ech}. In the latter, 
\(\bbR\)-invariant large-perturbation Seiberg-Witten equations on \(\bbR\times Y\) are
considered, with the role of \(*_3\nu\) here played by a contact form. This is
extended to the case of 4-dimensional exact symplectic cobordisms in
\cite{Arnold2}.  The generalization is straightforward, as it does not
require new techniques. 
Taubes's equivalence theorem has a sister version in \cite{LT}, which
asserts that (a different variant of) \(HM\) is isomorphic to
\(PFH\). Much of the proof for this theorem is similar to that in
\cite{T:ech}, except for the \(SW\Rightarrow Gr\) part. In both settings, a
crucial ``energy bound'' is required. In \cite{T:ech}, this 
hinges on certain spectral flow estimates, while in \cite{LT}, a more
topological argument is used, basing on 
harmonicity of \(\nu\). The latter is also the strategy of this
article; see Section \ref{sec:4}. In fact, what was proved in \cite{LT}  is the special case of our
Theorem \ref{thm:l-conv} when \(X\) is cylindrical and \(\nu\) is
nowhere-vanishing. Thus, one may expect that, like what happened in
the \(ECH\) case (Proposition 5.2 of \cite{Arnold2}), the results of this article can
also be obtained by simple modifications of the work in \cite{LT}. This
unfortunately turns out {\em not} to be the case. 

In the case of 
\cite{Arnold2}, the contact condition makes it possible to choose the
perturbation form on \(X\) so that it is identical with that in the
cylindrical case of \cite{T:ech}. In the setting of this article, the
condition that \(\nu\) be harmonic is rather rigid,  and the behavior of
the perturbation form on the ends of \(X\) is determined by that on
its compact piece, \(X_c\). In general, it does not agree with
the cylindrical case, but only asymptotic to it.  Surprisingly, this
small difference turns out to be a major obstacle, making the desired results 
substantially more elusive. See remarks near the end of Section
\ref{sec:top-energy} and in the beginning of Section \ref{sec:E_top-lower}. 

Another complication stems from the existence of zero locus of
\(\nu\). In \cite{Ts}, the harmonic form \(\omega\) is also allowed to vanish along a
1-dimensional submanifold of \(X\). While part of
the ingredients in our proof can be obtained from 
straightforward adaptation of their counterparts in \cite{Ts},  the
same existence of zero locus causes additional
problems in our setting. In \cite{Ts}, much of what was done in
\cite{T} easily carries over if one restricts attention to 
complement of small tubular neighborhoods of \(\omega^{-1}(0)\). This
is no longer possible in our setting, 
since several important steps here rely on global estimates where
contribution from tubular neighborhood of \(\nu^{-1}(0)\) can not be 
ignored. Lack of knowledge on the behavior of Seiberg-Witten
solutions on this region thus becomes an unavoidable problem. 

The purpose of \cite{KLT1}-\cite{KLT5}'s modification to the original program in
\cite{L} is precisely to sidestep the aforementioned issues. In the
particular setting of \cite{L}'s program, the zero locus in the
cylindrical \(\bbR\times M\) is explicitly known. Thus, one may
perform certain surgery operation along the zero locus of \(*_3\nu\),
and ``perturb'' \(\nu\) near the surgered part so that \(*_3\nu \) approximates
certain standard contact 3-form. This makes \(\nu\) nowhere-vanishing,
and also allows one to take advantage Taubes' proof of \(HM=ECH\),
which was not available when \cite{L} was written. 

The major part of \cite{KLT5} is dedicated to proving a 
variant of Proposition \ref{cor:F-positive} for the connecting sum
cobordisms. (Cf. Propositions 3.9-3.14 in \cite{KLT5}). Here one also
encounter the first issue above, arising from the asymptotic
behavior of the perturbation form. Because of this problem, the
results there only applies to very special types of manifolds equipped
with rather stringent geometric conditions, 
essentially only the specific cobordisms used in the article. For example, the perturbation form there is required to be
\(\bbR\)-invariant on the ends (cf. e.g. (2.12) of \cite{KLT5}), and
the pseudo-holomorphic submanifold \(P\) (in Proposition
\ref{cor:F-positive}'s notation) is a cylinder, near which the metric
is ``standard''. As mentioned before, this does not
hold in most cases. Therefore the entire Section 9 therein (about 60
pages) is
dedicated to constructing such metrics and forms case by case for the specific
cobordisms needed in \cite{KLT5}. Except for its Proposition 3.13,
where \(X=\bbR\times Y\) and the perturbation form is one that
interpolates the harmonic version of \cite{L} and the modified version
of \cite{KLT1}-\cite{KLT5}, all other positivity results in
\cite{KLT5} are special cases of Proposition \ref{cor:F-positive}
above, and its Proposition 3.13 says that as far as \cite{KLT5} part
of the program is
concerned, it makes no difference to either work along the original
plan in \cite{L} and or  its modified form in \cite{KLT1}-\cite{KLT5}.

Alternatively, an
avatar (2006-2008) of the present article also assumed that the
perturbation form is constant on the ends, and this is made
permissible by only requiring the form \(\nu\) to be Lipshitz, not
necessarily smooth. (The assumptions of  \cite{KM} is sufficiently
general 
to apply in this case). In this way, the proof therein was no more
difficult than that in \cite{LT}. The trade off is that one would have
to deal with Lipshitz almost complex structures and Lipshitz
holomorphic curves.

One should also note that in \cite{KM}, the perturbations used to
define the differentials in \(HM\) and the cobordism maps are also
assumed to be constant on the ends. Though this is in general not the
case when applying our theorem, we believe that the results in
\cite{KM} may be easily generalized to our setting, as our
perturbations decay to constant ones on the ends. One should be able
to do so by either going through the arguments in \cite{KM}, or use a
direct limit argument. 

\subsection{Convention, notations, etc.}\label{sec:convention}

Here is a list of notations and conventions frequently
used throughout this article.
\begin{enumerate}
\item \(W^\circ\) or \(\mathring{W}\) denotes the interior of a
  topological space \(W\). 
\item \(p.d.\) stands for ``Poincar\'e dual''.
\item Given an interval \(I\subset \bbR^{\geq 0}\),  we use \(\hat{Y}_{i,
    I}\) to   denote the subset \(\mathfrc{s}_i^{-1}I\) in the end \(\hat{Y}_i\),
  and \(\hat{Y}_{i, L}:=\hat{Y}_{i, [L, \infty)}\). The latter should
  be contrasted with  \[Y_{i:L}:=\mathfrc{s}_i^{-1}(L).\]
We shall often omit
  the index \(i\) when the statement refers to any of 
the ends of \(X\). The interval \(I\) above is allowed to be either
open or closed, or of the mixed type \([L, L')\) or \((L,
L']\). We use \(\hat{Y}_{i, [L, \infty]}\simeq [L, \infty]\times Y_i\supset
[L, \infty)\times Y_i\simeq \hat{Y}_{i, L}\) to denote the natural compactification of
\(\hat{Y}_{i, L}\) by adding a \(Y_{i:\infty}\simeq Y_i\) at
infinity. Correspondingly, \(\ov{X}=X_c\cup \bigcup_{i\in \grY}\hat{Y}_{i, [L,
  \infty]}\supset X\) denotes the compactification of \(X\) by adding
\(\bigcup_{i\in \grY}Y_{i:\infty}\) at infinity. Given a subspace \(M\subset
  \ov{X}\), \(\ov{M}\) denotes its closure in
  \(\ov{X}\).
\item \(\pi _2\co I\times Y\to Y\) denotes the projection to the
  second factor of the product.  

\item Given a bundle \(\rmV\)
  over \(M\), \(C^k_{A}(M;
  \rmV)\), \(L^p_{k, A}(M; \rmV)\) respectively denote the spaces of
  \(C^k\)- and \(L^p_k\)-sections associated to a connection \(A\) on
  \(\rmV\). Occasionally \(M\), \(V\), or the connection \(A\) is
  omitted from the notation when it is clear from the
  context.  The bundles appearing in this article are typically
  constructed out of \(T^*M\) and/or the spinor bundle over
  \(M\). When the connection \(A\) is omitted from the notation, e.g.  
  \(L^p_k(M)\), the implicit connection \(A\) that induced from the Levi-Civita connection on
  \(T^*M\), and/or the (previously specified) reference connection
  \(A_0\) on the spinor bundle. 
\item Given a bundle \(\rmV\) over \(\ov{X}\) with
\(\rmV_i:=\rmV\big|_{Y_{i: \infty}}\), a section \(q\in \Gamma
(\rmV\big|_{X})\) is said to have \(q_i\in \Gamma (\rmV_i)\) as a {\em
  \(Y_i\)-end limit} if  
\[
\|q-\pi _2^*q_{i}\|_{C^2(\hat{Y}_l)}+\|q-\pi _2^*q_{i}\|_{L^2_1(\hat{Y}_l)}\to 0 \quad\text{as \(l\to \infty\).}
\]
 \item 
Given a function \(\mathbf{L}\co \grY\to [0, \infty]\) with
\(\bfL(i)=L_i\), let \(X_{\mathbf{L}}\subset \ov{X}\) denote the subset
\(\{ x\, |\, \mathfrc{s}_i(x)\leq L_i  \, \forall i, \, \text{wherever
  \(\mathfrc{s}_i\) is defined}\}\). Given two functions \(\bfl, \mathbf{L}\co \grY\to [0, \infty]\) We write
\(\mathbf{l}<\mathbf{L}\) if \(X_{\mathbf{l}}\subset
X_{\mathbf{L}}^\circ\), and we write \(\mathbf{l}\leq\mathbf{L}\) if \(X_{\mathbf{l}}\subset
X_{\mathbf{L}}\). 
Also, let \(\hat{Y}_{[{\bf l}, {\bf L}]}:=\bigcup_{i\in \grY}
\hat{Y}_{i, [l_i, L_i]}\).


\item The notation \(X_\bullet\) will be used to denote any of the
  \(X_{\bf L}\), \(\hat{Y}_{i, I}\) above, or their interiors. Given \(X_\bullet\subset
  \bar{X}\), \(\ov{X_\bullet}\) denotes its closure in
  \(\bar{X}\). 
The {\em length}  of
  \(X_\bullet\), denoted by \(|X_\bullet |\), is given by:
  \(|X_\bullet|:=|\ov{X_\bullet}|\in [0, \infty]\) and 
\[
|X_\bullet|:=\begin{cases}
| I |=L-l & \text{when \(X_\bullet=\hat{Y}_I\), and \(I=[l, L]\).}\\
1+\sum_iL_i &\text{when \(X_\bullet=X_{\bf L}\), \({\bf L}(i)=:L_i\).}
\end{cases}
\]  
\item 
 Given an admissible pair \((X, \nu )\), let  \(X_\delta \subset X\) denote the subspace of \(X\) consisting of
  points whose distance to \(\nu ^{-1}(0)\) is greater than \(\delta
  \), and let \(Z_\delta :=X-X_\delta \subset X\).  Given 
  \(X_\bullet\subset \ov{X}\), \(X_{\bullet, \delta }:=X_\delta \cap
  X_\bullet\) and \(Z_{\bullet, \delta }:=Z_\delta \cap
  X_\bullet\). For a
  closed 3-manifold \(Y\) and a 1-form \(\theta \) over \(Y\), 
  \(Y_\delta \subset Y\) is similarly defined. 
\item \(\tau_L\co [l, \infty)\times  Y\to [l+L, \infty)\times Y\)
  denotes translation in the first factor. Slightly abusing the
  notation, we also sometimes write \(\tau_{-L}^*\) as
  \(\tau_L\). I.e., given a function or section \(\xi(s, p)\) on a half cylinder
  \(\hat{Y}\simeq \bbR^{\geq 0}\times Y\) and \(L\in \bbR\), \((\tau_L\xi)(s,
  p):=\xi(s-L, p) \) on \([L, \infty)\times Y\).
\item Throughout this article, a ``\(k\)-dimensional current''
  \(\tilde{C}\) on an
  \(n\)-manifold \(X\) refers to a linear functional on the space of 
compactly-supported smooth \((n-k)\)-forms on \(X\), where the latter space
is endowed with the topology defined by the sup-norm. The
{\em convergence} of currents refers to the convergence in weak\(^*\)
topology. Given a smooth \((n-k)\)-form \(\nu\), we sometimes use
\(\int_{\tilde{C}}\nu\) or \(\langle \td{C}, \nu \rangle \) to denote the value of \(\tilde{C}\) on
\(\nu\). 
\item A (\(k\)-dimensional) {\em subvariety} \(C\) of a manifold \(X\) is the image of an immersion
  \(u_0 \co C_0\to X\), where \(C_0\) is a (\(k\)-dimensional)
  manifold (possibly disconnected or empty) containing a
  countable set \(\Lambda \subset C_0\) such that
  \(u_0\Big|_{C_0-\Lambda }\) is an embedding.
The immersion \(u_0\) is said to be a {\em representing map} of
\(C\). 
  Given a \(k\)-form \(\vartheta\) on
  \(X\), \(\int_C\vartheta :=\int_{C_0}u_0^*\vartheta \). The {\em current associated to \(C\)}, often
  denoted \(\td{C}\), refers to that given by \(\vartheta \mapsto \int_{\td{C}}\vartheta :=\int_C\vartheta \). An {\em
  irreducible component} of \(C\) is the closure of a component of
\(C-u_0(\Lambda)\). \(C\) is {\em irreducible} if it has only one
irreducible component. 

\item A {\em weighted subvariety} \(\bfC\) in a manifold \(X\) with
{\em  underlying subvariety} \(C\) is a set of pairs \(\{(C_a, m_a)\}_a\), where
\(C_a\) are the irreducible components of \(C\), and \(m_a\in
\bbZ^+\) is the {\em weight} associated to the irreducible component \(C_a\). The {\em current} associated to \(\bfC\), denoted
\(\td{C}\), is \(\sum_am_a\td{C}_a\). Since the underlying subvariety \(C\)
and its associated current \(\td{C}\) together determine the
weighted subvariety \(\bfC\), it is alternatively denoted as a pair
\(\bfC=[C, \td{C}]\).

\item Let \(\bfC=\{(C_a, m_a)\}_a\) be  an oriented weighted
  subvariety, and \(P\) be  an oriented submanifold in the
  oriented manifold \(X\).  Let \(u_0\co C_0\to X\) be the immersion representing the
      underlying subvariety \(C=\bigcup_aC_a\) of \(\bfC\). Write
      \(C_0=\coprod _aC'_a\), where \(C'_a\) is the connected
      component which is the preimage of \(C_a\) under \(u_0\). Then \(P\)
      intersects  \(\bfC\) 
      {\em transversely} if  
      \(u_0\) intersects \(P\subset X\) transversely. 
      The {\em intersection number}
      \(\# \, (P\cap \bfC):=\sum_am_a\, \chi\big((u_0|_{C'_a})^{-1}P\big)\), where \(\chi
      \) denotes the Euler characteristic. 
      
\item The letters \(\zeta\), \(z\), \(\sO\), or any of their decorated
  forms (e.g. \(\zeta _i\), \(\zeta '\)) 
  customarily denote positive constants independent of \(r\), \(\delta \), and \((A,
  \Psi )\) or \((B, \Phi )\), and
  whose precise values are immaterial. Unless
  otherwise specified, they can be taken so as to increase with each
  appearance. Likewise, \(r_0\geq 1\) is a constant whose precise value may 
  change with appearances. It is taken so as to increase with each
  appearance.
\item Unless otherwise specified, \(\chi\co \bbR\to \bbR\)
  denotes a smooth cutoff function supported on \((-\infty,2]\), and equal 1
  on \((-\infty,1]\). 
\item \(B(x, R)\) or \(B_x(R)\) denotes a geodesic ball centered at \(x\) with
  radius \(R\).
\item \(*_3\) and \(*_4\) respectively denote Hodge dual on
  3-manifolds and 4-manifolds. Sometimes we also use \(*_Y\) to denote
  the Hodge dual on the specific manifold \(Y\).
\item Given a \(\Spin^c\) 4-manifold \(X\) with spinor bundle
  \(\bbS^+\oplus \bbS^-\), we use \(\Conn (\bbS^+)\), \(\Conn (\det
  \bbS)\) to denote respectively the space of \(\Spin^c\)-connections
  on \(\bbS^+\), and the space of unitary connections on \(\det
  \bbS^+\). Given a 3-manifold \(Y\) with spinor bundle \(\bbS\), let
  \(\Conn (\bbS)\), \(\Conn (\det \bbS)\) be similarly defined. A
  \(\Spin^c\) connection induces a unitary connection on its
  associated determinant line bundle, and therefore we have a map from
  \(\Conn (\bbS^+)\) to  \(\Conn (\det \bbS^+)\) (or \(\Conn (\bbS)\)
  to \(\Conn (\det \bbS)\) in the case of 3-manifolds). This map is an
  isomorphism when \(H^1\) of the 3- or 4-manifold in question has no
  2-torsion; in particular, \(\Conn (\bbS)\simeq\Conn (\det \bbS)\) in
  the case of 3-manifolds. 
Abusing notation, we often use the same notation (\(A\) in the case of
4-manifolds and \(B\) in the case of 3-manifolds) to denote either a \(\Spin^c\)
  connection or its induced connection on the associated determinant
  bundle. It should be clear from the context which is
  meant. 
E.g. \(A-A_0\) and \(F_A\) are
  \(i\bbR\)-valued differential forms; so the connection \(A\) in
  these cases must refer to a connection on \(\det \bbS^+\). The \(A\)
  appearing in covariant derivatives \(\nabla_A\) or Dirac operators
  \(\slp_A^+\) on a spinor-bundle necessarily refer to a \(\Spin^c\)-connection. 
In the same token, \(\Conn (X)\) (resp. \(\Conn (Y)\)) is used to
denote either  \(\Conn (\bbS^+)\) or  \(\Conn (\det \bbS^+)\) (resp.  \(\Conn (\bbS)\)
  or \(\Conn (\det \bbS)\)). 
\item We frequently use the shorthand \(\partial_s\) for
  \(\frac{d}{ds}\).
\item ``LHS'' and ``RHS'' respectively stand for ``left hand side''
  and ``right hand side''.
\item Given a function \(f\), \((f)_+\) denotes the function \(\max\, 
  (f, 0)\). Given a 2-form \(\mu\) on a 4-manifold, \(\mu^+\) denotes
  its self-dual part. 
\item \((A_r, \Psi_r)\) is used to denote a solution to the \(r\)-th
  version of the Seiberg-Witten equations; i.e. \(\grS_{\mu _r,
    \grp}(A_r, \Psi _r)=0\). However, in the context of
  a sequence \(\{(A_n, \Psi_n)\}_{n\in\bbZ^+}\), \((A_n, \Psi_n)\)
  is often used as a short hand for \((A_{r_n}, \Psi_{r_n})\), where
  \(\{r_n\}_{n\in \bbZ^+}\) is a sequence of positive numbers
  approaching \(\infty\).

\item The book \cite{T} contains four previously published
  articles. When referring to theorems or equations therein, we use
  roman numerical to label the specific article we refer to. 
  For example: Theorem I.1.2 stands for
  Theorem 1.2 of the first article in \cite{T}. 
\end{enumerate}

This article frequently refers to formulas in various literature, which
unfortunately use different conventions. For the reader's convenience,
we clarify some of their relations here. The Seiberg-Witten equations
in this article follow the convention of \cite{KM}. In Taubes' articles, \(F_A/2\) above is
replaced by \(F_A\). This results in a difference of factor 2 in many
expressions below from their analogs in Taubes' articles.
To sum up, 
\[\begin{split}\Psi & =\Psi_{KM}=\Psi_{PFH}/\sqrt{2}=\Psi_{har}/\sqrt{2};\\ 
\frac{i\mu}{4} &=irw_{f}|_{PFH}=-2w|_{KM};\\
\frac{i\mu^+}{4}& =\frac{ir\omega}{8}|_{har},
\end{split}
\]
where the first expressions in all three lines are in the notation used in
this article, and the subscripts \(KM\), \(PFH\), \(har\) refer
respectively to their counterparts in \cite{KM}, \cite{LT}, and
\cite{Ts}. 

The rest of this article is organized as follows: Sections 2 and 3
respectively contain prerequisites regarding \(Gr\) and \(SW\) sides of
the story in our context, in the sense of Taubes' \(SW=Gr\) in
\cite{T}. Section \ref{sec:synopsis} enumerates some key steps of
Taubes' arguments in the \(SW\Rightarrow Gr\) part of \cite{T}, which basically all subsequent literature in
this direction also go through, including the present article. It also
indicates where each step is carried out in the remainder of the
article. 
Section 4 is the crux of the proof, and contains a
preliminary version of the crucial ``energy bound''. The latter is
iteratively improved in Sections 5 and 6 to serve our
purposes. Sections 5 and 6 describe how
some of Taubes' original argument should be modified in our
context. The final section, Section 7, brings everything together to
complete the proofs of the main theorems. 

Lastly, a word of caution is in order: For the sake of brevity, we
often quote minute details from \cite{Ts} and \cite{T} directly
without preparation.  
The reader is thus strongly advised to keep copies of \cite{Ts} and \cite{T} at hand.

\subsubsection*{Acknowledgements.} The author's debt to Taubes is
abundantly evident.
This work is supported in part by
Hong Kong Research Grant Counsel grants GRF-401913, 14316516, 143055419, and in earlier
stages, by U.S. NSF grants.

\section{Preliminaries: the \(Gr\) side}
Some basic notions and facts about the right hand side of our variant
of \(SW\Rightarrow Gr\) are gathered in this section.

\subsection{t-orbits and t-curves}\label{sec:t-orbit}
Recall (weighted) subvarieties and related notions and conventions from Section
\ref{sec:convention}; in particular, items (11)-(14). The notion of t-orbits and t-curves are special
cases of weighted subvarieties. 

Let \(Y\) be a closed oriented riemannian 3-manifold and \(\theta\) be
a harmonic, regular Morse-Novikov 1-form on \(Y\).
\begin{defn}
 A {\em t-orbit} with respect to 
\((Y, \theta)\) is a (possibly empty) finite set \(\pmb{\gamma} =\{(\gamma_a, m_a)\}_a\) where:
\begin{itemize}
\item Each \(\gamma_a\subset Y-\theta ^{-1}(0)\) is
a finite-length flow line of the vector field dual to \(\theta\), and
\(\gamma _a\neq\gamma _b\) when \(a\neq b\). \(m_a\in \bbZ^+\). 
\item \(\partial \, (\bigcup_a \gamma_a)=\theta^{-1}(0)\) as oriented manifolds, 
\end{itemize}
\(\gamma_a\) is called a {\em constituent flow line} of
\(\pmb{\gamma}\), and \(m_a\) its {\em weight}.

A {\em multi-orbit} of \((Y, \theta)\) is a union
\(\cup_a\gamma_a\subset Y\), where \(\gamma_a\) are disjoint
finite-length flow lines of  the vector field dual to \(\theta\), and
the index set \(\{a\}\) has finite (possibly 0) elements.
\end{defn}

Let \((\scrX, \omega)\) be a (possibly noncompact) 
symplectic 4-manifold, and \(J\) be a compatible almost
complex structure. A subvariety \(C\subset \scrX\) is
said to be a {\em pseudo-holomorphic (or \(J\)-holomorphic) subvariety} if the representing
map \(u_0\co  C_0\to \scrX\)  of \(C\) satisfies then following additional
conditions: 
\(C_0\) (possibly non-compact, non-connected or empty) 
is \(C^\infty \) complex curve, and \(u_0\co  C_0\to \scrX\) is a
proper \(J\)-holomorphic map. A {\em weighted pseudo-holomorphic (or
  \(J\)-holomorphic) subvariety} in \(\scrX\) is a weighted subvariety
\(\{(C_a, m_a)\}_a\) such that \(C_a\) are irreducible \(J\)-holomorphic
subvarieties.

\label{sec:t-curve}
Let \(X\) be an oriented MCE, and \(\nu\) be an admissible form on
\(X\). Let \(J\) be the
almost complex structure on \[X_0:=X-\nu^{-1}(0)\]
determined by \(g\) and \(\omega=2\nu^+\). 



\begin{defn}\label{defn:t-curve}
A {\em t-curve} on \(X_\bullet\) is a weighted \(J\)-holomorphic subvariety 
\(\mathbf{C}=[C, \td{C}]\) 
in \(X_0\cap X_\bullet\) such that the intersection number of \(\mathbf{C}\) with each linking 2-sphere
  of \(\nu^{-1}(0)\) is one.
\end{defn}



Note that when \((X, \nu )=(\bbR\times Y, \pi _2^*(*_3\theta ))\) is
cylindrical, a translation (by \(\tau _L\)) of a
t-curve/pseudo-holomorphic subvariety gives another
t-curve/pseudo-holomorphic subvariety. Moreover, every multi-orbit \(\gamma \) on \((Y, \theta )\) corresponds
to a translation-invariant pseudo-holomorphic subvariety \(\bbR\times
\gamma \) in \((X, \nu
)\), and vice versa.  Such a translation-invariant pseudo-holomorphic
subvariety is said to be {\em constant}. The notion
of  {\em constant} t-curves in cylindrical admissible pairs is
similarly defined.   

\subsection{Topology on the space of t-curves}

Following Taubes  \cite{T99}, we define
\begin{defn} \label{def:geom-conv} 
Given two subsets \(S, S'\) in a riemannian manifold \(\scrX\), the {\em geometric distance}
between \(S, S'\) in \(\scrX\) is 
\[\dist_\scrX(S, S'):=\Big( \sup_{x\in S}\dist (x, S')+\sup_{x\in S'}\dist
    (x, S)\Big).\]
Let \((\scrX, \omega)\) be a (possibly noncompact) 
symplectic 4-manifold, and \(J\) be a compatible almost
complex structure. Let \(C_n\), \(n\in\bbZ^{\geq 0}\) be 
pseudo-holomorphic subvarieties in \(\scrX\) with
\(\int_{C_n}\omega<\infty\). The sequence \(\{C_n\}_{n\in
  \bbZ^+}\) is said to {\em converge geometrically} to \(C_0\) if the
following hold:
\begin{itemize}
\item  Regard each \(C_n\) and its irreducible components as
  2-dimensional rectifiable currents. Then \(\{C_n\}_{n\in\bbZ^+}\)
  converges weakly as currents to a current of the form \(\sum_a m_a C'_a\),
  where \(m_a\in\bbZ^+\), and \(C'_a\) are irreducible components of
  \(C_0\).
\item For any compact subset \(K\subset \scrX\),
  \(\lim_{n\to\infty}\dist_K(C_n|_K, C_0|_K)=0\).
\end{itemize}
\end{defn}

There is a version of Gromov compactness result for pseudo-holomorphic
subvarieties in 4-dimensional symplectic manifolds that hold
simultaneously for all genera:
\begin{prop}\label{thm:gromov-cpt}
(\cite{T99} Propositions 3.8) Let \((\scrX, \omega), J\) be as in
Definition \ref{def:geom-conv}. Suppose there is a countable exhaustion \(\{
\scrX_n\}_{n\in \bbZ^+}\) of \(\scrX\) by open subsets with compact closure,
\(\scrX_1\subset\cdots\scrX_n\subset \scrX_{n+1}\subset \cdots \subset \scrX\), and let \(\{(\omega_n, J_n)\}_n\)
be a corresponding sequence whose typical element \((\omega_n, J_n)\)
consists of a symplectic form \(\omega_n\) and a compatible
almost complex structure \(J_n\) on \(\scrX_n\). Suppose that \(\{(\omega_n,
J_n)\}_n\) converges to \((\omega, J)\) in the \(C^\infty\) topology. 
Let \(z_E>0\), and let \(\{C_n\}_{n\in \bbZ^+}\) be a sequence of \(J_n\)-holomorphic
subvarieties in \((\scrX_n, \omega_n)\) such that 
\(\int_{C_n}\omega<z_E\)
\(\forall n\). Then there is a \(J\)-holomorphic subvariety
\(C\subset \scrX\) with \(\int_C\omega<z_E\) and a subsequence of
\(\{C_n\}_n\) which converges geometrically to \(C\).
\end{prop}

Now let \(X\) be a 4-dimensional MCE and \(\nu\) be an admissible form
on \(X\). 
\begin{defn}
A sequence of t-curves \(\{\mathbf{C}_n\}_{n\in\bbZ^+}\) in \(X^{'a}\) is said to {\em
  (locally) converge} to a t-curve \(\mathbf{C}_0\) in \(X^{'a}\) if
the underlying subvarieties of \(\mathbf{C}_n\)
converge geometrically to the underlying subvariety of
\(\mathbf{C}_0\), and the associated currents of \(\mathbf{C}_n\)
converge to the associated current of \(\mathbf{C}_0\).
\end{defn}

The notion of {\em convergence of t-orbits} is similarly defined.  

\begin{defn}\label{def:asymp-cur}
Let \(i\in \grY_m\).  A t-curve
\(\mathbf{C}\) on \(X^{'a}\) is said to be  {\em asymptotic} on the \(Y_i\)-end to some t-orbit
\(\pmb{\gamma}_i\) with respect to the pair \((Y_i, \nu_i)\) if:
\begin{itemize}
\item[(i)] Let \(\gamma_i\) and \(C\) be respectively the underlying
  subvarieties of \(\pmb{\gamma}_i\) and \(\mathbf{C}\).
Given any sequence of positive numbers \(\{L_n\}_{n\in \bbR^+}\) with
  \(\lim_nL_n\to \infty\), 
  \(\{(\tau_{-L_n}C)|_{\hat{Y}_i}\}_n\) geometrically converges to \(
  [0, \infty)\times \gamma_i\) in \(\hat{Y}_i\). In fact, 
for all sufficiently small \(\varepsilon>0\), there exists an 
  \(R(\varepsilon)>0\), such that
  \(\dist_{\hat{Y}_i}(\tau_{-L}C,[0, \infty)\times 
  \gamma_i)<\varepsilon  \), \(\forall L>R(\varepsilon)\).
\item[(ii)] The sequence of currents \(\{\tau_{-L_n}\tilde{C}\}_n\) converges weakly to
\(\pi _2^*\tilde{\gamma}_i\), where \(\pi _2\co \hat{Y}_i=\bbR\times Y_i\to
Y_i\) is the projection. 
\end{itemize}
 The t-orbit \(\pmb{\gamma}_i\) above is said to be
the {\em limiting t-orbit} of \(\mathbf{C}\) on the
\(Y_i\)-end. In the case when \(X=\bbR\times Y\) is cylindrical, the
limiting t-orbit on \((+\infty)\)-end is called the {\em
  \(+\infty\)-limit} of \(\mathbf{C}\), and the limiting t-orbit on
the \((-\infty)\)-end, as a t-orbit in \(Y\), is called the {\em
  \(-\infty \)-limit} of \(\mathbf{C}\). If \(\bfC\) has a limiting
t-orbit on each of its Morse ends, \(\bfC\) is said to be {\em
  asymptotically constant}. 

Replace the t-curve \(\bfC\) and the limiting t-orbits \(\pmb{\gamma
}_i\) above respectively by a pseudo-holomorphic subvariety \(P\), and
\(\pp_i\), where \(\pp_i\) is a union of flow lines of the vector
field dual to \(*\nu _i\). The notion of \(P\) being {\em asymptotic to
\(\pp_i\) on the \(Y_i\)}, is similarly defined.  
\end{defn}
\subsection{\(\Spin^c\) structures}\label{sec:gr-spin}

\subsubsection*{(a) \(\Spin^c\)-structure of t-orbits} A t-orbit \(\pmb{\gamma}=\{(\gamma_i, m_i)\}_i\) is assigned a \(\Spin^c\)
structure as follows. Let  \(B_\theta\) be a small tubular neighborhood
of \(\theta^{-1}(0)\subset Y\). Let \[
[\pmb{\gamma}]:=[\tilde{\gamma}]=\sum_im_i[\gamma_i]\in H_1
(Y, \theta^{-1}(0);\mathbb{Z})\simeq H_1(Y, B_\theta;\mathbb{Z})\] 
and let \(E\) be a complex line bundle on \(Y-B_\theta\) with
\(c_1(E)=p.d.[\pmb{\gamma}]\). Let \(K^{-1}\) denote the subbundle
\(\ker\theta\) of \(TY|_{Y-B_\theta}\). 
This is an oriented plane bundle; choose a
complex structure on \(K^{-1}\) compatible with the orientation. 
Then
\(\mathbb{S}=E\oplus E\otimes  K^{-1}\) is trivial over \(\partial
(Y-B_\theta)\). Extend \(\mathbb{S}\) trivially over \(B_\theta\) to
get a complex rank 2 vector bundle on \(Y\), also denoted by
\(\mathbb{S}\). Then the {\em  \(\Spin^c\)-structure of \(\pmb{\gamma}\)},
denoted \(\mathfrak{s}(\pmb{\gamma})\), is that with \(\mathbb{S}\) as its
spinor bundle, and with \(E\), \(E\otimes K^{-1}\) being eigen-bundles
under the Clifford  action of \(*\theta/|\theta|\) corresponding to
eigenvalues \(-i\) and \(i\) respectively. 

Let \(S_{Y, *\theta}\subset
H_1(Y, \theta^{-1}(0);\mathbb{Z})\) denote the subset that maps to the
fundamental class of \(\theta^{-1}(0)\) under the connecting map
\(\delta\co H_1(Y, \theta^{-1}(0);\mathbb{Z})\to H_0\, (\theta^{-1}(0);
\mathbb{Z})\) in the relative exact sequence of the pair \((Y,
\theta^{-1}(0))\). Note that the relative homology class \([\pmb{\gamma}]\in H_1
(Y, \theta^{-1}(0);\mathbb{Z})\) of any t-orbit falls in
\(S_{Y, *\theta}\), and in fact the
procedure described above defines an isomorphism from \(S_{Y, *\theta}\) to
the space of \(\Spin^c\) structures on \(Y\), 
\[
\varsigma_{(Y,*\theta)}\co S_{Y, *\theta}\stackrel{\sim}{\to }\Spin^c(Y), 
\] 
as affines spaces over \(H_1(Y;\mathbb{Z})\simeq H^2(Y;
\mathbb{Z})\). 

Note that the assumption that \(\theta\) is regular (cf. Definition
\ref{def:reg_form}) implies that for any fixed \(\grs\in \Spin (Y)\), there
are finitely many t-orbits (with respect to \((Y, \theta )\)), \(\pmb{\gamma}\), with
\(\grs(\pmb{\gamma})=\pmb{\gamma}\). Denote this set by \(\bbP (Y,
\theta ; \grs)\). (\(Y\), \(\theta \), or both are sometimes omitted
from the notation when they are clear from the context. For example, \(\bbP (Y,
\theta ; \grs)=\bbP (Y,\grs)=\bbP (\grs)\).)

\subsubsection*{(b) \(\Spin^c\)-structure of t-curves}

Let \(H_{*, closed}\) denote the Borel-Moore homology\footnote{The
  Borel-Moore homology is also
  called the 
``closed homology'' or ``homology with closed support''.  See e.g. 
\cite{BT}~Remark 5.17, the Wiki article  \\ {\tt
  https://en.m.wikipedia.org/wiki/Borel-Moore\_homology\#CITEREF1960}
and references therein. It is dual to the cohomology of
compact support.}.  It is useful to note that in our
setting, 
\[\begin{split}
 H_{*,closed}\, (X;\bbK) 
 & \simeq H_{*, closed}\, (X^{'a}, \partial X^{'a};\bbK)\\
 & \simeq H_*(X_c, \partial X_c;\bbK);\\
H_{*,  closed}\, \big(X^{'a}, \partial X^{'a}\cup Z^{'a}_0;\bbK\big)
& \simeq H_{*,
  closed}\, \big(\ov{X^{'a}}, \partial \ov{X^{'a}}\cup Z_0^{'a}
;\bbK)\\
& \simeq H_*(\ov{X^{'a}_\delta}
, \partial \ov{X^{'a}_\delta };\bbK\big),\\
\end{split}
\]
where \(X^{'a}_\delta \), \(Z^{'a}_\delta \) are  as defined in Section
\ref{sec:convention} (9). 

A t-curve \(\mathbf{C}=\{(C_k, m_k)\}_k\) in \(X^{'a}\) determines a relative homology
class through its associated current
\[\begin{split}
[\mathbf{C}]& =[\tilde{C}]\in  H_{2, closed}\, (X^{'a}, \partial X^{'a}\cup Z^{'a}_\delta;
\mathbb{Z})\\
& \quad \simeq H_{2, closed}\, (X^{'a}, \partial X^{'a}\cup Z^{'a}_0; \mathbb{Z})\simeq
H_2(\ov{X^{'a}_\delta }, \partial \ov{X^{'a}_\delta }
;\mathbb{Z}).
\end{split}
\] 
The last condition in Definition \ref{defn:t-curve} implies that \([\mathbf{C}]\)
maps to the fundamental class \([Z_0^{'a}]\in H_{1,
  closed}\, (Z_0^{'a}, \partial Z_0^{'a};\mathbb{Z}^{'a})\) under the connecting map
\(\delta \) in the relative exact sequence
\[\begin{split}
\cdots  0\to  & H_{2, closed}\, (X^{'a}, \partial X^{'a};\mathbb{Z})\stackrel{p}{\to }H_{2, closed}\,
(X^{'a}, \partial X^{'a}\cup Z^{'a}_0;\mathbb{Z})\\
& \quad \quad \stackrel{\delta }{\to }H_{1,
  closed}\, \big (Z^{'a}_0, \partial  Z^{'a}_0;\mathbb{Z}^{'a}\big)\cdots. 
\end{split}
\]
Let  
\[
S_{X,  \nu}:=\delta ^{-1}[Z_0^{'a}]\subset H_{2, closed}\, (X^{'a}, \partial X^{'a}\cup Z^{'a}_0; \mathbb{Z}).\]

In parallel to the case with t-orbits, there is an isomorphism 
\[
\varsigma_{X, \nu}\co S_{X, \nu}\stackrel{\sim}{\to}\Spin^c(X)\simeq \Spin^c(X^{'a})
\] 
as affine spaces
  over the abelian group \[
\ker \delta =\im p\simeq H_{2, closed}\, (X^{'a}, \partial X^{'a}; \mathbb{Z})\simeq H^2(X;
  \mathbb{Z}).
\]
(In this article, we do not
 distinguish between the set
 of \(\Spin^c\)-structures and \(X\) and the set of
 \(\Spin^c\)-structures on \(X^{'a}\), as they are canonically isomorphic.) 
The inverse of \(\varsigma_{X, \nu}\) can be described more concretely
as follow: Given 
a \(\Spin^c\) structure \(\mathfrak{s}\) on \(X\), 
take
\(\bbS_X=\mathbb{S}=\mathbb{S}^+\oplus \mathbb{S}^-\) to be the associated spinor
bundle. Over \(X^{'a}_0\), the Clifford action of \(\nu^+/|\nu^+|\) splits \(\mathbb{S}^+|_{X_0}\) as a
direct sum \(E\oplus E\otimes   K^{-1}\) of eigenbundles, where \(E\)
is a complex line bundle 
corresponds to eigenvalue \(-i\). Then  \(\varsigma_{X,
  \nu}^{-1}(\mathfrak{s})\in H_{2, closed}\, (X^{'a}, \partial X^{'a}\cup Z_0;
\mathbb{Z})\) is the Poincar\'{e} dual of \(c_1(E)\in H^2(X^{'a}-\nu ^{-1}(0);\bbZ)\).

Notice that the definition of \(\Spin^c\) structures of t-orbits
and t-curves only depend on their associated currents. So are the
various notions of energy defined in the upcoming subsection.
Also, when \((X, \nu )\) is cylindrical, two t-curves related by
translation have the same \(\Spin^c\) structure.

\subsection{Energy of t-curves and chains of t-curves}\label{sec:energy}

Let \(X\) be an MCE and \(\nu\) an admissible form on \(X\). 
\begin{defn}
Given  an oriented 2-dimensional subvariety  \(C\) in \(X^{'a}\) and
\(X_\bullet\subset X^{'a}\), the {\em \(\nu\)-energy} of \(C|_{X_\bullet}\) is:\[
\scrF_\nu(C|_{X_\bullet}):=\int_{C|_{X_\bullet}}\nu.
\] 
When \(X_\bullet\) is compact, the following is called the
\(*\nu\)-energy of \(C|_{X_\bullet}\):
\[
\scrF^*_\nu(C|_{X_\bullet}):=\int_{C|_{X_\bullet}}*\nu. 
\]
Finally, setting \(\omega=2\nu^+\), the {\em \(\omega\)-energy} of \(C|_{X_\bullet}\) is: 
\[
\scrF_\omega(C|_{X_\bullet}):=\scrF_\nu
(C|_{X_\bullet})+\sup_{X_*\in \grX_1}\scrF^*_\nu (C|_{X_*}),
\]
where \(\grX_1=\{ X_*\, |\, X_*\subset X_\bullet,
|X_*|\leq 1\}\). (Cf. \textsection \ref{sec:convention} for the
definition of \(|X_\bullet|\)).
\end{defn}
When \(\omega^{-1}(0)=\nu^{-1}(0)=\emptyset\), the finite-energy condition above
agrees with that in  \cite{BEHWZ}. It also straightforwardly from the
definitions that in the cylindrical case, two t-curves that are
related by translation have the same \(\nu \)-energy, \(*\nu
\)-energy, and \(\omega\)-energy. 
\begin{rem}
Just as in \cite{BEHWZ}, a bound on \(\scrF_\omega\) is the
prerequisite for local compactness results, while a bound on
\(\scrF_\nu\) is used to establish global compactness (weak convergence to
broken-trajectories or ``buildings''). In the context of \cite{Ts}
where \(X\) is closed, \(\scrF_\omega\) is simply the LHS of its
Equation (7.10), which also give a bound on this in terms of the
\(\Spin^c\) structure of \(C\), the cohomology class \([\omega]\), and
the metric on \(X\).
\end{rem}

The following useful lemma is an immediate consequence of the
compatibility condition of \(J\).
\begin{lemma}
\label{lem:F-curve}
Under the assumptions of the preceding definition, if \(C\) is a
pseudo-holomorphic subvariety in \(X_\bullet\), then \(\scrF_\omega(C), \scrF_\nu(C),
\scrF_\nu^*(C)\geq 0\). Moreover, when \((X_\bullet =I\times Y,
\nu =\pi _2^*\ud{\nu })\) is cylindrical, \(\scrF_\nu(C)=0\)
iff  \(C\) is constant. 
Namely, \(C=I\times \gamma\), where 
\(\gamma\) is a union of flow lines
of the vector dual to \(*\ud{\nu}\) in \(Y\).
\end{lemma}

Given a t-curve \(\mathbf{C}=\{(C_k, m_k)\}_k\) in \(X_\bullet\), its {\em
  \(\nu\)-energy}, {\em \(*\nu\)-energy}, {\em \(\omega\)-energy} are
respectively defined to be \[\begin{split}
\scrF_\nu(\mathbf{C}):=\sum_k m_k \scrF_\nu (C_k|_{X_\bullet}),\\
\scrF^*_\nu(\mathbf{C}):=\sum_k m_k \scrF^*_\nu (C_k|_{X_\bullet}),\\
\scrF_\omega(\mathbf{C}):=\sum_k m_k \scrF_\omega(C_k|_{X_\bullet}).
\end{split}
\]

We will see in Proposition \ref{prop:t-curve-asymp} below that by the
preceding lemma, any t-curve
\(\mathbf{C}\) with finite
\(\omega\)-energy in \(X^{'a}\) is
asymptotic to a t-orbit \(\pmb{\gamma}_i\) on each of its Morse
end \(\hat{Y}_i\). 
It also follows from the preceding lemma and Proposition \ref{thm:gromov-cpt}  that 
\begin{prop}
Let \(X\) be an MCE and \(\nu \) be an admissible form on \(X\). Let
\(\{\mathbf{C}_n\}_{n\in \bbZ^+}\) be
a sequence of t-curves on \(X^{'a}\), 
with a uniform upper bound on
their \(\omega\)-energy. Then there is a t-curve
\(\mathbf{C}_0\) and a subsequence of \(\{\mathbf{C}_n\}_{n\in
  \bbZ^+}\) that locally converges to \(\mathbf{C}_0\). 
\end{prop}


\begin{defn}\label{def:chain}
Let \((X, \nu )\) be an admissible pair. 
\begin{itemize}
\item[(i)] Suppose the admissible pair \((X=\bbR\times Y, \nu)\) is
  cylindrical. (And thus \(X\) has no vanishing ends and \(X^{'a}=X''=X\).) Let \(\grs_Y\) denote the restriction of the \(\Spin^c\)
  structure \(\grs\) to \(\{s\}\times Y\), for any \(s\in \bbR\). For
  \(\grK\in \bbZ^{\geq 0}\), a (\(\grK\)-component)  
  {\em chain of t-curves} \(\grC\) on \(X^{'a}\) consists of:
  \begin{itemize}
\item   a sequence of \(\grK+1\) t-orbits with \(\Spin^c\)-structure
  \(\grs_Y\), \\ \(\{\pmb{\gamma }_0, \pmb{\gamma }_1,
  \ldots, \pmb{\gamma }_{\grK}\}\), which we call the {\em rest
    orbits} of \(\grC\); and 
\item a sequence (possibly empty) 
  of non-constant, asymptotically constant
  t-curves
  \(\mathbf{C}_k\) on \(X=\bbR\times Y\) with  \(\Spin^c\)-structure \(\grs\), with \(-\infty\)- and
  \(+\infty\)-limits respectively \(\pmb{\gamma}_{k-1}\) and 
  \(\pmb{\gamma}_{k}\). 
\end{itemize}
We say that \(\mathfrak{C}\) has \(\Spin^c\)-structure \(\grs\). The
{\em (\(-\infty\))-limit of \(\mathfrak{C}\)} is \(\pmb{\gamma}_{0}\);
the {\em \(+\infty\)-limit of \(\mathfrak{C}\)} is
\(\pmb{\gamma}_{\grK}\). The {\em \(k\)-th component} of \(\grC\) is
\(\bfC_k\). We typically write \(\grC\) as an ordered set, \(\grC=\{\bfC_k\}_{k=1}^\grK\), when
\(\grK>0\). Two chains of t-curves, \(\grC\), \(\grC'\),  
are said to be {\em equivalent} and
denoted \(\grC\sim_\tau\grC'\), if they have the same number of components
and the same rest orbits, and when \(\grK>0\), each pair of corresponding components
\(\bfC_k\), \(\bfC'_k\) are related by
translation. I.e. \(\bfC'_k=\tau _L\bfC_k\) for a certain \(L\in \bbR\). 
\item[(ii)] Suppose the admissible pair \((X, \nu)\) is not
  cylindrical. Denote the \(\Spin^c\) structure on \(X\) by \(\grs\),
  and the restriction \(\grs|_{\mathfrc{s}_i^{-1}(0)}=:\grs_i\) for
  each \(i\in \mathfrak{Y}_X\). A
  {\em chain of t-curves} on \(X^{'a}\) is a pair \(\mathfrak{C}=\{\mathbf{C}_0,
  \{\mathfrak{C}_i\}_{i\in \mathfrak{Y}_m}\}\) such that:
\begin{itemize}
\item[(1)] \(\mathbf{C}_0\) is an asymptotically constant t-curve with \(\Spin^c\)-structure
  \(\grs\) on \(X^{'a}\). For each \(i\in \grY_m\), let \(\pmb{\gamma}_{0,i}\) denote the limiting
  t-orbit of \(\mathbf{C}_0\) on the \(Y_i\)-end.
\item[(2)] Let \((Y_i, \nu_i)\) denote the ending pair of the end
  \(\hat{Y}_i\subset X\). For each \(i\in \grY_m\), \(\mathfrak{C}_i\) is a chain of 
  t-curves of \(\Spin^c\)-structure \(\hat{\grs}_i\) on the cylindrical admissible pair
  \((\mathbb{R}\times Y_i, \hat{\nu}_i)\), where \(\hat{\grs}_i\),
  \(\hat{\nu}_i\) are respectively the pull-back of \(\grs_i\) and
  \(\nu_i\) under the projection \(\bbR\times Y_i\to Y_i\).
\item[(3)] For each \(i\in \grY_m\), the \(-\infty\)-limit of \(\mathfrak{C}_i\) is \(\pmb{\gamma}_{0,i}\).
\end{itemize}
We say that the chain \(\mathfrak{C}\) has \(\Spin^c\)-structure
\(\grs\). The rest orbits of \(\grC_i\) are said to be the rest orbits of
\(\grC\) (in the \(\hat{Y}_i\)-end). Two chains of t-curves,
\(\grC=\{\bfC_0, \{\grC_i\}_{i\in \grY_m}\}\), \(\grC'=\{\bfC'_0, \{\grC'_i\}_{i\in \grY_m}\}\), are said to be {\em equivalent} and
denoted \(\grC\sim_\tau\grC'\), if \(\bfC_0=\bfC'_0\) and \(\grC_i\sim_\tau\grC'_i\) \(\forall i\in \grY_m\). 
\end{itemize}
\end{defn}
An equivalence class of chains of t-curves  is an analog of what is called an ``(unparametrized) broken
trajectory'' or ``broken \(X\)-trajectory''  in \cite{KM}'s Definitions 16.1.2 and 24.6.1.

The \(\nu\)-energy and \(\omega\)-energy of chains of 
t-curves on \((X^{'a}, \nu)\) are defined as follows: When \((X, \nu)\) is
cylindrical and \(\grC=\{\mathbf{C}_k\}_k\) has at least one component, 
\[
\scrF_\nu (\grC):=\sum_k\scrF_\nu(\mathbf{C}_k); \quad \scrF_\omega(\grC):=\sum_k\scrF_\omega(\mathbf{C}_k).
\]
If \(\grC\) has 0 components and a single rest orbit \(\pmb{\gamma}\), set \(\scrF_\nu
(\grC):=\scrF_\nu (\bbR\times\pmb{\gamma })=0\) and \(\scrF_\omega(\grC):=\scrF_\omega(\bbR\times\pmb{\gamma })\).

When \((X,\nu)\) is non-cylindrical and 
 \(\mathfrak{C}=(\mathbf{C}_0,
  \{\mathfrak{C}_i\}_{i\in \mathfrak{Y}_m})\), 
\[
\scrF_\nu (\grC):=\sum_k\scrF_\nu(\mathbf{C}_0)+\sum_i\scrF_\nu (\grC_i); \quad \scrF_\omega(\grC):=\sum_k\scrF_\omega(\mathbf{C}_0)+\sum_i\scrF_\omega(\grC_i).
\]
Note that two equivalent chains of t-curves have the same \(\nu
\)-energy and \(\omega\)-energy.

\begin{defn}\label{def:w-t-conv}
Let \(\Gamma =\{r_n\}_{n\in \bbZ^+}\subset [1, \infty)\) be a sequence of strictly
increasing numbers such that \(r_n\to \infty\) as \(n\to \infty\), and
\(\{(A_r, \Psi _r)\}_r\in \Gamma \) be a corresponding sequence of
admissible solutions to the Seiberg-Witten equation \(\grS_{\mu _r,
  \hat{\grp}_r}(A_r, \Psi _r)=0\), where \(\mu _r\), \(\hat{\grp}_r\)
are as in the statement of Theorem \ref{thm:l-conv}. We say that the
sequence \(\{(A_r, \Psi _r)\}_{r\in \Gamma }\) {\em weakly t-converges}
to a chain of t-curves \(\grC\) if the following holds: 
\begin{itemize}
\item[(i)] When \((X, \mu_r, \hat{\grp}_r)\) is cylindrical: If
  \(\grC\) has 0 components, namely, it has a single rest orbit
  \(\pmb{\gamma }\), then \(\{(A_r, \Psi _r)\}_{r\in \Gamma }\)  t-converges
to the constant t-curve \(\bbR\times \pmb{\gamma }\). If \(\grC\) has
at least one component, write \(\mathfrak{C}=
\{\mathbf{C}_k\}_{k=1}^\grK\).  There is 
a sequence \(\{{\bf L}_r=(L_{1, r}, \cdots, L_{\grK,r})\}_{r\in\Gamma}\) in
\(\bbR^\grK\), with \(\lim_{r\to\infty}(L_{k_1, r}-L_{k_2, r})=\infty\)
for any fixed pair \(k_1, k_2\) with \(k_1>k_2\), such that  
\(\{\tau_{-L_{k, r}}(A_r, \Psi_r)\}_{r\in\Gamma}\) t-converges to
\(\mathbf{C}_k\) for each \(k\). (Note that the admissibility
assumption on \(\nu \) implies that in this case \(X\) has no
vanishing ends, and \(X^{'a}=X\).)

\item[(ii)] When \((X, \mu_r, \hat{\grp}_r)\) is not cylindrical, write
\(\mathfrak{C}=(\mathbf{C}_0, \{\mathfrak{C}_{i}\}_{i\in
  \mathfrak{Y}_m}),
\mathfrak{C}_i=\{\mathbf{C}_{i,k}\}_{k=1}^{\grK_i}\) in the notation
of Definition \ref{def:chain} (ii). 
There exists a sequence \(\{{\bf L}_r^i=(L^i_{1, r}, \cdots, L^i_{\grK_i,
  r})\}_{r\in \Gamma}\) in
\((\bbR^{+})^{\grK_i}\) for each \(i\in \mathfrak{Y}_m\), 
with \(\lim_{r\to\infty}(L^i_{k_1, r}-L^i_{k_2, r})=\infty\) for any fixed pair \(k_1, k_2\) with \(k_1>k_2\),
such that 
\(\{\tau_{-L^i_{k, r}}(A_r, \Psi_r)\Big|_{\hat{Y}_i}\}_{r\in \Gamma}\) t-converges to
\(\mathbf{C}_{i, k}\), and \(\{(A_r, \Psi_r)\}_{r\in \Gamma}\) t-converges to
\(\mathbf{C}_0\).
\end{itemize}
\end{defn}
Note that if a sequence \(\{(A_r, \Psi _r)\}_r\) weakly t-converges to
a chain of t-curves \(\grC\), and \(\grC\sim_\tau\grC'\), then \(\{(A_r,
\Psi _r)\}_r\) also weakly t-converges to \(\grC'\). It therefore
makes sense to talk about weakly t-convergence to equivalence classes
of chains of t-curves.

\subsection{Relative homology classes of (chains of) t-curves}\label{sec:h-class}

Let \((X, \nu )\) be an admissible pair as before. 

Given two asymptotically constant t-curves \(\mathbf{C}, \mathbf{C}'\) in \(X^{'a}\) with the
same \(\Spin^c\)-structure and limiting t-orbits on all its Morse ends,
\([\mathbf{C'}-\mathbf{C}]:=[\tilde{C}'-\tilde{C}]\) defines a class in
\(H_2(X^{'a}, \partial X^{'a};\bbZ)\). The
two t-curves \(\mathbf{C}'\), \(\mathbf{C}\) are said to be of
the same  {\em relative homology class} if \([\mathbf{C'}-\mathbf{C}]=0\).
The precise value of \(a\) does not matter, as all \(H_2(X^{'a},
\partial X^{'a};\bbZ)\) are mutually canonically isomorphic, and the 
relative homology classes of two \(\bfC, \bfC'\) in \(X^{'a}\) and in
\(X^{'b}\), \(b\leq a\) are identified by such a canonical
isomorphism.  
We now describe the set of relative homology classes for (chains of) 
t-curves in more detail. 
Let \(\{(Y_i, \nu_i)\}_i\) be the set of ending pairs of \((X, \nu)\), and
let \(\grs\) denote the \(\Spin^c\)-structure on \(X\). Let 
\(\grs_i:=\grs|_{\mathfrc{s}_i^{-1}(1)}\). 
Consider the following
commutative diagram of relative exact sequences, and recall from
Section \ref{sec:gr-spin} (a) and (b) the
definitions of the maps and spaces \(\varsigma _{X, \nu }\co S_{X, \nu}\to \Spin^c(X)\), \(\varsigma _{Y_i, \nu _i}\co S_{Y_i, \nu
  _i}\to \Spin^c(Y_i)\) respectively via the second and the third
columns of this diagram:

{\footnotesize
\begin{equation}\label{eq:CD}
\minCDarrowwidth3pt
\begin{CD}
@ VVV @ VVV @ VVV \\
\cdots H_{2}(X^{'a}, \partial X^{'a}) @> j_\bbK>> H_{2,cl}\, (X^{'a}, \partial X^{'a})
@>\partial >> \bigoplus \limits_{i\in \grY_m} H_{1}(Y_i)\cdots\\
@ VVV @ VpVV @ VVV \\
\cdots H_2\, (X^{'a}, \partial X^{'a}\cup Z^{'a}_0)@>>> H_{2, cl}\, (X^{'a}, \partial X^{'a}\cup Z^{'a}_0) @>\partial >>\bigoplus \limits_{i\in \grY_m} 
H_{1}(Y_i, \nu_i^{-1}(0))\cdots\\ 
@ VVV @ V\delta  VV @ V VV \\
\cdots H_1(Z^{'a}_0, \partial Z^{'a}_0) @>>> H_{1, cl}\, (Z^{'a}_0, \partial
Z^{'a}_0)
@>>> \bigoplus \limits_{i\in \grY_m} H_0\, (\nu_i^{-1}(0))\cdots\\
@ VVV @ VVV @ VVV \\
\end{CD}
\end{equation}}
where the homologies are all with coefficient \(\bbK\), and \(H_{*,
  cl}:=H_{*, closed}\) denote the Borel-Moore homology.

Let \(\pmb{\gamma}_i\) be
respectively t-orbits of \((Y_i, \nu_i)\) with
\(\Spin^c\)-structure \(\grs_i\). The relative homology class of a
t-curve \(\mathbf{C}\) with \(\Spin^c\)-structure \(\grs\) and
with limiting t-orbits \(\pmb{\gamma}_i\) takes value in the
space  \(\scrH \big((X^{'a}, \nu, \grs),\{\tilde{\gamma}_i\}_i\big)\) defined below.

For each Morse ending pair \((Y_i, \nu_i)\) \(i\in \grY_m\), define the space \(\scrZ_{(Y_i,
  \nu_i, \grs_i)}\) to be the space of all 1-dimensional integral currents
\(\tilde{\gamma}_i\)  on
  \(Y_i\) that represents the class \(\varsigma_{(Y_i,
    \nu_i)}^{-1}(\grs_i)\in S_{Y_i, \nu
    _i}\subset H_1(Y_i, \nu_i^{-1}(0);\bbZ)\). 

Given \(\grs\in \Spin^c(X)\) and \(\tilde{\gamma }_i\in \scrZ_{(Y_i,
  \nu_i, \grs_i)}\) for each \(i\in \grY_m\),  define 
\[
\scrH \big((X^{'a}, \nu, \grs),\{\tilde{\gamma}_i\}_{i\in
  \grY_m}\big)=\{\tilde{C}\}/\sim,
\]
where \(\tilde{C}\) is a 2-dimensional integral current on
\(X^{'a}\) satisfying the following conditions. 
\begin{itemize}
\item For each \(i\in \grY_m\),  \(\tilde{C}\) is asymptotic to \(\tilde{\gamma}_i\) on the
\(\hat{Y}_i\)-end in the sense of item (ii) in Definition
\ref{def:asymp-cur}, and the limit \(\lim_{l\to
  \infty}\int_{\td{C}|_{\hY_{i,l}}} \pi _2^*h<\infty\) exists for any
\(h\in \Omega ^2(Y_i)\). 
\item \(\tilde{C}\) represents the element in \(\varsigma _{X, \nu
  }^{-1}(\grs )\in  S_{X, \nu }\subset H_{2, closed}\, (X^{'a},
 \partial X^{'a}\cup Z_0^{'a};\bbZ)\). 
\end{itemize}
Meanwhile, we write \(\tilde{C}\sim \tilde{C}'\) for two
such currents \(\tilde{C}\), \(\tilde{C}'\) if
\(\tilde{C}-\tilde{C}'\)  is the boundary of an integral current that is
asymptotic to \(0\) on all Morse ends of \(X\). 

It follows from the preceding definition and Equation (\ref{eq:CD})
that the set \(\scrH \big((X^{'a},
\nu, \grs),\{\tilde{\gamma}_i\}_{i\in \grY_m}\big)\) is an affine space over the abelian
group \[
\scrH_X:=\ker j_\bbZ\subset H_2(X^{'a}, \partial X^{'a}; \bbZ).\] 
An element in \(\scrH
\big((X^{'a},\nu, \grs),\{\tilde{\gamma}_i\}_i\big)\) is said to be a (\(\bbZ\)-) {\em relative
homology class in \((X^{'a}, \nu, \grs)\) (rel
\(\{\tilde{\gamma}_i\}_i\))}. Later on, we sometimes omit \(\grs\)
from the notation when the \(\Spin^c\) structure is fixed. 

\begin{rem}
Assertions (c) and (d) of Theorem \ref{thm:g-conv} postulates that \(b^1(Y_i)=0\)
for all vanishing ends \(i\in \grY_v\).
(In particular, this holds if \((X, \nu )\)
has no vanishing ends; especially, when \((X, \nu )\) is cylindrical.) This
constraint  implies that 
\(H_2(\partial X^{'a};\bbZ)=\bigoplus_{i\in
  \grY_v}H_2(Y_i;\bbZ)=0\), and hence 
\(H_2(X^{'a}, \partial X^{'a};\bbZ)\simeq H_2(X^{'a};\bbZ)\simeq H_2(X;\bbZ)\). 
It also implies that the \(\nu\)-energy of a chain \(\grC\)  depends
only on its
relative homology class and its limiting t-orbits. 
\end{rem}

When \((X^{'a}, \nu)=(X, \nu)=(\bbR\times Y_i, \pi_2^*\nu_i)\) is cylindrical,
\(\scrH_X=H_2(X; \bbZ)=H_2(Y_i;\bbZ)\). We use \(\scrH(Y_i, \nu
_i,\grs_i; \tilde{\gamma }_i, \tilde{\gamma }_i')\) to denote \(\scrH((\bbR\times Y_i, \pi _2^*\nu _i,
\hat{\grs}_i), \{\tilde{\gamma }_i, \tilde{\gamma }_i'\})\), where
\(\tilde{\gamma }_i, \tilde{\gamma }_i'\in \scrZ_{Y_i, \nu _i,
  \grs_i}\) are respectively the limiting currents at the \(-\infty\)-
and the \(+\infty\)-ends of \(X\). 
By concatenation, there is a composition map that sends every pair
\((a, b)\in \scrH(Y_i, \nu
_i,\grs_i; \tilde{\gamma }_i^-, \tilde{\gamma }_i)\times \scrH(Y_i, \nu
_i,\grs_i; \tilde{\gamma }_i, \tilde{\gamma }_i^+)\) to an element
\(b*a\in \scrH(Y_i, \nu
_i,\grs_i; \tilde{\gamma }_i^-, \tilde{\gamma }_i^+)\). Also, every element \(\{h_i\}_i\) in
\(\prod_{i\in\grY_m}\scrH(Y_i, \nu _i, \grs_i; \tilde{\gamma }_i,
\tilde{\gamma }_i')\) defines an isomorphism 
\[
\op{c}_{\{h_i\}_i}\co \scrH \big((X^{'a},
\nu, \grs),\{\tilde{\gamma}_i\}_{i\in \grY_m})\to \scrH \big((X^{'a}, \nu,
\grs),\{\tilde{\gamma}'_i\}_{i\in \grY_m}).
\]

Let \(\mathfrak{C}=\{\mathbf{C}_1, \ldots,
\mathbf{C}_{\grK}\}\) be a chain of t-curves on \((\bbR\times Y_i, \pi
_2^*\nu _i)\) with \(\grK>0\). 
The {\em relative homology class of \(\grC\)} is
the element
\[
[\grC]:=
[\tilde{C}_\grK]* \cdots *[\tilde{C}_2]*[\tilde{C}_1]. 
\]
When \(\grC\) is a 0-component chain of t-curves with a single rest
orbit \(\pmb{\gamma }\), we set the relative homology class of
\(\grC\) to be 
the element  \([\bbR\times \pmb{\gamma }]\in \scrH(Y_i, \nu
_i,\grs_i; \tilde{\gamma }, \tilde{\gamma })\), which maps to the
element  \(0\in H_2(Y_i;\bbZ)\) under the canonical isomorphism \(\scrH(Y_i, \nu
_i,\grs_i; \tilde{\gamma }, \tilde{\gamma })\simeq H_2(Y_i;\bbZ)\). 

Suppose the admissible pair \((X, \nu)\) is not cylindrical. 
Given a  chain of
t-curves \(\mathfrak{C}=(\mathbf{C}_0,
\{\mathfrak{C}_i\}_{i\in\grY_m})\) on \(X^{'a}\), 
the {\em relative homology class} of \(\grC\) is defined to be 
\[
[\grC]=\op{c}_{\{[\grC_i]\}_i} ([\bfC_0]). 
\]
It follows directly from the definition that equivalent chains of
t-curves have the same relative homology class. 

What follows are some observations for future reference. When the currents \(\tilde{\gamma }_i\), \(\tilde{C}\) above are not
required to be integral, the same strategy can be employed to define
the real-variants of the spaces  \(\scrZ_{(Y_i,
  \nu_i, \grs_i)}\), \(\scrH \big( (X^{'a}, \nu,
\grs),\{\tilde{\gamma}_i\}_i\big)\) and other related notions. These are respectively denoted by  \(\scrZ^\bbR_{(Y_i,
  \nu_i, \grs_i)}\), \(\scrH ^\bbR\big((X^{'a}, \nu,
\grs),\{\tilde{\gamma}_i\}_i\big)\), etc. The set of real relative
homology classes, \(\scrH ^\bbR\big((X^{'a}, \nu,
\grs),\{\tilde{\gamma}_i\}_i\big)\), is now an affine space under
\(\scrH^\bbR_X:=\ker j_\bbR\subset H_2(X^{'a}, \partial X^{'a}; \bbR)\). Note that \(\ker
j_\bbZ=\scrH_X\) is torsion-free; so we may and will often identify it with the
integral lattice in \(\ker j_\bbR=\scrH_X^\bbR\). 
Similarly, When \(\{\tilde{\gamma}_i\}_i \) consists of integral
currents, \(\scrH  \big((X^{'a}, \nu,
\grs),\{\tilde{\gamma}_i\}_i\big)\) embeds in \(\scrH ^\bbR\big((X^{'a}, \nu,
\grs),\{\tilde{\gamma}_i\}_i\big)\) as an orbit under the \(\scrH_X\subset
\scrH_X^\bbR\)-action. \(\bbR\)-coefficient versions of  composition
maps \(*\) and \(\op{c}_{\{h_i\}_i}\) are similarly defined and have 
  similar properties. 
\begin{lemma}\label{gr:hR-isom}
 Fix \(i\in \grY_m\). For any \(\tilde{\gamma }_i, \tilde{\gamma }_i'\in \scrZ^\bbR_{(Y_i,
  \nu_i, \grs_i)}\), the metric on \(Y_i\) determines an isomorphism
\(I_\scrH\) from \(\scrH^\bbR(Y_i, \nu
_i,\grs_i; \tilde{\gamma }_i, \tilde{\gamma }_i')\) to
\(H_2(Y_i;\bbR)=H_2(\bbR\times Y_i;\bbR)=\scrH_{\bbR\times
  Y}^\bbR\) as affine spaces under \(H_2(Y_i;\bbR)\). In particular,
all  the affine spaces  \(\scrH^\bbR(Y_i, \nu
_i,\grs_i; \tilde{\gamma }_i, \tilde{\gamma }_i')\) are identified 
via these isomorphisms. By concatenation (the real versions of
\(\op{c}_{\{h_i\}_i}\) above, which we shall denote by the same
notation), for any pair \(\{\tilde{\gamma}_i\}_i,
,\{\tilde{\gamma}'_i\}_i\), the metric on \(X\) determines an
isomorphism between the corresponding affine spaces \(\scrH ^\bbR\big((X^{'a}, \nu,
\grs),\{\tilde{\gamma}_i\}_i\big)\), \(\scrH ^\bbR\big((X^{'a}, \nu,
\grs),\{\tilde{\gamma}'_i\}_i\big)\)  as
well. Moreover, given  \(\tilde{\gamma }_i, \tilde{\gamma }_i',
\td{\gamma }_i'' \in \scrZ^\bbR_{(Y_i,
  \nu_i, \grs_i)}\) and \([\td{C}]\in \scrH^\bbR(Y_i, \nu
_i,\grs_i; \tilde{\gamma }_i, \tilde{\gamma }_i')\),
\([\td{C}']\in\scrH^\bbR(Y_i, \nu_i,\grs_i; \tilde{\gamma }'_i,
\tilde{\gamma }_i'')\),
\begin{equation}\label{I_H-comp}
  I_\scrH([\td{C}]*[\td{C'}])=I_\scrH([\td{C}])+I_\scrH([\td{C'}]) .
\end{equation}
\end{lemma}
\pf
Let \(\tilde{C}\) be a 2-current on \(\bbR\times Y_i\) representing and
element \([\tilde{C}]\in \scrH^\bbR(Y_i, \nu
_i,\grs_i; \tilde{\gamma }_i, \tilde{\gamma }_i')\), and let \(h\) be
a harmonic 2-form on \(Y_i\) representing an element in
\([h]\in H^2(Y_i;\bbR)\). The integral of \(\pi ^*_2h\) over \(\tilde{C}\) is
a real number depending only on the class \([\tilde{C}]\); so
in this way we have a map 
\[
I_\scrH\co \scrH^\bbR(Y_i, \nu
_i,\grs_i; \tilde{\gamma }_i, \tilde{\gamma }_i')\to \Hom 
(H^2(Y_i;\bbR); \bbR)\simeq H_2(Y_i;\bbR). 
\]
It is straightforward to verify that the \(H_2(Y_i;\bbR)\)-action on
the space \(\scrH^\bbR(Y_i, \nu
_i,\grs_i; \tilde{\gamma }_i, \tilde{\gamma }_i')\) intertwines with
that on \(H_2(Y_i;\bbR)\), as well as the composition rule (\ref{I_H-comp}). 
\epf

The \(\bbZ\)-coefficients version of statement in the preceding lemma
does not hold. However, since when \(\tilde{\gamma }_i, \tilde{\gamma }_i'\in \scrZ_{(Y_i,
  \nu_i, \grs_i)}\subset  \scrZ^\bbR_{(Y_i,
  \nu_i, \grs_i)}\), the space \(\scrH(Y_i, \nu
_i,\grs_i; \tilde{\gamma }_i, \tilde{\gamma }_i')\) embeds in \(\scrH^\bbR(Y_i, \nu
_i,\grs_i; \tilde{\gamma }_i, \tilde{\gamma }_i')\) as an orbit under
the \(H_2(Y_i;\bbZ)\subset H_2(Y_i;\bbR)\) action, the 
isomorphism \(I_\scrH\) above also identifies it as an
\(H_2(Y_i;\bbZ)\)-orbit in \(H_2(Y_i;\bbR)\); that is, an element
\(t_\scrH (\tilde{\gamma }_i, \tilde{\gamma }_i')\) in
the orbit space \(\bbT_{Y_i}:=H_2(Y_i;\bbR)/H_2(Y_i;\bbZ)\). Equip the
torus 
\(\bbT_{Y_i}\) with the
standard Euclidean metric. This map,
\[
t_\scrH\co \scrZ_{(Y_i,
  \nu_i, \grs_i)}\times \scrZ_{(Y_i,
  \nu_i, \grs_i)}\to \bbT_{Y_i} ,\] 
is continuous with respect to the
current topology on \(\scrZ_{(Y_i,
  \nu_i, \grs_i)}\), and is linear with respect to \(\tilde{\gamma
}_i'-\tilde{\gamma }_i\in \scrZ(Y_i)\), where \(\scrZ(Y_i)\) is the
space of 
exact integral 1-currents on \(Y_i\). Thus, one has: 
\begin{lemma}\label{gr:h-isom}
Suppose \(\tilde{\gamma }_i, \tilde{\gamma }_i'\in \scrZ_{(Y_i,\nu_i,
  \grs_i)}\) are sufficiently close in the sense that \(t_\scrH
(\tilde{\gamma }_i, \tilde{\gamma }_i')\) falls in the ball \(B_{0}(1/2)\subset
\bbT_{Y_i}\). Then there is a distinguished element in \(\scrH(Y_i, \nu
_i,\grs_i; \tilde{\gamma }_i, \tilde{\gamma }_i')\).   The notion of
``sufficiently close'' depends on the metric on \(Y_i\), but the
distinguished element is independent of the
metric on \(Y_i\). 

Suppose a pair \(\{\tilde{\gamma}_i\}_i,
,\{\tilde{\gamma}'_i\}_i\in \prod_{i\in \grY_m}\scrZ_{(Y_i,\nu_i,
  \grs_i)}\) are sufficiently close in the sense described above, and
let \(h_i\in \scrH(Y_i, \nu
_i,\grs_i; \tilde{\gamma }_i, \tilde{\gamma }_i')\) denote the
aforementioned distinguished element. Then the concatenation map 
\(\op{c}_{\{h_i\}_i}\)  defines a canonical isomorphism of affine spaces
under \(\scrH_X\) from  \(\scrH \big((X^{'a}, \nu,
\grs),\{\tilde{\gamma}_i\}_i\big)\) to \(\scrH \big((X^{'a}, \nu,
\grs),\{\tilde{\gamma}'_i\}_i\big)\). 
\end{lemma}

\subsection{Some basic facts about t-curves}

Analogs of the results below are well-known in the context of
symplectic field theory.
\begin{prop}
\label{prop:t-curve-asymp}
Let \((X, \nu)\) be an admissible pair, and let \(\{(Y_i,
\nu_i)\}_{i\in\grY}\) be the  ending pairs. A t-curve
\(\mathbf{C}=[C, \td{C}]\) with finite \(\omega\)-energy in \(X^{'a}\) is asymptotic to some t-orbit
\(\pmb{\gamma}_i=[\gamma _i, \td{\gamma }_i]\) in the pair \((Y_i, \nu_i)\), \(\forall i\in
\grY_m\), in the following sense:
\begin{itemize}
\item[(i)] 
Given any sequence of positive numbers \(\{L_n\}_{n\in \bbZ^+}\) with
  \(\lim_nL_n\to \infty\), 
  \(\{(\tau_{-L_n}C)|_{\hat{Y}_i}\}_n\) geometrically converges to
  \([0, \infty)\times\gamma_i\) in \(\hat{Y}_i\simeq[0, \infty)\times Y_i\). In fact, 
for all sufficiently small \(\varepsilon>0\), there exists an 
  \(R(\varepsilon)>0\), such that \(\dist_{\hat{Y}_i}(\tau_{-L}C,
  [0, \infty)\times\gamma_i)<\varepsilon  \), \(\forall L>R(\varepsilon)\).
\item[(ii)] The sequence of currents \(\{\tau_{-L_n}\tilde{C}\}_n\) converges weakly to
\(\pi_2^*\tilde{\gamma}_i\), where \(\pi_2\co \hat{Y}_i\simeq[0,
\infty)\times Y_i\to
Y_i\) is the projection. 
\end{itemize}
\end{prop}

\pf  We first show that the statements (i) and (ii) hold for a
subsequence of \(\{L_n\}_n\).

The finite \(\omega\)-energy assumption in the definition of
t-curves  implies that
\(\forall i\in \grY_m\), 
\begin{equation}\label{E-w-bdd}
\sup_{L> 0}\int_{C|_{\hat{Y}_{i, [L, L+1]}}}\omega < z_i\quad
\text{for some $z_i>0$.}
\end{equation} 
Thus, by the admissibility condition on $\nu$ we may apply 
Proposition \ref{thm:gromov-cpt} to find a subsequence \(\{R_n\}_{n\in\bbZ^+}\) of \(\{L_n\}_n\)
with \(\lim_{n\to \infty}R_n=\infty\), such that 
\(\{\tau_{-R_{n}}C\}_n\) converges geometrically over
\([0,1]\times (Y_i\backslash \nu_i^{-1}(0))\) to a $J_i$-holomorphic
subvariety $C_i$. Here, $J_i$ is the almost complex structure
determined by the metric and the symplectic form $\nu_i+ds
\wedge*_3\nu_i$ on \([0,1]\times (Y_i\backslash \nu_i^{-1}(0))\). Because of (\ref{E-w-bdd}) and the fact that
\(\scrF_\nu(C|_{\hat{Y}_{i}})<\infty\), we know
that \(\int_{C_i}\omega< z_i\), and \(\int_{C_i}\nu=0\). By Lemma
\ref{lem:F-curve} and the conditions that define a t-curve, 
we know that \(C_i=[0,1]\times \gamma_i\), where \(\gamma_i\) is
the underlying subvariety of a certain t-orbit. (The weights on
closed orbits will be determined later).

Observe that statement (a) below implies statement (b), which in turn
implies statement (c). 
\begin{itemize}
\item[(a)] the geometric convergence of \(\{\tau_{-R_{n}}C\}_n\) over
\([0,1]\times (Y_i\backslash \nu_i^{-1}(0))\) to \([0,1]\times \gamma_i\).
\item[(b)] the geometric convergence of \(\{\tau_{-R_{n}}C\}_n\) over
\([0,1]\times Y_i\) to \([0,1]\times \gamma_i\).
\item[(c)] the geometric convergence of (a subsequence of)
  \(\{\tau_{-R_{n}}C\}_n\) over \(\bbR^+\times Y_i\) to \(\bbR\times
  \gamma_i\). We use the same notation $\{\tau_{-R_n}C\}_n$ to denote
  the subsequence.
\end{itemize}
The statement (c) follows from statement (b) by a diagonalization
argument. To see why (a) implies (b), let $U_\epsilon\subset Y_i$
denote the set of points with distance less or equal $\epsilon$ from
$\nu_i^{-1}(0)$. By (a), for any small \(\varepsilon >0\), \(\exists
\, N(\varepsilon)\) such that \[\dist_{[0,1]\times (Y_i\backslash   U_{\varepsilon/2})}(\tau_{-R_{n}}C,
[0,1]\times \gamma_i)<\varepsilon/2\quad \forall n>N(\varepsilon).\] 
On the
other hand, \(\dist_{[0,1]\times U_{\varepsilon/2}}(\tau_{-R_{n}}C,
[0,1]\times \gamma_i)<\varepsilon/2\). Together, this means
\(\dist_{[0,1]\times Y_i}(\tau_{-R_{n}}C, [0,1]\times 
\gamma_i)<\varepsilon\) \(\forall n>N(\varepsilon)\). 

We have previously established statement (a) for
\(\{\tau_{-R_n}C\}_n\). Thus, the above observation establishes
assertion (i) in the statement of the Proposition for \(\{\tau_{-R_n}C\}_n\).
Next we establish assertion (ii) for \(\{\tau_{-R_n}\mathbf{C}\}_n\).
Define the current \(\tilde{\gamma}_i\) such that
\(\pi_2^*\tilde{\gamma}_i\) is the (weak) limit of the currents
\(\tau_{-R_{n}}\tilde{C}\) over \([0,1]\times (Y_i\backslash\nu_i^{-1}(0))\). The currents \(\tau_{-R_{n}}\tilde{C}\), \(\pi_2^*\tilde{\gamma}_i\) over
\([0,1]\times (Y_i\backslash\nu_i^{-1}(0))\) extend over \([0,1]\times
Y_i\), since \([0,1]\times \nu_i^{-1}(0)\) is a co-dimension 3
submanifold 
in \([0,1]\times Y_i\). Set \(\pmb{\gamma}_i\) to be the t-orbit
with underlying subvariety \(\gamma_i\), and with associated current \(\tilde{\gamma}_i\).

To complete the proof of this proposition, it now suffices 
to show that \([0,1]\times \gamma_i\) is the geometric limit
of \(\{\tau_{-L_n}C|_{\hat{Y}_{i, [0,1]}}\}_n\)
for {\em any} sequence \(\{L_n\}_n\) with \(\lim_{n\to \infty
 }L_n=\infty\). Suppose the contrary. Then there is a sequence
\(\{L_n\}_n\), \(\lim_{n\to\infty}L_n=\infty\) such that \(\forall n\),
\(\tau_{-L_n}C|_{\hat{Y}_{i,[0,1]}}\) has distance larger than
\(\varepsilon\) from \([0,1]\times \gamma_i\) for some small \(\varepsilon >0\). The same argument via Theorem
\ref{thm:gromov-cpt} as above shows that there is a subsequence
\(\{L'_n\}_n\) of \(\{L_n\}_n\), such that \(\lim_{n\to\infty}L_n'\to
\infty\) and \(\{\tau_{-L'_n}C|_{\hat{Y}_{i,[0,1]}}\}_n\)
converging geometrically to \([0,1]\times \gamma_i'\). Here, \(\gamma_i'\) is the
underlying subvariety of a different t-orbit, and \(\dist_{Y_i}
(\gamma_i, \gamma_i')>\varepsilon\). Without loss of generality, we
may choose \(\{L_n'\}_n\) such that \(L'_n>3R_n\). By Proposition \ref{thm:gromov-cpt}
 again, \(\{\tau_{-(2R_n)}C\}_n\) converges
geometrically to a pseudo-holomorphic variety  \(C'\subset
\bbR\times Y_i\) asymptotic to \(\bbR\times \gamma_i\) and
\(\bbR\times \gamma_i'\) respectively on
the \(-\infty\)- and \(+\infty\)-ends. However, \(\int_{C'}\nu=0\),
and thus by Lemma \ref{lem:F-curve}, \(C'=\bbR\times \gamma_i\) and
\(\gamma_i=\gamma '_i\), contradicting our assumption on \(\gamma_i'\). 
\epf

We say that a t-curve in a cylindrical admissible pair
\((\bbR\times Y_i, \pi_2^*\nu_i)\) is {\em constant} if it is of the form
\(\bbR\times \pmb{\gamma}\) for some t-orbit \(\pmb{\gamma}\) of
\((Y_i, \nu_i)\). 
\begin{lemma}\label{lem:hbar}
Let \((Y, \grs_Y)\) be a closed \(\Spin^c\) 3-manifold. Let \(g\) denote a
riemannian metric on \(Y\), and let \((X, \nu)\) denote a
cylindrical admissible pair with \(X=\bbR\times Y\). Let \(\grs\) be the \(\Spin^c\) structure on
\(X=\bbR\times Y\) induced by \(\grs_Y\). 
There exists a
constant \(\hbar >0\) such that any nonconstant t-curve
\(\mathbf{C}\) of \(\Spin^c\)-structure \(\grs\) in \((X, \nu)\) has \(\scrF_\nu (\mathbf{C})>\hbar\).
\end{lemma}
\pf
By assumption, \(\nu=\pi_2^*\nu_Y\), where \(\nu _Y\) is Hodge dual to
a regular Morse-Novikov 1-form. There are finitely many t-orbits
 with a fixed \(\Spin^c\) structure in \(Y\) (cf. e.g. \cite{HL}).
 Then without loss  of generality, we may restrict to t-curves asymptotic to fixed
 t-orbits \(\pmb{\gamma}_\pm\) respectively in the \(\pm\infty\)-ends. 

Suppose the lemma is false and there is a 
sequence of nontrivial t-curves \(\{\mathbf{C}_n\}_n\) asymptotic to
t-orbits \(\pmb{\gamma}_\pm\) on the \(\pm\infty\)-ends, such
that \(\scrF_\nu(\mathbf{C}_n)\to 0\). We claim that \(\forall \varepsilon >0\),
\(\exists N\) such that \(\forall n >N\), \(C_n\) is of distance less than \(\varepsilon\) from
\(\bbR\times \gamma\) for some t-orbit
\(\pmb{\gamma}\). (\(C_n\) and \(\gamma\) respectively denote the
underlying subvariety of \(\mathbf{C}_n\) and \(\pmb{\gamma}\)). It would then follow that \(C_n\) has \(\gamma\) as
both the \(+\infty\)- and \(-\infty\) ends, and \(\pi_2 (C_n)\)
represents a trivial element in \(H_2(Y;\bbZ)\), since
\(H_2(N_\varepsilon; \bbZ)\) is trivial. Here, \(N_\varepsilon\)
denotes the tubular neighborhood of \(\bbR\times \gamma\) consisting
of points of distance less than \(\varepsilon\) from \(\bbR\times
\gamma\). Thus,
\(\scrF_\nu(C_n)=0\) and we saw from Lemma \ref{lem:F-curve}
that \(C_n\) must be constant, contradicting our assumption. 

To prove the claim, again suppose the contrary, that there exists of
sequence of t-curves \(\{\mathbf{C}_n\}_n \), \(C_n\not\subset
N_\varepsilon\) \(\forall n\) for any t-orbit \(\pmb{\gamma}\)
with the prescribed $\Spin^c$ structure. 
Let \(\Gamma\subset Y\) denote the union of all underlying
subvarieties of such t-orbits and let
\(N^\Gamma_\varepsilon\subset X\) be the set consisting of points with
distance less than \(\varepsilon\) from \(\bbR\times\Gamma\). Any pseudo-holomorphic
variety \(C_n\) from the above sequence is not included in \(N_\varepsilon ^\Gamma\).
We may then choose \(L_n\in \bbR\)
such that \((\tau_{L_n}C_n)\cap (\{0\}\times Y)\) contains a point with
distance larger than \(\varepsilon \) from \(\Gamma\). By Theorem \ref{thm:gromov-cpt},
\(\{\tau_{L_n}C_n\}_n\) geometrically converges to a
nonconstant t-curve \(\mathbf{C}\) with
\(\scrF_\nu(\mathbf{C})=0\), again contradicting Lemma \ref{lem:F-curve}. 
\epf

\section{Preliminaries: the \(SW\) side}

This section consists of a minimal review of Seiberg-Witten theory
in the context of MCE's as well as some setups and definitions. 

Recall the definitions of \(\Conn (M)\) and related notions in Section
\ref{sec:convention} (19), where \(M\) is a \(\Spin^c\) 3- or 4-manifold,
with spinor bundle 
\(\bbS\) in the 3-dimensional case or 
\(\bbS^+\oplus \bbS^-\) in the 4-dimensional case. Let \(\scrC(M)\) denote \(\Conn (M)\times \Gamma
(\bbS)\) in the 3-dimensional case, and \(\Conn (M)\times \Gamma
(\bbS^+)\) in the 4-dimensional case. This is said to be the
(Seiberg-Witten) {\em configuration space}, and an element of \(\scrC(M)\)
is said to be a (Seiberg-Witten) {\em configuration}. The {\em gauge group} over
\(M\) is \(\scrG=\scrG (M)=C^\infty(M; S^1)\). It acts on the space \(\scrC(M)\)
by gauge transformations: Recall in particular that an element \(u\in C^\infty(M; S^1)\) defines
an automorphism on \(\Conn (\det \bbS^+)\times \Gamma (\bbS^+)\) (or
\(\Conn (\det \bbS)\times \Gamma (\bbS)\)) by 
\(
u\cdot (A, \Psi )=(A-2u^{-1}du, u\cdot \Psi ). 
\)
We use \(\scrB(M)\) to denote the quotient space
\(\scrB(M):=\scrC(M)/\scrG\). This is called the quotient
configuration space, and its elements quotient configurations (or
gauge equivalence classes of configurations). The gauge equivalence of
a configuration \((A, \Psi )\) is denoted by \([(A, \Psi )]\).


\subsection{Setup and Assumptions}\label{sec:adm-SW}

Let \(Y\) be a closed connected \(\Spin^c\) 3-manifold. 
We first consider the cylindrical case \(X=\bbR\times Y\), or more
generally, consider \(Z=I\times Y\subset X\) for an interval \(I\subset \bbR\).  Let \(\bbS_X=\bbS^+\oplus \bbS^-\) denote spinor bundle
over \(X\), and \(s\in\bbR\) be an affine coordinate for the first
factor of \(X=\bbR\times Y\). Let \(\mu _Y\)
be a closed 2-form on \(Y\), and set  \(\mu=\pi^*_2\mu_Y\), \(\pi_2\co
\bbR\times Y\to Y\) being the projection to the second factor. Let
\(\grq\in \scrP(Y)\), where \(\scrP (Y)\) denote \(Y\)'s version of the Banach space of
large tame perturbatons \(\scrP\) introduced in \cite{KM}'s Theorem 11.6.1 and Definition 11.6.3. 
As
explained in II.4.3 of \cite{KM}, the Clifford action \(\rho(\partial_s)\)
defines an isomorphism \(\bbS^+\simeq\bbS^-\), which in turn is 
identified with the pull back of the spinor bundle \(\bbS=\bbS_Y\) on
\(Y\). Meanwhile, a \(\Spin^c\)-connection \(A\) on \(\bbS^+\) may be
written (via a gauge tranformation) in temporal gauge as a path of \(\Spin^c\) connections
\(B(s)\) on \(\bbS\) by way of Equation (4.8) of \cite{KM}. 
In this manner, an \((A, \Psi)\in \scrC(X) 
\) in temporal gauge on the cylinder \(\bbR\times
Y\) can alternatively be expressed as a path \((B(s), \Phi(s))\) in
\(\scrC(Y)\). By way of  the latter interpretation of
\((A, \Psi)\), each \(\grq\in \scrP(Y)\) determines a 4-dimensional
perturbation denoted by  \(\hat{\grq}\). (Cf. Equation (10.2)
of \cite{KM} and the subsequent discussion.)  
Also in terms of the aforementioned 
temporal gauge 
expression, the Seiberg-Witten equation \(\grS_{\mu, \hat{\grq}}(A,
\Psi)=0\)  can be written as a gradient flow equation: \[\big(\frac{1}{2}\frac{dB}{ds},
\frac{d\Phi}{ds}\big)=-\grF_{\mu_Y}(B, \Phi)-\grq(B, \Phi).\] 
In the above, \(\grF_{\mu_Y}\) is as defined in 
(\ref{eq:SW3}), regarded as a section of the tangent bundle over
\(\scrC(Y)\). The term \(\grq\) appearing in the preceding equation is
understood as another section of this bundle: As in \cite{KM}, the same notation is used to denote both an element
\(\grq\in \scrP(Y)\) and its  image under the map \(\grQ\) in
\cite{KM}'s Theorem 11.6.1. 
Write \[
  \grF_{\mu_Y, \grq}=\grF_{\mu_Y}+\grq.\]
This is the
formal gradient of the (perturbed) Chern-Simons-Dirac functional on \(\scrC(Y)\): 
\begin{equation}\label{eq:CSD_q}
\op{CSD}_{\mu_Y, \grq}(B, \Phi) := \op{CSD}_{\mu_Y}(B, \Phi)+f_\grq (B,
\Phi),
\end{equation}
where
\begin{equation}\label{eq:CSD}\begin{split}
\op{CSD}_{\mu_Y}(B, \Phi) &:=-\frac{1}{8}\int_Y (B-B_0)\wedge
(F_B+F_{B_0})+\frac{1}{2}\int_Y \langle \slp_B \Phi, \Phi\rangle \,
d \op{vol} \\ 
&\qquad -\frac{i}{8}\int_Y \mu_Y\wedge (B-B_0). \\
\end{split}
\end{equation}

In the above, \(B_0\in \scrC(Y)\) is a reference connection, and $f_\grq$ is the
$\bbR$-valued, gauge 
invariant function on the space
\(\scrC(Y)\) 
with formal gradient \(\grq\),
and with \(f_\grq(B_0, 0)=0\). Given \((B, \Phi)\in \scrC(Y)
\), we use
\((\hat{B}, \hat{\Phi})\) to denote the element in \(\scrC(\bbR\times
Y) 
\) corresponding to the constant path at \((B, \Phi)\). 
In particular, a solution \((B, \Phi)\)
to the 3-dimensional Seiberg-Witten equation \(\grF_{\mu_Y, \grq}(B,
\Phi)=0\) gives rise to an
\(\bbR\)-invariant Seiberg-Witten solution: \(\grS_{\mu,
  \hat{\grq}}(\hat{B}, \hat{\Phi})=0\) on \(\bbR\times Y\). 

Note that \(\op{CSD}_{\mu_Y}\) is in general not gauge invariant, but is
invariant under the identity component of the gauge group, denoted \(\scrG_0\subset C^\infty(Y,
S^1)\). It also depends implicitly on the choice of the reference
connection \(B_0\). It is written as \(\op{CSD}^{B_0}_{\mu_Y}\) when we
wish to emphasize this dependence. Similarly for \(\op{CSD}_{\mu_Y,
  \grq}\).  A further superscript \(Y\) is added, e.g. 
\(\op{CSD}_{\mu_Y}^{Y, B_0}=\op{CSD}_{\mu_Y}\) when needs for specifying the
3-manifold \(Y\)  arise. 

Now let \(X\) be a general \(\Spin^c\) MCE with spinor bundle \(\bbS_X=\bbS^+\oplus \bbS^-\), and let 
\(\bbS_i\) 
denote the spinor
bundle on the ending \(\Spin^c\) 3-manifold \(Y_i\). 
{\em Fix a reference connection} 
\(A_0\in \Conn (X) \) such that 
\begin{equation}\label{eq:A_0}
A_0=\hat{B}_{0,i}
\end{equation}
for a reference
connection \(B_{0,i}\in \Conn (Y_i)\)
on each end \(\hat{Y}_i\subset X\).

\begin{defn}\label{def:adm-SW}
An {\em admissible} element \((A, \Psi)\) of \(\scrC(X)\) 
is one with the following property:
 There exists \((B_i, \Phi_i)\in\scrC(Y_i)
\) such that  \((A,
\Psi)-(\hat{B_i}, \hat{\Phi}_i)\) is in \(L^2_1(T^*\hat{Y}_i)\times
L^2_{1, A_0}(\bbS^+)\) over each end \(\hat{Y}_i\). 
In this case,
\((B_i, \Phi_i)\) is said to be the {\em \(Y_i\)-end limit} of \((A,
\Psi)\).
\end{defn}

All elements \((A, \Psi)\) of \(\scrC(X)\) are assumed to be
admissible in this article. 

\begin{remarks}
(a) Note that the preceding definition of admissible \((A, \Psi)\) does not depend on the
choice of the reference connection \(A_0\). 

(b) When \((A, \Psi
)\) is a  Seiberg-Witten solution, the definition
of  ``\(Y_i\)-end limits of  \((A, \Psi)\)'' above is consistent with the
notion  of \(Y_i\)-end limits of general sections over \(X\) given in
Section \ref{sec:convention}, by way
of a well-known elliptic bootstrapping argument; see e.g. \cite{KM}. 
\end{remarks}

This article concerns admissible solutions \((A, \Psi)=(A_r, \Psi
_r)\) to the Seiberg-Witten equation \(\grS_{\mu _r,
  \hat{\grp}_r}(A_r, \Psi _r)=0\), where \(\mu _r=r\nu +w_r\) and
\(\hat{\grp}_r\) are as described in Theorem \ref{thm:l-conv}. 
The terms  \(w_r\) and 
\(\hat{\grp}_r\) are regarded as ``additional  perturbations'' to the
dominating perturbation \(r\nu \) in the Seiberg-Witten equation;
they are ultimately inconsequential to
the \(SW\Rightarrow Gr\) story told in this article. They are chosen
to satisfy the conditions spelled out in Assumption \ref{assume}
below. These conditions are meant to both 
simplify the arguments yet remain sufficiently general
for the purposes of Floer theory (e.g. transversality and compactness). The subscripts \(r\)
are sometimes dropped from \(w_r\) and \(\hat{\grp}_r\), as their
dependence on \(r\) is immaterial. 

In order to describe the assumptions as well as for future reference, we list some basic facts on \cite{KM}'s nonlocal perturbations: 
\begin{facts} \label{list:scrP}
{\bf 1)} Given \(Y\), there exists  a constant \(m\) such that 
\begin{equation}\label{eq:q-bdd}
\|\grq(B, \Phi)\|_{L^2}\leq
m\|\grq\|_\scrP\,( \|\Phi\|_{L^2}+1).
\end{equation}
for every \(\grq\in \scrP (Y)\) and \((B, \Phi )\in \scrC(Y)\). Here, 
\(\|\grq\|_\scrP\) denote the norm of \(\grq\) as an element in the
Banach space \(\scrP(Y)\). (Cf. Equation (11.18) of \cite{KM}.) Given an compact interval
\(I\subset \bbR\), let \(Z\) denote the corresponding \(\Spin^c\)
cylindrical manifold \(I\times Y\). Let \(\bbS^+\) denote the spinor
bundle over \(Z\) and fix a reference connection \(A_0\in \Conn
(\bbS^+)\). Then there exists a continuous function over \(\bbR\),
denoted 
\(\hat{m}\), such that 
\begin{equation}\label{eq:q-hat-bdd}
\|\hat{\grq}(A, \Psi)\|_{L^2_{1, A}}\leq  \|\grq\|_\scrP  \, \hat{m}\,
(\|(A-A_0,
  \Psi)\|_{L^2_{1, A_0}}).
\end{equation}
(See Definition 10.5.1 and Theorem 11.6.1 of \cite{KM}. Note that
there is a misprint in Theorem 11.6.1 (iv). By a Sobolev inequality,
(\ref{eq:q-hat-bdd}) is equivalent to \(\|\hat{\grq}(A, \Psi)\|_{L^2_{1, A}}\leq  \|\grq\|_\scrP  \, \hat{\rmm}\,
(\|(A-A_0,
  \Psi)\|_{L^2_{1, A}})\) for another continuous function \(\hat{\rmm}\).)

{\bf 2)} Let \(Z\), \(\grq\) be as before, and let the Banach manifolds
\(\scrC_k(Z)\) and the Banach bundles 
  \(\scrV_k(Z)\to \scrC_k(Z)\)  be as in \cite{KM}'s Definition 10.1.1. (The
  subscript \(k\) signifies \(L^2_k\)-completion.)  For \(k\geq 1\),
  the map \(\grQ\) from \(\scrP(Y)\times \scrC_k(Z)\) to \(\scrV_k(Z)\) given
  by \((\grq, c)\mapsto \hat{\grq}(c)\) is a continuous map of Banach
  manifolds which vanishes over \(\{0\}\times \scrC_k(Z)\subset
  \scrP(Y)\times \scrC_k(Z)\). For \(k\geq 2\), it is a smooth map. 

{\bf 3)} Let \(X\) be a \(\Spin^c\) MCE with ending 3-manifolds
  \(Y_i\). A nonlocal perturbation \(\hat{\grp}\) on \(X\), as defined
  in Equation (24.2) of \cite{KM}, is of the following
  form: For every \(i\in \grY\), there exists a compact interval \(\bbI_i=[\grl_i,
  \grl_i']\subset\bbR\) and \(\grq_i, \grp'_i\in \scrP(Y_i)\) such that
  \(\hat{\grp}\) vanishes over \(X_{\pmb{\grl}}\); it restricts to
  \(\hat{\grq}_i\) over \(\hat{Y}_{i, \grl_i'}\), and over
  \(\hat{Y}_{i, \bbI_i}\) it takes the form 
\begin{equation}\label{eq:hatp}
\hat{\grp}=\chi_i(\mathfrc{s}_i)\hat{\grq}_i+\lambda_i(\mathfrc{s}_i)\hat{\grp}'_i, 
\end{equation}
where  \(\chi_i\) is a
  smooth cutoff function supported on \((\grl_i, \infty)\) and equals
  \(1\) as \(\mathfrc{s}_i>\grl'_i\); \(\lambda_i\) is a smooth bump
  function supported on \(\bbI_i\) and equals 1 on a strictly smaller
  interval \(\bbI'_i\subset \bbI_i\). As usual, \(\pmb{\grl}\) above denotes a
  function from \(\grY\) to \(\bbR^{\geq 0}\) given by \(i\mapsto
  \grl_i\). 
\end{facts}

\begin{assumption}\label{assume}
Let \(X\) be a \(\Spin^c\) MCE with ending 3-manifolds
  \(Y_i\), \(i\in \grY\).  Assume that \(w_r\), \(\hat{\grp}_r\) satisfy the following conditions:
\begin{itemize}
\item[(1)] \(w_r\) is a smooth closed 2-form depending smoothly on
\(r\), and \(\|w_r\|_{C^{k+3}}< \varsigma_w\) for some \(k\geq 3\) and an \(r\)-independent
positive constant \(\varsigma_w\). 
\item[(2)]  \(w_r\) are in the cohomology class \(4\pi  c_1(\grs_X)\) for all \(r\);
\item[(3)] Over each end \(\hat{Y}_i\subset X\), \(w_r=\pi_2^*w_{i, r}\) for  a closed 2-forms \(w_{i, r}\) on
  \(Y_i\) .  If \(\hat{Y}_i\) is a
  vanishing end, the  closed form \(w_{i,r}\) is independent of
   of \(r\) and is written as \(w_i\). 
\item[(4)]  \(\hat{\grp}_r=\hat{\grp}\) is of the type given
  in Equation (24.2) of \cite{KM}. Adopt the notations from
 Item 3) of  Facts \ref{list:scrP}, with the ingredients for constructing
  \(\hat{\grp}\) (i.e. \(\grq_i\), \(\grp'_i\), \(\grl_i\),
  \(\grl_i'\), \(\chi _i\), \(\lambda _i\)) understood to 
  (implicitly) depend on \(r\). (A subscript \(r\) is added,
  e.g. \(\grq_i=\grq_{i,r}\), when the \(r\)-dependence is
  emphasized.) We require: 
\begin{itemize}
\item[a)] \(\grq_i=0=\grp'_i\)  when \(\hat{Y}_i\) is a Morse
    end. In this case we formally set \(\grl_i=\grl_i'=\infty\) and
    \(\bbI_i=\emptyset\). 
\item[b)]  The triples \(\bbI_i\), \(\chi _i\), \(\lambda _i\)
  corresponding to different \(i\in \grY_v\) are related by
  translations. Namely, they take the form of 
\[
\bbI_i=\tau _{\grl_i}\bbI_0; \quad \chi _i=\tau _{\grl_i}^*\chi
_0; \quad \lambda _i=\tau _{\grl_i}^*\lambda _0, 
\] 
where \(\bbI_0:=[0,\zeta _0]\), \(\chi _0, \lambda _0\) are all independent of both \(r\) and \(i\). 
\item[c)]  \(\grl_i>\ul_i^+\) 
\(\forall  i\in \grY_v\) and \(r\),  where \(\ul_i^+\in \bbR^+\) is as in Definition \ref{def:adm}. In
 particular, \(\hat{\grp}\) vanishes over \(X''\).
\item[d)]  Abusing notation, let \(\|\hat{\grp}\|_{\scrP}:=\sum_{i\in
  \grY_v}\big(\|\grp'_i\|_{\scrP}+\|\grq_i\|_{\scrP}\big)\). There is an \(r\)-independent constant \(\zzz_\grp<1/8\)
  such that 
\[
\|\hat{\grp}\|_{\scrP}=\sum_{i\in
  \grY_v}\big(\|\grp'_i\|_{\scrP}+\|\grq_i\|_{\scrP}\big)<\zzz_\grp
\quad \forall r. 
\]
\end{itemize}

\item[(5)]  Let \(k\) be  as in Item (1) above, and  let \(\bbI_\epsilon :=(-\epsilon , \zeta
_0+\epsilon )\supset\bbI_0\) be a slightly larger open interval containing
\(\bbI_0\). Denote \(\upsilon _i(r):=\|\mu_r\|_{C^{k+3}(\hat{Y}_{i, \grl_i})}\). 
Fix an \(r\)-independent constant \(\zeta _\grp<\zeta
_w/8\). With the constants and functions \(\zeta _0\), \(\lambda _0\), \(\chi _0\), \(\grl_i\) from Assumption \ref{assume}
4b)-4c) fixed, \(\hat{\grp}\) is determined by \(\grq_i\) and
\(\grp'_i\). We imposes a stronger, \(r\)-dependent constraint
(cf. the upcoming remark) on their 
sizes (controlled by
\(\|\hat{\grp}_r\|_\scrP\)) than that from Item 4d) above, 
such that the following holds
for every \(i\in \grY_v\) and \(r\geq 1\): Fix \(i\in \grY_v\) and drop the
subscript \(i\) for the rest of this item. Write \(Z:=\bbI_0\times Y\); \(Z_\epsilon
:=\bbI_\epsilon \times Y\), and let 
\(\hat{\grx}_r:=\chi_0\hat{\grq}_r+\lambda _0\hat{\grp}'_r=\tau
_{\grl_r}^*\hat{\grp}_r\). Then give \(r\geq 1\): 
\begin{equation}\label{zeta_p2a}
 \|\hat{\grx}_r(A, \Psi )\|_{C^k_{A}(Z)}<\zeta _\grp
\end{equation}
for any \((A, \Psi )\in \scrC(Z_\epsilon )\) that satisfies:
\begin{itemize}
\item[(i)] \(\|F_A\|^2_{L^2(Z_\epsilon )}+\|\nabla_A\Psi
  \|^2_{L^2(Z_\epsilon )}+\|\Psi \|^4_{L^4(Z_\epsilon )}\leq 2^8
  \upsilon (r)^2\);
\item[(ii)]  Let \(\grW\) denote the space of \(C^{k+3}\) closed
  2-forms on \(Z_\epsilon \) cohomologous to \(\pi _2^*w\).  Then \(\grS_{\rmw, \hat{\grx}_r}(A, \Psi )=0\) for a \(\rmw\in
  \grW\)  with  \(\|\rmw\|_{C^{k+3}(Z_\epsilon )}\leq \upsilon (r)\). 
\end{itemize}
Likewise, 
\begin{equation}\label{zeta_p2b}
 \|\hat{\grq}_r(A, \Psi )\|_{C^k_{A}(Z)}<\zeta _\grp
\end{equation}
for any \((A, \Psi )\in \scrC(Z_\epsilon )\) that satisfies both item
(i) above and 
\begin{itemize}
\item[(ii')]
  \(\grS_{\rmw, \hat{\grq}_r}(A, \Psi )=0\) for a \(\rmw\in
  \grW\)  with
  \(\|\rmw\|_{C^{k+3}(Z_\epsilon )}\leq \upsilon (r)\). 
\end{itemize}

\end{itemize}
\end{assumption}

\begin{remarks}
{\bf (a)} Item (5) in the preceding assumption is invoked only in Section
\ref{sec:pt-est} to make the contribution from \(\hat{\grp}\)
ignorable while performing pointwise estimates. It may likely be
removed or weakened with extra work, appealing to the weak maximum principle
in lieu of the strong-maximum-principle  arguments in Section
\ref{sec:pt-est}. 
The stronger upper bound on
\(\|\hat{\grp_r}\|_{\scrP}\)  mentioned in Item (5) may also be
weakened by imposing stronger
(\(r\)-dependent) lower bounds on \(\grl_i\) than those required by
Item 4c). 

{\bf (b)} Here is a more detailed description of the stronger bound on
\(\hat{\grp}_r\) that is 
required by Item (5) to guarantee the conditions (\ref{zeta_p2b}) and
(\ref{zeta_p2a}). 
Fix \(r\geq 1\). Given \(\grq\in \scrP(Y)\), let \(\scrM (\rmw; \grq)\)
denote the space of gauge equivalences of \((A, \Psi )\) satisfying
Condition (i) of Item (5) and the Seiberg-Witten equation
\(\grS_{\rmw, \hat{\grq}}(A, \Psi )=0\). Given \(\grq, \grp'\in
\scrP(Y)\), let \(\scrM (\rmw; \grq, \grp')\)
denote the space of gauge equivalences of \((A, \Psi )\) satisfying
Condition (i) of Item (5) and the Seiberg-Witten equation
\(\grS_{\rmw, \hat{\grx}}(A, \Psi )=0\), where
\(\hat{\grx}:=\chi_0\hat{\grq}+\lambda _0\hat{\grp}'\). Then both the space 
\(\td{\scrM}_r:=\bigcup_{\rmw, \|\rmw\|_{C^{k+3}}\leq \upsilon (r)}\bigcup_{\grq,
  \|\grq\|_\scrP\leq \zzz_\grp} \scrM (\rmw; \grq)\) and the space \(\td{\scrM}'_r:=\bigcup_{\rmw, \|\rmw\|_{C^{k+3}}\leq \upsilon (r)}\bigcup_{\grq,
  \|\grq\|_\scrP\leq \zzz_\grp} \bigcup_{\grp',
  \|\grp'\|_\scrP\leq \zzz_\grp} \scrM (\rmw; \grq, \grp')\) are 
represented by compact subspaces in \(\scrC_{k+3}(Z)\). (To see this,  repeat the elliptic
bootstrapping arguments in \cite{KM}'s proofs of Theorems 10.7.1 and
5.2.1, with  the condition (i) of Item (5) playing the role of \cite{KM}'s
(10.15), and noting the following:  First, the
\(C^{k+3}\)-bound on \(\rmw\) in conditions (ii) and
(ii') implies that such \(\rmw\) lies in a compact subspace in the
space of \(L^2_j\) closed 2-forms in cohomology class \(\pi ^*_2w\),
\(\forall \, 0\leq j<k+3\).  Second, instead of appealing to Condition
(ii) of Definition 10.5.1 as in \cite{KM}'s proof of its Theorems
10.7.1, one makes use of its parametrized generalization 
established in \cite{KM} Theorems 11.6.1 (iii)-(iv).)  
Let \(\scrP(Y;\delta )\subset \scrP(Y)\) denote the ball in
\(\scrP(Y)\) of radius \(\delta\) centered at 0, and let 
\(N_\scrV(\delta )\)  denote the tubular neighborhood of the
zero-section of \(\scrV_{k+3}(Z_\epsilon )\) whose intersection with
the fiber of \(\scrV_{k+3}(Z_\epsilon )\) over any given \((A, \Psi )\in
\scrC_{k+3}(Z_\epsilon )\), namely \(L^2_{k+3, A}\), is a ball of
radius \(\delta \). By the 
Sobolev embedding theorem, there is a constant \(\zeta >0\) such that
\(\|\eta\|_{C^{k}_A(Z)}\leq\zeta  \|\eta\|_{L^2_{k+3, A}(Z_\epsilon
  )}\). According
to Fact \ref{list:scrP} (2), there is an \(\zzz'(r)\in
\bbR^+\) such that \(\scrP(Y;\zzz')\times (\td{\scrM}_r\cup
\td{\scrM}'_r)\subset \grQ^{-1}N_\scrV(\zeta ^{-1}\zeta _\grp)\). We
require that \(2\|\hat{\grp}_r\|_\scrP\) is  smaller than every 
\(i\in \grY_v\)'s version of \(\zzz'(r)\). 

{\bf (c)} 
See Remark \ref{rmk:CSD-change} and \cite{KLT1}-\cite{KLT5} for motivations of the
condition Item (2) in the preceding assumption.  
\end{remarks}

\begin{remarks}\label{rem:M-disconn}
{\bf (a)} The notions  \(\scrC(M)\), \(\scrB(M)\) are defined for possibly
disconnected 3- or 4-manifold \(M\). The formula (\ref{eq:CSD}) for
\(\op{CSD}_\mu \) also applies when the 3-manifold \(Y\) is
disconnected. 
Though in \cite{KM} the Banach space of large tame perturbations \(\scrP\) is defined
for connected 3-manifolds, it generalizes readily to the disconnected
case: Suppose \(M=\bigcup_k M_k\) is a closed \(\Spin^c\)
3-manifold with connected components \(M_k\). Take 
\(\scrP(M):=\prod_k\scrP(M_k)\). Then each \(\grq=(\grq_1, \ldots, \grq_k, \ldots)\in
\prod_k\scrP(M_k)=\scrP(M)\) defines a section of the tangent bundle \(T\scrC(M)=  \boxtimes _k
 T\scrC(M_k)\) over \(\scrC(M)=\prod_k\scrC(M_k)\):  given
\(\grc =(\grc_1,\ldots, \grc_k,
\ldots) \) in \\ \(\prod_k\scrC(M_k)=\scrC(M)\),
\(\grq(\grc):=(\grq_1(\grc_1), ,\ldots, \grq_k(\grc_k),
\ldots) \in \boxtimes_kT\scrC(M_k)=T\scrC(M)\). As before, given a
reference connection \(B_0\in \scrC(M)\), the function \(f_\grq\co
\scrC(M)=\prod_k\scrC(M_k)\to \bbR\) is defined to be the unique
gauge-invariant function with formal gradient \(\grq\) and with
\(f_\grq(B_0, 0)=0\). (So \(f_\grq(\grc)=\sum_k
f_{\grq_k}(\grc_k)\).) The properties (1)-(2) in Fact \ref{list:scrP}
remain valid for such generalized \(\grq\) when \(M\) is compact. 
In this way, the definitions of \(\grF_{\mu ,
  \grq}\), \(\op{CSD}_{\mu , \grq}\) both generalize to the case
of disconnected manifolds. 

{\bf (b)} Let \(X\) be a MCE with ending 3-manifolds \(Y_i\), \(i\in
\grY\) and let \(\mu \) be a closed 2-form on
\(X\) with the closed 2-form \(\mu _i\) on \(Y\) as it \(Y_i\)-end
limit for each \(i\). Let  \(\hat{\grp}\) be a nonlocal perturbation of the type
given in Equation (24.2) of \cite{KM}. For our applications, the
relevant disconnected 3-manifold
\(M=\bigcup_k M_k\) usually consists of slices of \(X\) of the form
\(M_k=Y_{i_k:  s_k}\). When the 3-manifold \(M\) under discussion is
clear from the context, we often abbreviate the
restriction of \(\mu \) to \(M\) as \(\mu \big|=\mu \big|_M\). For any such
\(M\), the nonlocal perturbation \(\hat{\grp}\) on \(X\) also induces
an element in the Banach space of large tame perturbations, denoted
\(\hat{\grp}\big|_M\) or simply \(\hat{\grp}\big|\in \scrP(M)\):  In
the notation of (\ref{eq:hatp}), \( \hat{\grp}|_{Y_{i:s}}=\chi_i(s)\hat{\grq}_i+\lambda_i(s)\hat{\grp}'_i\in
\scrP(Y_i)\). To make the notations less cumbersome, when \(M\) is
understood from the context, 
\(f_{\hat{\grp}|}\), \(\op{CSD}_{\mu |, \hat{\grp}|}^M\) are also 
abbreviated respectively as \(f_{\hat{\grp}}\),  \(\op{CSD}_{\mu ,
  \hat{\grp}}\). 
\end{remarks}

\subsection{\(Y_i\)-end limits of admissible Seiberg-Witten
  solutions}\label{sec:end-limit}

We begin with some general definitions and observations.
  \begin{defn}
Let \(Y\) be a closed \(\Spin^c\) 3-manifold and \(B_0\in \Conn 
(Y)\) be a reference connection as in
(\ref{eq:CSD}). A connection \(B\in \Conn (Y)\) is said to be in a
{\em Coulomb gauge} (with respect to \(B_0\)) if the following holds: 
\begin{equation}\label{eq:deltab0}
d^*(\delta B)=0, \quad \text{where \(\delta B=B-B_0\).} 
\end{equation}
An element \((B, \Phi )\in \scrC(Y)\) is  said to be in a
{\em Coulomb gauge} (with respect to \(B_0\)) if \(B\) is in a Coulomb
gauge.

Fix an orthonormal basis \(\{[h_k]\}_k\) of
\(H^1(Y;\bbZ)\), where the latter space is equipped with the inner product
given by the Riemannian metric on \(Y\).  Let \(h_k\in \Omega^1(Y)\) denote the
harmonic representative of \([h_k]\). We say that a \(B\in \Conn (Y)\)
or a \((B, \Phi )\in \scrC(Y)\) is in 
the {\em normalized Coulomb gauge} (with
respect to \(B_0\) and the basis \(\{[h_k]\}_k\)) if in addition to
(\ref{eq:deltab0}), it also satisfies 
\begin{equation}\label{eq:deltab1}
\text{\(0\leq \Big\langle h_k, \frac{i}{4\pi }\delta
  B\Big\rangle<1\) \(\forall k\)}.
\end{equation}
\end{defn}


For each
gauge equivalence class \([(B, \Phi)]\in \big(\Conn (\bbS_Y)\times
\Gamma(\bbS_Y)\big)/\scrG\), there is a unique representative in the
normalized Coulomb gauge. We usually denote this particular
representative as \([(B,
\Phi)]_c\in \scrC(Y)\).

Working in the normalized Coulomb gauge has the following advantage: 
If \(B\) (or \((B, \Phi )\) is in a normalized Coulomb gauge, then by a G\aa rding inequality,
\begin{equation}\label{eq:garding}
\|\delta B\|_{L^2}\leq \zeta_b \big(\|(d+d^*)\delta B\|_{L^2}+1\big)=\zeta_b \big(\|F_B-F_{B_0})\|_{L^2}+1\big)
\end{equation}
for a positive constant \(\zeta_b\) depending only on the metric on
\(Y\). We shall need the following simple consequence of the preceding observations: 
\begin{lemma}\label{lem:f_q}
Let \(Y\) be closed \(\Spin^c\) 3-manifold \(Y\), and fix a reference
connection \(B_0\in \Conn (Y)\). Recall that \(f_\grq\) denotes the
gauge-invariant function on \(\scrC(Y)\) with formal gradient
\(\grq\in\scrP(Y)\) and with \(f_\grq (B_0, 0)=0\).  
There is a constant \(\zeta \) depending only on the metric of \(Y\) such that 
such that \(\forall\)  \(\grq\in \scrP(Y)\) and \((B, \Phi )\in
\scrC(Y)\), 
\[
|f_\grq (B, \Phi )|\leq \zeta \|\grq\|_\scrP\, \big(
\|F_B-F_{B_0}\|^2_{L^2}+\|\Phi \|^2_{L^2}+1\big).\quad
\]
\end{lemma}
\pf Since \(f_\grq\) is gauge-invariant, it suffices to consider the
case when \((B, \Phi )\) is in the normalized Coulomb gauge. Write
\(\delta B:=B-B_0\) as before and let 
\((B_t, \Phi _t):=(B_0+t\delta B, t\Phi )\) for \(t\in [0, 1]\). Then
according to (\ref{eq:garding}) and (\ref{eq:q-bdd}), 
\begin{equation}\label{ineq:f_q}
  \begin{split}
|f_\grq (B, \Phi )|& \leq \|(\delta B, \Phi )\|_{L^2}\|\grq (B_t, \Phi
_t)\|_{L^2}\quad \text{for a certain \(t\in [0,1]\)}\\
& \leq  m\|\grq\|_\scrP\big(\zeta_b\|F_B-F_{B_0}\|_{L^2}+\zeta_b+\|\Phi \|_{L^2}\big)(\|\Phi
_t\|_{L^2}+1)\\
& \leq \zeta \|\grq\|_\scrP\, \big( \|F_B-F_{B_0}\|^2_{L^2}+\|\Phi \|^2_{L^2}+1\big).
\end{split}
\end{equation}
\epf

Returning now to the setting of this article: Let \((X, \nu)\) be an admissible pair and let \((Y_i, \nu_i)\) denote
the limiting pair of \((X, \nu)\) of the \(Y_i\)-end. Let \(\mu_r\),
\(\hat{\grp}_r\) be as in the statement of Theorem \ref{thm:l-conv}, and
write \(\mu_{i,r}=r\nu_i+w_{i,r}\in \Gamma (Y_i)\) for the
\(Y_i\)-end limit of \(\mu_r\). As previously mentioned, we frequently drop the subscript \(r\)
when \(Y_i\) is a vanishing end, since in this case
\(\mu_{i,r}=w_i\) is independent of \(r\). 
Let \((A, \Psi)=(A_r, \Psi_r)\) be an admissible solution to the Seiberg-Witten
equations \(\grS_{\mu_r, \hat{\grp}_r}(A, \Psi)=0\). The standard
  compactness/properness results in
  Seiberg-Witten theory (see e.g. \textsection 10.7 in \cite{KM}) imply that
  the \(Y_i\)-end limit of \((A, \Psi)\) is a solution to the
  3-dimensional Seiberg-Witten equations \(\grF_{\mu_{i,r},
    \grq_i}(B_i, \Phi_i)=0\). Moreover, \((A, \Psi)\big|_{Y_{i:L}}\)
  converges to \((B_i, \Phi_i)\) in \(C^k(Y_i)\) topology as \(L\to \infty\). 
In the remainder of this subsection,  we gather some facts about
such 3-dimensional Seiberg-Witten
solutions that will be of use later. In what follows, all  statements
concern only the pair of 3-manifold and harmonic 2-form
\((Y_i, \nu _i)\), the closed 2-form and the nonlocal perturbation,
\(w_{i, r}\) and \(\grq_{i,r}\), over \(Y_i\); they are not necessarily ending pairs and \(Y_i\)-end
limits of  a 4-dimensional  \((X, \nu )\), \(w_r\), 
\(\hat{\grp}_r\). The expression \(i\in \grY_v\) or \(i\in \grY_m\)
only serves to indicate whether \(\nu _i\equiv 0\) or
otherwise. 
\(w_{i,r}\) and \(\grq_{i,r}\) are assumed to satisfy the 
uniform bounds \(\|w_{i,r}\|_{C^{k+3}}\leq \varsigma _w\) and
\(\|\grq_{i, r}\|_\scrP\leq \zzz_\grp\) for an integer \(k\geq 2\).  For simplicity we also
assume that the cohomology class \([w_{i,r}]\) is independent of
\(r\), and that 
\(\grq_{i,r}=0\) when \(i\in \grY_m\), as in Assumption
\ref{assume}. (The preceding conditions hold if  \(w_{i,r}\) and \(\grq_{i,r}\) arise as
the \(Y_i\)-end limits of certain 4-dimensional \(w_r\) and
\(\hat{\grp}_r\) that satisfy Assumption \ref{assume}.) The constants
mentioned below will depend on \(w_{i,r}\), \(\grq_{i,r}\) only
through  \(\varsigma_w\), \(\zzz_\grp\), and \([w_{i,r}]\), not on the precise form of
\(w_{i,r}\), \(\grq_{i,r}\). 

The next lemma is more or less standard, and is usually the first step
in typical proofs of compactness results in Seiberg-Witten theory. A short sketch of proof
is nevertheless provided, both in order to keep track of the \(r\)-dependence of
various coefficients, and because some formulas appearing therein 
will be used again later. 

\begin{lemma}\label{lem:3d-Phi}
Let \((Y_i, \nu _i)\),  \(\grq_i=\grq_{i,r}\),
\(i\in \grY\) \(\varsigma _w\), \(\zzz_\grp\) be as described above, and define 
\[
\rma_i:=
\begin{cases} 0 &\text{when \(i\in \grY_v\);} \\
1 & \text{ when \(i\in \grY_m\). }
\end{cases}
\]
Let \(\mu _{i,r}=r\nu_i+w_{i,r}\) and let \((B_i, \Phi_i)=(B_{i,r}, \Phi_{i,r})\) be a solution to the 3-dimensional
Seiberg-Witten equation \(\grF_{\mu_{i,r}, \grq_i}(B_i,
\Phi_i)=0\). 
  There exist an \(r\)-independent constants \(\zeta_i\in  \bbR^+\)
  such that 
\begin{equation}\label{eq:3d-Phi}
\|\Phi_i\|_{L^\infty}\leq\zeta_i \, r^{\rma_i/2}.  
\end{equation}
Let \(k\geq 2\) be as in Assumption \ref{assume}. In the case when \(i\in \grY_v\) and \((B_i,
\Phi_i)\) is in the normalized Coulomb gauge with respect to a
reference connection \(B_{0,i}\),  \(\eta:=(\delta
B_i=B_i-B_{0,i}, \Phi_i)\) satisfies \[
\|\eta\|_{L^2_{k+4,
    B_{0,i}}}+\|\eta\|_{C^{k+2}_{
    B_{0,i}}}+\|\Phi _i\|_{L^4}\leq \zeta'_i\] 
for some
\(r\)-independent positive constant
\(\zeta'_i\). The constant \(\zeta_i\) is 
determined by the metric 
on \(Y_i\) and the norms of the
perturbations via \(\|\nu_i\|_{C^2}\), \(\|w_{i,r}\|_{C^2}\leq \varsigma _w\), and
\(\|\grq_i\|_{\scrP}\leq \zzz_\grp\).  
The constant \(\zeta '_i\) is determined by
all the above, plus \(B_{0,i}\). 
\end{lemma}
\pf The proof hinges on the 3-dimensional Weitzenb\"{o}ck
formula 
\begin{equation}\label{eq:Weit3d}
\slp_{B_i}^2\Phi_i=\nabla _{B_i}^*\nabla_{B_i}\Phi_i+\frac{R_g}{4}
\Phi_i+\frac{1}{2}\rho(F_{B_i})\Phi_i,
\end{equation}
where \(R_g\) denotes the scalar curvature. 
Applying this formula together with (\ref{eq:q-bdd}) to the identity 
\(\|\grF_{\mu _{i,r}}(B_i, \Phi _i)\|_{L^2}^2=\|\grq_i(B_i, \Phi
_i)\|_{L^2}^2 \), one has: 
\[
\begin{split}
&\frac{1}{4}\|F_{B_i}\|^2_{L^2}+\|\nabla_{B_i}\Phi_i\|^2_{L^2}+\big\|\frac{i}{4}\rho(\mu_{i,r})-(\Phi_i\Phi_i^*)_0\big\|^2_{L^2}\\
& \qquad +\int_{Y_i}\frac{R_g}{4}|\Phi_i|^2-\frac{i}{4}\int_{\hat{Y}_{[l,L]}}F_{B_i}\wedge
*\mu _{i,r}\leq 4m^2\|\grq_i\|_\scrP^2\, \big(1+\|\Phi _i\|^2_{L^2}\big).\\
\end{split}
\] 
Judicious applications of the triangle inequality then leads to 
\[\begin{split}
\frac{1}{8}\|F_{B_i}\|^2_{L^2} +\|\nabla_{B_i}\Phi_i\|^2_{L^2} &+\frac{1}{2}\big\|\frac{i}{4}\rho(\mu_{i,r})-(\Phi_i\Phi_i^*)_0\big\|^2_{L^2}\\
& \leq\zeta
\big(\|\mu
_{i,r}\|_{C^0}^2+\|R_g\|_{C^0}^2+m^4\|\grq_i\|_\scrP^4\big)+\zeta
', 
\end{split}
\]
where \(\zeta \), \(\zeta '\) are constants independent of
anything. (For more details, see  (\ref{f-trick}) and  subsequent discussions
for 4-dimensional analogs.) 

According to Assumption \ref{assume}, when \(i\in \grY_v\), \(\|\mu
_{i,r}\|_{C^0}=\|w_{i,r}\|_{C^0}\leq \varsigma_w\) and
\(\|\grq_i\|_\scrP\leq \zzz_\grp\), and so it follows that 
\[
\|F_{B_i}\|^2_{L^2} +\|\nabla_{B_i}\Phi_i\|^2_{L^2} +\|\Phi
_i\|_{L^4}^4\leq \zeta _v,
\]
where \(\zeta _v\) is an \(r\)-independent constant determined by the
metric, \(\varsigma_w\) and \(\zzz_\grp\). With the preceding bound at
hand, a standard elliptic bootstrapping argument (see e.g. the proofs
of \cite{KM}'s Lemma 5.1.2 and Theorem 10.7.1 for the 4-dimensional version) then leads to 
\[
\|\eta\|_{L^2_{k+4,
    B_{0,i}}}+\|\Phi_i\|_{C^{k+2}_{B_{0,i}}}\leq \zeta _v',
\]
where \(\zeta _v'\) is a constant determined by the
metric, \(B_0\), \(\varsigma _w\) and \(\zzz_\grp\). Together with the
preceding bound on \(\|\Phi
_i\|_{L^4}\), this verified the assertions of the lemma in the \(i\in
\grY_v\) case.

In the case when  \(i\in \grY_m\), \(\grq_i=0\) and \(\|\mu_{i,r}\|_{C^2}\leq r\|\nu _i\|_{C^2}+\varsigma_w\). 
The preceding argument provides an \(r\)-dependent bound on \(\|\Phi
_i\|_{L^\infty}\) that grows faster with \(r\) than what we
need. Meanwhile, in this case the right hand side of (\ref{eq:Weit3d})
vanishes by the second line of the Seiberg-Witten equation \(\grF_{\mu
  _{i, r}} (B_i, \Phi _i)=0\). Combining this resulting identity with
  the first line of the  Seiberg-Witten equation, one has a pointwise
  inequality of the form 
\[
d^*d|\Phi_i|^2+|\nabla_{B_i}\Phi_i|^2+\frac{|\Phi_i|^4}{2}\leq \zeta
r|\Phi_i|^2, 
\]
where \(\zeta \) is a constant depending only on the metric, \(\|\nu
_i\|_{C^2}\), and \(\varsigma _w\). This permits
a straightforward application of the strong maximum principle, which
leads to the asserted bound on \(\|\Phi_i\|_{L^\infty}\). \epf


The next lemma is an ingredient for the proof of Theorem
\ref{thm:strong-t}.  

 Suppose \(i\in \grY_m\). Let \(\sigma(x)=\sigma _{Y}(x)\) denote the distance from \(x\in Y_i\) to
\(\nu_i^{-1}(0)\), and let \(Y_{i,\delta}
  \subset Y_i\)
  denote the set of all points with distance greater or equal to
  \(\delta\) from \(\nu_i^{-1}(0)\). Over \(Y_i-\nu_i^{-1}(0)\), the Clifford action of \(\nu_i\)
  splits the spinor bundle \(\bbS\) into a direct sum of
  eigen-bundles, 
  \(\bbS= E\oplus E\otimes K^{-1}\), where \(E\) to be the summand
  corresponding to the eigenvalue \(-i|\nu_i|\), and 
\(K^{-1}=\ker (*\nu _i)\subset TY_i\) is a subbundle of oriented
2-planes over \(Y_i\). The Clifford action by \(\nu _i/|\nu _i|^{-1}\)
endows \(K^{-1} \) with a complex structure. Write
  \[
\phi_i:=(2r)^{-1/2}\Phi_i=(\alpha, \beta)
\] in terms of the
  preceding 
  decomposition. Let \(\nabla_{B_i}\alpha\), \(\nabla_{B_i}\beta\)
  respectively denote the covariant derivatives of \(\alpha\in
  \Gamma(E)\) and \(\beta\in \Gamma(E\otimes K^{-1})\) with respect to
  the connections induced fron \(B_i\) and the Levi-Civita connection.

Let \begin{equation}\label{def:sigma-3d}
\tilde{\sigma}:=\sigma\chi(\sigma)+1-\chi(\sigma).
\end{equation}
Note that
our assumption that \(\nu_i\) has nondegenerate zeros implies that
there exists a positive constant \(\zeta_i\) satisfying 
\begin{equation}\label{bdd:nu-3d}
\zeta_i^{-1}\ts\leq |\nu_i|\leq \zeta_i\ts.
\end{equation}

\begin{lemma}\label{lem:3d-ptws}
Adopt the assumptions and notations in Lemma \ref{lem:3d-Phi}, and 
suppose that \(i\in \grY_m\). There is a positive constant \(r_0>1\) such that for all \(r\geq
r_0\), the following holds: Let \((B_i, \Phi_i)\) be as in the
previous lemma, and \(i\in \grY_m\).  Then:
\begin{itemize}
\item[(1)] There exists an
\(r\)-independent positive constants \(z'\), \(z''\) such that over \(Y_i\),
\begin{equation}\label{eq:Phi-ptws}
\begin{split}
|\phi_i|^2 &\leq |\nu_i|+z r^{-1/3};\\
|\phi_i|^2 &\leq |\nu_i|+z' r^{-1}(\sigma ^{-2}+1).
\end{split}
\end{equation}
\item[(2)] There exist positive constants \(\textsc{o}>100\), \(c\),
  \(c'\), \(\zeta'\) and
  \(\zeta\) that are independent of \(r\) and \((B, \Phi )\), such that 
  for \(\delta_0:=\sO r^{-1/3}\), 
the following pointwise estimates
  hold on \(Y_{i,\delta_0}\):
\begin{eqnarray}
\label{ineq:3d-beta}
&& |\beta|^2\leq c \ts^{-3}r^{-1} (
|\nu |-|\alpha|^2)+\zeta \ts^{-5} r^{-2};\\
&&
|\beta|^2\leq c' \ts^{-3}r^{-1} (
|\nu |-|\phi_i|^2)+\zeta '\ts^{-5} r^{-2}.
\label{ineq:3d-beta'}
\end{eqnarray}
\end{itemize}
All the constants above (\(\zeta\), \(z\), \(c\), \(\sO\) \(r_0\) etc.) are determined
by \(\|\nu_i\|_{C^1}\), \(\varsigma_w\), and the metric on \(Y_i\).
\end{lemma}

That \(r\) is no less than the constant \(r_0\) above will be a standing {\em assumption throughout} the rest
of this article. 

\pf 
The estimates
(\ref{eq:Phi-ptws}) are  the
3-dimensional analogs of Lemma 3.2 in \cite{Ts} and follow from the
same proof.  Similarly, a straightforward modification of the proof of
Proposition 3.1 in \cite{Ts}, with  the
dependence of constants therein taken into account, yields the
following: 

There exist positive constants \(r_0\), \(\textsc{o}\), \(c\),
  \(c'\), \(\zeta'\) and
  \(\zeta\), which depend only on by \(\|\nu_i\|_{C^1}\),
  \(\varsigma_w\), and the metric on \(Y_i\), 
 such that for any \(r\geq r_0\) and 
 \(\delta\geq\sO r^{-1/3}/2=\delta_0/2\), 
the following pointwise estimate
  holds on \(Y_{i,\delta}\):
\[
 |\beta|^2 \leq \frac{c}{8} \delta^{-3}r^{-1} (
|\nu |-|\alpha|^2)+\zeta \,( 2\delta)^{-5} r^{-2}.
\]
Now, for a fixed \(x\in Y_{i,\delta_0}\), we can apply the
preceding formulae for any choice of 
\(\delta\geq\delta_0/2\), and thus obtaining infinitely many bounds on
the value \(|\beta(x)|^2\). 
 The inequality (\ref{ineq:3d-beta}) is given by the
(\(x\)-dependent) choice by setting \(\delta=\ts/2\). (Note that for
any \(x\in Y_{i,\delta_0}\), \(\ts/2\geq\delta_0/2\) and so it is an
admissible choice of \(\delta\)). 

The inequality (\ref{ineq:3d-beta'}) is a direct
consequence of (\ref{ineq:3d-beta}). \epf


\begin{lemma}\label{lem:F-L_1}
Adopt the assumptions and notations in Lemma \ref{lem:3d-Phi}, and 
when \(i\in \grY_m\), suppose  in addition that \(r\geq r_0\), where  \(r_0>1\) is  as
in Lemma \ref{lem:3d-ptws}.
 There exist an \(r\)-independent constant \(\zzz_i\in \bbR^+\)
  such that 
\[
\|F_{B_i}\|_{L^1(Y_i)}\leq \zzz_i.
\]
In the case when \(i\in \grY_v\),  the constant
\(\zzz_i\) above is determined by the metric on \(Y_i\),
\(\varsigma_w\),
 and \(\|\grq_i\|_{\scrP}\leq \zzz_\grp\). In the case when \(i\in \grY_m\), the constant
\(\zzz_i\) above is  determined by the metric on \(Y_i\),
\(\|\nu_i\|_{C^1}\), \([w_{i,r}]\), 
\(\varsigma_w\), and the \(\Spin^c\) structure \(\grs_i\) through the quantity \[\wp_i:=
\frac{\pi}{2}
c(\grs_i)\cdot[*\nu_i].
\]
\end{lemma}
\pf 
When \(i\in \grY_v\), by the 3-dimensional Seiberg-Witten
equation \(\grF_{\mu_{i,r, \grq_i}}(B_i, \Phi_i)=0\), one may bound 
\(\int_{Y_i}|F_{B_i}|\leq \zeta_1 (\|\Phi_i\|_{L^\infty}^2+\varsigma_w)+\zeta_2
\|\grq_i(B_i, \Phi_i)\|_{L^2}\), where \(\zeta _1, \zeta _2\) depend
only on the metric of \(Y_i\). The asserted \(L^1\)-bound on
\(F_{B_i}\) then follows from Lemma \ref{lem:3d-Phi} and
(\ref{eq:q-bdd}).

Suppose now that \(i\in \grY_m\). The proof in this case takes four
steps. 

\subsubsection*{\it Step 1.} This step derives some pointwise bounds on
\(|F_{B_i}|\). By direct computation: 
\begin{eqnarray}
|2F_{B_i}-iw_{i,r}|& =&r\big((|\nu_i|-|\alpha |^2+|\beta |^2)^2+4|\alpha |^2\,
|\beta|^2\big)^{1/2}\nonumber\\
& =& r \big((|\nu_i|-|\phi_i|^2)^2+4|\nu_i|\, |\beta|^2\big)^{1/2}\nonumber\\
&\leq   & 
r\Big||\nu_i|-|\phi_i|^2\Big|+2 r|\nu_i|^{1/2} |\beta| \label{17.1}\\
&\leq   & r\Big||\nu_i|-|\phi_i|^2\Big|+\zeta
'r^2\ts^3|\beta|^2+\zeta''\ts^{-2}\quad \text{ over \(Y_i\),}\label{17.2}
\end{eqnarray}
where $\zeta '$, $\zeta ^{''}$ are constants depending only on $\nu_i$. Thus, for $r\geq r_0$, one has: 
\begin{equation}\label{3d-F}
 |F_{B_i}|\leq \begin{cases}\zeta_f\, r\ts+\zeta'_f \, \ts^{-2} &\text{ over
     \(Y_i\) by (\ref{17.1}) and (\ref{eq:Phi-ptws});}\\
\zeta_0 r\Big|
|\nu_i|-|\phi_i|^2\Big|+\zeta'_0\ts^{-2}&\text{ over \(Y_{i, \delta
    _0}\) by (\ref{17.2}) and (\ref{ineq:3d-beta'}).}
\end{cases}
\end{equation}
In the above, $r_0$ and \(\delta_0\) are as in Lemma
\ref{lem:3d-ptws}, and $\zeta _f$, $\zeta_f'$, $\zeta _0$, $\zeta_0'$
are constants depending only on $\nu_i$, $\varsigma _w$, and the metric. 

\subsubsection*{\it Step 2.} Use the first line of 
Seiberg-Witten equation \(\grF_{\mu_{i,r}}(B_i, \Phi_i)=0\) to write: 
\[
\begin{split}
\frac{r}{2} & \int_{Y_i}\big(|\nu_i|\, (|\nu_i|-|\alpha|^2+|\beta|^2)\big) 
 \\& \quad  =\int_{Y_i}(iF_{B_i}-w_i/2)\wedge*\nu_i\\&\quad =\big(2\pi
c(\grs_i)-[w_i]/2\big)\cdot[*\nu_i]=4\wp_i-[w_i]\cdot [*\nu_i]/2.
\end{split}
\]

 Since \(r\int_{Y_i}\big(|\nu_i|(|\phi_i|^2-|\nu_i|)_+\big)\leq \zeta\) by
Lemma \ref{lem:3d-ptws} (1), the preceding formula implies that
\begin{equation}\label{eq:3d-E_t}
r\int_{Y_i}|\nu_i|\Big|
|\nu_i|-|\phi_i|^2\Big|+2r\int_{Y_i}|\nu_i||\beta|^2\leq \zeta_c,
\end{equation}
where \(\zeta _c\) depends only on the metric, \(\nu _i\), \(\varsigma
_w\), \([w_i]\), and \(\wp_i\). 


\subsubsection*{\it Step 3.} While (\ref{eq:3d-E_t}) is useful for bounding
\(r\int_{Y_{i, \rho }}\Big|
|\nu_i|-|\phi_i|^2\Big|\) for \(\rho \geq\delta _0\), the bound is
hardly efficient when \(\rho \) is small. This step derives a better
integral bound for \(r\Big|
|\nu_i|-|\phi_i|^2\Big|\) near the zeros of \(\nu _0\). 

By assumption, the zero locus \(\nu _i^{-1}(0)\) is nondegenerate;
therefore it consists of finitely many points \(\nu
_i^{-1}(0)=\bigcup_pp\), and there exists  
a small positive \(\rho _0<1\) and a function \(f\) over \(Y_i-Y_{i,
  \rho _0}\), such that \(Y_i-Y_{i, \rho _0}\) consists of mutually
disjoint balls \(B_p(\rho _0)\), and over \(Y_i-Y_{i, \rho _0}=\bigcup_pB_p(\rho _0)\),
\(df=*\nu_i\) and \(|f|\leq \frac{|\nu |\sigma }{2}+O(\sigma^3)\). Recall also from
(\ref{bdd:nu-3d}) that \(\zeta_i\sigma\geq |\nu_i|\geq\zeta
_i^{-1}\sigma\) over \(Y_i-Y_{i, \rho _0}\). 

Let \(p\in\nu_i^{-1}(0)\), and let \(\scrW_p(\rho ):=\displaystyle\frac{r}{2}\displaystyle\int_{B_p(\rho)}  \Big(|\nu_i|\Big||\nu_i|-|\phi_i|^2\Big|\Big)\). 
The first line of the
Seiberg-Witten equation \(\grF_{\mu _r}(B_i, \Phi _i)=0\), together with the second line of
(\ref{3d-F}), implies that for  \(\delta _0\leq \rho \leq\rho _0\) and \(r\geq r_0\),  
\begin{equation}\label{22}
\begin{split}
\scrW_p(\rho )-\zeta _1\rho ^2&\leq \frac{r}{2}\int_{B_p(\rho)}  \big(|\nu_i|(|\nu_i|-|\phi_i|^2+2|\beta|^2)\big)\\
&  =\int_{\partial B_p(\rho)}f\,( iF_{B_i}-w_{i,r}/2)\\
& \leq \frac{\rho (1+\zeta '\rho )}{2}\frac{d}{d\rho }\scrW_p(R)+\zeta_2 \rho^2,\\ 
\end{split}
\end{equation}
where \(\zeta _1\), \(\zeta _2\), \(\zeta '\)
are constants depending only on the metric, \(\nu_i\) and \(\varsigma
_w\). Rewrite the preceding differential inequality as:
\[
  \frac{d}{d\rho }\big(e^h \scrW_p(R)\big)\geq-\zeta \rho , \quad
  \text{where \(h:=-2\ln \big(\frac{\rho }{1+\zeta '\rho }\big)\),} 
  \]
  and integrate from \(\rho \) to \(\rho _0\). One has:
  \[
    \rho ^{-2}_0\scrW_p(\rho _0)-\rho ^{-2}\scrW_p(\rho )\geq-\zeta _3(\rho _0^2-\rho ^2).
  \]
  Thus, from (\ref{eq:3d-E_t}) one has
  \begin{equation}\label{22.1}
    \scrW_p(\rho )\leq \zeta _4\rho ^2.
  \end{equation}

let \(N\in \bbZ^{\geq 0}\) be the largest  integer such that \(\delta _0':=2^{-N+1}\rho
_0\geq\delta _0\). (Therefore \(2\delta _0\geq\delta _0'\geq \delta
_0\).) Write \(A_n:=A (2^{-n+1}\rho _0, 2^{-n}\rho
_0)\). Then it follows from (\ref{22.1}) that there is a constant
\(\zeta ''\) depending only on  the metric, \(\nu_i\) and \(\varsigma
_w\), such that for all \(n, r,\) satisfying \(0\leq n\leq N\) and
\(r\geq r_0\), \(r\int_{A_n}\Big| |\nu_i|-|\phi_i|^2\Big|\leq\zeta  ''\,
2^{-n}\). Consequently 
\begin{equation}\label{3d-monotone}
r\int_{Y_{i,
    \delta _0'}-Y_{i, \rho _0}}\Big| |\nu_i|-|\phi_i|^2\Big|\leq\zeta_b
\end{equation}
for a constant \(\zeta _b\) depending only on  the metric, \(\nu_i\) and \(\varsigma
_w\).

\subsubsection*{\it Step 4.} 
Divide \(Y_i\) into the
three regions: \(Y_i-Y_{i,\delta _0'}\), \(Y_{i, \delta _0'}-Y_{i, \rho
  _0}\), and \(Y_{i,\rho _0}\). Use the first line of the pointwise
bounds (\ref{3d-F}) 
to estimate the \(L^1\) norm of \(F_{B_i}\) on the first region. Use the
second line of the pointwise bounds (\ref{3d-F}) for both the second and the third region;
and combine with the integral bounds (\ref{3d-monotone}) and (\ref{eq:3d-E_t})
respectively for the second and the third region. We have 
\[\begin{split}
\int_{Y_i}|F_{B_i}|& =\int_{Y_i-Y_{i, \delta _0'}}|F_{B_i}|+\int_{Y_{i,
    \delta _0'}-Y_{i, \rho _0'}}|F_{B_i}|+\int_{Y_{i, \delta _0'}}|F_{B_i}|
\leq 
\zzz_i, 
\end{split}
\]
where  
\(\zzz_i\) is a constant depending only on \(\nu _i\), \(\varsigma _w\), and the
metric. \epf

\begin{lemma}\label{lem:CSD-est}
Adopt the assumptions and notations in Lemma \ref{lem:3d-Phi}. When
\(i\in \grY_m\), assume in addition that \(r\geq
r_0\), where \(r_0\) is as in Lemma \ref{lem:F-L_1}. 
There 
 exist  \(r\)-independent constants \(\zeta_i,  \zeta _i', z'_i\in \bbR^+\) with the
following significance: 

Suppose that \((B_i, \Phi_i)=[(B_i, \Phi_i)]_c\) 
is  {\em in the
normalized Coulomb gauge} with respect to \(B_{0,i}\).  Then: 
\begin{itemize}
\item[(a)] One has
the following bounds for \(\delta B_i:=B_i-B_{i,0}\):
\[\text{ \(\|\delta B_i\|_{L^1(Y_i)}\leq \zeta_i\) and \(\|\delta
  B_i\|_{L^\infty(Y_i)}\leq \zeta'_i r^{2\ra_i/3}.\) }
\]

\item[(b)] The values of various versions of \(\op{CSD}\) functionals
  at \((B_i, \Phi _i)\) are
  bounded as follows: 
\[
\begin{split} |\op{CSD}_{0}& (B_i, \Phi_i)| +|\op{CSD}_{w_{i,r}}(B_i,
  \Phi_i)|+|\op{CSD}_{w_{i,r}, \grq_i}(B_i, \Phi_i)|\\
& +r^{-\ra_i/3}|\op{CSD}_{\mu
    _{i,r}}(B_i, \Phi_i)|+r^{-\ra_i/3}|\op{CSD}_{r\nu_i}(B_i, \Phi_i)|\leq z_i'r^{2\ra_i/3}.
\end{split}
\]
\end{itemize}

In the case when \(i\in \grY_v\),  the constants
\(\zeta _i, \zeta'_i, z'_i\) above are determined by the metric on \(Y_i\),
\(\varsigma_w\),
 and \(\|\grq_i\|_{\scrP}\leq\zzz_\grp\). In the case when \(i\in \grY_m\), the constants
\(\zeta _i, \zeta'_i, z'_i\) above are  determined by the metric on \(Y_i\),
\(\|\nu_i\|_{C^1}\), 
\(\varsigma_w\), \([w_i]\), and \(\Spin^c\) structure \(\grs_i\) through the quantity \(\wp_i\).
\end{lemma}

\pf 
{\bf (a)} When \(i\in\grY_v\), the asserted \(L^\infty\)-bound on \(\delta
B_i\) follows from Lemma \ref{lem:3d-Phi}. The asserted \(L^1\)-bound on \(\delta  B_i\)
follows from Lemma \ref{lem:3d-Phi} together with 
the fact that \(\|\delta  B_i\|_{L^1}\leq \zeta \|\delta
B_i\|_{L^2}\), where \(\zeta \) is a constant depending only on the
metric of \(Y_i\). 

When \(i\in \grY_m\), combine the previous lemma and 
 (\ref{eq:deltab1}) to get: 
\[\begin{split}
\|\delta B_i\|_{L^1(Y_i)}& \leq \zeta
\|G\|_{L^1(Y_i)}(\|F_{B_i}\|_{L^1(Y_i)}+1)\\
& \leq \zeta _i,
\end{split}
\]
where \(G\) denotes the Green's function for \(d+d^*\) on \(Y_i\). In
the above, \(\zeta \) and \(\zeta _i\) are \(r\)-independent constants:
\(\zeta \) depends only on the metric of \(Y_i\);
\(\zeta _i\) is  determined by the metric on \(Y_i\),
\(\|\nu_i\|_{C^1}\), \([w_{i,r}]\), 
\(\varsigma_w\), and \(\wp_i\). 

To get the  
asserted \(L^\infty\) bound of \(\delta  B_i\), notice first that
Lemma \ref{lem:3d-Phi} together with the Seiberg-Witten equation
\(\grF_{\mu _{i,r}}(B_i, \Phi _i)=0\) imply that 
\begin{equation}\label{3d-F-inf}
  \|F_B\|_{L^\infty}\leq \zeta  r^{\ra_i},
\end{equation}
where \(\zeta \) is a constant depending only on \(\|\nu
_i\|_{L^\infty}\) and \(\varsigma _w\). 
The assumption that \((B_i, \Phi _i)\) is in the normalized Coulomb
gauge implies that \(\delta  B_i=b_h+b\), where \(b_h\) is harmonic
and \(b\) is coexact. Moreover,
\begin{equation}\label{bdd:b_h}
  \|b_h\|_{L^\infty(Y_i)}\leq \zeta
  _h,
\end{equation}
where \(\zeta _h\) only depends on the first betti number of
\(Y_i\). Meanwhile, using the fact that the Green's function satisfies
\(|G(x,y)|\leq \zeta _g\dist (x, y)^{-2}\),  one has for a given \(x\in Y_i\),
\[\begin{split}
  |b(x)| & \leq \zeta _g'\int_{Y_i} \dist (x, \cdot )^{-2}
  \big(|F_{B_i}|+|F_{B_{i,0}}|\big)\\
     &\leq \zeta _g' \int_{B_x ( \rho )}\dist (x,
     \cdot)^{-2}|F_{B_i}|\, +\zeta _g' \int_{Y_i-B_x ( \rho )}\dist (x, \cdot)^{-2}|F_{B_i}| \, +\zeta '.
  \end{split}
\]
Using (\ref{3d-F-inf}) for the first integral in the last line above,
and applying Lemma \ref{lem:F-L_1} to the second integral therein, we
have:
\[
|b(x) |\leq \zeta _1 r^{\ra_i}\rho +\zeta _2\rho ^{-2}+\zeta ', 
  \]
where \(\zeta _1\) only depends on the metric, \(\|\nu
_i\|_{L^\infty}\) and \(\varsigma _w\), the constant \(\zeta _2\) only depends on
the metric and the constant \(\zzz_i\) in Lemma \ref{lem:F-L_1}, and
\(\zeta '\) depends only on the metric and the choice of \(B_{i,0}\).
Now take \(\rho =r^{\ra_i/3}\) in the preceding inequality, and
combine with (\ref{bdd:b_h}), we have \(\|\delta
B_i\|_{L^\infty(Y_i)}\leq \zeta _i'r^{2\ra_i/3}\) as asserted.

{\bf (b)} When \(i\in \grY_v\), the asserted bounds on
\(\op{CSD}_0=\op{CSD}_{r\nu _i}\) and \(\op{CSD}_{w_i}=\op{CSD}_{\mu
  _{i,r}}\) follows directly from Lemma \ref{lem:3d-Phi}.  The bound
on \(\op{CSD}_{w_i, \grq_i}\) follows from the preceding bound on
\(\op{CSD}_{w_i}\), together with a bound on \(|f_{\grq_i}(B_i, \Phi
_i)|\). This bound is provided by combining 
Lemma \ref{lem:f_q} and Lemma \ref{lem:3d-Phi}.

Now suppose \(i\in \grY_m\). Use the 3-dimensional Seiberg-Witten equation, \(\grF_{\mu_{i,r}}(B_i,
\Phi_i)=0\),  Part (a) of the lemma and Lemma \ref{lem:F-L_1}
to get: 
\begin{equation}\label{eq:c1}
\begin{split}
|\op{CSD}_{\mu
  } & (B_i, \Phi_i)|  =
\Big|-\frac{1}{8}\int_Y \delta B_i\wedge
\big(F_{B_{0,i}}+F_{B_i}+i\mu\Big)\Big| \\
&\leq  \frac{1}{8}\|\delta
B_i\|_{L^\infty(Y_i)}\|F_{B_{0,i}}+F_{B_i}\|_{L^1(Y_i)}
+\frac{1}{8}\|\delta B_i\|_{L^1(Y_i)} \|\mu\|_{L^\infty(Y_i)}\\
&\leq \zeta _1 r^{2\ra_i/3}+\zeta _2 \|\mu\|_{L^\infty(Y_i)}, 
\end{split}
\end{equation}
where \(\zeta _1\), \(\zeta _2\) are constants determined by the
constants \(\zeta _i\), \(\zeta _i'\) from Part (a), and the constant
\(\zzz_i\) in Lemma \ref{lem:F-L_1}. 
Setting   \(\mu\) above to be \(0\), \(w_{i,r}\), \(\mu_{i, r}\), and
\(r\nu_i\),  we arrive at  the asserted bounds for the values of \(\op{CSD}_0\), \(\op{CSD}_{w_i}=\op{CSD}_{w_i,\grq_i}\), \(\op{CSD}_{\mu_{i,r}}\) and \(\op{CSD}_{r\nu _i}\) respectively. \epf


 \begin{rem}\label{rmk:CSD-change}
(1) Recall that any  \((B_i, \Phi_i)\in \scrC(Y_i)\) 
may be written as \((B_i, \Phi_i)=u_i\cdot [(B_i, \Phi
_i)]_c\) for a gauge transformation \(u_i\in \scrG\). 
The values of various \(\op{CSD}\) functionals at \((B_i, \Phi _i)\)
may be expressed in terms of those at \([(B_i, \Phi
_i)]_c\) and  the cohomology class of \(u_i\), \([u_i]:=\frac{1}{2\pi}[-i u_i^{-1}
du_i]\in H^1(Y_i; \bbZ)\), via the following identity: 
\[\begin{split}
& \op{CSD}_{\mu, \grq_i}(B_i, \Phi_i)-\op{CSD}_{\mu, \grq_i}([(B_i,
\Phi_i)]_c)\\
& \qquad =\op{CSD}_{\mu}(B_i, \Phi_i)-\op{CSD}_{\mu}([(B_i,
\Phi_i)]_c)\\
& \qquad =2\pi\big(\pi  c_1(\grs_i)-\frac{[\mu]}{4}\big)\cdot [u_i], 
\end{split}
\]
In particular, when \([w_{i,r}]=4\pi  c_1(\grs_i)\), 
the bound for \(|\op{CSD}_{w_{i,r}}(B_i,
  \Phi_i)|+|\op{CSD}_{w_{i,r}, \grq_i}(B_i, \Phi_i)|\) given in Part
  (b) of the preceding lemma holds for \((B_i, \Phi _i)\) {\em in
    arbitrary gauge}. 

 (2) Via a more involved argument, the \(L^\infty\) bound for \(\delta  B_i\) in the preceding lemma
 may be improved as \(\|\delta  B_i\|_{L^\infty (Y_i)}\leq \zeta '_i
 r^{\ra_i/2}\), and accordingly, the bounds for the CSD functionals
 may be improved as  \(|\op{CSD}_0(B_i,
  \Phi_i)|+|\op{CSD}_{w_{i,r}}(B_i,
  \Phi_i)|+|\op{CSD}_{w_{i,r}, \grq_i}(B_i, \Phi_i)|\leq z_i'
  r^{\ra_i/2}\).  This argument makes use of 3-dimensional analogs
  (which are simpler and stronger) of Lemma \ref{lem:loc-mod} and
  Proposition \ref{prop:exp-decay} below. 
\end{rem}

\subsection{t-convergence of large-\(r\) Seiberg-Witten solutions}\label{sec:t-conv}

Let \((\scrX,\omega)\) be a symplectic 4-manifold
with a riemannian metric
\(g\) and an almost complex structure \(J\) so that \(|\omega |\, g=\omega (\cdot,
J\cdot)\). Let \(\mu_r^+\) be defined in terms of \(\omega\) as in
Equation (\ref{eq:mu-pert}). Fix a \(\Spin^c\)-structure on \(\scrX\) and
let \(\bbS=\bbS^+\oplus\bbS^-\) be the associated spinor bundle. The
Clifford action of \(\omega\) splits \(\bbS^+\) into a direct sum of
eigen-bundles, \(E\oplus E\otimes K^{-1}\), where \(E\) corresponds to
the eigenvalue \(-|\omega|i\) 
and \(K^{-1}\) is the anti-canonical
bundle. Given a \(\Spin^c\) connection \(A\) on \(\bbS^+\), let \(A^E\in
\Conn (E)\) denote
the connection that it induces on \(E\). Conversely,
with the riemannian metric on \(\scrX\) fixed, an  \(A^E\in
\Conn (E)\) uniquely determines
a \(\Spin^c\) connection on \(\bbS^+\), and we have an isomorphism 
\(\Conn (\bbS^+)\simeq\Conn (E)\). 

\begin{defn}\label{defn:t-convergence}
Let \((\scrX, \omega)\), \(J\) be as above. Let \(\mathbf{C}:=[C,
\td{C}]\) be a weighted \(J\)-holomorphic subvariety in \(\scrX\).  
Let \(\scrX_1\subset\cdots\scrX_n\subset \scrX_{n+1}\subset \cdots
\subset \scrX\) be a countable exhaustion of \(\scrX\) by open subsets with compact closure,
and let \(\{(A_n, \Psi_n)\}_n\) be a corresponding sequence in
\(\Conn (\bbS^+)\times \Gamma (\bbS^+)\), where \((A_n, \Psi_n)\) is
defined over \(\scrX_n\). This sequence is said to {\em
  t-converge} to \(\mathbf{C}\) if the
following hold for any compact subset \(\mathcal{K}\subset X-\omega^{-1}(0)\), 
\begin{itemize}
\item[(1)] \label{eq:dist-conv} 
Let \(\Psi^E_n\) denote the \(E\)-component of \(\Psi_n\). Then
\((\Psi^E_n)^{-1}(0)\) are closed sets of Hausdorff dimension 2, and \(\lim_{n\to \infty}
  \dist_{\mathcal{K}} \big(C\cap \mathcal{K}, (\Psi^E_n)^{-1}(0)\cap
  \mathcal{K}\big)=0\);
\item[(2)] \(\{\frac{i}{2\pi }F_{A_n^E}\big|_\mathcal{K}\}_n\) converges weakly as
  currents to \(\tilde{C}|_\mathcal{K}\).
\end{itemize}
\end{defn}

Note that the definition of t-convergence above only depends on the
gauge equivalence class of \((A_n,\Psi_n)\). Therefore, we shall often
refer to the ``t-convergence'' of gauge equivalence classes of configurations. 

The notion of t-convergence is similarly defined in the
3-dimensional case: Let \(M\) be a \(\Spin^c\) 3-manifold and
\(\theta\) be a nowhere vanishing harmonic 1-form on \(M\). Replace
\(\omega\) in the preceding definition by \(*\theta\), and replace \(C_k\) by
orbits of the dual vector to \(\theta\). Let \(\Psi_n\) be sections of
the spinor bundle \(\bbS\) instead, and assume that
\((\Psi^E_n)^{-1}(0)\) now have Hausdorff dimension 1. Here, \(E\)
denotes the first summand in \(\bbS=E\oplus E\otimes K^{-1}\), where
\(K^{-1}\) is the 2-plane bundle \(\ker\theta\subset TM\), equipped
with a complex structure given by Clifford multiplication by \(|\theta
|^{-1}\theta\). 

In either the 3-dimensional case \((M, \theta )\) or the 4-dimensional
case \((\scrX, \omega)\) described above,  \(A^K\) will denote the
connection on \(K^{-1}\) induced from the Levi-Civita connection.
For our applications, \((\scrX, \omega)\)  typically takes the form 
\(((X_\bullet \cap X_0)^\circ, 2\nu ^+)\) for an admissible pair \((X,
\nu)\). (Recall from Section \ref{sec:convention} the definition of
\(X_0:=X-\nu ^{-1}(0)\) and \(X_\bullet\).)  In the 3-dimensional case, typically 
\(M=Y-\theta^{-1}(0)\) for a harmonic Morse 1-form on a closed
\(\Spin^c\) 3-manifold \(Y\). In the former case, we say that a
sequence of configurations on \(X\) {\em t-converges to a t-curve}
\(\mathbf{C}\) in \(X_\bullet\) if it t-converges on the interior of
\(X_\bullet \cap X_0\). In the latter case, 
a sequence of configurations on \(Y\) is said to {\em t-converge to a
  t-orbit \(\pmb{\gamma}\)} if it t-converges on \(Y-\theta^{-1}(0)\).
In the 3-dimensional case, the assumption that the zero locus of
\(\theta \) is nondegenerate implies the existence of a constant
\(\zeta \) such that  \(|F_{A^K}|\leq\zeta  \ts^{-2}\), and therefore
\(\|F_{A^K}\|_{L^1(Y)}\leq \zeta '\) is finite. In the 4-dimensional
case, the assumption that 
\(\nu^{-1}(0)\) is non-degenerate, together with the asymptotic
condition of \(\nu \), imply a similar pointwise bound for
\(|F_{A^K}|\) and that \(\|F_{A^K}\|_{L^1(X_\bullet)}\leq \zeta '\,
|X_\bullet|\) for a constant \(\zeta '\) depending only on the metric and \(\nu
\). (Cf. Section \ref{sec:convention} for the definition of \(|X_\bullet|\).)

What follows is a restatement of Theorem \ref{thm:strong-t} in weaker form. 
\begin{prop}\label{prop:t-conv3d}
Adopt the notations and assumption in Theorem \ref{thm:strong-t}.
Then there exists a t-orbit \(\pmb{\gamma}\) on \(Y-\theta^{-1}(0)\) and a subsequence
of the gauge equivalence classes \(\{[(B_r, \Phi_r)]\}_r\) which t-converges to \(\pmb{\gamma}\). 

Furthermore, \(
\int_{\tilde{\gamma}}\theta\geq 0\) and equals 
\(\frac{1}{2}(c_1(\grs)\cdot[\theta]-\zeta_\theta)\), where
\[
\zeta_\theta=\int_{Y}\frac{i}{2\pi}F_{A^K}\wedge\theta,
\]
and \(A^{K}\) is the connection on \(K^{-1}\) induced by the
Levi-Civita connection. In
particular, sequences \(\{(B_r,
\Phi_r)\}_r\) satisfying the assumption of Theorem \ref{thm:strong-t}
exist only when \(c_1(\grs)\cdot[\theta]\geq
\zeta_\theta\). 
\end{prop}
\pf As mentioned in the introduction, the claim about t-convergence is
a consequence of \cite{Ts}'s Theorem \ref{thm:T1}. Alternatively, it can
be proved by going through the proof in \cite{Ts}, which
simplifies significantly in the 3-dimensional setting. (See \cite{LT}
Section 3 for the case when \(\theta^{-1}(0)=\emptyset \)). As the constituent
flow lines of \(\pmb{\gamma}\) are flow lines of \(-\check{\theta}\),
\(\int_{\tilde{\gamma}}\theta\geq 0\). By the t-convergence we also
have
\[\int_{\tilde{\gamma}}\theta=\lim_{r\to\infty}\int_{Y}\frac{i}{2\pi}F_{B^E_r}\wedge\theta=\frac{1}{2}(c_1(\grs)\cdot[\theta]-\zeta_\theta),\]
where \(B^E_r\) is the connection on \(E\) induced from \(B_r\).  \epf

\begin{remarks}\label{rmk:thurston}
Note the cohomological nature of \(\zeta_\theta
\). Remove small radius \(\delta \) balls centered at points in
\(\theta^{-1}(0)\), and denote the resulting 3-manifold with boundary
\(Y_\delta\). Then \(\frac{i}{2\pi}F_{A^K}\) defines a class \(c_1(K^{-1})\in
H^2(Y_\delta;\bbZ)\), while \(\theta\) defines a class \([\theta]\) in
\(H^1(Y_\delta, \partial Y_\delta; \bbR)\simeq
H_2(Y_\delta;\bbR)\). \(\zeta_\theta \) is then the pairing \(
\zeta_\theta =c_1(K^{-1})\cdot [\theta]\). In the case when \(\theta
=df\) for a circle-valued Morse function \(f\co Y\to \bbR/\bbZ\),
\(\zeta _\theta \) is the maximal  Euler characteristic among the regular
level surfaces of \(f\). Note that the latter is invariant under
homotopy of \(f\). Thus, \(-\zeta _\theta \geq \|\theta \|_T\). In
fact, \(-\zeta _\theta \) often coincides with \(\|\theta \|_T\), such
as in the case when \(\theta \) has no zeros. We believe that when
\([\theta ]\) is an integral primitive class with its Poincar\'e dual
represented by a connected embedded minimal \(\chi _-\) surface, \(-\zeta _\theta =
\|\theta \|_T\). In the above, \(\chi _-\) is as defined in \cite{Th}.
\end{remarks}

\subsection{Strong t-convergence }\label{sec:strong-t}

Fix a closed connected \(\Spin^c\) 3-manifold \((Y, \grs)\) and let \(\bbS\)
denote its associated
\(\Spin^c\) bundle. Let \(\Conn=\Conn (\bbS)=:\Conn (Y)\) denote the
space of \(\Spin^c\)-connections on \(\bbS\), and let
\(\scrC:=\{(A, \Psi)\}=\Conn (Y)\times \Gamma(Y, \bbS)\). Let
\(\scrG = C^\infty (Y, S^1)\) be the gauge group. We use the short
hand 
\(\scrB=\scrB_Y=\scrB (Y)=\scrC/\scrG\). 
Let  \(
\Pi\co \scrC\to \Conn 
\)
be the projection map taking \((B, \Phi
)\in \scrC\) to 
\(B\in \Conn\). It induces a map, also denoted \(\Pi \), from \(\scrB\) to 
\(\Conn /\scrG\). 

\begin{defn}
A sequence in \(\scrB(Y)\) is said to be {\em strongly t-convergent} if
it is t-convergent and its image under \(\Pi\) in \(\Conn (Y)/\scrG(Y)\) is
convergent in the ``current topology'', to be specified below.
\end{defn} 

The aforementioned topology is weaker than the customary \(C^\infty\)
topology, and requires knowledge of some basic structure of
\(\Conn /\scrG\) to describe. To begin, note that \(\pi _0(\scrG)=H^1(Y;\bbZ)\) is the
fundamental group of both \(\scrB\) and \(\Conn /\scrG\).  Let
\(\scrG_0\subset \scrG\) denote the identity component. The
universal cover of \(\Conn /\scrG\), \(\Conn /\scrG_0\), is an affine space under \(
\Omega^1(Y)/ (d\Omega^0(Y))\). The latter  has a bundle structure
by way of the following short exact sequence:
\begin{equation}\label{eq:ex}
0\to H^1(Y;\bbR)\to \Omega^1(Y)/ (d\Omega^0(Y))\stackrel{d}{\to }
d\Omega^1(Y)\to 0. 
\end{equation}
The deck transformations by \(H^1(Y;\bbZ))\) on \(\Conn /\scrG_0\) preserves the
fibers, \(H^1(Y;\bbR)\), 
of this fibration; it acts freely on \(H^1(Y;\bbR)\) by
addition. (Recall that \(H^1(Y;\bbZ)\) is torsion-free,
cf. e.g. \cite{D}'s Exercise 7.22, and we may and will  identify it  with the
integral lattice in \(H^1(Y;\bbR)\).) Consequently, \(\Conn /\scrG\) inherits the structure
of a torus bundle: 
\begin{equation}\label{eq:CD1}
\begin{CD}
H^1(Y;\bbR)/H^1(Y;\bbZ)@>>>\Conn /\scrG\\
@. @V\op{d} VV\\
@. \op{B}^2(\grs),
\end{CD}
\end{equation}
where \(\op{B}^2(\grs)\) is the affine space under \(\op{B}^2:=d\Omega^1(Y)\)
consisting of all closed 2-forms with cohomology class
\(c_1(\grs)/2\in H^2(Y; \bbZ)\) modulo torsions, and \(\op{d}\) is the
map sending the gauge-equivalence class \([B]\in \Conn /\scrG \)
to the closed 2-form \(\frac{i}{4\pi }F_B\).  Note that by
Poincar\'e duality, the space \(H^1(Y;\bbR)/H^1(Y;\bbZ)\) is canonically
isomorphic to  the
torus denoted as \(\bbT_{Y}\) in Section
\ref{sec:h-class}; and we shall denote it by the same notation. 
With the preliminaries out of the way, the {\em current topology
on \(\Conn /\scrG\)} is that induced by the current topology on the
base space, \(\op{B}^2(\grs)\), of (\ref{eq:CD1}), together with the
standard topology on the fiber, \(\bbT_Y\), as a
\(b^1(Y)\)-dimensional torus. Let \(\op{C}(Y)\),
\(\tilde{\op{C}}(Y)\) respectively denote the completion
of \(\Conn /\scrG \) and \(\Conn /\scrG _0\) (with \(C^\infty\)
topology) with respect to the current topology, and let
\(\scrZ^\bbR(\grs)\supset\op{B}^2(\grs)\) denote the space of closed
1-currents with cohomology class
\(c_1(\grs)/2\in H^2(Y; \bbZ)\) modulo torsions. This space is an
affine space under \(\scrZ^\bbR(Y)\), the space of exact
1-currents. (Recall the 
notation convention from Section \ref{sec:h-class}.) We have fibrations:
\begin{equation}\label{CD:1.5}
\begin{CD}
\bbT_Y@>>>\op{C}(Y)\\
@. @V\op{d} VV\\
@. \scrZ^\bbR(\grs),
\end{CD} \qquad \begin{CD}
H^1(Y;\bbR)@>>>\tilde{\op{C}}(Y)\\
@. @V\op{d} VV\\
@. \scrZ^\bbR(\grs).
\end{CD}
\end{equation}
When \((Y_i, \nu_i
)\) is a Morse end, \(\nu _i\) together with the metric on \(Y_i\)
determines an isomorphism from  \(\scrZ^\bbR(\grs)\) to
\(\scrZ^\bbR_{(Y_i, \nu _i, \grs_i)}\) as affine spaces under \(\scrZ^\bbR(Y)\), given by
\(\frac{i}{4\pi } F_B\mapsto \frac{i}{2\pi }F_{B^E}\).

\subsubsection*{\it Proof of Theorem \ref{thm:strong-t}.} As explained in the
paragraph prior to statement of the Theorem, by the 3-dimensional
version of Theorem \ref{thm:T1}, there is a subsequence of \(\{(B_r,
\Phi_r)\}_r\) that t-converges to a t-orbit \(\pmb{\gamma}\). We
show that this subsequence is in fact strongly t-convergent. 

Choose \(b^1:=b^1(Y)\) smooth closed
2-forms \(\{\xi_1, \ldots , \xi_{b^1}\}\) whose cohomology classes
form a basis of \(H^2(Y;\bbR)\), and denote by \(\{\xi_k^*\}_{k=1}^{b^1}\) the
dual basis for \(H^1(Y; \bbR)\). Define a splitting \(\Omega^1(Y)/d\Omega ^0(Y)\to
H^1(Y;\bbR)\) of the exact sequence (\ref{eq:ex}) given by 
\[
a\mapsto \sum_k(\int_Y a\wedge\xi_k) \, \xi_k^*.
\] Since \((B_r, \Phi_r)\) are t-convergent, their curvature currents 
already converge in \(\op{B}^2(\grs)\). Thus, by (\ref{eq:ex}) it suffices to
check that there exists \(u_r\in \scrG\) such that \(B_r'=u_r\cdot
B_r\) satisfies 
\begin{equation}\label{eq:conv}
i\int_Y (B_r'- B_0)\wedge\xi_k\quad  \text{converges as \(r\to\infty\)}
\end{equation} 
for all \(k\). By Stokes' theorem and the convergence
of \(F_{B_r}\), for this purpose the choice of these 1-forms may be arbitrary.
By the 3-dimensional version of
Proposition 6.1 in \cite{Ts}, (\ref{eq:conv})
can easily shown to be true if the
2-forms \(\xi_k\) are chosen to be supported on small neighborhoods of
1-cycles in \(Y\), which are disjoint from  \(\alpha^{-1}_r(0)\) and
\(\nu^{-1}(0)\) for all sufficiently large \(r\).
\epf

\subsection{Relative homotopy classes and relative homology classes}\label{sec:SW-class}
Let \(X\) be \(\Spin^c\) MCE, and use \(\bbS=\bbS_X=\bbS^+\oplus \bbS^-\) to denote its
associated spinor bundle. For each end \(\hat{Y}_i\) of
\(X\), fix a \(\grc_i\in\scrB_{Y_i}\). Denote by
\(\scrB_X(\{\grc_i\}_{i\in \grY})=\{(A, \Psi)\}/\sim\) the quotient
configuration space, where \((A, \Psi)\) is a
pair consisting of a \(\Spin^c\) connection on \(\bbS^+\) and a
\(\Psi\in \Gamma (\bbS^+)\), such that \((A, \Psi)\) is asymptotic to a
representative of \(\grc_i\) on the \(\hat{Y}_i\)-end. Here, \((A',
\Psi')\sim (A, \Psi)\) if they are related by a gauge transformation,
i.e., 
\((A', \Psi')=u\cdot (A, \Psi)\) for some \(u\in C^\infty(X, S^1)\). Given an
element in \(\scrB_X(\{\grc_i\}_i)\), its {\em relative homotopy
  class} refers to the element in \(\pi_0\scrB_X(\{\grc_i\}_i)\) given
by the path component it is in. In the case when \(X=\bbR\times Y\)
and \(Y\) is a \(\Spin^c\) closed 3-manifold, given \(\grc_-,
\grc_+\in \scrB_Y\), we use \(\pi _Y(\grc_-, \grc_+)\) to denote the
set of relative homotopy classes of quotient configurations with
\(-\infty\)-limit \(\grc_-\) and \(+\infty\)-limit \(\grc_+\). This is
a torsor under the group \(\pi_1(\scrB_Y)=H^1(Y;\bbZ)\). More
generally, given any \(X_\bullet\subset \ov{X}\), let \[
  \grY_{X_\bullet}:=\pi_0(\partial\ov{X_\bullet})\quad \text{and}
  \quad \scrB_{\partial\ov{X_\bullet}}:=\prod_{i\in
    \grY_{X_\bullet}}\scrB_{Y_i}.\]
Given \(\{\grc_i\}_{i\in
  \grY_{X_\bullet}}\in \scrB_{\partial\ov{X_\bullet}}\), let
\[\scrB_{X_\bullet}(\{\grc_i\}_{i\in\grY_{X_\bullet}})=\{(A, \Psi
  )\}/\sim\quad \text{and}\quad 
\pi_0\scrB_{X_\bullet}(\{\grc_i\}_{i\in
  \grY_{X_\bullet}})\] be respectively the direct generalization
of \(\scrB_{X}(\{\grc_i\})\) and  \(\pi_0\scrB_{X}(\{\grc_i\})\), with 
\((A, \Psi )\) now required to  either be asymptotic to or restricts to a representative of
\(\grc_i\), depending on whether \(Y_i\subset \partial\ov{X_\bullet}\) is a
ending 3-manifold of \(X\) or lies in the interior of \(X\). The set
\(\pi_0\scrB_{X_\bullet}(\{\grc_i\})\) is a torsor under
the group \[\pi_{X_\bullet}:=H^1(\partial \ov{X_\bullet}; \bbZ)/\im (i^*).\] Here, \(i^*\co
H^1(\ov{X_\bullet}; \bbZ)\to H^1(\partial\ov{X_\bullet}; \bbZ)\) is the map which forms part
of the commutative diagram of long exact sequences below:

{\footnotesize
\begin{equation}\label{CD3}
\minCDarrowwidth13pt\begin{CD}
\cdots H_3(\ov{X_\bullet}, \partial \ov{X_\bullet};\bbZ) @>>> H_2(\partial\ov{X_\bullet}; \bbZ) @>i
>> H_2(X_\bullet; \bbZ)@>j >> H_{2} (\ov{X_\bullet}, \partial\ov{X_\bullet};\bbZ)\cdots\\
@V\iota _{PD}VV@V\iota _{PD}VV @V\iota _{PD}VV@V\iota _{PD}VV\\
\cdots H^1(X_\bullet;\bbZ) @> i^* >> H^1(\partial\ov{X_\bullet}; \bbZ) @>d >>
H^2(\ov{X_\bullet}, \partial \ov{X_\bullet};\bbZ) @>j^*>>
H^2(X_\bullet;\bbZ) \cdots.
\end{CD}
\end{equation}}

The rows above are parts of relative exact sequences of the pair
\((\ov{X_\bullet}, \partial \ov{X_\bullet})\), and the vertical maps
\(\iota _{PD}\) are the Poincar\'{e}
duality maps.
As a result of the above commutative diagram, we have a chain of
canonical isomorphisms: 
\[\begin{split}
\pi_{X_\bullet} & =H^1(\partial \ov{X_\bullet}; \bbZ)/\ker d\simeq \im(d)\\
& \quad =\ker
j^*\simeq \ker j\subset H_2(X_\bullet;\bbZ),
\end{split}
\] 
where the first isomorphism
\(H^1(\partial \ov{X_\bullet}; \bbZ)/\ker d\simeq \im (d)\) is
induced from \(d\), and the second isomorphism \(\ker
j^*\simeq \ker j\) is the inverse of the 
Poincar\'{e} duality map. 
\begin{rem}\label{rem:b_1=0}
Note that when \(X_\bullet= X^{'a}\) and
(\ref{b_1=0}) holds, the long exact sequence in first row in
(\ref{CD3}) agrees with the \(\bbK=\bbZ\)-version of the long exact
sequence in the first row in 
 (\ref{eq:CD}) (but with each showing a different portion of the long exact
 sequence). Thus, recalling the defintion of \(\scrH_X\) from Section
 \ref{sec:h-class}, from (\ref{CD3}) we have a canonical
 isomorphism \[
   \op{h}\co \pi_X\stackrel{\sim}{\to}\scrH_X\quad \text{
     under the assumption (\ref{b_1=0}).}\]
 More generally, applying
 (\ref{CD3}) to the case \(X_\bullet= X^{'a}\) yields a map, \(\op{h}\), from \(\pi _X\simeq
 \pi _{X^{'a}}\simeq\im (d)\subset H^2(\ov{X^{'a}}, \partial \ov{X^{'a}};\bbZ)\)
 to \(\scrH_X\subset H_2(X^{'a}, \partial X^{'a};\bbZ)\) via
 \(j'\circ\iota_{PD}^{-1}\), where \(j'\co H_2(X^{'a};\bbZ)\to
 H_2(X^{'a}, \partial X^{'a};\bbZ)\) is part of the relative long exact
 sequence of the pair \((X^{'a}, \partial X^{'a})\). 
\end{rem}

Parallel to the relative homology classes discussed in Section
\ref{sec:h-class}, concatenation defines maps among sets of relative
homotopy classes in analogy to the operations \(*\),
\(\op{c}_{\{h_i\}_i}\) in Section
\ref{sec:h-class}, and we use the same notations for them: Assigned to
every \((c, c')\in \pi _Y(\grc_-, \grc)\times \pi _Y(\grc, \grc_+)\)
is an element \(c'*c\in \pi _Y(\grc_-, \grc_+)\). (This is nothing but
the composition in the fundamental groupoid of \(\scrB_Y\). In
general, we use  \(\pi (M; a,
b)\) to denote set of morphisms from \(a\) to \(b\) in the fundamental
groupoid of a topological space \(M\). Namely, it is the set of relative homotopy classes of paths in \(M\)
from \(a\) to \(b\).) 
 Fix 
\(X_\bullet\subset \ov{X}\) and  
  two arbitrary elements \(\grc=\{\grc_i\}_{i\in
  \grY_{X_\bullet}}, \grc'=\{\grc_i'\}_{i\in
  \grY_{X_\bullet}}\) in 
\(\scrB_{\partial\ov{X_\bullet}}\). Then corresponding to every  \(c\in
\pi _{\partial \ov{X_\bullet}}(\grc, \grc'):=\prod_{i\in
  \grY_{X_\bullet}}\pi _{Y_i}(\grc_i, \grc_i')\), there is an
isomorphism (the concatenation map): 
\begin{equation}\label{map:concant}
\op{c}_{c}\co \pi _0 \scrB_{X_\bullet}(\grc)\to \pi _0
\scrB_{X_\bullet}(\grc'). 
\end{equation}
This is used to make sense of Condition (2) 
in the statement of Theorem \ref{thm:l-conv}. 
\begin{lemma}\label{rem:rel_class}
(a) Given any \(i\in \grY_m\) and 
a strongly converging sequence \(\{\grc_{i, r}\}_r\), there is an
\(r_0>1\) with the following significance: There is a  
distinguished element \(o_{Y_i}(\grc_{i,r}, \grc_{i,r'}) \) in \(\pi _{Y_i}(\grc_{i, r}, \grc_{i, r'})\)
for any pairs of  \(r\), \(r'>r_0\). The number \(r_0\) depends on the
metric on \(Y_i\), but the distinguished element is independent of
it. 

(b) Let \(\{\pmb{\grc}_{ r}|\, \pmb{\grc}_r\in \Pi _{i\in
  \grY}\scrB_{Y_i}\}_r\) be a sequence satisfying condition (1) of
Theorem \ref{thm:l-conv}. Namely, \(\pmb{\grc}_r=(\grc_{i, r})_i\) is
such that \(\grc_{i,r}=\grc_i\in \scrB_{Y_i}\) is independent of \(r\)
when \(\hat{Y}_i\) is a vanishing end,
and \(\{\grc_{i, r}\}_r\subset\scrB_{Y_i}\) is a strongly t-converging sequence when
\(\hat{Y}_i\) is a Morse end. Let  \(r_0>1\) be such that for any
pair of \(r, r'>r_0\), and any \(i\in \grY_m\), the distinguished element \(o_{Y_i}(\grc_{i,r}, \grc_{i,r'})\) from
item (a) is defined. Set \(o_{Y_i}(\grc_{i,r}, \grc_{i,r'}):=1\in \pi _1(\scrB_{Y_i})=\pi _{Y_i}(\grc_{i, r},
\grc_{i, r'})\) when \(i\in \grY_v\). Let \(c(r, r'):=\big(o_{Y_i}(\grc_{i,r}, \grc_{i,r'})\big)_i
\in  \prod_{i\in \grY}\pi
_{Y_i}(\grc_{i,r}, \grc_{i,r'})\). Then the concatenation map \(\op{c}_{c(r,
  r')}\) defines a canonical isomorphism from 
\(\pi_0\scrB_X(\pmb{\grc}_{r})\) to 
\(\pi_0\scrB_X(\pmb{\grc}_{r'})\). 
\end{lemma}
\pf Item (b) in the statement of the lemma is a direct consequence of
Item (a). Item (a) follows from Item (a) of the next lemma, together with the
following observations: 

(i) Recall the map \(\Pi \co \scrB_Y\to
\Conn /\scrG\subset \op{C}(Y)\)
from the previous subsection.  For every pair \(\grc, \grc'\in
\scrB_Y\), 
The map \(\Pi \) induces an isomorphism 
\[
\Pi_*\co \pi_Y (\grc, \grc')\to \pi \, (\Conn /\scrG; \Pi \grc, \Pi
\grc')=\pi \, (\op{C}(Y); \Pi \grc, \Pi
\grc').\] 
Thus, to assign a distinguished element in \(\pi_Y (\grc, \grc')\), it suffices to assign distinguished elements in \(\pi \, (\Conn /\scrG; \Pi \grc, \Pi
\grc')\) or \(\pi
(\op{C}(Y); \Pi \grc, \Pi \grc')\). 

(ii) By definition, \(\{\grc_r\}_r\subset \scrB_Y\) is strongly
convergent iff \(\{\Pi \grc_r\}_r\subset \Conn  /\scrG\subset
\op{C}(Y) \) converges
in the current topology. \epf

The upcoming lemma is the \(SW\) counterpart of Lemma \ref{gr:h-isom} on
the \(Gr\) side. In preparation for stating it, we introduce a
counterpart of the map \(t_\scrH\) in Section \ref{sec:h-class}. 

Let \(\tilde{c}\co \bbR\to \Conn  /\scrG\) be a path from \(b\) to \(b'\)
in \(\Conn /\scrG\). This corresponds to an element \([A]\) in
\(\Conn \, (X; b, b')/C^\infty(X, S^1)\), where \(X=\bbR\times Y\),
and 
\(\Conn \, (X; b,
b')\subset \Conn \, (X)\) denotes the space of \(\Spin^c\) connections
on \(\bbS^+\) asymptotic respectively to representatives of \(b,
b'\) on its \(-\infty\)- and \(+\infty\)-end. We extend the space
\(\Conn \, (X; b, b')/C^\infty(X, S^1)\) to a space, denoted  \(\op{C} (X;b,
b')\), that  corresponds to the space of paths in \(\op{C}(Y)\) from \(b\)
to \(b'\), where \(b, b'\in \op{C}(Y)\supset \Conn /\scrG\).  Let \(A\) be a unitary
connection on \(\det \bbS^+\) representing \([A]\), and let \(h\) be a
harmonic 2-form on \(Y\) representing an element \([h]\) in \(H^2(Y;\bbR)\). The
integral \(\int_X \frac{i F_A}{4\pi }\wedge \pi _2^* h\) is finite and
depends only on the relative homotopy class of \(\tilde{c}\), \([h]\), and the
metric on \(Y\).  In this way, the metric on \(Y\) determines a map
\[
\jmath_c\co \pi \, (\Conn /\scrG; b, b')\to H^1(Y;\bbR)=\Hom 
(H^2(Y;\bbR), \bbR).
\]
In fact, \(\jmath_c\) extends to a map from \(\pi \, (\op{C}(Y); b, b')\) to
\(H^1(Y;\bbR)\), where \(b, b'\in \op{C}(Y)\). The latter map is also
denoted by \(\jmath_c\), and may be regarded as a SW analog of the map
\(I_\scrH\) introduced in Lemma \ref{gr:hR-isom}. A more precise
description of the relation between \(\jmath_c\) and \(I_\scrH\) will
be given prior to Lemma \ref{lem:htpy}. 
It follows straightforwardly from construction that \(\jmath_c\)
intertwines with the \(\pi_1(\Conn /\scrG)=\pi _1(\op{C}(Y))\simeq H^1(Y;\bbZ)\)-action on \(\pi \, (\Conn /\scrG;
b, b')\) (or more generally \( \pi _1(\op{C}(Y), b, b')\)) and the \(H^1(Y;\bbZ)\subset
H^1(Y;\bbR)\)-action on \(H^1(Y;\bbR)\). Thus, with the metric on
\(Y\) fixed, one has
a canonical isomorphism from \(\pi \, (\Conn /\scrG;
b, b')\) (or more generally \(\pi _1(\op{C}(Y), b, b')\)) to an orbit of \(H^1(Y;\bbZ)\) in \(H^1(Y;\bbR)\). Let
\(t_{\op{C}}(b, b')\) denote the corresponding element in the orbit
space, \(H^1(Y;\bbR)/H^1(Y;\bbZ)=\bbT_Y\). Summarizing, we have a map 
\[
  t_{\op{C}}\co
  \op{C}(Y)\times   \op{C}(Y)\to\bbT_Y.  
\]
This map is an analog of  the map \(t_\scrH\) defined after  Lemma
\ref{gr:hR-isom} in Section \ref{sec:h-class}. 
It is continuous with respect to the current topology on \(
\Conn /\scrG\) and \(\op{C}(Y)\). 
Furthermore, recall from (\ref{eq:CD1}) and (\ref{CD:1.5})
the structure of  \(
\Conn /\scrG\), \(\op{C}(Y)\) as \(\bbT_Y\)-bundles, \(t_{\op{C}}\)
is \(\bbT_Y\)-equivariant both with respect to the \(\bbT_Y\)-action
on the right factor of \(\op{C}(Y)\times
\op{C}(Y)\), and the \(\bbT_Y^{op}\)-action on the left factor of \(\op{C}(Y)\times
\op{C}(Y)\) (or \(( \Conn /\scrG)\times
(\Conn /\scrG)\)). (\(\bbT_Y^{op}\) above refers to the inverse
  action of the torus group.) It maps the diagonal to the identity
  element \(0\in
  \bbT_Y\). 
\begin{lemma}\label{sw:h-isom}
(a) Suppose \(b, b'\in \op{C}(Y)\) are sufficiently close in the sense that \(t_{\op{C}}
(b, b')\) falls in the ball \(B_{0}(1/2)\subset
\bbT_{Y}\). Then there is a distinguished element \(\td{o}_Y(b, b')\) in \(\pi
(\op{C}(Y);b, b')\).  This distinguished element is independent of the
metric on \(Y\), though the notion of being  ``sufficiently close'' does. 

(b) Fix \(X_\bullet \subset \ov{X}\). 
Suppose a pair \(\{b_i\}_i , \{b'_i\}_i\in \prod_{i\in
  \grY_{X_\bullet}}\op{C}(Y_i)\) are sufficiently close in the sense
described above. 
 Given any \(\grc_i\in \Pi ^{-1} b_i\),
\(\grc_i'\in \Pi ^{-1} b_i'\), let \(o_{Y_i}(\grc_i, \grc_i'):=(\Pi _*)^{-1}\td{o}_{Y_i}(b_i, b'_i)\in \pi _{Y_i}
(\grc_i, \grc_i')\). Then the concatenation map 
\(\op{c}_{\{o_{Y_i}(\grc_i, \grc_i')\}_i}\)  defines a canonical isomorphism from  \(\pi
_0\scrB_{X_\bullet} (\{\grc_i\}_i)\) to \(\pi
_0\scrB_{X_\bullet} (\{\grc_i'\}_i)\) as affine spaces under \(\pi _{X_\bullet}\). 
\end{lemma}
\pf Item (a) follows directly from the discussion preceding to
the statement of the lemma. Item (b) generalizes Item (b) of the
Lemma \ref{rem:rel_class}, and follows from Item (a) above together
with Part (i) in the proof of Lemma \ref{rem:rel_class}. 
\epf

\begin{rem}\label{rem:rel-htpy-convention} 
The arguments in the proof of the preceding lemma 
also establishes the following: A choice of basis for \(H^1(Y; \bbZ)\) gives a
way of (simultaneously) identifying all sets of relative
homotopy classes \(\pi _Y (\grc, \grc')\),  \(\grc, \grc'\in \scrB_Y\),
with \(H^1(Y;\bbZ)\) (as affine spaces). (Recall that such a choice
(cf. (\ref{eq:deltab1})) is
required to define the normalized Coulomb gauge, and has been fixed
implicitly in this article.) To see this, observe that as explained in Part (i) in the proof of
Lemma \ref{rem:rel_class}, it suffices to identify sets of relative
homotopy classes  \(\pi    (\op{C}(Y); b, b')\), 
with \(H^1(Y;\bbZ)\) (as affine spaces) for every pair \(b, b'\in
\op{C}(Y)\). 
This is equivalent to
(consistently) choosing a base point in \(\pi   (\op{C}(Y); b,
b')\) for each pair \(b, b'\in
\op{C}(Y)\).  To do so, note
that a choice of basis for \(H^1(Y;\bbZ)\) defines an isomorphism
\[
i_h \co \bbR^{b^1} \stackrel{\sim}{\to} H^1(Y;\bbR). 
\]
Recall also that \(t_{\op{C}}(b, b')\) corresponds to an isomorphism from \(\pi    (\op{C}(Y); b,
b')\) to an orbit of
\(H^1(Y;\bbZ)\) in \(H^1(Y;\bbR)\). We define the
base point in \(\pi    (\op{C}(Y); b,b')\) to be the unique element 
whose image under this isomorphism lies in \(i_h\, ([0,1)^{\times
  b^1})\subset H^1(Y;\bbR)\). Denote this base point by \(\td{o}_Y(b,
b')\), and given \(\grc, \grc'\in \scrB_Y\), let
\[
o_{Y}(\grc, \grc'):=(\Pi _*)^{-1}\td{o}_Y(\Pi \grc, \Pi \grc')\in
\pi _Y (\grc, \grc')\]
denote the corresponding base point in \(\pi _Y (\grc, \grc')\).
The argument for Item (b) of Lemma
\ref{sw:h-isom} then implies that once a choice of basis is fixed for
every \(H^1(Y_i; \bbZ )\),  then for any given \(X_\bullet\), \(\pi
_0\scrB_{X_\bullet} (\{\grc_i\}_i)\), \(\pi
_0\scrB_{X_\bullet} (\{\grc_i'\}_i)\) are identified for every pair
\(\{\grc_i\}_i, \{\grc_i'\}_i\) in \(\scrB_{\partial\overline{X_\bullet}}\). 
\end{rem}
\begin{rem}\label{rem:based-rel-htpy}
Fix a reference connections \(A_0\) on \(\bbS^+\) as in
(\ref{eq:A_0}). Then for every given \(X_\bullet\), \((A_0,
0)|_{X_\bullet}\) determines a base
point in \(\pi_0(\scrB_{X_\bullet}(\{\grc_{0,i}\}_i))\), where \(i\in
\grY_{X_\bullet}\), and \(\grc_{0,i}\in \scrB_{Y_i}\) denotes the
restriction of \((A_0,
0)\) to the \(i\)-th connected component of
\(\partial\ov{X_\bullet}\). This in turn defines an
isomorphism 
\[
\imath_{A_0}\co \pi_0(\scrB_{X_\bullet}(\{\grc_{0,i}\}_i))\stackrel{\sim}{\to }\pi
_{X_\bullet}
\]
as affine spaces. 

Thus, a choice of \(A_0\) determines, for any given \(X_\bullet\subset \ov{X}\),  a way of simultaneously identifying \(\pi_0(\scrB_{X_\bullet}(\{\grc_i\}_i))\simeq \pi
_{X_\bullet }\) for all \(\{\grc_i\}_i\in \scrB_{\partial\overline{X_\bullet}}\): Given any \(\{\grc_i\}_i\in \prod_{i\in
  \grY_{X_\bullet}}\scrB_{Y_i}\), Let \(h_{A_0}\) denote the
isomorphism
\[
h_{A_0}:=\imath_{A_0}\circ  \op{c}_{\{o_{Y_i}(\grc_i,
  \grc_{0,i})\}_{i\in \grY_{X_\bullet}}}\co
\pi_0(\scrB_{X_\bullet}(\{\grc_{i}\}_{i\in \grY_{X_\bullet}}))\stackrel{\sim}{\to }\pi
_{X_\bullet}.
\]
When \(X_\bullet\) is of the form \(\hat{Y}_{i, I}\), this isomorphism
agrees with the isomorphism \(\pi _{Y_i}(\grc,
\grc')\stackrel{\sim}{\to } H^1(Y_i;\bbZ)\) introduced in  Remark
\ref{rem:rel-htpy-convention} under our assumptions on \(A_0\).
\end{rem}

We next compare the relative homology classes of
t-curves and the relative homotopy classes of Seiberg-Witten quotient 
configurations. Recall from
Section \ref{sec:h-class} the definition of \(\scrH^\bbR(X^{'a}, \nu , \grs,
\{\tilde{\gamma}_i\})\) and various other notions. Let \(c=[(A,
\Psi)]\in \scrB_X(\{\grc_i\}_i)\), with \(\grc_i=[(B_i, \Phi_i)]\in
\scrB_{Y_i}\). Assign to \(c\) the 2-current \(\tilde{c}=\frac{i}{2\pi} F_{A^E}\) on \(X^{'a}\). This
current determines a class in \(\scrH^\bbR(X^{'a}, \nu , \grs,
\{\tilde{\grc}_i\}_{i\in \grY_m})\), where
\(\tilde{\grc}_i=\frac{i}{2\pi} F_{B_i^E}\in\scrZ(Y_i, \nu _i,
\grs_i)\subset
\scrZ^\bbR(Y_i, \nu _i, \grs_i)\). This class only depends on the
relative homotopy class of \(c\), and so in this way we have a map
\[
\grh'  \co \pi _0\scrB_X(\{\grc_i\}_{i\in \grY})\to \scrH^\bbR((X^{'a}, \nu , \grs),
\{\tilde{\grc}_i\}_{i\in \grY_m}). 
\] 
This map intertwines with the \(\pi _X\)-action on \(\pi
_0\scrB_X(\{\grc_i\}_{i\in\grY})\), and the \(\op{h}(\pi _X)\subset
\scrH^\bbR_X\)-action on \(\scrH^\bbR((X^{'a}, \nu , \grs),
\{\tilde{\grc}_i\}_{i\in\grY_m})\). (The map \(\op{h}\co \pi _X\to \scrH_X\) is as defined in Remark
\ref{rem:b_1=0}.) It is also natural with respect to
the concatenation maps on both sides. In the special case when \((X, \nu
)=(\bbR\times Y_i, \pi _2^*\nu _i)\), \(i\in\grY_m\),  is cylindrical,
the map 
\(\grh' \co \pi _{Y_i} (\grc, \grc')\to \scrH^\bbR(Y_i, \nu _i, \grs_i; \tilde{\grc},
\tilde{\grc}')\) by construction factors through a map \[\ud{\grh}\co \pi \,
 (\op{C}(Y_i), \Pi \grc, \Pi \grc')\to \scrH^\bbR(Y_i, \nu _i, \grs_i; \tilde{\grc},
\tilde{\grc}'):\]  \(\grh'=\ud{\grh}\circ \Pi
_*\). It follows from construction that in this case, 
\[
  \jmath_c=\iota_{PD}\circ I_{\scrH}\circ\underline{\grh}, 
  \]
  where \(\iota_{PD}\) again denotes the Poincar\'e  map.

  Next, consider the case when (\ref{b_1=0}) holds. 
  Recall  from Remark \ref{rem:b_1=0} that
\(\pi_X\stackrel{\rmh}{\simeq }\scrH_X\) in this case.
\begin{lemma}\label{lem:htpy}
Assume that all vanishing ends of  \((X, \nu)\) has
zero first betti numbers. (See  (\ref{b_1=0}).) For each \(i\in \grY_m\), let
\(\{\grc_{i, r}\}_r\in \scrB_{Y_i}\) be a sequence that strongly
t-converges to a t-orbit \(\pmb{\gamma }_i\), and for each \(i\in
\grY_v\), let \(\grc_{i,r}=\grc_i\in \scrB_{Y_i}\) be independent of
\(r\). 
Let  \(r_0>1\) be as in Lemma \ref{rem:rel_class}. 
Then:

(a) For all \(r>r_0\), there is a canonical
isomorphism \(\grh\) from the set \(\pi_0\scrB_X(\{\grc_{i,r}\}_{i\in
  \grY})
\) to an
\(\scrH_X\)-orbit in \(\scrH^\bbR\big((X^{'a}, \nu),\{\tilde{\gamma }_i\}_i\big)\) as affine spaces under
\(\pi_X\simeq \scrH_X\).

(b) Let \(\{(A_r, \Psi _r)\}_r\) be as in the statement of  Theorem
\ref{thm:l-conv}, and suppose that for each \(i\in \grY_m\), the sequence \(\{\grc_{i,r}\}_r\), 
consists of  the \(Y_i\)-end limits of a sequence of admissible solutions
 \(\{(A_r, \Psi _r)\}_r\). 

Then the
image of the map \(\grh\) is \(\scrH\big((X^{'a}, \nu),\{\tilde{\gamma
}_i\}_i\big)\subset \scrH^\bbR\big((X^{'a}, \nu),\{\tilde{\gamma
}_i\}_i\big)\). Thus, in this case \(\grh\) is a canonical
isomorphism from \(\pi_0\scrB_X(\{\grc_{i,r}\}_i)\) to \(\scrH\big((X^{'a},
\nu),\{\tilde{\gamma }_i\}_i\big)\) as affine spaces under \(\pi_X\simeq \scrH_X\).
\end{lemma}
\pf   By construction, we have the following commutative diagram under
the assumptions of the lemma: For any \(r, r'>r_0\)
{\footnotesize
\begin{equation}\label{CD:grh}
\xymatrix@=18pt{
   \pi_0\scrB_X(\{\grc_{i,r'}\}_{i\in
     \grY})\ar@{->}[r]^{}\ar@{->}[d]^{\grh'}&
   \pi_0\scrB_X(\{\grc_{i,r}\}_{i\in \grY}) \ar@{->}[d]^{\grh'}
   \\
\scrH ^\bbR\big((X^{'a}, \nu),\{\tilde{\grc}_{i,r'}\}_{i\in \grY_m}\big)
\ar@{->}[r]^{} & \scrH ^\bbR\big((X^{'a},
\nu),\{\tilde{\grc}_{i,r}\}_{i\in \grY_m}\big) \ar@{->}[r]^{i_\infty} & \scrH ^\bbR\big((X^{'a}, \nu),\{\tilde{\gamma }_i\}_i\big).
 }
\end{equation}}
where all the maps are morphisms of affine spaces under \(\pi
_X\simeq\scrH_X\subset \scrH_X^\bbR\); the horizontal map in the top
row is the canonical isomorphism from Lemma \ref{rem:rel_class}; the
horizontal maps in the bottom
row are the canonical isomorphisms from Lemma \ref{gr:hR-isom}.
The composition  \(i_\infty\circ \grh'\) maps
 \(\pi_0\scrB_X(\{\grc_{i,r}\}_{i\in \grY})\) to another orbit of the
 \(\scrH_X\)-action in \(\scrH ^\bbR\big((X^{'a}, \nu),\{\tilde{\gamma
 }_i\}_i\big)\). The canonical isomorphism claimed in Item (a) of the
 statement of the present lemma is defined to be this composition map:
 \(\grh:=i_\infty\circ \grh'\).

 The proof of Item (b) will follows as a by-product of the proof of
 Theorem \ref{thm:g-conv} (c), and will be deferred to Section
 \ref{sec:g-conv:a} 
 \epf

\subsection{Taubes' proof for \(SW\Rightarrow Gr\): A synopsis}\label{sec:synopsis}

The proof of Theorem \ref{thm:l-conv} follows the outline of Taube's
arguments in \cite{T}. To serve as a roadmap for the remainder of this
article, a brief summary of the ingredients in \cite{T} is
provided here for reader's convenience. For each step listed below, we
indicate where it takes place in  \cite{T}, \cite{Ts}, and
the present article, under their respectively contexts.

Let \((X, \nu)\) be an admissible pair and let 
 \((A, \Psi)=(A_r, \Psi_r)\) be as in the statement of Theorem
 \ref{thm:l-conv}. Write \[
\Psi=(r/2)^{1/2} \psi, \quad \text{and \(\psi=(\alpha, \beta)\)}\]
 according to the decomposition \(\bbS^+=E\oplus E\otimes K^{-1}\) on
 \(X^{'a}\). The \(\Spin^c\)-connection corresponding to \(A\) induces a connection
on \(E\), and together with the Levi-Civita connection, also a
connection on \(E\otimes K^{-1}\). We use \(\nabla_A\alpha\),
\(\nabla_A\beta\) to denote the covariant derivative with respect to
the aforementioned induced connections. 

Taubes' proof proceeds with the steps listed below: 
\begin{itemize}
\item[(1)] Obtaining pointwise estimates of \(|\beta|, |F_A|, |\nabla_A\alpha|,
|\nabla_A\beta|\) in terms of the ``energy
density'' \(r|\varpi|:=r||\nu|-|\alpha|^2|\).
The estimates show that these quantities are interesting where \(\alpha\) is
small. Cf. \cite{Ts} Sections 3 b)-e), \cite{T} Section I.2. 
Section \ref{sec:pt-est} in this article contains the corresponding estimates
on a MCE.

\item[(2)] Obtaining an integral bound for \(r|\varpi|\) in terms of constants depending
only on the \(\Spin^c\) structure and the metric. This is basically
equivalent to an integral bound on
\(\frac{iF_A}{2\pi}\wedge\nu\). Cf. \cite{Ts} Section 3 a) and
references therein. The
corresponding results are establishes in \textsection \ref{sec:int-est} below. 

\item[(3)] Establishing a monotonicity formula for \(\mathcal{W}_B\),
  which is  a certain notion of  ``energy'' 
over a ball \(B\) in \(X\). 
Cf. \cite{Ts} Section 4 and \cite{T} Section I.3 and the beginning of
Section \ref{sec:mono} below.
Roughly speaking, \(\scrW_B\) is the integral of the ``energy
density'' \(r|\varpi|\) over \(B\). Note that this notion of energy 
is however of a
different nature from other, more typical notions of energy for
Seiberg-Witten theory on MCE, such as
the Chern-Simons-Dirac functional, or more generally, the topological
energy \(\scrE_{top}\) introduced in \cite{KM}. This \(\scrW_B\) should be regarded
as an analogue of the notion of energy on the Gromov side, namely the
area of the holomorphic curves. As a result of the monotonicity
formula, one obtains  an \(r\)-independent bound on the
2-dimensional Hausdorff measure of \(\alpha^{-1}(0)\), and hence also
on that of the t-curve that \(\alpha^{-1}(0)\) geometrically converges to as
\(r\to \infty\). The relevant monotonicity formula in our context is
given in Section \ref{sec:mono} below.

\item[(4)] Rescaling that \((A, \psi)\) over a ball of radius \(O(r^{-1/2})\)
yields an approximate solution to the version of the 
Seiberg-Witten equation (\ref{eq:SW}) on the Euclidean space \(X=\bbR^4=\{(x_1, \cdots, x_4)\}\) and
\(\mu^+=\frac{1}{2}(dx_1\wedge dx_2+dx_3\wedge dx_4)\). The behavior of
Seiberg-Witten solutions for such \((X, \mu^+)\) are 
well-understood in terms of vortex solutions on \(\bbR^2\), and these
give the local models
for \((A, \psi)\). In particular, combining with Items (1)--(3) above, this
implies the exponential decay of  \(|\beta|, |F_{A^E}|, |\nabla_A\alpha|,
|\nabla_A\beta|\) away from \(\alpha^{-1}(0)\). This step appeared in
Section 6 of \cite{Ts} and Section I.4 of \cite{T}. This
adapts easily to our context; see Sections \ref{sec:exp} and
\ref{sec:local_m} below.

\item[(5)] Regarding \(\frac{i}{2\pi}F_{A^E}\) as a 2-current , Items (1) and (2) above imply that this
current is bounded independently of \(r\). Thus, Alaoglu's theorem
(cf. e.g. \cite{Rd} Theorem 10.6.17) implies that these currents
converge to a current \(\mathcal{F}\) in weak\(^*\)
topology. Furthermore, via (3) and (4) it is shown that the support of
\(\mathcal{F}\), denoted temporarily by \(C\), is a closed space of 2-dimensional measure,
and is the \(r\to\infty\) limit of \(\alpha^{-1}(0)\) in the sense of
(\ref{eq:dist-conv}). Cf. the first half of Section I.5 of \cite{T};
Sections 7(a) and (b) of \cite{Ts}, and Section \ref{sec:l-conv} in
this article.

\item[(6)] Using Item (4), it is shown that the current \(\mathcal{F}\) defines a
``positive cohomology assignment'', namely, the map from the set of a
certain kind of generic disks in \(X\) (the ``admissible disks'') to \(\bbZ\), whose
value is positive for admissible pseudo-holomorphic disks which
intersect \(C\). Taubes showed that this fact guarantees that \(C\) is
a pseudo-holomorphic curve. Cf. Section 6 of \cite{T} and Sections
7(c)-(e) of \cite{Ts}. This part of Taubes' argument applies directly in our setting with
no need for modification.
\end{itemize}

\section{\(L^2_{1, loc}\)-bounds and integral estimates}\label{sec:4}

The main results of this section come in two groups. 
 The first group of  main results, Propositions
\ref{prop:SW-L2-bdd} and  \ref{est-L^2_1}, provide 
\(L^2_1\)-bounds on solutions \((A_r, \Psi _r)\) to the Seiberg-Witten
equation \(\grS_{\mu _r, \hat{\grp}}(A_r, \Psi _r)=0\), with careful
control on the growth rate of these bounds on \(r\). While (for each
fixed \(r\)) such
\(L^2_1\) bounds usually serve as the starting point of  typical
proofs compactness results in Seiberg-Witten
theory (cf. e.g. \cite{KM}'s Corollary 10.6.2 and Theorem 5.2.1
(ii)), the usual arguments to derive them (e.g. those in  \cite{KM}),
do not provide sufficient bounds on the \(r\)-growth rate needed for t-convergence. 
The second group of main results, Lemmas \ref{co:E-omega-bdd3} and
\ref{T:lem3.1}, supply the type of ``energy bounds'' required to begin the proof
of \(t\)-convergence results. This corresponds to Step (2) of the list
in Section \ref{sec:synopsis}.



\subsection { \(\scrE_{top}\) and  \(\scrE_{anal}\): the topological energy
  and the analytic energy}\label{sec:top-energy}

We begin with some general definitions and observations. 
Similar to what was done in \cite{KM} (Definition 4.5.4 therein), we introduce two
Floer-theoretic notion of ``energy'', \(\scrE_{anal}(A, \Psi)\) and
\(\scrE_{top}(A,\Psi)\) below. They are related when \((A, \Psi )\)
satisfies a Seiberg-Witten equation. (The subscripts ``anal'' and ``top'' signify
``analytical'' and ``topological'' respectively.) As is typical in
Floer theories, \(\scrE_{anal}\) is useful for estimating
\(L^2_1\)-norms, while \(\scrE_{top}\) depends only on the end points
of the Floer trajectory and its relative homotopy class. In the
cylindrical case, 
\(\scrE_{top}\) is the difference in the values of the CSD functional at end points
(cf.  \textsection \ref{sec:adm-SW} for the definition of CSD functionals).

Let \(X\) be a \(\Spin^c\) MCE with ending 3-manifolds \(Y_i\), and let \(\mu \)
be a closed 2-form on \(X\) that has a closed 2-form \(\mu _i\) as the 
\(Y_i\)-end limit for each \(i\in \grY\), in the sense defined in
Section \ref{sec:convention}.  Recall also from Section \ref{sec:convention} the
definition of \(X_\bullet\). (Note in particular that the boundary
components of \(X_\bullet\) consist of slices of the form \(Y_{i:l}\)).
\begin{defn}
Fix a compact \(X_\bullet\subset X\). Given \((A, \Psi)\in \scrC (X_\bullet)\), let 
\begin{equation}\label{eq:E-an}
\begin{split}
\scrE_{anal}^\mu(X_{\bullet})(A, \Psi)= &\frac{1}{4}\int_{X_{\bullet}}|F_A|^2+\int_{X_{\bullet}}|\nabla_A\Psi|^2+\int_{X_{\bullet}}\Big|\frac{i}{4}\rho(\mu^+)-(\Psi\Psi^*)_0\Big|^2\\
&\qquad 
+\int_{X_{\bullet}}\frac{R_g}{4}|\Psi|^2-\frac{i}{4}\int_{X_{\bullet}}F_A\wedge*_4\mu,
\\
\end{split}\end{equation}
where \(R_g\) denotes the scalar curvature.  

Over a slice
\(Y_{i:s}:=\mathfrc{s}_i^{-1}(s)\simeq Y_i\) of \(\bar{X}\), let \((B,
\Phi)=(B(s),
\Phi(s))=(A, \Psi )\big|_{Y_{i:s}}\) denote
the restriction of \((A, \Psi)\) to \(\mathfrc{s}_i^{-1}(s)\subset
X\). 
We define 
\begin{equation}\label{eq:E-top}
\scrE_{top}^\mu(X_{\bullet})(A, \Psi)=\frac{1}{4}\int_{X_{\bullet}}F_A\wedge
F_A-\int_{\partial \ov{X_\bullet}}\langle \Phi,
\slp_B\Phi\rangle+\frac{i}{4}\int_{X_{\bullet}}F_A\wedge \mu. 
\end{equation}

We then define
\begin{equation}\label{def:E_top nonlocal}
\scrE_{top}^{\mu, \hat{\grp}}(X_\bullet)(A,
\Psi)=\scrE_{top}^\mu (X_{\bullet})(A,
\Psi)-2\, f_{\hat{\grp}}((A, \Psi )|_{\partial\overline{X_\bullet}}).
\end{equation}
\end{defn}

The preceding definition of \(\scrE_{top}^\mu \) and \(\scrE^\mu
_{anal}\) is motivated by the following identity: 
\begin{equation}\label{eq:E-top=anal}
\|\mathfrak{S}_{\mu_r}(A, \Psi) \|_{L^2(X_{\bullet})}^2=\mathcal{
E}_{anal}^{\mu_r}(X_{\bullet})(A, \Psi)-{\mathcal E}_{top}^{\mu_r}(X_{\bullet})(A, \Psi). \quad
\end{equation}
In particular, when \((A, \Psi )\in \scrC(X_\bullet)\) 
is a solution to the
Seiberg-Witten equation \(\grS_{\mu, \hat{\grp}}(A, \Psi)=0\), one
has:  
\[
\scrE_{anal}^{\mu}(X_\bullet)(A,
  \Psi)
=\scrE_{top}^{\mu,    \hat{\grp}}(X_\bullet)(A, \Psi)+\|\hat{\grp}(A, \Psi )\|_{L^2(X_\bullet)}^2.
\]
This identity is used in Part 1 of the next subsection to bound the square
terms in the first line of (\ref{eq:E-an}) in terms of
\(\scrE_{top}\).  What follows are some observations demonstrating the
topological nature of \(\scrE_{top}\).

Let \(A_0\), \(B_{0,i}\) be the reference connections from 
(\ref{eq:A_0}). Recall the definitions of CSD functionals from (\ref{eq:CSD_q}),
(\ref{eq:CSD}),  and  Remark \ref{rem:M-disconn}. It is often convenient to re-express  the topological
energies \(\scrE_{top}^\mu \), \(\scrE_{top}^{\mu , \hat{\grp}}\) in
terms of the CSD functionals:  
By the Stokes' theorem, 
\begin{equation}\label{eq:E-top-q}
\begin{split}
 \scrE_{top}^{\mu}(X_\bullet ) (A, \Psi) &=\frac{1}{4}\int_{X_{\bullet}}F_{A_0}\wedge
  F_{A_0}+  \frac{i}{4}\int_{X_{\bullet}}F_{A_0}\wedge\mu
  -2\op{CSD}_{\mu}^{\partial \ov{X_\bullet }}(B, \Phi);\\
\scrE_{top}^{\mu, \hat{\grp}}(X_\bullet ) (A, \Psi) &=\frac{1}{4}\int_{X_{\bullet}}F_{A_0}\wedge
  F_{A_0}+  \frac{i}{4}\int_{X_{\bullet}}F_{A_0}\wedge\mu_r
  -2\op{CSD}_{\mu, \hat{\grp}}^{\partial \ov{X_\bullet }}(B, \Phi).\\
\end{split}
\end{equation}

Given the formulas (\ref{eq:E-top}), (\ref{def:E_top nonlocal}),
\(\scrE_{top}^\mu \) and \(\scrE_{top}^{\mu , \hat{\grp}}\) are by
definition gauge-invariant. Moreover, \(\scrE_{top}^\mu (X_\bullet)(A, \Psi
)\), \(\scrE_{top}^{\mu , \hat{\grp}}(X_\bullet)(A, \Psi
)\) are 
can be computed from the relative homotopy class of
\([(A, \Psi )]\) via the following explicit formulas: 
Recalling that \([(B, \Phi )]_c\in \scrC(Y)\) denotes the representative of \([(B,
\Phi )]\in \scrB(Y)\)
in the normalized Coulomb gauge, one may re-express (\ref{eq:E-top-q}) as
\begin{equation}\label{eq:E-top-h}
\begin{split}
& \scrE_{top}^{\mu}(X_\bullet ) (A, \Psi)=\frac{1}{4}\int_{X_{\bullet}}F_{A_0}\wedge
  F_{A_0}+  \frac{i}{4}\int_{X_{\bullet}}F_{A_0}\wedge\mu\\
& \qquad \quad   -2\op{CSD}_{\mu}^{\partial \ov{X_\bullet }}([(B,
\Phi)]_c)-2i^*\big(
  (\pi  c_1(\grs_X|_{X_\bullet})-\frac{\mu |_{X_\bullet}}{4}\big)\cdot
  h\,(A, \Psi );\\
& \scrE_{top}^{\mu, \hat{\grp}}(X_\bullet ) (A, \Psi)=\frac{1}{4}\int_{X_{\bullet}}F_{A_0}\wedge
  F_{A_0}+  \frac{i}{4}\int_{X_{\bullet}}F_{A_0}\wedge\mu\\
& \qquad \quad   -2\op{CSD}_{\mu, \hat{\grp}}^{\partial \ov{X_\bullet }}([(B,
\Phi)]_c)-2i^*\big(
  (\pi  c_1(\grs_X|_{X_\bullet})-\frac{\mu |_{X_\bullet}}{4}\big)\cdot
  h\,(A, \Psi ), 
\end{split}
\end{equation}
where \(h\, (A, \Psi )\in H^1(\partial \ov{X_\bullet}; \bbZ)/\im i^*\simeq \pi
_{X_\bullet}\) denotes the image of the relative homology class of
\([(A, \Psi )]\) under the map \(h_{A_0}\) from Remark \ref{rem:based-rel-htpy}, and \(i^* \co
H^*(\ov{X_\bullet}; \bbZ)\to H^*(\partial \ov{X_\bullet };\bbZ)\) is part of the
relative long exact sequence of the pair \((\ov{X_\bullet}, \partial\ov{
X_\bullet})\). Note that the pairing between the modules \(H^1(\partial \ov{X_\bullet};
\bbZ)/\im i^*\simeq H_2(\partial \ov{X_\bullet};
\bbZ)/\ker i\) and  \(\im (i^*\co H^2(\ov{X_\bullet}; \bbZ)\to H^2(\partial \ov{X_\bullet};\bbZ))\) is
well-defined because the submodule \( \im (i^*\co H^2(\ov{X_\bullet}; \bbZ)\to H^2(\partial\ov{
X_\bullet};\bbZ))\) pairs trivially with \( \ker (i\co H_2(\partial \ov{X_\bullet}; \bbZ)\to H_2(\ov{
X_\bullet};\bbZ))\). 

To summarize:
\begin{lemma}
Let \(X\), \(\mu \), \(\hat{\grp}\) be as before. Given arbitrary
(possibly non-compact) \(X_\bullet \subset X\) and an admissible \((A,
\Psi)\in \scrC(X_\bullet)\), the topological
energy  \(\scrE_{top}^\mu(X_\bullet )(A, \Psi)\) is
well-defined. 
Moreover, their values depend only on the
\(\Spin^c\)-structure of \(X\), \([(A, \Psi )|_{\partial
  \ov{X_\bullet}}]\in \scrB(\partial
  \ov{X_\bullet})\), and the relative
homotopy class of \([(A, \Psi)]\).
\end{lemma}

Moreover, combining  
(\ref{eq:E-top-h}) 
and Lemmas \ref{lem:CSD-est}, \ref{lem:htpy}, one 
has: 
\begin{lemma}\label{lem:E_topX}
Adopt the assumptions and notations in the statement of Theorem
\ref{thm:l-conv}. (In particular, recall that  \((A_r, \Psi_r)\) are
Seiberg-Witten solutions \(\grS_{\mu _r, \hat{\grp}}(A_r, \Psi _r)=0\)
with \(Y_i\)-end limits proscribed by \(\{\pmb{\gamma
}_i\}_{i\in\grY_m}\), \(\{(B_i, \Phi _i)\}_{i\in \grY_v}\), 
and with fixed relative homotopy class 
\(\mathfrc{h}\).) 
Then there exist  a  constant
\(\smE, \in\bbR^+\) depending only on \(\|\nu \|_{C^1}\),
\(\varsigma_w\), and:
\BTitem\label{dep:rel-htpy}
\item the \(\Spin^c\)-structure on
\(X\), 
\item the relative homotopy class \(\mathfrc{h}\), or equivalently
  \(\grh(\mathfrc{h})\) (and hence also 
implicitly on  \(\{\pmb{\gamma
}_i\}_{i\in\grY_m}\)), 
\item the cohomology class \([\nu ]\),  
\ETitem
such that for
all \(r\geq r_0\) (\(r_0\) being larger or equal to that in Lemma \ref{lem:CSD-est}), 
\begin{equation}\label{eq:CSD-est}
|\scrE_{top}^{\mu_r}(X)(A_r, \Psi_r)|\leq  \smE r,\quad\text{and}
\quad 
|\scrE_{top}^{\mu_r, \hat{\grp}}(X)(A_r, \Psi_r)|\leq \smE \, r.
\end{equation}
\end{lemma}

The next lemma says more about the coefficient \(\smE\) above. 

\begin{lemma}
(a) Let \(Y_i\) be a Morse end. 
Suppose \(\{(B_r, \Phi _r)\}_r\) is a sequence of 
Seiberg-Witten solutions \(\grF_{\mu _{i,r}}(B_r, \Phi _r)=0\) that
strongly t-converges to \(\pmb{\gamma }\). Then the
limit \[\lim_{r\to \infty}(r^{-1} \op{CSD}_{\mu _{i,r}}([B_r, \Phi
_r]_c))=\lim_{r\to \infty}(r^{-1} \op{CSD}_{\mu _{i,r}, \grq_i}([B_r, \Phi _r]_c))\]
exists, and 
\begin{equation}\label{ineq:CSD}
\Big|\lim_{r\to \infty}(r^{-1} \op{CSD}_{\mu _{i,r}}([B_r, \Phi
_r)]_c)\Big|\leq \frac{\pi }{2}\sqrt{b^1(Y_i)}\, |[\nu _i]|. 
\end{equation}

(b)  Adopt the notations and assumptions of Theorem \ref{thm:l-conv}. Then the
limit 
\[
\bbE:=\lim_{r\to \infty}(r^{-1} \scrE_{top}^{\mu _r, \hat{\grp}}(A_r, \Psi
_r))
\]
exists, and equals \(\lim_{r\to \infty}(r^{-1} \scrE_{top}^{\mu _r}(A_r, \Psi
_r))\). It is determined by the items listed in (\ref{dep:rel-htpy})
via the formula (\ref{eq:bbE}) below.  
\end{lemma}
\pf {\em (a):} Let \(\grc_r\) denote the gauge equivalence class of \((B_r,
\Phi _r)\). 
Recall the map \(\jmath_c\co \pi \, (\Conn /\scrG;  [B_0], \Pi
\grc_r)\to H^1(Y_i;\bbR)\simeq H_2(Y_i; \bbR)\) from Section \ref{sec:SW-class} and let
\(\td{c}_0(\grc_r)\in \pi \, (\Conn /\scrG; [B_0], \Pi \grc_r)\) be the element represented
by the path \(s\mapsto B_0+(1-\chi (s))(B_r-B_0)\) on \(\Conn (Y_i)\). Then
\begin{equation}\label{ineq:CSD_r}
  \begin{split}
r^{-1} \op{CSD}_{\mu _{i,r}}([B_r, \Phi _r]_c)& =r^{-1} \op{CSD}_{\mu
  _{i,r}, \grq_i}([B_r, \Phi _r]_c)\\
& = r^{-1} \op{CSD}_{w_{i,r}}(\grc_r)
-\frac{i}{8}\int_{Y_i} ([B_r]_c-B_0)\wedge\nu _i\\
& =r^{-1} \op{CSD}_{w_{i,r}}(\grc_r)
-\frac{\pi }{2}[\nu _i]\cdot \jmath_c (\td{c}_0(\grc_r)).\\
\end{split}
\end{equation}
The assumption that \((B_r, \Phi _r)\) strongly converges implies that
the limit
\(\lim_{r\to\infty}\jmath_c (\td{c}_0(\grc_r))\) exists, and its norm is bounded
by \(\sqrt{b^1(Y_i)}\). Denote this limit by 
\(
\jmath_h (\pmb{\gamma })\in 
H_2(Y_i;\bbR)\).  
Together with Lemma
\ref{lem:CSD-est}, this implies that for \(i\in \grY_m\), 
\[
\lim_{r\to \infty}(r^{-1} \op{CSD}_{\mu _{i,r}}([B_r, \Phi
_r]_c))=-\frac{\pi }{2}[\nu _i]\cdot \jmath_h (\pmb{\gamma }), 
\]
and hence Assertion (a) of the lemma. 

{\em (b):} By (\ref{eq:E-top-h}), 
\[
\begin{split}
& r^{-1} \scrE_{top}^{\mu _r}(A_r, \Psi
_r)\\
& \qquad =\frac{1}{4r}\int_{X}F_{A_0}\wedge
  F_{A_0}+  \frac{i}{4r}\int_{X}F_{A_0}\wedge\mu_r\\
& \qquad \qquad -2r^{-1}\sum_{i\in \grY}\op{CSD}_{\mu_{i,r}}^{Y_i}([(B_{i,r},
\Phi_{i,r})]_c)+i^*[\nu ]\cdot h_{A_0}(\fh)/2;\\
& r^{-1} \scrE_{top}^{\mu _r, \hat{\grp}}(A_r, \Psi
_r) \\
& \qquad =\frac{1}{4r}\int_{X}F_{A_0}\wedge
  F_{A_0}+  \frac{i}{4r}\int_{X}F_{A_0}\wedge\mu_r\\
& \qquad \qquad -2r^{-1}\sum_{i\in \grY}\op{CSD}_{\mu_{i,r}, \grq_i}^{Y_i}([(B_{i,r},
\Phi_{i,r})]_c)+i^*[\nu ]\cdot h_{A_0}(\fh)/2, 
\end{split}
\]
where \(h_{A_0}\) is the map defined in Remark \ref{rem:based-rel-htpy}.

Combined with Assertion (a) of the lemma and Lemma
\ref{lem:CSD-est}, this gives
\begin{equation}\label{eq:bbE}
\begin{split}
\bbE & :=\lim_{r\to \infty}\big(r^{-1} \scrE_{top}^{\mu _r, \hat{\grp}}(A_r, \Psi
_r)\big)=\lim_{r\to \infty}\big(r^{-1} \scrE_{top}^{\mu _r}(A_r, \Psi
_r)\big)\\ 
&  = \frac{i}{4}\int_{X}F_{A_0}\wedge\nu 
+\pi \sum_{i\in \grY_m}[\nu _i]\cdot \jmath_h (\pmb{\gamma }_i)+ i^*[\nu ]\cdot h_{A_0}(\fh)/2. 
\end{split}
\end{equation}
\epf

\begin{remarks}\label{rmk:E-dep}
(a)  By Theorem \ref{thm:strong-t} and the fact that there are finitely
  many t-orbits with a fixed \(\Spin^c\) structure, the bound
  (\ref{ineq:CSD}) holds for any sequence of solutions to the 3-dimensional
  Seiberg-Witten equations \(\grF_{\mu _{i,r}}(B_r, \Phi _r)=0\).

(b) By Item (b) of the previous lemma, the constant \(\smE\) in Lemma
\ref{lem:E_topX} may be chosen to depend only on
(\ref{dep:rel-htpy}) (though \(r_0\) may still depend on \(\|\nu
\|_{C^1}\) and \(\varsigma_w\)). 
\end{remarks}

To obtain the first group of results mentioned in the beginning of the
present section, we need a generalization of Lemma \ref{lem:E_topX} to general
\(X_\bullet\subset X\), with counterparts of the coefficient \(\smE\) in
(\ref{eq:CSD-est}) independent of both \(r\) and \(X_\bullet\). This
is much more difficult to achieve, mainly due to the fact that the
perturbation form \(\nu \) is not translation-invariant on the ends of
\(X\); cf. the second paragraph of Section \ref{sec:literature}. In
fact, instead of bounding \(\scrE_{top}^{\mu _r}\) and
\(\scrE_{top}^{\mu _r, \hat{\grp}}\), we find it more convenient to
work with a modified version of them, which agree with them in the
case when \(\nu \) is translation-invariant on the ends.

First, on each
\(Y_{i:s}\subset X-X_c\),  let \(\nu _+=\nu_+(s)\), \(\mu _+=\mu_+(s)\) respectively
denote the \(s\)-dependent closed 2-forms on \(Y_i\):
\[
\nu _+(s):=2\nu ^+|_{Y_{i:s}}; \quad \mu _+(s)=2(\mu _r^+)|_{Y_{i:s}}=r\nu _++w_{i,r}
\]
Modifying (\ref{eq:E-top-q}), we set
\begin{equation}\label{def:E'}
\begin{split}
\scrE_{top}^{' \mu_r}(X_{\bullet})(A, \Psi) & :=
 \scrE_{top}^{\mu_r}(X_\bullet ) (A, \Psi)+\frac{i}{4}\int_{\partial \ov{X_\bullet}
 }(A-A_0)\wedge (*_4\mu_r)\\
&=\frac{1}{4}\int_{X_{\bullet}}F_{A_0}\wedge
  F_{A_0}+  \frac{i}{4}\int_{X_{\bullet}}F_{A_0}\wedge\mu_r -2\op{CSD}_{\mu_+}^{\partial \ov{X_\bullet }}(B, \Phi).\\
\end{split}
\end{equation}
Define \(\scrE_{top}^{' \mu_r,\hat{\grp}}(X_{\bullet})\) similarly by
replacing the term \(\scrE_{top}^{\mu_r}(X_\bullet )\) above with 
\(\scrE_{top}^{\mu_r, \hat{\grp}}(X_\bullet )\). 
Note that \[
  \scrE_{top}^{'\mu_r}(X)= \scrE_{top}^{\mu_r}(X);\qquad
\scrE_{top}^{'\mu_r,\hat{\grp}}(X)= \scrE_{top}^{\mu_r,
  \hat{\grp}}(X).\]
Therefore the bounds in (\ref{eq:CSD-est}) hold
for \(\scrE_{top}^{'\mu_r}(X)\),
\(\scrE_{top}^{'\mu_r,\hat{\grp}}(X)\) as well.

A preliminary version of the aforementioned generalization of Lemma
\ref{lem:E_topX}  is stated in terms of \(\scrE_{top}^{'\mu _r}\) and
\(\scrE_{top}^{'\mu _r, \hat{\grp}}\) as follows:

\begin{lemma}\label{lem:Etop-bdd1}
Let \((A, \Psi )=(A_r, \Psi _r)\) be an admissible solution to the Seiberg-Witten
equation \(\grS_{\mu _r,
  \hat{\grp}}(A_r, \Psi _r)=0\) on \(X\) that satisfies in addition 
\begin{equation}\label{assume:EtopX-ubdd}
  \text{either (a)} \,\, \scrE_{top}^{\mu_r}(X)(A_r,
    \Psi_r)\leq \smE \, r\quad \text{ or (b)}\, \,  \scrE_{top}^{\mu_r, \hat{\grp}}(X)(A_r,
    \Psi_r)\leq \smE \, r
\end{equation}
for
a positive constant \(\smE\) independent of \(r\) and \((A, \Psi
)\). (In particular, according to Lemma \ref{lem:E_topX} and Remarks \ref{rmk:E-dep}, this
holds for those \((A_r, \Psi _r)\) from the statement of Theorem
\ref{thm:l-conv}, with \(\smE\) determined by (\ref{dep:rel-htpy}) via
(\ref{eq:bbE}).)  
Then there exist constants
\(\zeta>0\) \(\zeta'> 0\),  a function \(\pmb{\hatl}\co \grY\to \bbR^+\) depending only
on \(\nu \),  and an \(r_0>8\) depending only on the constants
\(\smE\), \(\grl'_i\), \(i\in \grY_v\), 
and \(\nu \), such that for all \(r\geq r_0\) and 
\(\bfl_r:=(\ln r)\, \pmb{\hatl}\),
 \begin{equation}\label{eq:E_top-M1}
\begin{split}
\mathrm{(i) } & \begin{cases}
-\zeta 'r\ln r\leq \scrE_{top}^{' \mu_r }(X_\bullet)(A_r, \Psi_r)\leq
r(\smE+\zeta) &\text{ assuming (\ref{assume:EtopX-ubdd}) (a)}\\
-\zeta '_pr\ln r\leq \scrE_{top}^{' \mu_r, \hat{\grp}}(X_\bullet)(A_r, \Psi_r)\leq
r(\smE+\zeta_p)  &\text{assuming (\ref{assume:EtopX-ubdd}) (b)}
\end{cases}
\\& \qquad \qquad \text{\(\forall X_\bullet\supset X_{\bfl_r}\)};\\
\mathrm{(ii) } & -\zeta '_5 r\leq \scrE_{top}^{' \mu_r} (X_\bullet)(A_r, \Psi_r) =\scrE_{top}^{' \mu_r, \hat{\grp}}(X_\bullet)(A_r, \Psi_r)\leq r(\smE+ \zeta    _e\ln r)\\
& \qquad \qquad  \text{\rm \(\forall X_\bullet\subset X-\mathring{X}_{\bfl_r}\)}.
\end{split}
\end{equation}
The positive constants 
\(\zeta \), \(\zeta '\), \(\zeta '_5\), \(\zeta _e\)  above depend only on the metric, the \(\Spin^c\)
structure, \(\nu \), \(\varsigma _w\), \(B_0\), and the constants \(\zzz_i\)
in Lemma \ref{lem:F-L_1}. \(\zeta _p\), \(\zeta '_p\) are positive
constants that depend only on the preceding
list of parameters, together with the constant \(\zzz_p\) in
Assumption \ref{assume}. 
In particular, these constants as well as \(\pmb{\hatl}\) and \(r_0\) above
 are independent of \(X_\bullet\) and \(r\). 
\end{lemma}
A proof of the preceding lemma will be given in Section
\ref{sec:pf-E_top-bdd0}. The  undesirable factors of \(\ln r\) in
(\ref{eq:E_top-M1}) appear due to the previously-mentioned trouble
with non-translation invariant \(\nu \) at the ends
(cf. (\ref{bdd:CSD-lower}) and remarks that follow). They 
 will eventually be removed in Section \ref{sec:improved}. (See Proposition
\ref{prop:integral-est2}.)  


Looking ahead, the definitions and lemmas above are relevant to the first group of
results mentioned in the beginning of this section in the following
manner: The typical first step towards such results is to use the
relation between analytical and topological energies to bound  \(L^2(X_\bullet)\)-norms of 
gauge invariant terms such as \(F_A\) and \(\nabla_A\Psi \) in terms
of the topological energy over \(X_\bullet\). (In our context, this appears as  Lemmas \ref{lem:E-top-bdd0} and
\ref{E-top-bdd-ends} below.) Upper bounds on the topology energy over \(X_\bullet\), such
as those given in Lemma \ref{lem:Etop-bdd1} and Proposition
\ref{prop:integral-est2} in our context, then give rise to upper
bounds on the aforementioned squares of \(L^2(X_\bullet)\)-norms of 
gauge invariant terms. The latter are used to obtain \(L^2_{1,
  A}(X_\bullet)\) bounds on Seiberg-Witten solutions. (This is done in
Section \ref{sec:L^_1-bdd} in our context.)

In existing literature,\((X, \nu )\) is cylindrical on the ends, and upper
bounds on the topological energy over \(X_\bullet\) follows directly
from a corresponding bound over \(X\). (The analog of  Lemma
\ref{lem:E_topX} in our context.) This is due to the fact that the \(\op{CSD}\)
functional is decreasing on a Floer trajectory, which can be
interpreted as the gradient flow line of the \(\op{CSD}\)
functional. In the more general
setting of ours,  local upper bounds for \(\scrE_{top}^\mu\) likewise follows from
a lower bound on \(\scrE_{top}^\mu\) for \(X_\bullet\) that are 
contained in the cylindrical ends of \(X\); see Section
\ref{sec:E_top-lower}. This lower bound makes use of an
interpretation of the Floer trajectory as a gradient flow line of a
{\em time-dependent} \(\op{CSD}\)-functional; see Section
\ref{sec:SW-grad_flow} below.



\subsection{  \(L^2_{local}\) bounds on gauge invariant terms in terms of
 (modified) topological energies 
}\label{sec:4.2}

Let \((A_r, \Psi_r)=(A, \Psi)\) be an admissible solution to the Seiberg-Witten
equation \(\grS_{\mu _t, \hat{\grp}}(A_r, \Psi _r)=0\) with
\(Y_i\)-end limit \((B_i, \Phi _i)\), 
and let \(A_0\) \(\{B_{0,i}\}_i\) be the reference
connections fixed previously. 
Recall the definition of \(|X_\bullet|\) from Section
\ref{sec:convention} and let 
\[
|X_\bullet|_1:=\min \{ 1, |X_\bullet|\}. 
\]
Recall also the definitions of \(X^{'a}\), \(X''\) from Definition
\ref{def:adm}, and note that \(\hat{\grp}(A,
\Psi)\) is supported on the vanishing ends; in fact, on
\(X-X''\supset X-X^{'a}\). Let
\[
X_{\bullet, v}:=X_\bullet \cap (X-X'')\qquad \text{and}\qquad X_{\bullet
  ,m}:=X_\bullet\cap X''. 
\]

In this subsection, we show that:
\begin{lemma}\label{lem:E-top-bdd0}
Let \(X_\bullet\subset X\) be compact 
and let \(r\geq 1\). 
Then there exist
  positive constants \(\zeta _0\), \(\zeta '\), \(\zeta '_p\), \(\zeta ''\), \(\zeta
  ''_p\) that depend only on 
\BTitem\label{parameters}
\item the metric, 
\item the   \(\Spin^c\)-structure on \(X\), 
\item the cohomology class of \(\nu \), 
\item the constant \(\varsigma _w\)  in Assumption \ref{assume}, 
\item the constant \(\zzz_\grp\) in  Assumption \ref{assume}, 
\ETitem 
 such that the following hold: 
 \begin{equation}\label{eq:E-top-bdd0}
\begin{split}
\text{(a)} & \quad
\frac{1}{8}\int_{X_{\bullet}}|F_A|^2+\int_{X_{\bullet}}|\nabla_A\Psi|^2+\frac{1}{2}\int_{X_{\bullet}}\big|\frac{i}{4}\rho(\mu_r^+)-(\Psi\Psi^*)_0\big|^2\\
&\qquad \quad 
\leq \scrE_{top}^{' \mu_r }(X_{\bullet})(A, \Psi)+r\, (\zeta _0|
X'_{\bullet,m}|+\zeta '|X_{\bullet}|_1)+\zeta''|X_{\bullet}|;\\
\text{(b)} &\quad \frac{1}{16}\int_{X_{\bullet}}|F_A|^2+\int_{X_{\bullet}}|\nabla_A\Psi|^2+\frac{1}{4}\int_{X_{\bullet}}\big|\frac{i}{4}\rho(\mu_r^+)-(\Psi\Psi^*)_0\big|^2\\
&\qquad \quad 
\leq \scrE_{top}^{' \mu_r,  \hat{\grp}}(X_{\bullet})(A, \Psi)+r\, (\zeta _0|
X'_{\bullet,m}|+\zeta '_p|X_{\bullet}|_1)+\zeta''_p |X_{\bullet}|. 
\end{split}
\end{equation}
\end{lemma}

\pf 
To begin, combine 
(\ref{eq:E-top=anal}) with the Seiberg-Witten equation \\ \(\grG_{\mu_r,
  \hat{\grp}}(A_r, \Psi_r)=0\) to get: 
\begin{equation}\label{X_c-energy}
\begin{split}
& \frac{1}{4}\int_{X_{\bullet}}|F_A|^2+\int_{X_{\bullet}}|\nabla_A\Psi|^2+\int_{X_{\bullet}}\big|\frac{i}{4}\rho(\mu_r^+)-(\Psi\Psi^*)_0\big|^2\\
&\quad 
=-\int_{X_{\bullet}}\frac{R_g}{4}|\Psi|^2+\frac{i}{4}\int_{X_{\bullet}}F_A\wedge*_4\mu_r+\scrE_{top}^{\mu_r}
(X_{\bullet})(A,\Psi)\\
& \qquad \qquad +\|\hat{\grp}(A, \Psi)\|_{L^2(X_{\bullet})}^2\\
& \quad =-\int_{X_{\bullet}}\frac{R_g}{4}|\Psi|^2+\|\hat{\grp}(A,
\Psi)\|_{L^2(X_{\bullet})}^2+\frac{i}{4}\int_{X_{\bullet}}F_A\wedge
*_4w_r\\
&\qquad \quad +\frac{ir}{4}\int_{X_{\bullet
    }}F_{A_0}\wedge*_4\nu
+\scrE_{top}^{' \mu_r }
(X_{\bullet})(A, \Psi).
\end{split}
\end{equation}
The terms  in the third and fourth line above are bounded in Steps
(1)-(3) below. 

(1)  The terms \(-\int_{X_{\bullet}}\frac{R_g}{4}|\Psi|^2\), \( \|\hat{\grp}(A,
\Psi)\|_{L^2(X_{\bullet})}^2\) are bounded by  the same general trick.
Observe that for any real
valued function \(\grf\),
\[\begin{split}
0& \leq
\Big|\frac{i}{4}\rho(\mu^+_r)-\big(1-\frac{\grf}{2|\Psi|^2}\big)(\Psi\Psi^*)_0\Big|^2\\
& \leq\Big|\frac{i}{4}\rho(\mu^+_r)-(\Psi\Psi^*)_0\Big|^2+\frac{\grf^2}{4}-\grf|\Psi|^2+\frac{|\grf|}{8}|\mu_r|
\end{split}\]
where \(\Psi\neq 0\). 
So, for general \(\Psi \),
\begin{equation}\label{f-trick}
\grf|\Psi|^2\leq \Big|\frac{i}{4}\rho(\mu^+_r)-(\Psi\Psi^*)_0\Big|^2+\frac{\grf^2}{4}+\frac{|\grf|}{8}|\mu_r|.
\end{equation}
Taking  \(\grf=\pm R_g\)  in (\ref{f-trick}) then gives us 
\begin{equation}\label{bdd(1)}
\begin{split}
  & \Big|-\int_{X_{\bullet}}\frac{R_g}{4}|\Psi|^2\Big| \\
  & \qquad \leq
\frac{1}{4}\Big\|\frac{i}{4}\rho(\mu^+_r)-(\Psi\Psi^*)_0\Big\|^2_{L^2(X_\bullet)}+\int_{X_{\bullet}}\frac{R_g^2}{16}+\int_{X_{\bullet}}\frac{|R_g|}{32}|\mu_r|\\
& \qquad \leq
\frac{1}{4}\Big\|\frac{i}{4}\rho(\mu^+_r)-(\Psi\Psi^*)_0\Big\|^2_{L^2(X_\bullet)}+\zeta
_g|X_\bullet|+\zeta _h \, r(|X_{\bullet,m}|+|X_{\bullet}|_1). 
\end{split}
\end{equation}
In the above, \(\zeta _g\), \(\zeta _h\) are positive constants
depending only on the metric, \(\varsigma_w\), and \(\nu \). To go from the
first line to the second line, we also made use
of the well-known fact that \(\nu \) exponentially decays to \(\nu _i\) on the
\(Y_i\)-end for each \(i\in \grY\) (cf. (\ref{eq:xi-exp})). 

Meanwhile by  (\ref{eq:q-bdd}), 
\[\begin{split}
    & \|\hat{\grp}(A, \Psi )\|_{L^2(X_\bullet)}^2 
   \leq z\big( \|
\Psi\|_{L^2(X_{\bullet, v})}^2+|X_{\bullet, v}|\big),
\end{split}\]
where \(z\) is a positive constant depending only on the constant
\(\zzz_\grp\) and the metric on the vanishing ends. Taking
\(\grf=4z\) in (\ref{f-trick}) then gives us 
\begin{equation}\label{bdd(1')}\begin{split}
    & \|\hat{\grp}(A, \Psi )\|_{L^2(X_\bullet)}^2 \\
    & \qquad \leq
\frac{1}{4}\Big\|\frac{i}{4}\rho(\mu^+_r)-(\Psi\Psi^*)_0\Big\|^2_{L^2(X_{\bullet,
    v})}+z_1|X_{\bullet, v}|_1+ z_2 r\int_{X_{\bullet, v}} |\nu |\\
&\qquad  \leq \frac{1}{4}\Big\|\frac{i}{4}\rho(\mu^+_r)-(\Psi\Psi^*)_0\Big\|^2_{L^2(X_{\bullet,
    v})}+z_1|X_{\bullet, v}|_1+ z_3 r|X_{\bullet,v}|_1, 
\end{split}
\end{equation}
noting that \(|\nu |\) decays exponentially to 0 on vanishing
ends. (For a more precise statement, see (\ref{eq:xi-exp}).) In
the above, \(z_i\) are positive constants that depend only on \(\nu\),
the metrics, and the constants \(\varsigma _w\), \(\zzz_\grp \).

(2) The term \(\frac{i}{4}\int_{X_{\bullet}}F_A\wedge *_4w_r \) in
(\ref{X_c-energy}) is bounded as follows. 
\begin{equation}\label{bdd(2)}
\frac{i}{4}\int_{X_{\bullet}}F_A\wedge *_4w_r\leq
\frac{1}{16}\|F_A\|_{L^2(X_\bullet)}^2+\|w_r\|_{L^2(X_\bullet)}^2\leq \frac{1}{16}\|F_A\|_{L^2(X_\bullet)}^2+z_4|X_\bullet|,
\end{equation}
where \(z_4\) depends only on the constant \(\varsigma _w\).

(3) The term \(\frac{ir}{4}\int_{X_{\bullet}}F_{A_0}\wedge*_4\nu\) in
(\ref{X_c-energy}) is bounded as follows. Recall that by assumption,
\(F_{A_0}\) is the pull-back of \(F_{B_{0, i}}\) on the ends
\(\hat{Y}_i\), and  that \(\wp_i=0\) for vanishing ends. Thus, 
\begin{equation}\label{bdd(3)}
\begin{split}
& \frac{ir}{4}\int_{X_{\bullet}}F_{A_0}\wedge*_4\nu \\
& \qquad =\frac{ir}{4}\int_{X_{\bullet}\cap
X_c}F_{A_0}\wedge*_4\nu+\sum_{i\in \grY}
r\wp_i|X_{\bullet}\cap\hat{Y}_i|\\
& \qquad \qquad \quad +r\sum_{i\in \grY}
\int_{X_{\bullet}\cap\hat{Y}_i}F_{A_0}\wedge*_4(\nu-\nu _i)\\
& \qquad = \frac{ir}{4}\int_{X_{\bullet}\cap
X_c}F_{A_0}\wedge*_4\nu+\sum_{i\in \grY_m}
r\wp_i|X_{\bullet}\cap\hat{Y}_i|+r\, \zeta _5 \sum_{i\in \grY}|X_{\bullet}\cap\hat{Y}_i|_1\\
& \qquad \leq r\, (z_5 |X_{\bullet, m}|+\zeta '_5). 
\end{split}
\end{equation}
Here, \(\zeta _5\) and \(\zeta _5'\) are positive constants  depending only on the choice of
\(A_0\) and \(\nu \), while \(z_5\) is a  positive constant  depending only on the choice of
\(A_0\), the \(\Spin^c\)-structure, and \(\nu \). In estimating the
last term in the second line above, we also used the exponential decay
of \(\nu \) on the ends of \(X\). (Cf. (\ref{eq:xi-exp}) for a precise statement.)

Combining the bounds (\ref{bdd(1)}), (\ref{bdd(1')}), (\ref{bdd(2)}),
(\ref{bdd(3)}) with (\ref{X_c-energy}) and re-arranging, we arrive at
item (a) of  (\ref{eq:E-top-bdd0}). The inequality in item (b) 
follows from item (a) when \(X_\bullet\subset X^{'a}\). By the additivity of
\(\scrE_{top}^{' \mu_r ,\hat{\grp}}\), it remains to verify the
inequality on the vanishing ends. The proof in this case is similar and is deferred
to the next subsection. 
\epf 

Recall the definition of \(\wp_i\) from Lemma \ref{lem:F-L_1}. For each \(i\in \grY\), set 
\[
\wp_i^+:=
\begin{cases}
\wp_i+\frac{\pi }{2}\|*\nu _i\|_T& \text{when \(\hat{Y}_i\) is a Morse
  end};\\
0 & \text{when \(\hat{Y}_i\) is a vanishing end.}
\end{cases}
\]
Note that according to Proposition \ref{prop:t-conv3d} and Remarks
\ref{rmk:thurston}, 
\(\wp_i^+\geq \frac{\pi }{2}(\|*\nu _i\|_T+\zeta _{*\nu _i})\) in our
context, with the right hand side, \(\frac{\pi }{2}(\|*\nu
_i\|_T+\zeta _{*\nu _i})\) vanishing in many cases, such as in the
context of \cite{LT}.

When \(X_\bullet=\hat{Y}_{i, [l, L]}\), the bound
(\ref{eq:E-top-bdd0}) (a) can be refined as follows. 

\begin{lemma}\label{E-top-bdd-ends}
  Let \((A, \Psi )=(A_r, \Psi _r)\) be as in Lemma
  \ref{lem:E-top-bdd0}. 
When \(X_\bullet=\hat{Y}_{i, [l, L]}\), we have: 
\begin{equation}\label{eq:E-top-bdd2}
\begin{split}
& \frac{1}{8}\int_{X_{\bullet}}|F_A|^2+\int_{X_{\bullet}}|\nabla_A\Psi|^2+\frac{1}{2}\int_{X_{\bullet}}\Big|\frac{i}{4}\rho(\mu_r^+)-(\Psi\Psi^*)_0\Big|^2\\
&\qquad \quad 
\leq \scrE_{top}^{' \mu_r }(X_{\bullet})(A, \Psi)+r\, ( \wp_i^++z_i)
|X_{\bullet}|+rz'_i|X_\bullet|_1+\zeta ''|X_\bullet|,
\end{split}
\end{equation}
where \(z_i\) is a non-negative constant depending only on the metric
on \(Y_i\), and \(z'_i\) is a positive constant that only depends on the
metric on \(X\), the cohomology class of \(\nu \), and the constant \(\varsigma _w\) in Assumption
\ref{assume}. (In particular, \(z_i\) and \(z'_i\) are independent of \(r\) and the \(\Spin^c\)
structure.) Moreover, the constant \(z_i\) has the following
properties: \(z_i=0=\wp^++z_i\) when \(\hat{Y}_i\) is a vanishing
end. When \(\hat{Y}_i\) is a Morse end and \(Y_i\) is irreducible and atoroidal, then for any \(\epsilon >0\), there exists a metric
on \(Y_i\) so that its corresponding \(z_i\) is less than \(\epsilon
\). Meanwhile, \(\zeta ''\) is a positive constant that depends only
on the metric and \(\varsigma_w\). 

\end{lemma}
\pf
To begin, recall  \cite{KM:sc}'s
Theorem 2. Reformulated as \cite{BD}'s Theorem 5.4, it asserts that
when \(Y_i\) is irreducible, and atoroidal, 
\begin{equation}\label{eq:KM}
\|*\nu _i\|_T=\frac{1}{4\pi }\inf_h \{\|R_h\|_h \, \|\theta _h\|_h\},
\end{equation}
where \(h\) ranges through all Riemannian metrics on \(Y_i\). In the above, \(\|\cdot\|_h\)
denotes \(L^2\)-norm with respect to \(h\), and \(\theta _h\) denotes
the harmonic representative of the cohomology class \([*\nu _i]\) with
respect to the metric \(h\). \(R_h\) denotes the scalar curvature on
\(Y_i\) corresponding to \(h\). Recall also that \(\|*\nu
_i\|_T=\|[*\nu _i]\|_T\) is the Thurston semi-norm of \([*\nu _i]\in
H^1(Y_i;\bbR)\). 

By assumption, over \(\hat{Y}_i\simeq [0, \infty)\times Y_i\), \(g\) is 
the product metric of the standard affine metric on \(\bbR\) and a
metric \(g_i\) on \(Y_i\). Thus, \(R_g|_{Y_{i: s}}=R_{g_i}\) for any
\(s\geq 0\). Combined with the asymptotic behavior of
\(\nu \) (described explicitly in (\ref{eq:dnu}) below), (\ref{eq:KM}) implies that for any
\(s\geq 0\),
\[
\|R_g\|_{L^2(Y_{i: s})}\|\nu|_{Y_{i: s}}\|_{L^2(Y_{i: s})}\leq 4\pi  \|*\nu _i\|_T+z_i+\grz_i(s),
\] 
where \(z_i=z_{Y_i, \nu _i}:=\|R_{g_i}\|_{L^2(Y_i)}\|\nu _i\|_{L^2(Y_i)}-4\pi \|*\nu
_i\|_T\).  According to (\ref{eq:KM}), \(z_i\) has the properties
asserted in the statement of the lemma. 
\( \grz_i(s)\) is
a positive function on \([0, \infty)\) that depends only on \(\nu \)
and exponentially decays to \(0\) as \(s\to \infty\). (Explicitly, it
is a multiple of the right hand side of (\ref{eq:xi-exp}).)
Use the preceding inequality to replace (\ref{bdd(1)}) in this context by: 
\begin{equation}\label{bdd(1-)} 
\begin{split}
& \Big|-\int_{X_{\bullet}}\frac{R_g}{4}|\Psi|^2\Big|\\
& \quad \leq
\frac{1}{4}\Big\|\frac{i}{4}\rho(\mu^+_r)-(\Psi\Psi^*)_0\Big\|^2_{L^2(X_\bullet)}+\int_{X_{\bullet}}\frac{R_g^2}{16}+\int_{X_{\bullet}}\frac{|R_g|}{32}|\mu_r|\\
& \quad \leq
\frac{1}{4}\Big\|\frac{i}{4}\rho(\mu^+_r)-(\Psi\Psi^*)_0\Big\|^2_{L^2(X_\bullet)}+\zeta
_1|X_\bullet|\\
& \qquad \qquad +\frac{r}{32}\big(\int _l^Lds \|R_g\|_{L^2(Y_{i:
    s})}\|\nu|_{Y_{i: s}}\|_{L^2(Y_{i: s})}\big)\\
& \quad \leq
\frac{1}{4}\Big\|\frac{i}{4}\rho(\mu^+_r)-(\Psi\Psi^*)_0\Big\|^2_{L^2(X_\bullet)}+\zeta
_1|X_\bullet|\\
& \qquad \qquad + r\, \big( \frac{\pi }{2}\|*\nu
_i\|_T+z_i\big)|X_\bullet|+ r\zeta _2|X_\bullet|_1, 
\end{split}
\end{equation}
where \(\zeta _1\), \(\zeta _2\) are positive constants depending only on \(R_g\)
and \(\varsigma_w\). (In particular, they are independent of \(r\) and
\(\grs_i\).) 

Meanwhile, replace (\ref{bdd(3)}) in this context 
by
\begin{equation}\label{bdd(3)'}
\frac{ir}{4}\int_{X_{\bullet}}F_{A_0}\wedge*_4\nu  \leq
r\wp_i|X_{\bullet}|+r\zeta _5|X_\bullet|_1,
\end{equation}
\(\zeta _5\) being a positive constant depending only on \(\nu \). 
Combine the bounds (\ref{bdd(1-)}), (\ref{bdd(2)}),
(\ref{bdd(3)'}) with (\ref{X_c-energy}) and re-arranging, we arrive at
(\ref{eq:E-top-bdd2}). 
\epf

\begin{remarks}
Since both \(r\) and \(|X_\bullet|\) can be arbitrarily large,
applying (\ref{eq:E-top-bdd2}) to the case when \(X=\bbR\times Y\) and
\((A_r, \Psi _r)=(\hat{B}_r, \hat{\Phi }_r)\) where \((B_r, \Phi _r)\)
is a sequence of solutions to the 3-dimensional Seiberg-Witten
equation \(\grF_{\mu _r}(B_r, \Phi _r)=0\) from Theorem
\ref{thm:strong-t} implies that \(c_1(\grs)\cdot [*\nu ]+\|*\nu
\|_T\geq -\frac{2z_{Y, \nu }}{\pi }\) when the
sequence of \((B_r, \Phi _r)\) in Theorem \ref{thm:strong-t} exists.
This bound is comparable to  \cite{KM}'s
Propositions 40.1.1 and 40.1.3, where certain non-existence results of
Seiberg-Witten solutions 
are obtained under similar constraints on \(c_1(\grs)\).  In
particular, in the case when \(b_1(Y)>0\),  the
aforementioned propositions in \cite{KM} were used to prove that for irreducible
\(Y\), the dual
Thurston polytope is the convex hull of the  ``Seiberg-Witten
basic classes'' (cf. \cite{KM:sc} Theorem 1 and \cite{KM} Theorem
41.5.2).
\end{remarks}

\subsection{Seiberg-Witten solutions as gradient flow lines of
  time-dependent CSD functionals}\label{sec:SW-grad_flow}

We begin with some preliminary observations on the behavior of the form \(\nu \)
on the ends, and introduce some notations along the way. In
particular, the previously imentioned exponential decay of
\(\nu \) on the ends of \(X\) is made precise. 

Consider an end \(\hat{Y}_i\)
of \(X\). 
Write the harmonic 2-form
\(\nu\) as \[\nu=\underline{\nu}+ds \wedge *_Y\grv \ \quad \text{on
  \(\hat{Y}_i\),}
\]
where \(\underline{\nu}\), \(\grv\) are \(s\)-dependent 2-forms on \(Y:=Y_i\),
\(*_Y\) denotes the 3-dimensional Hodge dual on \(Y_i\), and \(s=\mathfrc{s}_i\). The harmonicity of \(\nu\) implies
that both \(\underline{\nu}\) and \(\grv\) are closed 2-forms. 
Moreover, \(\grv\) and \(\underline{\nu}-\nu_i\) are \(s\)-dependent
exact 2-forms on \(Y_{i:s}=\mathfrc{s}_i^{-1}(s)\), and satisfy 
\begin{equation}\label{eq:dnu}
\partial_s(\underline{\nu}-\nu_i)=-d_Y*_Y\grv\quad
  \text{and}\quad \partial_s\grv=-d_Y*_Y(\underline{\nu}-\nu_i).
\end{equation}
Here, \(d_Y\) denotes the
exterior derivative on \(Y_i\).

Let 
\[
\xi_\nu:=\underline{\nu}-\nu_i+\grv;\quad 
\xi_\nu':=\underline{\nu}-\nu_i-\grv,
\] 
so that in this notation, \(\nu _+=\nu _i+\xi_\nu \) on
\(\hat{Y}_i\). \(\xi_\nu\) and \(\xi_\nu '\) are exact
2-forms satisfying 
\begin{equation}\label{DE:xi}
\partial _s\xi_\nu=-d_Y*_Y\xi_\nu; \quad \partial
_s\xi_\nu'=d_Y*_Y\xi_\nu'.
\end{equation}
As \(d_Y*_Y\) is a self-dual operator on
the space of \(L^2\) exact 2-forms and \(Y_i\) is closed, it has
discrete real eigenvalues \(\{\kappa_j^{(i)}\}_j\) with
\(\min_j|\kappa_j^{(i)}|=:\kappa_i>0\). Thus, there exists a constant \(\zeta_i>0\)
such that 
\begin{equation}\label{eq:xi-exp}
\begin{split}
& \|\xi_\nu\|_{C^k(Y_{i:s})}+\|\ud{\nu}-\nu_i\|_{C^k(Y_{i:s})}+\|\grv\|_{C^k(Y_{i:s})}\leq
\zeta_ie^{-\kappa _is}; \quad \text{and via (\ref{DE:xi}), }\\
& \|\nu-\pi _2^*\nu_i\|_{C^k(Y_{i,s})}\leq \, \zeta_ie^{-\kappa _is}.
\end{split}
\end{equation}

The (\(s\)-dependent) form \(\nu _+\) 
on \(Y_i\), introduced previously in Section \ref{sec:top-energy}, is 
expressed in this part's notation as: 
\[
\begin{split}
\nu_+ & =\ud{\nu}+\grv. \\
& =\nu _i+\xi _\nu .
\end{split}
\]

Returning now to the task of interpreting of Seiberg-Witten solutions
on the ends of \(X\), we choose to work in a gauge such that \((A,
\Psi)|_{\hat{Y}_i}\) is in the temporal gauge over each end \(\hat{Y}_i\)
of \(X\). Namely, over the ends \(X-X_c\) we can write 
\begin{equation}\label{temp-gauge}
(A, \Psi)=(\partial_s+B,
\Phi),
\end{equation}
where \((B(s), \Phi(s))\), \(s\geq 0\), is a path in \(\Conn (Y_i)\times
\Gamma (\bbS_{Y_i})\). Recall that the reference connection \(A_0\) is
chosen such that
over  the ends of \(X\), \(A_0=\partial_s +B_0\), where
\(B_0=B_{0,i}\) on the \(Y_i\)-end.

Restricting to an end \(\hat{Y}_i\) of \(X\),  
the 4-dimensional
Seiberg-Witten equation \(\grS_{\mu_r, \hat{\grp}}(A, \Psi)=0\) 
may be re-expressed in terms of \((B(s), \Phi(s))\) over \(\hat{Y}_i\)
as 
\begin{equation}\label{t-dep-flow}
\big(\frac{1}{2}\partial_sB,
\partial _s\Phi\big)+\mathfrak{F}_{\mu_+, \hat{\grp}}(B, \Phi)=0,
\end{equation}
where \(\hat{\grp}=\hat{\grp}(s)=\chi_i(s)\grq_i+\lambda_i(s)\grp'_i\) is regarded as a path of tame
perturbations as in (\ref{eq:hatp}). 

Square both sides of the previous equation and use integration by parts plus
(\ref{eq:dnu}) 
to get:
\begin{equation}\label{eq:flow-energy}\begin{split} 
& \|\partial_sB\|_{L^2(\hat{Y}_{[l, L]})}^2 +4\|\partial_s\Phi\|^2_{L^2(\hat{Y}_{[l, L]})}+4\|\mathfrak{F}_{\mu_+,
      \hat{\grp}} (B, \Phi)\|_{L^2(\hat{Y}_{[l,
    L]})}^2\\
& \quad\quad   = 8\op{CSD}_{\mu_+(l), \hat{\grp}(l)}^{Y}(B,
\Phi)|_{Y_l}-8\op{CSD}_{\mu_+(L), \hat{\grp}(L)}^{Y}(B,
\Phi)|_{Y_L}\\
& \qquad\quad 
-ir\int_{\hat{Y}_{[l,L]}}ds\,(  *_Y\xi_\nu)\wedge(F_B-F_{B_0})
-8\int_{\hat{Y}_{[l,L]}}(\partial_sf_{\hat{\grp}(s)})(B, \Phi)
 \,ds 
\end{split}
\end{equation}
for any interval \([l, L]\subset [0, \infty]\). (Here \(Y\) stands for
any of the ending 3-manifolds \(Y_i\).)

In what follows we use the preceding formula in two ways:  (i) In the
remainder of this subsection, we use the upper
bound on the right hand side of (\ref{eq:flow-energy}) to get an upper bound on the square terms
on the left hand side; in particular, this leads to the completion of
the proof of 
 (\ref{eq:E-top-bdd0}) (b). (ii) In the subsequent subsection, we use the non-negativity of the left hand
 side to obtain a lower bound of the first two terms on the right hand
 side, and hence a lower
bound on \(\scrE^{' \mu_r, \hat{\grp}}_{top}(\hat{Y}_{[l,L]}) (A,\Psi)\), 
generalizing the positivity results in the cylindrical case mentioned
towards the end  of  Section \ref{sec:top-energy}. 

\subsubsection*{\it Verifying the remaining case of (\ref{eq:E-top-bdd0}) (b). }
As observed in the previous subsection, it remains to verify the claimed inequality for
the case when 
\(X_\bullet \) lies in vanishing ends. 
Nevertheless, assume for the moment that \(X_\bullet=\hat{Y}_{[l,L]}\),
\(\hat{Y}=\hat{Y}_i\) is a general end of \(X_\bullet\). (Either
vanishing or Morse.)

To begin, recall from Assumption \ref{assume} (3) and use the 
3-dimensional Weitzenb\"ock formula (\ref{eq:Weit3d}) to get 
\[
\begin{split}
& \frac{1}{4}\|\partial_sB\|_{L^2(\hat{Y}_{[l, L]})}^2
+\|\partial_s\Phi\|^2_{L^2(\hat{Y}_{[l, L]})}+\|\mathfrak{F}_{\mu_+} (B, \Phi)\|_{L^2(\hat{Y}_{[l,
    L]})}^2\\
&\quad 
=\frac{1}{4}\int_{\hat{Y}_{[l,L]}}|F_A|^2+\int_{\hat{Y}_{[l,L]}}|\nabla_A\Psi|^2+\int_{\hat{Y}_{[l,L]}}\big|\frac{i}{4}\rho(\mu_r^+)-(\Psi\Psi^*)_0\big|^2\\
& \qquad +\int_{\hat{Y}_{[l,L]}}\frac{R_g}{4}|\Psi|^2-\frac{i}{4}\int_{\hat{Y}_{[l,L]}}F_A\wedge
*_4w_r-\frac{ir}{4}\int_{\hat{Y}_{[l,L]}} ds \, F_B\wedge
(*_Y\nu_+).\\
\end{split}
\] 
Inserting both the preceding formula and (\ref{eq:flow-energy}) into 
the following inequality 
\[\begin{split}
& \frac{1}{4}\|\partial_sB\|_{L^2(\hat{Y}_{[l, L]})}^2
+\|\partial_s\Phi\|^2_{L^2(\hat{Y}_{[l, L]})}+\|\mathfrak{F}_{\mu_+} (B, \Phi)\|_{L^2(\hat{Y}_{[l,
    L]})}^2\\
& \qquad \leq 
\frac{1}{4}\|\partial_sB\|_{L^2(\hat{Y}_{[l, L]})}^2 +\|\partial_s\Phi\|^2_{L^2(\hat{Y}_{[l, L]})}+\|\mathfrak{F}_{\mu_+,
      \hat{\grp}} (B, \Phi)\|_{L^2(\hat{Y}_{[l,
        L]})}^2\\
    &\qquad \qquad\qquad  + \| \hat{\grp} (A, \Psi )\|_{L^2(\hat{Y}_{[l,
    L]})}^2,
\end{split}
\]
we have 
\begin{equation}\label{bdd:sq}
\begin{split}
&
\frac{1}{4}\int_{\hat{Y}_{[l,L]}}|F_A|^2+\int_{\hat{Y}_{[l,L]}}|\nabla_A\Psi|^2+\int_{\hat{Y}_{[l,L]}}\big|\frac{i}{4}\rho(\mu_r^+)-(\Psi\Psi^*)_0\big|^2\\
& \quad \leq  -2 \op{CSD}_{\mu _+, \hat{\grp}}^{\partial
  \hat{Y}_{[l,L]}}(B, \Phi ) +r\wp_i (L-l)
-2\int_{\hat{Y}_{[l,L]}}(\partial_s\, f_{\hat{\grp}(s)})(B, \Phi)\\
& \qquad -\int_{\hat{Y}_{[l,L]}}\frac{R_g}{4}|\Psi|^2+\frac{i}{4}\int_{\hat{Y}_{[l,L]}}F_A\wedge
*_4w_r+\frac{ir}{4}\int_{\hat{Y}_{[l,L]}}ds\,(  *_Y\xi_\nu )\wedge
F_{B_0}\\
& \qquad + \|
\hat{\grp} (A, \Psi)\|_{L^2(\hat{Y}_{[l,
    L]})}^2 .
\end{split}
\end{equation}
Argue as in  (\ref{bdd(1-)}), (\ref{bdd(1')}), (\ref{bdd(2)}), (\ref{bdd(3)})  to bound the
fourth to the last terms on the right hand side of the preceding
formula and rearranging, we get: 
\begin{equation}\label{energy-v0}
\begin{split}
&
\frac{1}{8}\int_{\hat{Y}_{[l,L]}}|F_A|^2+\int_{\hat{Y}_{[l,L]}}|\nabla_A\Psi|^2+\frac{1}{2}\int_{\hat{Y}_{[l,L]}}\big|\frac{i}{4}\rho(\mu_r^+)-(\Psi\Psi^*)_0\big|^2\\
& \quad \leq  -2 \op{CSD}_{\mu _+, \hat{\grp}}^{\partial
  \hat{Y}_{[l,L]}}(B, \Phi ) +r (\wp_i^++z_i)|\hat{Y}_{[l,L]}|+\zeta
r|\hat{Y}_{[l,L]}|_1+\zeta '|\hat{Y}_{[l,L]}|\\
& \quad \quad -2\int_{\hat{Y}_{[l,L]}}(\partial_s\, f_{\hat{\grp}(s)})(B, \Phi) , \\
\end{split}
\end{equation}
where \(\zeta , \zeta '\) are positive constants that depend only on
the choice of reference connection \(A_0\), the metric, \(\nu \), and the constants
\(\varsigma _w\), 
\(\zzz_\grp\). Note that over 
\(\hat{Y}_{[l,L]}\), 
\begin{equation}\label{diff:CSD-E}\begin{split}
& -2 \op{CSD}_{\nu _+, \hat{\grp}}^{\partial
  \hat{Y}_{[l,L]}}(B, \Phi )-\scrE^{' \mu_r  ,\hat{\grp}}(
\hat{Y}_{[l,L]}) (A, \Psi )\\
& \quad =-\frac{i}{4}\int_{\hat{Y}_{[l,L]}}
F_{B_0}\wedge \mu _r
 \leq \zeta
_e r\, e^{-\kappa _il} |\hat{Y}_{[l,L]}|_1 , 
\end{split}
\end{equation}
where  \(\zeta _e\) is a positive constant depending
only on \(B_0\) and \(\nu \). 
Thus, the bound (\ref{energy-v0}) is equivalent to 
\begin{equation}\label{energy-v}
\begin{split}
&
\frac{1}{8}\int_{\hat{Y}_{[l,L]}}|F_A|^2+\int_{\hat{Y}_{[l,L]}}|\nabla_A\Psi|^2+\frac{1}{2}\int_{\hat{Y}_{[l,L]}}\big|\frac{i}{4}\rho(\mu_r^+)-(\Psi\Psi^*)_0\big|^2\\
& \quad \leq  \scrE^{' \mu_r  \hat{\grp}}(
\hat{Y}_{[l,L]}) (A, \Psi ) +r (\wp_i^++z_i)|\hat{Y}_{[l,L]}|+\zeta
_or|\hat{Y}_{[l,L]}|_1+\zeta '_o|\hat{Y}_{[l,L]}|\\
& \quad \quad -2\int_{\hat{Y}_{[l,L]}}(\partial_sf_{\hat{\grp}})(B, \Phi), \\
\end{split}
\end{equation}
where \(\zeta _o, \zeta '_o\) are likewise positive constants that depend only on
the choice of reference connection \(A_0\), the metric, \(\nu \), and the constants
\(\varsigma _w\), \(\zzz_\grp\). 


We now estimate the last term above,
\(-2\int_{\hat{Y}_{[l,L]}}(\partial_s\, f_{\hat{\grp}(s)})(B,
\Phi)\). This vanishes unless \(\hat{Y}_i\) is a vanishing end; so we assume
from now on that \(\hat{Y}=\hat{Y}_i\) is a vanishing end. 
According to (\ref{eq:hatp}) (and recalling the definitions and
assumptions on  \(\chi _i\), \(\lambda _i\) there), 
\[
(\partial_sf_{\hat{\grp}})(B,
  \Phi)=\chi '_i(s)\,f_{\grq_i}(B, \Phi)+\lambda
  '_i (s) \, f_{\grp'}(B, \Phi).
\] 
This is supported on \(\hat{Y}_{[\grl_i,
    \grl_i']}\), and by our assumption on the cutoff functions \(\chi
  _i\) and \(\lambda _i\), its absolute value is bounded by
  \(|f_{\grq_i}(B(s), \Phi(s))|+|f_{\grp'}(B(s), \Phi(s))|\). At each \(s\in [\grl_i, \grl_i']\),
  \(|f_{\grq_i}(B(s), \Phi(s))|\) and \(|f_{\grp'}(B(s), \Phi(s))|\) can
  be bounded via a variant of Lemma \ref{lem:f_q}: Applying the
  triangle inequality differently in the last step of
  (\ref{ineq:f_q}) and recalling Assumption \ref{assume} 4d), one may arrange so that: 
\[\begin{split}
   &  |f_{\grq_i}(B(s), \Phi(s))|+|f_{\grp'}(B(s), \Phi(s))|\\
& \qquad  \leq \frac{1}{32}\|F_B-F_{B_0}\|_{L^2(Y_{i: s})}^2+\zeta _2'\big(  \|\Phi
\|_{L^2(Y_{i: s})}^2+1\big), 
\end{split}
\]
where \(\zeta '_2\) is a positive constants depending only on \(\zzz_\grp\) 
and the metric on \(Y_i\). Integrating the preceding
inequality over \(I:=[l,L]\cap [\grl_i, \grl_i']\subset \bbR\), and
appealing again to (\ref{f-trick}), we get 
\begin{equation}\label{bdd:f-ps}
\begin{split}
& \Big|-2\int_{\hat{Y}_{[l,L]}}(\partial_sf_{\hat{\grp}})(B,
  \Phi)\Big|\\
&\quad  \leq \frac{1}{16}\|F_B-F_{B_0}\|_{L^2(\hat{Y}_{I})}^2+\zeta '_3\big(\|\Psi
  \|^2_{L^2(\hat{Y}_{I})}+ |\hat{Y}_{I}|\big)\\
& \quad \leq
\frac{1}{16}\|F_B-F_{B_0}\|_{L^2(\hat{Y}_{I})}^2+\frac{1}{4}\Big\|\frac{i}{4}\rho(\mu^+_r)-(\Psi\Psi^*)_0\Big\|^2_{L^2(\hat{Y}_I)}+\zeta
'_4r|\hat{Y}_{I}|, 
\end{split}
\end{equation}
where \(\zeta _3'\) is a positive constant depending only on
\(\zzz_\grp\) 
and the metric on \(Y_i\), and \(\zeta _4'\) is a
positive constant depending only on \(\nu \), and \(\varsigma _w\),
\(\zzz_\grp\), 
and the metric. 
Insert the preceding estimate into (\ref{energy-v}), 
we arrive at \((\ref{eq:E-top-bdd0})\) (b). \epf

\subsection{Lower bounds for \(\scrE^{'\mu_r,  \hat{\grp}}_{top}(\hat{Y}_{[l,L]})\)}\label{sec:E_top-lower}
We need a lower bound that
is independent of \(\hat{Y}_{[l,L]}\), \([l, L]\subset[0,\infty]\)
being arbitrary. When \(\hat{Y}=\hat{Y}_i\) is a vanishing end,
suppose in addition that \(l\geq \grl_i'\), 
so that the last term of
(\ref{eq:flow-energy}) vanishes on such \(\hat{Y}_{[l,L]}\). Then, for
an admissible solution \((A,
\Psi )=(A_r, \Psi _r)\) to \(\grS_{\mu _r, \hat{\grp}}(A_r, \Psi
_r)=0\), we have:
 \begin{equation}\label{bdd:CSD-lower}
\begin{split} 
& -2\op{CSD}_{\mu_+, \hat{\grp}}^{\partial \hat{Y}_{[l,L]}}(B,\Phi)\\
&\quad \quad =\frac{1}{2}\|\partial_sB\|_{L^2(\hat{Y}_{i, l})}^2
+2\|\partial_s\Phi\|^2_{L^2(\hat{Y}_{i,
    l})}-\frac{ir}{4}\int_{\hat{Y}_{[l,L]}}ds\,(
*_Y\xi_\nu)\wedge(F_B-F_{B_0})\\
& \quad \quad \geq -\frac{ir}{4}\int_{\hat{Y}_{[l,L]}}ds\,(
*_Y\xi_\nu)\wedge(F_B-F_{B_0}). 
\end{split}
\end{equation}
Ideally, \(1/r\) times the last term above should be bounded
independently of both
\(r\) and \(\hat{Y}_{[l,L]}\). Unfortunately, this is far from
straightforward. Note that similar terms do not appear in the setting
of \cite{LT} or \cite{Arnold2}. What follows is a first attempt to
tackle this trouble-making term. 

To proceed, use (\ref{eq:flow-energy}), (\ref{DE:xi}), (\ref{eq:xi-exp}),
Lemma \ref{lem:F-L_1} to get that for \(r\geq r_0\), 
\begin{equation}\label{eq:76a}
\begin{split} 
& \frac{1}{2}\|\partial_sB\|_{L^2(\hat{Y}_{i, l})}^2 +2\|\partial_s\Phi\|^2_{L^2(\hat{Y}_{i, l})}
\\
& \quad = -2\op{CSD}_{\mu_+, \hat{\grp}}^{\partial \hat{Y}_{i,[l,\infty]}}(B, \Phi ) 
+\frac{ir}{4}\int_{\hat{Y}_{i, l}}ds\,(  *_Y\xi_\nu)\wedge(F_{B_i}-F_{B_0})
\\
& \qquad\quad +\frac{ir}{4}\int_{\hat{Y}_{i, l}}ds\,(  *_Y\xi_\nu)\wedge(F_{B}-F_{B_i})
\end{split}
\end{equation}
The last term above may be re-expressed as follows:
\begin{equation}\label{eq:76b}\begin{split}
& \frac{ir}{4}\int_{\hat{Y}_{i,
    l}}ds\,(*_Y\xi_\nu)\wedge(F_{B}-F_{B_i})\\
&\qquad \qquad =\frac{ir}{4}\int_{\hat{Y}_{i,
    l}}ds\,(  *_Y\xi_\nu)\wedge d_Y(B-B_i)\\
& \qquad \qquad =\frac{ir}{4}\int_l^\infty ds\int_{Y_{i,
    s}}d(  *_Y\xi_\nu)\wedge (B-B_i)\\
&\qquad \qquad =-\frac{ir}{4}\int_l^\infty ds\int_{Y_{i,
    s}}d(  *_Y\xi_\nu)\wedge \int_s^\infty \partial_tB(t) \,dt\\
&\qquad \qquad =\frac{ir}{4}\int_l^\infty ds\int_{Y_{i,
    s}}(  \partial_s\xi_\nu)\wedge \int_s^\infty \partial_tB(t) \,dt.\\
  \end{split}
\end{equation}
The second equality above uses the Stokes' theorem. The third equality
uses the facts that \(\int_s^\infty \partial_tB(t)\,
dt=B_i-B(s)\). The last equality follows from (\ref{DE:xi}).

The last expression may in turn be estimated as follows:
\begin{equation}\label{eq:76c}
  \begin{split}
  &  \Big| \frac{ir}{4}\int_l^\infty  ds\,\int_{Y_{i:
    s}}( \partial_s\xi_\nu)\wedge\big(\int_s^\infty\partial_tB(t)\, 
dt\big)\Big| \\
&\qquad \qquad \leq \zeta  r\int_l^\infty ds\,  e^{-\kappa _i
  s}\int_s^\infty dt\int_{Y_i}|\partial_t B(t) |\\
&\qquad \qquad\leq \zeta  'r\int_l^\infty ds \, e^{-\kappa _i
  s}\int_s^\infty dt\, \|\partial_t B(t) \|_{L^2(Y_{i,t})}\\
&\qquad \qquad\leq \zeta  ''r\int_l^\infty ds \, e^{-\kappa _i
  s} \|\partial_s B \|_{L^2(Y_{i,s})}\\
 &\qquad \qquad\leq \frac{1}{4}\|\partial_s B \|^2_{L^2(\hat{Y}_{i,l})}+\zeta _2 r^2
 e^{-2\kappa _i l}. 
  \end{split}   
\end{equation}
The third inequality above uses the fact that \((B_i, \Phi _i)\) is by the
assumption (0) in Theorem \ref{thm:l-conv} regular, and therefore
\(\|\partial_t B(t) \|_{L^2(Y_{i,t})}\) decays exponentially to
\(0\).

Combining (\ref{eq:76a}), (\ref{eq:76b}), (\ref{eq:76c}), one has 
\[\begin{split}
    & \frac{1}{4}\|\partial_sB\|_{L^2(\hat{Y}_{i, l})}^2 +2\|\partial_s\Phi\|^2_{L^2(\hat{Y}_{i, l})}
\\
& \quad \leq -2\op{CSD}_{\mu_+, \hat{\grp}}^{\partial \hat{Y}_{i,[l,\infty]}}(B, \Phi
) +\zeta  _ir\,  e^{-\kappa _i l} +\zeta '_i r^2e^{-2\kappa _i l},
\end{split}
\]
where \(\zeta _i\), \(\zeta _i'\) are positive constants; \(\zeta _i\) depends only on \(\nu \), \(B_0\), the constant \(\zzz_i\)
in Lemma \ref{lem:F-L_1}; \(\zeta _i'\) depends only on \(\nu\). 
Rearranging, we get 
\[
-2\op{CSD}_{\mu_+, \hat{\grp}}^{\partial \hat{Y}_{i,[l, \infty]}}(B, \Phi
) \geq- \zeta _i'' (r^2e^{-2\kappa _i l}+1), 
\]
and by (\ref{diff:CSD-E}), \(\scrE^{'\mu _r}_{top}(\hat{Y}_{i,l})(A,
\Psi )=\scrE^{'\mu_r, \hat{\grp}}_{top}(\hat{Y}_{i,l})(A, \Psi )\) satisfies
a similar inequality: 
\[
\scrE^{'\mu_r, \hat{\grp}}_{top}(\hat{Y}_{i,l})(A, \Psi )\geq - \ud{\zeta }_i
(r^2e^{-2\kappa _i l}+1). 
\]
In particular, setting \(\hatl_i:=\kappa _i^{-1}\), and choosing
\(r_0\) to be greater than 
\(\max_{i\in \grY_v}\exp \, (4\grl_i'/\hatl_i)\), we have that for \(r\geq
r_0\),
\begin{equation}\label{bdd:E_l-lower}
\begin{split}
 -2\op{CSD}_{\mu_+, \hat{\grp}}^{\partial Y_{i,[l, \infty]}}(B, \Phi
) & \geq -2\, \zeta _i'' r \qquad \forall l\geq \big(\frac{\ln r}{2}\big)\hatl_i;\\
\scrE^{'\mu _r}_{top}(\hat{Y}_{i,l})(A,
\Psi )= \scrE^{'\mu_r, \hat{\grp}}_{top}(\hat{Y}_{i,l})(A, \Psi )& \geq - 2\ud{\zeta
}_i r\qquad \forall l\geq \big(\frac{\ln r}{2}\big)\hatl_i
\end{split}
\end{equation}
since \(r^2e^{-2\kappa _i l}\leq r\) when \(l\geq \big(\frac{\ln
  r}{2}\big)\hatl_i\). 
In the above, the positive constants \(\zeta _i'', \ud{\zeta
}_i >0\) depend only on \(\nu \), \(B_0\), and the constant \(\zzz_i\)
in Lemma \ref{lem:F-L_1}. 

Let \(\pmb{\hatl}\co
\grY\to \bbR^+\) be the function given by
\(\pmb{\hatl}(i)=(\hatl_i)_i\), and write \(\bfl_r:=(\ln r)\, 
\pmb{\hatl}\). 
Under the assumption (\ref{assume:EtopX-ubdd}) (b), the preceding pair of
inequalities then implies that for \(r\geq r_0\) and \( \forall {\bf
  l}\geq \pmb{\hatl}_r/2\) with \(\bfl(i)=:l_i\), 
\begin{equation}\label{bdd:E-up-i}
\begin{split}
\scrE^{'\mu_r, \hat{\grp}}_{top}(X_{\bf l})(A, \Psi
)& =\scrE^{'\mu_r, \hat{\grp}}_{top}(X)(A, \Psi
)-\sum_{i\in \grY}\scrE^{'\mu_r, \hat{\grp}}_{top}(\hat{Y}_{i,l_i})(A, \Psi )\\
& \leq  r \, (\smE+2\sum_{i\in \grY}\ud{\zeta
}_i ) \quad 
\text{assuming (\ref{assume:EtopX-ubdd}) (b); similarly, }
\\
\scrE^{'\mu_r}_{top}(X_{\bf l})(A, \Psi
) & \leq  r \, (\smE +2\sum_{i\in \grY}\ud{\zeta
}_i ) \quad 
\text{assuming (\ref{assume:EtopX-ubdd}) (a).}  
\end{split}
\end{equation}

Combining the preceding upper bounds for
\(\scrE^{'\mu_r, \hat{\grp}}_{top}(X_{\bf l})\),
\(\scrE^{'\mu_r}_{top}(X_{\bf l})\) with (\ref{eq:E-top-bdd0}), we
have for any \(r\geq r_0\) and \( \forall {\bf
  l}\geq \pmb{\hatl}_r/2\) with \(\bfl(i)=:l_i\), 
\begin{equation}\label{bdd:E_Xl}
\begin{split}
  & 
  \frac{1}{8}\int_{X_{\bfl}}|F_A|^2+\int_{X_{\bfl}}|\nabla_A\Psi|^2+\frac{1}{2}\int_{X_{\bfl}}\big|\frac{i}{4}\rho(\mu_r^+)-(\Psi\Psi^*)_0\big|^2\\
&\qquad \quad 
\leq r \, (\smE+\zeta )+\zeta '\sum_{i\in \grY}r^{\ra_i}l_i \quad
\text{assuming (\ref{assume:EtopX-ubdd}) (a)};\\
& \frac{1}{16}\int_{X_{\bfl}}|F_A|^2+\int_{X_{\bfl}}|\nabla_A\Psi|^2+\frac{1}{4}\int_{X_{\bfl}}\big|\frac{i}{4}\rho(\mu_r^+)-(\Psi\Psi^*)_0\big|^2\\
 &\qquad \quad 
 \leq r \, (\smE+\zeta _p)+\zeta '_p\sum_{i\in \grY}r^{\ra_i}l_i\quad
\text{assuming (\ref{assume:EtopX-ubdd}) (b).}
\end{split}
\end{equation}
(As in Lemma \ref{lem:CSD-est}, \(\ra_i=1\) when \(\hat{Y}_i\) is a
Morse end, and \(\ra_i=0\) when it is a vanishing end.) 
In particular, given \(i\in \grY\), by taking \(\bfl\) above to be
such that \(\bfl(j)=(\frac{\ln r}{2})\,  \hatl_j\) when
\(j\neq i\) and \(\bfl(i)=l\geq (\frac{\ln r}{2})\, \hatl_i\), we have 
\[\begin{split}
 \frac{1}{16}\int_{\hat{Y}_{i, [l-1, l]}}|F_A|^2& <
 \frac{1}{16}\int_{X_{\bfl}}|F_A|^2\\
& \leq r \, (\smE+\zeta _1\ln
 r)+\zeta _2 r^{\ra_i} \big(l-\frac{\ln r}{2}\hatl_i\big)\qquad
 \forall\, l\geq \frac{\ln r}{2}\hatl_i. 
\end{split}
\] 
The positive
constants \(\zeta _1, \zeta _2\), as well as \(\zeta \), \(\zeta '\) in
(\ref{bdd:E_Xl}), depend only on the metric, the \(\Spin^c\)
structure, \(\nu \), \(\varsigma _w\) \(B_0\), and \(\max_{j\in \grY}
\zzz_j\),  \(\zzz_j\) being the constants 
in Lemma \ref{lem:F-L_1}. In particular, they are
independent of \(r\), \(l\), and \(i\). 

Combining the preceding bound with (\ref{bdd:CSD-lower}), and
choosing \(r_0\) to be sufficiently large (depending on \(\smE\)), we have for all
\(l\geq\frac{\ln r}{2}\, \hatl_i\)
and any \(L\geq l\), 
\[\begin{split}
& -2\op{CSD}_{\mu_+, \hat{\grp}}^{\partial \hat{Y}_{i, [l,L]}}(B,\Phi)\\
& \qquad \qquad\geq -\zeta r\sum_{n=0}^{\infty} e^{-\kappa
  _i(l+n)}\|F_B-F_{B_0}\|_{L^2(\hat{Y}_{i, [l+n, l+n+1]})}\\
& \qquad \qquad\geq -\sum_{n=0}^{\infty}e^{-\kappa
  _i(l+n)}\big( \zeta 'r^2+r \, (\smE+\zeta _1\ln
 r)+\zeta _2 r^{\ra_i} (l+n)\big)\\
& \qquad \qquad\geq-  \zeta _3 \, r^2 e^{-\kappa
  _il}-\zeta _4\, r e^{-\kappa
  _il/2}\quad \text{when \(r\geq r_0\).} 
\end{split}
\]
In particular, for any \(r\geq r_0\), \(l\geq (\ln
  r) \, \hatl_i\), and \(L\geq l\), we have 
\begin{equation}\label{bdd:E-lower-ii}
\begin{split}
-2\op{CSD}_{\mu_+, \hat{\grp}}^{\partial \hat{Y}_{i,
    [l,L]}}(B,\Phi)& \geq -\zeta _5 \, r\quad \text{and equivalently (by (\ref{diff:CSD-E})), }\\
\scrE^{' \mu _r}_{top} (\hat{Y}_{i,
    [l,L]})(A, \Psi )& =\scrE^{' \mu _r, \hat{\grp}}_{top} (\hat{Y}_{i,
    [l,L]})(A, \Psi ) \geq  -\zeta '_5 \, r.\quad 
\end{split}
\end{equation}
The positive constants \(\zeta _3, \zeta _4, \zeta _5, \zeta _5'\)
above depend only on the metric, the \(\Spin^c\)
structure, \(\nu \), \(\varsigma _w\), \(B_0\), and \(\max_{j\in \grY}
\zzz_j\).

\subsection{Lemma \ref{lem:Etop-bdd1} and some of its variants
}\label{sec:pf-E_top-bdd0}

In this part we prove Lemma \ref{lem:Etop-bdd1}
and some of its variants.

\subsubsection*{\it Proof of Lemma \ref{lem:Etop-bdd1}.}
The rightmost inequalities  in  (\ref{eq:E_top-M1}) (i) and the left
most 
inequalies in (\ref{eq:E_top-M1}) (ii) follow respectively from 
(\ref{bdd:E-up-i}) and 
(\ref{bdd:E-lower-ii}). To get the remaining two inequalities in
(\ref{eq:E_top-M1}), we need a lower bound on \(\scrE^{'\mu_r}_{top}(X_{\bf l})(A, \Psi
)\) or \(\scrE^{'\mu_r, \hat{\grp}}_{top}(X_{\bf l})(A, \Psi
)\). Here is preliminary bound from (\ref{eq:E-top-bdd0}): for all \({\bf l}\co \grY\to
\bbR^+\), \(\bfl(i)=:l_i\).  
\[\begin{split}
& \scrE^{'\mu_r, \hat{\grp}}_{top}(X_{\bf l})(A, \Psi
)\geq - \zeta _p (r+ \sum_{i\in \grY}r^{\ra_i}l_i);\\
& \scrE^{'\mu_r}_{top}(X_{\bf l})(A, \Psi
)\geq - \zeta   (r+\sum_{i\in \grY}r^{\ra_i}l_i),
\end{split}
\]
where \(\zeta\), \(\zeta_p\) are positive constants depending only on
the parameters listed in (\ref{parameters}). 
Combined with (\ref{bdd:E-lower-ii}), this implies:  For any \(r\geq r_0\), 
\(X_\bfl\), \(\bfl\geq \bfl_r\) and \(l_i:=\bfl(i)\) (possibly
\(\infty\)), 
\begin{equation}\label{bdd-lower-i}
\begin{split}
\scrE^{'\mu_r}_{top}(X_{\bf l})(A, \Psi
)& =\scrE^{'\mu_r}_{top}(X_{\bfl _r})(A, \Psi
)+\sum_{i\in \grY}\scrE^{'\mu_r}_{top}(\hat{Y}_{i, [\bfl_r(i), l]})\\
& \geq - \zeta 'r\ln r;\quad \text{similarly,} \\
\scrE^{'\mu_r, \hat{\grp}}_{top}(X_{\bf l})(A, \Psi
)&\geq  - \zeta _p'r\ln r, 
\end{split}
\end{equation}
where \(\zeta'\), \(\zeta'_p\) are positive constants depending only on
the parameters listed in (\ref{parameters}), 
\(B_0\), and the constants   \(\zzz_i\)  
in Lemma \ref{lem:F-L_1}. This is precisely what was asserted in the
leftmost inequalities in (\ref{eq:E_top-M1}) (i).  
Meanwhile, combining (\ref{bdd-lower-i}) and (\ref{bdd:E-up-i}),
one has: For any \(\hat{Y}_{i, [l, L]}\subset X-\mathring{X}_{\bfl
  _r}\), \(L\) possibly \(\infty\), we have: 
\begin{equation}\label{eq:4.38+}
\begin{split}
\scrE^{'\mu_r}_{top}(\hat{Y}_{i, [l, L]})(A, \Psi
)& =\scrE^{'\mu_r}_{top}(X_{\bfL_i})(A, \Psi
)-\scrE^{'\mu_r}_{top}(X_{\bfl_i})(A, \Psi
)\\
& \leq r \, (\smE + \zeta _e\ln r), \\
\end{split}
\end{equation}
where \(\bfl_i, \bfL_i\co \grY\to \bbR^+\) are such that
\(\bfl_i(j)=\bfL_i(j)=\bfl_r(j)\) for \(j\neq i\), and
\(\bfl_i(i)=l\); \(\bfL_i(i)=L\). In the above, 
\(\zeta '\) depends only  on
the parameters listed in (\ref{parameters}), 
\(B_0\), and the constants   \(\zzz_i\)  
in Lemma \ref{lem:F-L_1}. The preceding inequality leads directly to
the second inequality in (\ref{eq:E_top-M1}) (ii). \epf

A straightforward corollary of  Lemma \ref{lem:Etop-bdd1} is: 
\begin{cor}\label{cor:E_top(X')}
Adopt the assumptions and notation of Lemma \ref{lem:Etop-bdd1}. There is a  positive constant \(\zeta \) independent of
\(r\), \(X_\bullet\), and \((A, \Psi )\) such that 
\[\begin{split}
\text{\rm (i) }  & \scrE^{'\mu_r}_{top}(X_\bullet)(A, \Psi )\geq -\zeta  \, r\quad \text{for all
  \(r\geq r_0\) and \(X_\bullet\subset X-X''=X_v\); }\\
  & \qquad \text{consequently}\\
\text{\rm (ii) }  & \scrE^{'\mu_r}_{top}(X'')(A, \Psi )\leq r\, (\smE +\zeta );
\end{split}
\]
similarly for \(\scrE^{'\mu_r,\hat{\grp}}_{top}(X_\bullet)(A, \Psi
)\). 
\end{cor}
\pf We shall prove only the statement for
\(\scrE^{'\mu_r}_{top}(X_\bullet)(A, \Psi )\), since the proof for
\(\scrE^{'\mu_r, \hat{\grp}}_{top}(X_\bullet)(A, \Psi )\) is the same.

According to (\ref{eq:E-top-bdd0}),  for all
\(X_\bullet\subset X_v\)
\[
\scrE^{'\mu_r}_{top}(X_\bullet)(A, \Psi )\geq -(\zeta  'r+\zeta_0
|X_\bullet|).  
\]
It follows that for any \(X_\bullet \subset X_v\), 
\[
\begin{split}
& \scrE^{'\mu_r}_{top}(X_\bullet)(A, \Psi )\\
& \quad =\scrE^{'\mu_r}_{top}(X_\bullet\cap
X_{\bfl_r})(A, \Psi )+\scrE^{'\mu_r}_{top}(X_\bullet- X_{\bfl_r})(A, \Psi )\\ &\quad  \geq-(\zeta  'r+\zeta '_0 \ln r)-\zeta '' r\geq -\zeta  r.
\end{split}
\]
In the above, we applied the preceding inequality  to bound the first term in
the second line, and used the first inequality in (\ref{eq:E_top-M1})
(ii) to bound the second term. This proves Inequality (i)
asserted by the corollary. To obtain Inequality (ii), simply
write \(\scrE^{'\mu_r}_{top}(X'')(A, \Psi )=\scrE^{'\mu_r}_{top}(X)(A, \Psi
)-\scrE^{'\mu_r}_{top}(X-X'')(A, \Psi )\), then combine with the bounds from (i)
and (\ref{assume:EtopX-ubdd}). 
\epf 

For future reference, note that the arguments above  in fact establishes the following
generalization of Lemma \ref{lem:Etop-bdd1}: 
\begin{lemma}
\label{lem:Etop-bddf}
Adopt the assumptions and notation of Lemma \ref{lem:Etop-bdd1}.
Suppose 
\(\grt(r)\) is a function from \([r_0, \infty)\)
to \((1, \infty)\) such that the lower bound on
\(\scrE^{'\mu_r}_{top}(X_\bullet)\) in (\ref{bdd:E-lower-ii}) holds for any
given 
\(r\geq r_0\) and \(X_\bullet\subset X'-X_{\grt(r)
  \pmb{\hatl},m}\). That is, assume that 
\begin{equation}\label{assume:t(r)}
\scrE^{'\mu_r}_{top} (X_\bullet)(A, \Psi ) \geq  -\zeta '_5 \, r\quad \forall X_\bullet\subset X'-X_{\grt(r)
  \pmb{\hatl},m}, r\geq r_0.
\end{equation}
(In particular, it follows from Lemma \ref{lem:Etop-bdd1} that this holds
with \(\grt(r)=\ln r\).)
Then there exist positive constants \(\zeta \), \(\zeta _1\), \(\zeta_2\), \(\zeta _3\) 
such that for any 
 \(r\geq r_0\), 
 \begin{equation}\label{eq:E_top-Mf}
\begin{split}
\text{\rm (i) } & 
-r(\zeta _1+\zeta _2 \grt(r))\leq \scrE_{top}^{'\mu_r}(X_\bullet)(A,
\Psi )\leq r\, (\smE+\zeta) 
\quad \text{\(\forall X_\bullet\supset X_{\grt(r) \pmb{\hatl},m}\)};\\
\text{\rm (ii) } & -\zeta _1 r\leq \scrE_{top}^{'\mu_r}(X_\bullet)(A, \Psi
) \leq r\, (\smE+
\zeta    _2\grt (r)+\zeta _3)
\quad  \text{ \(\forall X_\bullet\subset X-\mathring{X}_{\grt(r)\pmb{\hatl},m}\)}.
\end{split}
\end{equation}
The positive constants 
\(\zeta \), \(\zeta _1\), \(\zeta_2\), \(\zeta _3\)  above depend only on the metric, the \(\Spin^c\)
structure, \(\nu \), \(\varsigma _w\), \(B_0\), and the constants \(\zzz_i\)
in Lemma \ref{lem:F-L_1}. 
In particular, these constants 
are independent of \(X_\bullet\) and \(r\).

A similar statement holds for \(\scrE_{top}^{'\mu_r,
  \hat{\grp}}(X_\bullet)(A, \Psi ) \), with the constants in the
inequalities depending additionally on \(\zzz_\grp\). 
\end{lemma}
\pf We shall again prove only the inequalities for
\(\scrE_{top}^{'\mu_r}(X_\bullet)(A, \Psi )\), since the proof for \(\scrE_{top}^{'\mu_r,
  \hat{\grp}}(X_\bullet)(A, \Psi ) \) is identical.

The first inequality in (\ref{eq:E_top-Mf}) (ii) follows from the
assumption (\ref{assume:t(r)}) and Inequality (i) in Corollary \ref{cor:E_top(X')}. 
Replacing (\ref{bdd:E-lower-ii}) and (\ref{bdd:E_l-lower}) by this inequality, the arguments in
the previous subsection and in the proof of Lemma \ref{lem:Etop-bdd1} above can be
repeated to replace (\ref{bdd:E-up-i}), (\ref{bdd-lower-i}) (\ref{eq:4.38+})
with their respective  companion versions. These are respectively the
second inequality in (\ref{eq:E_top-Mf}) (i), the first inequality in
(\ref{eq:E_top-Mf}) (i), and  the
second inequality in (\ref{eq:E_top-Mf}) (ii). 
\epf

\subsection{\(L^2_{1, loc /A}\)-bounds for \((A-A_0, \Psi)\)}\label{sec:L^_1-bdd}

It follows from the discussion in the preceding subsection that: 
\begin{prop}\label{prop:SW-L2-bdd}
Adopt the notation and assumptions of Lemma \ref{lem:Etop-bdd1}. 
Let  \(\zeta _0\), \(\zeta ''\), \(\zeta ''_p\)
be the constants from  Lemma
\ref{lem:E-top-bdd0},
Then there exist constants 
\(\zeta _1\), \(\zeta _1'\)
such that the following hold for any \(r\geq
r_0\), any  \((A, \Psi )\), and any compact \(X_\bullet\subset X\)
satisfying either \(\partial X_\bullet \subset
X-\mathring{X}_{\bfl_r}\) or \(X_\bullet\subset X_{\bfl_r}\): 
\begin{equation}\label{eq:L^2_1}
\begin{split}
  &   \frac{1}{8}\int_{X_{\bullet}}|F_A|^2+\int_{X_{\bullet}}|\nabla_A\Psi|^2+\frac{1}{2}\int_{X_{\bullet}}\big|\frac{i}{4}\rho(\mu_r^+)-(\Psi\Psi^*)_0\big|^2\\
&\qquad 
\leq r \, \big(\zeta _0|X_{\bullet,m}|+\smE+\zeta _1\ln r
  \big) +\zeta ''|X_\bullet| \quad \text{assuming 
    (\ref{assume:EtopX-ubdd}) (a)};
\\
& \frac{1}{16}\int_{X_{\bullet}}|F_A|^2+\int_{X_{\bullet}}|\nabla_A\Psi|^2+\frac{1}{4}\int_{X_{\bullet}}\big|\frac{i}{4}\rho(\mu_r^+)-(\Psi\Psi^*)_0\big|^2\\
&\qquad 
\leq r \, \big(\zeta _0|X_{\bullet,m}|+\smE +\zeta '_1\ln r \big)
+\zeta ''_p|X_\bullet| \quad \text{assuming 
    (\ref{assume:EtopX-ubdd}) (b)}.
\end{split}
\end{equation}
The constants \(\zeta _1\), \(\zeta _1'\) above depend only on the parameters listed in
(\ref{parameters}), together with \(B_0\) and the constants \(\zzz_i\) in
Lemma \ref{lem:F-L_1}.
In particular, it is independent of \(r\) and
\(X_\bullet\). 

Moreover, if the assumption (\ref{assume:t(r)}) in Lemma \ref{lem:Etop-bddf} holds, then
the statement above holds with 
all appearances of \(\ln r \) replaced by \(\grt(r)\). 
\end{prop}
\pf When \(X_\bullet \) is such that 
\(\partial X_\bullet \subset X-\mathring{X}_{\bfl_r}\), the claim of
the proposition 
follows directly from Lemmas \ref{lem:Etop-bdd1}, \ref{lem:Etop-bddf} and
\ref{lem:E-top-bdd0}. The case when \(X_\bullet \subset X_{\bfl_r}\) follows from the preceding
case,  together with the observation that when \(X_\bullet\subset
X_{\bfl_r}\), 
\[\begin{split}
& \frac{1}{8}\int_{X_{\bullet}}|F_A|^2+\int_{X_{\bullet}}|\nabla_A\Psi|^2+\frac{1}{2}\int_{X_{\bullet}}\big|\frac{i}{4}\rho(\mu_r^+)-(\Psi\Psi^*)_0\big|^2\\
&\qquad \leq \frac{1}{8}\int_{X_{\bfl_r}}|F_A|^2+\int_{X_{\bfl_r}}|\nabla_A\Psi|^2+\frac{1}{2}\int_{X_{\bfl_r}}\big|\frac{i}{4}\rho(\mu_r^+)-(\Psi\Psi^*)_0\big|^2;\\
\end{split}
\]
\epf

As a consequence,

\begin{prop}\label{est-L^2_1}
  Adopt the assumptions and notation of Lemma \ref{lem:Etop-bdd1}, 
  and recall the reference connection  \(A_0\) from (\ref{eq:A_0}). 
  Fix a compact  \(X_\bullet \subset X\) and an \(r\geq r_0\), \(r_0\) being
as in the previous proposition.  Let \(u_r\in
C^\infty(X_\bullet, S^1)\) be such that  \(
u_r\cdot  (A_r, \Psi_r)\) is in the gauge specified in
\cite{KM}'s Equations (5.2) and (5.3). Then  there exist positive
constants \(\zeta_2\), \(\zeta _2'\) depending only on the parameters listed in
(\ref{parameters}), such that 
\begin{equation}\label{est:L^2_1}\begin{split}
& \|(u_r\cdot A_r-A_0, \Psi_r)\|^2_{L_{1,
    A_r}^2(X_{\bullet})} 
\leq \zeta_2
r \, \big(|X_{\bullet,m}|+\smE +\ln r \big)|X_\bullet|
+\zeta _2' |X_\bullet|^2
\end{split}
\end{equation}
Again, if the assumption (\ref{assume:t(r)}) in Lemma \ref{lem:Etop-bddf} holds, then
the statement above holds with 
all appearances of \(\ln r \) replaced by \(\grt(r)\). 
\end{prop}
\pf 
Combining (\ref{eq:L^2_1}) with standard elliptic estimates
(cf. e.g. pp. 101-104 of \cite{KM} and noting that the lowest absolute
value of non-zero eigenvalues of \(*d+d*\) on \(X_\bullet\) has a
lower bound propotional to \(|X_\bullet|^{-1}\)),
one may find a positive
constant \(\zeta \) depending only on the metric on \(X\) such that 
\[
\|u_r\cdot A_r-A_0\|^2_{L_{1}^2(X_{\bullet})} +\|\nabla_A\Psi\|^2_{L^2(X_\bullet)}\leq
\zeta |X_\bullet|\cdot (\text{RHS of (\ref{eq:L^2_1})  (b)}). 
\]
Meanwhile, a combination of (\ref{eq:L^2_1}) with (\ref{f-trick})
(with \(\grf\) set to be 1) implies that 
\[\begin{split}
\|\Psi\|^2_{L^2(X_\bullet)}& \leq 4\cdot (\text{RHS of
  (\ref{eq:L^2_1})  (b)})+ |X_\bullet|+\frac{1}{2}\int_{X_\bullet}|\mu
_r|\\
& \leq \zeta_3
r \, \big(|X_{\bullet,m}|+\smE +\ln r \big)
+\zeta _3' |X_\bullet|, 
\end{split}
\]
where \(\zeta _3, \zeta _3'\) only on the parameters listed in
(\ref{parameters}).  Together with the preceding inequality, we arrive
at (\ref{est:L^2_1}). 

\epf 

\begin{rem}\label{rem:CNgauge}
The gauge transformation \(u_r\) in (\ref{est:L^2_1}) depends on \(X_\bullet\). Specifically, it is
determined, via Equations (5.2) and (5.3) in \cite{KM}, by \(A_0\) together with a choice of \(\{q_1, q_2, \ldots,
q_{b^1(X_\bullet)}\}\), where each \(q_k\), \(k\in \{1, \ldots, b^1(X_\bullet)\}\) is
a closed 3-form supported in a compact set in the interior of
\(X_\bullet\) and \(\{[q_k]\}_k\subset H^3(X_\bullet , \partial
X_\bullet;\bbR)\) forms a basis of \(H^3(X_\bullet , \partial
X_\bullet;\bbZ)/\text{Tors}\). As in \cite{KM}, we say that a \((A'_r, \Psi
'_r)=u_r\cdot (A_r, \Psi _r)\) or \(A'_r=u_r\cdot A_r\) is in a {\em Coulomb-Neumann} gauge
(with respect to \(A_0\)) if it satisfies Equation (5.2) of
\cite{KM}. We say that it is in the {\em normalized Coulomb-Neumann} gauge
(with respect to \(A_0\) and \(\{q_k\}_{k=1}^{b^1(X_\bullet)}\)) if it
satisfies both Equations (5.2) and (5.3) in \cite{KM}. Note that if 
\((A'_r, \Psi _r')\) is in a Coulomb-Neumann gauge with respect to
\(A_0\), then \((B', \Phi ')=(A'_r, \Psi _r')|_{\partial X_\bullet}\)
is in a Coulomb gauge with respect to \(B_0=A_0|_{\partial
  X_\bullet}\). 
In the case
when \(X_\bullet\) is of the form \(\hat{Y}_{[l, L]}\), for some
\(Y=Y_i\), \(H^3(X_\bullet , \partial
X_\bullet;\bbZ)\simeq H^2(Y;\bbZ)\), and we choose \(\{q_k\}_k\) so
that  \(q_k=d(\chi _+(s) *_3 h_k)\), with \(h_k\), \(k\in \{1, \ldots,
b^1(Y)\}\) being the harmonic
1-forms on \(Y\) appearing in (\ref{eq:deltab1}), and \(\chi _+\co [l,
L]\to [0,1]\) is a smooth non-decreasing function that is 0 on a
neighborhood of \(\{l\}\subset [l,L]\) and is 1 on a neighborhood of
\(\{L\}\subset [l,L]\). If \((A_r', \Psi _r')\) is in the normalized
Coulomb-Neumann gauge with respect to such \(\{q_k\}_k\) (and
\(A_0\)), then \((A_r', \Psi _r')|_{\hat{Y}_{:L}}\) is in the
  normalized Coulomb gauge with respect to \(\{h_k\}_k\) (and
  \(B_0=A_0|_{\hat{Y}_{:L}}\)). 
\end{rem}

\subsection{Some integral estimates}\label{sec:int-est}

This subsection contains some integral estimates which will be
instrumental to the pointwise estimates in next section. 

\begin{lemma}\label{co:E-omega-bdd}
Adopt the assumptions and notation of Lemma \ref{lem:Etop-bdd1} 
and recall from (\ref{eq:A_0}) that \(A_0=\hat{B}_{0,i}\) on the
\(\hat{Y}_i\)-end. 
Fix \(r\geq r_0\) and \(L_i\geq 1\). Let \((B, \Phi)\) denote
 the restriction of \((A_r, \Psi_r)\) to \(Y_{i:L_i}\subset X\) and 
 let \((B', \Phi ')\) be the representative of the gauge equivalence
 class \([(B, \Phi )]\) in the normalized Coulomb gauge.  Let
 \(\hat{\grp}| (B, \Phi ) \) denote the restriction of \(\hat{\grp}(A,
\Psi )\) to \(Y_{i:L_i}\subset X\).  
 Then 
 there exist positive constants \(\zeta', \zeta _0, \zeta _w, \zeta '_\grp\) such that
\begin{equation}\label{eq:CSD-bdd}
\begin{split}
{\rm (1)} & \, \|(B'-B_{0,i}, \Phi)\|_{L^2_{1/2,B}(Y_{i:L_i})}^2 \leq \zeta' r\, (\smE +\ln r);\\
{\rm (2)} &\, |\op{CSD}^{Y_i}_0(B', \Phi')| \leq \zeta_0 \, r\, (\smE +\ln r);\\
& \,   |\op{CSD}^{Y_i}_{w_r}(B', \Phi')| =|\op{CSD}^{Y_i}_{w_r}(B, \Phi)|
\leq \zeta_w  \, r\, (\smE +\ln r);\\
{\rm (3)} & \,   |\, f_{\hat{\grp}|} (B, \Phi)\, | \leq   \zeta '_\grp
\, r\, (\smE +\ln r)
\end{split}\end{equation}
The positive
constants \(\zeta', \zeta _0, \zeta _w, \zeta '_\grp\)  depend only on the parameters listed in
(\ref{parameters}), together with \(B_{0,i}\) and the constants \(\zzz_i\) in
Lemma \ref{lem:F-L_1}. As before,  the factors of 
\(\ln r\) in (\ref{eq:CSD-bdd}) can be replaced by \(\grt(r)\)
 if the assumption of Lemma \ref{lem:Etop-bddf}
 holds.
\end{lemma}
\pf 
{\bf (1)}: Let \((A'_r, \Psi '_r)=u_r\cdot (A_r, \Psi _r)\) be in
the normalize Coulomb-Neumann gauge on \(\hat{Y}_{i,
    [L_i-1, L_i]}\), with respect to the choices specified in Remark
  \ref{rem:CNgauge}.  As remarked above, with such choices, \((A'_r, \Psi
  '_r)|_{Y_{i:L_i}}=(B', \Phi ')\).  Note that \[\|(B'-B_{0,i}, \Phi)\|_{L^2_{1/2, B}(Y_{i:
    L_i})}\leq \zeta_b \|(A_r'-A_0, \Psi'_r)\|_{L^2_{1, A}(\hat{Y}_{i,
    [L_i-1, L_i]})},\] then appeal to Proposition \ref{est-L^2_1}. 

{\bf (2)}: 
This is a
direct consequence of {\bf (1)} above.

{\bf (3)}: Recall the proof of Lemma \ref{lem:f_q} and the notation therein. From the first
inequality in (\ref{ineq:f_q}), (\ref{eq:q-bdd}) and {\bf (1)} above, we have: 
\[  \begin{split}
|f_{\hat{\grp}|}(B, \Phi )|& \leq \zeta _\grp\|(\delta B, \Phi )\|_{L^2}(\|\Phi_t\|_{L^2}+1)\quad \text{for a certain \(t\in [0,1]\)}\\
& \leq   \zeta '_\grp  \, r\, (\smE +\ln r),
\end{split}\]
where \(\delta  B:=B'-B_{0,i}\). \epf

\begin{lemma}\label{co:E-omega-bdd3}
Adopt the assumptions and notation of Lemma \ref{lem:Etop-bdd1}.
Let \(r_0\), 
\(\grt(r)\) be as in Lemma \ref{lem:Etop-bddf}. Fix \(r\geq r_0\). 
Then there are positive constants \(\zeta\), \(\zeta _h\) independent of \(r\),
and \(X_\bullet\) such that 
\begin{equation}\label{eq:E-omega-bdd3}
\Big| \int_{X_{\bullet}}i F_{A_r}\wedge \omega \Big| \leq
\zeta _h\,(\smE +\ln r)+\zeta \,|X_{\bullet,m}| \quad \text{for compact
                                       \(X_\bullet\subset X\)}
\end{equation}
The constants \(\zeta , \zeta_h\) above only depend on the parameters listed in (\ref{parameters}).

Moreover, if the assumption (\ref{assume:t(r)}) in Lemma \ref{lem:Etop-bddf} holds, then
the statement above holds with 
the factor of \(\ln r \) replaced by \(\grt(r)\). 
\end{lemma}
\pf
Write 
\begin{equation}\label{eq:F-omega}
\begin{split}
  & \int_{X_{\bullet}}i F_A\wedge \omega \\
  &\qquad \quad = i\int_{X_{\bullet}}
F_{A}\wedge*\nu      +i\int_{X_{\bullet}}
F_{A}\wedge\nu     \\
& \qquad \quad =i \int_{X_\bullet}
F_{A_0}\wedge*\nu     
+  
 4r^{-1}\scrE_{top}^{'\mu_r}(X_{\bullet})(A,
 \Psi)-4r^{-1}\scrE_{top}^{w_r, \hat{\grp}}(X_{\bullet})(A, \Psi).\\
 & \qquad \quad =i \int_{X_\bullet}
F_{A_0}\wedge*\nu     
+  
 4r^{-1}\scrE_{top}^{'\mu_r, \hat{\grp}}(X_{\bullet})(A,
 \Psi)-4r^{-1}\scrE_{top}^{w_r, \hat{\grp}}(X_{\bullet})(A, \Psi).\\
\end{split}
\end{equation}
Consider the second equality of (\ref{eq:F-omega}) in the case of
(\ref{assume:EtopX-ubdd})(a),  the second equality in the case of
(\ref{assume:EtopX-ubdd}) (b). 
We bound the right hand side of these formulas term by term.  The
first term on the right hand side in both cases is
bounded via (\ref{bdd(3)}) as: 
\[
\Big| i \int_{X_\bullet}
F_{A_0}\wedge*\nu     \Big| \leq \zeta _1\, |X_{\bullet,m}|+\zeta _1'
\]
For the second term, use (\ref{eq:E_top-M1}) 
to obtain:
\[\begin{cases}
\big| 4r^{-1}\scrE_{top}^{'\mu_r}(X_{\bullet})(A,
\Psi)\big|\leq 4 \smE +\zeta _2\ln r &\text{assuming
  (\ref{assume:EtopX-ubdd}) (a)};\\
\big| 4r^{-1}\scrE_{top}^{'\mu_r, \hat{\grp}}(X_{\bullet})(A,
\Psi)\big|\leq 4 \smE +\zeta '_2\ln r &\text{assuming
  (\ref{assume:EtopX-ubdd}) (b)}.\\
\end{cases}
\]
To bound the last term, 
 note that by our
 assumptions on \(A_0\) and \(w_r\), 
\[
\begin{split}
 \scrE_{top}^{w_r}(X_{\bullet}) (A, \Psi) & = \frac{1}{4}\int_{X_{\bullet}}F_{A_0}\wedge
  (F_{A_0}+iw_r)-2\, \op{CSD}^{\partial X_\bullet}_{w_r}(B,
  \Phi)\\
&= \frac{1}{4}\int_{X_{\bullet}\cap X_c}F_{A_0}\wedge
  (F_{A_0}+iw_r)-2\, \op{CSD}^{\partial X_\bullet}_{w_r}(B,
  \Phi);\\
  \scrE_{top}^{w_r, \hat{\grp}}(X_{\bullet}) (A, \Psi) &= \frac{1}{4}\int_{X_{\bullet}\cap X_c}F_{A_0}\wedge
  (F_{A_0}+iw_r)-2\, \op{CSD}^{\partial X_\bullet}_{w_r,
    \hat{\grp}|}(B,
  \Phi). 
\end{split}
\]
Combined with (\ref{eq:CSD-bdd}) (2) and (3), this gives  
\begin{equation}\label{eq:E_w}\begin{split}
4r^{-1}\big| \scrE_{top}^{w_r}(X_{\bullet}) (A, \Psi) \big| 
&\leq \zeta_3 \, (\smE+\ln r)\qquad \text{assuming
  (\ref{assume:EtopX-ubdd}) (a)};\\
4r^{-1}\big| \scrE_{top}^{w_r, \hat{\grp}}(X_{\bullet}) (A, \Psi) \big| 
&\leq \zeta'_3 \, (\smE+\ln r)\qquad \text{assuming
  (\ref{assume:EtopX-ubdd}) (b)}. 
\end{split}
\end{equation}

Plugging all these back to (\ref{eq:F-omega}), we have the asserted
inequality (\ref{eq:E-omega-bdd3}).

The assertion regarding the general case assuming the
condition of Lemma \ref{lem:Etop-bddf} follows from the same argument, with
the role of Lemma \ref{lem:Etop-bdd1} above played by Lemma 
\ref{lem:Etop-bddf} instead.
\epf

The next lemma is an analog of \cite{Ts}'s Lemma 3.1. Let \(\psi:=(r/2)^{-1/2}\Psi\).
\begin{lemma}\label{T:lem3.1}
  Adopt the assumptions and notation of Lemma \ref{lem:Etop-bdd1} 
and let \(X_\bullet\subset X\) be arbitrary. Then
there exists a positive constant \(\zeta\) independent of \(r\) and
\(\bullet\), such that 
\begin{gather*}
 r \int_{X_{\bullet}}(|\nu|-|\psi|^2)^2\leq
\zeta ( \smE +\ln r+|X_{\bullet,m}|).\\
 \end{gather*}
If in addition, \(X_\bullet\subset X''\), then 
\[
r\int_{X_{\bullet}}|\nu |\Big||\nu |-|\psi|^2\Big|\leq \zeta (\smE +\ln r+|X_\bullet|)
\]
for a positive constant \(\zeta'\) independent of \(r\) and
\(\bullet\). The constants \(\zeta , \zeta '\) above only depend on the parameters listed in (\ref{parameters}).
Moreover, the factor of 
\(\ln r\) above can be replaced by \(\grt(r)\) 
if the assumption of Lemma \ref{lem:Etop-bddf}
holds. 
\end{lemma}
\pf
Using the first line of the Seiberg-Witten equation, the fact that \(\omega\) is
self-dual, our assumption on \(\mu_r\), and (\ref{eq:E-omega-bdd3}), one has:
\begin{equation}\label{T:3.3}\begin{split}
 \frac{r}{2}\int_{X_{\bullet}}|\nu |(|\nu|-|\psi|^2-\zeta'
r^{-1})& \leq \int _{X_\bullet}i F_A\wedge\omega \\
& \leq \zeta_h (\smE+\ln r)+\zeta  |X_{\bullet,m}|).
\end{split}
\end{equation}
Meanwhile, arguing as Equation (3.4) of \cite{Ts}, one has:
\begin{equation}\label{T:eq3.4}\begin{split}
2^{-1} d^*d |\psi|^2+|\nabla_A\psi|^2+4^{-1}
r|\psi|^2(|\psi|^2-|\nu |)\leq \zeta_2 |\psi|^2.
\end{split}\end{equation}
Hence, with \(u:=|\psi|^2-|\nu|\),
\[
2^{-1} d^*d |\psi|^2+4^{-1}r|\nu| u+4^{-1} ru^2\leq \zeta_2 |\psi|^2.
\]
Integrating this and using (\ref{T:3.3}) and Proposition \ref{est-L^2_1}, one has
\[\begin{split}
r \int_{X_{\bullet}} u^2 &\leq
\zeta_4\, (|X_{\bullet,m}|+\ln r )+2\Big|\int_{\partial X_\bullet} \partial_s |\psi|^2\Big|.\\
\end{split}\]
Meanwhile, using the Seiberg-Witten equation and item (1) of
(\ref{eq:CSD-bdd}), 
 \[\Big|\int_{\partial X_\bullet} \partial_s |\psi|^2\Big|\leq
\zeta_4\, \|\psi\|_{L^2_{1/2}(\partial X_\bullet)}^2\leq \zeta_5\, \ln
r.\] 
These together imply the first inequality asserted in the Lemma.

To derive the second inequality in the lemma, follow Taubes' argument in
\cite{Ts}. The harmonicity of \(\nu \) implies that over \(X''\), 
\begin{equation}\label{eq:omega-ineq}
-\zeta |\nu|\leq d^*d|\nu|+|\nu|^{-1}|\nabla\nu |^2\leq
\zeta |\nu|,
\end{equation}
and the analog of Equation (3.9) in \cite{Ts} reads: 
\[
\begin{split}
r \int_{X_{\bullet}} |\nu| u_+ &\leq \zeta  \int_{X_{\bullet}}
(|\nu|^{-1} |\nabla\nu|^2+\zeta )+\zeta'\Big|\int_{\partial X_{\bullet}}
\partial_s u\Big|\\
& \leq \zeta''\Big(|X_\bullet|+\Big|\int_{\partial
  X_\bullet} \partial_s |\psi|^2\Big|+\Big|\int_{\partial
  X_\bullet} \partial_s |\nu|\Big|\Big)\\
&\leq \zeta_3\, (|X_\bullet|+ \ln r),
\end{split}
\]
where \(u_+:=\max \, (u, 0)\). Combine this with (\ref{T:3.3}).

The assertion regarding the general case assuming the
condition of Lemma \ref{lem:Etop-bddf} follows from the same
argument.
\epf

\section{A priori pointwise estimates}\label{sec:pt-est}

This section consists largely of refinements and extension of the
pointwise estimates in Section 3 of \cite{Ts}, which is in turn based
on Section I.2 in \cite{T}. Familiarity with
these references is assumed. We begin with some preliminaries. 


First, note that by assumption,  there exist positive constants \(\zzz_v\geq 1\),
\(\zzz_v '\geq 1\) that depends only on \(\nu \),  such that
 \begin{equation}\label{eq:v-sig}
\begin{split}
 (\zzz'_{v})^{-1} \leq  & \inf_{x\in X''\cap \nu^{-1}(0)}|\nabla \nu|\leq \sup_{x\in X''\cap \nu^{-1}(0)}
|\nabla\nu|\leq \zzz'_v; \\ 
 & \zzz^{-1}_v \tilde{\sigma} \leq |\nu|\leq
 \zzz_v \, \tilde{\sigma}, \quad  \text{
 over \(X''\), }
\end{split}
\end{equation}
where \(\ts\) is a function on \(X''\) defined as follows: 
Suppose that \(|\nu|^{-1}(0)\neq\emptyset\). Let \(\sigma(\cdot)\) denote the distance function to
\(\nu ^{-1}(0)\) on \(X''\), and set 
\[
\tilde{\sigma}:=\chi(\sigma)\, \sigma+(1-\chi(\sigma)).
\]
When \(\nu ^{-1}(0)=\emptyset\), let \(\sigma =\infty\) and
\(\td{\sigma }=1\). 


Let \(\gamma_a\) be a smooth
cutoff function on \(X\) that equals \(1\) on \(X^{'a}\) and agrees
with \(\chi (\mathfrc{s}_j-a+1)\) over each vanishing end \(j\in
\grY_v\). 

Recall also that  \(Z^{'a}:=X^{'a}\, \cap \nu^{-1}(0)\). Given 
\(\delta>0\), 
\[\text{\(X_\delta^{'a}:=\{x\, |\, x\in X^{'a}, \sigma (x)\geq\delta\}\), and \(Z_\delta^{'a}:=X^{'a}-X^{'a}_\delta\).}
\] 
In the case when \(a=0\) (resp. \(a=10\)), the spaces \(Z^{'a}\),
\(Z_\delta^{'a}\), \(X_\delta^{'a}\) introduced above are
alternatively denoted by \(Z'\),
\(Z'_\delta\), \(X'_\delta\) (resp. \(Z''\),
\(Z_\delta''\), \(X_\delta''\)). We use \(\scrX^a_\delta \subset
X^{'a}_\delta \) to denote a manifold with smooth boundary obtained by
``rounding corners''. In particular, \(\partial \scrX^a_\delta \) is smooth
and \(\dist (\partial\scrX^a_\delta , \partial X^{'a}_\delta )\ll
2^{-8}\). Let 
\(\gamma _{a,\delta }\) be a smooth cutoff function on \(X\) that is
supported on \(\scrX^{a+1}_{\delta }\) and equals 1 on
\(X^{'a}_{2\delta }\), such that \(\|\gamma _{a, \delta }-\gamma _a \,
(1-\chi (\sigma /\delta ))\|_{C^2}\ll 2^{-8}\).

\begin{remarks} 
(1) As in  \cite{Ts}, the pointwise estimates provided in this section
typically 
hold over domains of the form \(X_\delta \) (or more generally,
\(X^{'a}_\delta \)). They are \(\delta\)-dependent, and constants appearing
the relevant inequalities usually depend on \(\delta\).  In \cite{Ts}, the dependence
of the constants on \(\delta\) are left unspecified. For our
purpose, this dependence is important and therefore shall be made explicit
below. As mentioned in Section \ref{sec:convention}, in what follows the
notation \(\zeta\) and its decorated variants such as \(\zeta '\),
\(\zeta _i\), are reserved for constants
independent of \(\delta\) (and also independent of \(r\) and \((A_r, \Psi _r)\)).

(2) Typically, we  improve the pointwise estimates in Section 3 of \cite{Ts} by replacing
factors of constants \(\delta^{-1}\) therein by the function 
 \(\ts^{-1}\leq \delta^{-1}\).   This is often made possible with the
 help of the following  observation:  Given \(\rmm>0\), there are
 constants \(\zeta _\rmm\), \(\zeta _\rmm'\) depending only on \(\nu \), the
 metric, and \(\rmm\), such that 
\begin{equation}\label{ts-loBdd}
d^*d (\ts^{-k})+\frac{r|\nu|}{\rmm}(\ts^{-k})\geq\zeta _\rmm\, 
(-\zeta'_\rmm\, \ts^{-k-2}+r\ts^{-k+1} )>\zeta _\rmm\, r
\ts^{-k+1}/2
\end{equation}
over \(X''_\delta \) when \(r\delta^3>\zeta '_\rmm\).  (In what follows, \(\rmm\) is typically
taken to be \(2^k\), \(k=1, 2, 3, 4\).) This enables one to replace
\(\delta^{-k}\) by \(\ts^{-k}\) as comparison functions in various
comparison principle arguments. 
\end{remarks}



The pointwise estimates in this section are made simpler thanks to the
next lemma, which motivated the introduction of Item (5) in Assumption \ref{assume}.

\begin{lemma}
Let \((X, \nu )\) be an admissible pair, and let \(w_r\)
and \(\hat{\grp}_r\) respectively be a \(r\)-parametrized family of
closed 2-forms and nonlocal
perturbations on \(X\) satisfying Assumption
\ref{assume}. Let  \(\mu
_r=r\nu+w_r\) as before, and let \((A, \Psi)=(A_r, \Psi_r)\) be as in Lemma
\ref{lem:Etop-bdd1}. Then for any given \(\smE>0\),  there is a constant
\(r_\smE\geq 8\) such that
\begin{equation}\label{eq:nonlocalptwise}
\|\hat{\grp}_r\, (A, \Psi)\|_{C^k_{A}(X)}\leq \zeta _\grp 
\end{equation}
for all \(r\geq r_\smE\) and admissible solutions \((A, \Psi )=(A_r,
\Psi _r)\) to the Seiberg-Witten equation
\(\grS_{\mu _r, \hat{\grp}_r}(A, \Psi )=0\) satisfying the energy
bound (\ref{assume:EtopX-ubdd}). 
\end{lemma}

\pf Invoke Assumption \ref{assume} (5). If \(X\) is cylindrical, then
\(\grY_v=\emptyset\) by our assumption on \(X\), and so
\(\hat{\grp}_r=\hat{\grq}_r\equiv 0\). When \(X\) is non-cylindrical, \(\hat{\grp}(A, \Psi )\) is
supported on \(\bigcup_{i\in \grY_v}\hat{Y}_{i, [\grl_i, \grl_i']}\),
it suffices to examine it on each \(\hat{Y}_{i, [\grl_i,
  \grl_i']}\). In the present setting, \(\upsilon (r)=\zeta
r\). Meanwhile, applying Proposition \ref{prop:SW-L2-bdd} to \(X_\bullet=\hat{Y}_{i, [\grl_i,
  \grl_i']}\), we see that there is an \(r_\smE\geq \max (2^8 \smE, 8)\) such
that Conditions (i) and (ii) in Assumption \ref{assume}
(5) holds for  \((A, \Psi )=(A_r, \Psi _r)\) when \(r\geq
r_\smE\). The assertion (\ref{eq:nonlocalptwise}) follows directly
from Assumption \ref{assume} (5) and the properties of
\(\hat{\grp}_r\) and \(w_r\) prescribed in items (4) and (1) in
Assumption \ref{assume}. 
\epf

We shall apply he preceding lemma to those  \((A_r, \Psi _r)\) from the statement of Theorem
\ref{thm:l-conv}, with the constant \(\smE\) taken to be that
given by (\ref{eq:CSD-est}).
The value of each occurrence of \(r_0\)
in the rest of this article will be taken to be larger or equal to all
its predecessors and the version of 
\(r_\smE\) corresponding to this value of \(\smE\).

Throughout the rest of this section, we tacitly invoke 
the bound (\ref{eq:nonlocalptwise}) to omit terms arising from the
nonlocal perturbation \(\hat{\grp}\) by adjusting the coefficients in
the inequalities.

Write \(\Psi=(r/2)^{1/2}\psi\), and write \(\psi=(\alpha, \beta)\) with
respect to the decomposition \(\bbS^+\simeq E\oplus E\otimes
K^{-1}\) over \(X-\nu^{-1}(0)\). Throughout this section, \((A,
\Psi)=(A_r, \Psi_r)\) is an admissible solution to the Seiberg-Witten equation
\(\grS_{\mu_r, \hat{\grp}}(A_r, \Psi_r)=0\) satisfying the assumptions
of  Lemma \ref{lem:Etop-bdd1}.

\subsection{Estimates for \(|\psi|^2\)}\label{assume6}

With our assumption on \(\nu\), \(w_r\), and \(\hat{\grp}\),  an 
\(L^\infty\)-estimate on \(\psi\) over \(X\) may be established easily.

\begin{lemma}\label{lem:4dPsi:l-infty}
Let \((A, \Psi)\), \(\psi\) be as described immediately preceding this
subsection. Then 
\begin{equation}\label{psi-infty}
\|\psi\|_{L^\infty(X)}\leq \zeta_\infty \quad \text{over \(X\)},
\end{equation}
where \(\zeta_\infty\) is a positive constant depending only on \(\sup_X
R_g\), \(\sup_X |\nu|\), and the constants \(\varsigma_w\),
\(\zeta _\grp\) from Assumption \ref{assume}.  
\end{lemma}
\pf
We argue as in the
   proof of the Morse-end case of Lemma
\ref{lem:3d-Phi}.

By the first line of the Seiberg-Witten equation
\(\grS_{\mu_r,\hat{\grp}}(A, \Psi)=0\), one has \[\langle \psi,
\bar{\slp}^-_A\bar{\slp}^+_A\psi\rangle=-r^{-1/2}\langle \psi ,
\hat{\grp}(A, \Psi )\rangle.\] 
It then follows from the Weitzenb\"ock formula, the
second line of the Seiberg-Witten equation and (\ref{eq:nonlocalptwise}) that 
\begin{equation}\label{ineq:psi}\begin{split}
\frac{1}{2}d^*d|\psi|^2 & +|\nabla_A\psi|^2+\frac{r}{4}
|\psi|^2(|\psi|^2-|\nu+r^{-1}w_r|)+\frac{R_g}{4}|\psi|^2\\ & \leq
\zeta_1|\psi|^2+\zeta_2\, r^{-1},
\end{split}
\end{equation}
where \(\zeta _1\), \(\zeta _2\) are positive constants depending only
on \(\zeta _\grp\). 
The smooth function \(|\psi|^2\) must have a maximum at a certain
point \(x_M\in X\), or it is
bounded by \(r^{-1}\|\Phi_i\|_{L^\infty(Y_i)}\) for certain \(i\), where \((B_i,
\Phi_i)\) is the \(Y_i\)-end limit of \((A, \Psi)\). In the former
case, consider the previous inequality at \(x_M\) and rearranging to
get
\[
|\psi(x_M)|^2\big(|\psi( x_M)|^2-|\nu(x_M)|^2\big)\leq
\zeta_3 |\psi(x_M)|^2+\zeta_4, 
\]
where \(\zeta _3\), \(\zeta _4\) are positive constants depending only
on \(\zeta _\grp\), \(\varsigma_w\), and \(\sup_X
R_g\). 
Hence, \(|\psi(x_M)|^2\leq \zeta_5\) for a positive constant \(\zeta
_5\) depending only
on \(\zeta _\grp\), \(\varsigma_w\), and \(\sup_X
R_g\), and \(\sup_X |\nu|\). In the latter case, invoke
Lemma \ref{lem:3d-Phi}. Either way, Equation (\ref{psi-infty}) holds. 
\epf

Over \(X^{'9}\), a better pointwise bound on \(|\psi|^2\) may be obtained from the
 \(L^\infty\) bound in Lemma \ref{lem:4dPsi:l-infty}.

\begin{prop}\label{T:lem3.2}
There is a constant \(\zeta\) depending only on
the metric, \(\nu \), \(\zeta
_\grp\), and \(\varsigma_w\). 
such that over \(X^{'a}\), \(0\leq a\leq 9\), 
\begin{eqnarray}
(a) & |\psi|^2 &\leq |\nu|+\zeta r^{-1/3};\nonumber\\
(b) & |\psi|^2 &\leq |\nu|+\zeta r^{-1}(\sigma^{-2}+1). \label{ineq:psi^2}
\end{eqnarray}
\end{prop}
\pf This is an analog of Lemma 3.2 of \cite{Ts}\footnote{Equation (3.11)
  of \cite{Ts} contains some errors/typos which are easy to
  fix.}.

From Equation (\ref{ineq:psi}) and Lemma \ref{lem:4dPsi:l-infty} we have
\begin{equation}\label{u-ineq1}
\frac{1}{2}d^*d|\psi|^2+|\nabla_A\psi|^2+\frac{r}{4}
|\psi|^2(|\psi|^2-|\nu|)\leq \zeta_0, 
\end{equation}
where the constant \(\zeta _0>0\) depends on the parameters listed in
the previous lemma. 
Let \(\mathfrc{v}:=\gamma_a |\nu|+(1-\gamma_a)(\zeta_\infty^2+1)\),
\(\zeta_\infty\) being the constant from the
previous lemma. Then by (\ref{eq:omega-ineq}) 
\[
-d^*d \mathfrc{v}\leq
\zeta_1\gamma_a|\nu|^{-1}|\nabla\nu|^2+\zeta_2\gamma_{a+1}, 
\]
where \(\zeta _1\), \(\zeta _2\) depend only on \(\nu \) and \(\zeta
_\infty\).

Combine the above two inequalities, and setting
\(u:=|\psi|^2-\mathfrc{v}\), we have the following analog of (3.12) in
\cite{Ts}:
\begin{equation}\label{ineq:u}
2^{-1}d^*du+4^{-1} r|\nu| u\leq \zeta'(\gamma_a\, \tilde{\sigma}^{-1}+\gamma_{a+1}),
\end{equation}
where \(\zeta '\) depends only on the metric, \(\nu \), \(\zeta
_\grp\), and \(\varsigma_w\).

Set \(\chi_R:=\chi(r^{1/3}\sigma/R)\), and argue similarly to
(3.14)-(3.15) of \cite{Ts}.
We may find positive constants \(\zeta_3,
\zeta_4\) such that with \begin{equation}
\hat{u}:=u+\zeta_3\chi_1\sigma-\zeta_4\, 
r^{-1/3}\quad \text{on \(X''\)},
\end{equation}
one has
\[
2^{-1} d^*d\hat{u}+4^{-1}r|\nu|\hat{u}\leq 0\quad \text{on \(X''\)}.
\]
Suppose \(\hat{u}\) has a maximum in the interior of \(X''\), then the
above inequality implies that \(\hat{u}<0\) at this maximum. Otherwise,
\(\sup_{x\in X''}\hat{u}(x)\) appears in \(\partial X''\) or as a
value of
\(\hat{u}_i:=r^{-1}|\Phi_i|^2-|\nu_i|u+\zeta_3\chi_1\sigma-\zeta_4\, 
r^{-1/3}\) for one of the Morse ends \(\hat{Y}_i\). In the first case, 
\(\hat{u}\leq 0\) by our choice of \(\mathfrc{v}\), \(\zeta_3\), and
\(\zeta_4\). In the second case, this value is also nonpositive by
(\ref{eq:Phi-ptws}). (Adjust the values of \(\zeta_3, \zeta_4\) if necessary).
Equation (\ref{ineq:psi^2}) (a) now follows. 

For (\ref{ineq:psi^2}) (b), argue as in \cite{Ts} using
(\ref{ineq:u}) to find positive constants \(\zeta_5, \zeta_6\)
so that for \(\check{u}:=u-\zeta_6\, r^{-1}(\tilde{\sigma}^{-2}+1)\), \[2^{-1}
d^*d\check{u}+4^{-1}r|\nu|\check{u}\leq 0
\]
over the region \(\{x\, |\, \sigma(x)\geq
\zeta_5\, r^{-1/3}, \, x\in X''\}\). Now apply the maximum principle type
arguments 
over this region as in the proof of  (\ref{ineq:psi^2}) (a), using
 (\ref{ineq:psi^2}) (a)  to ensure that \(\hat{u}\leq 0\) on the boundary of this region. 

When \(|\nu|^{-1}(0)=\emptyset\), simply set
\(\mathfrc{v}=|\nu|\). Then \(-d^*d\mathfrc{v}\leq \zeta\) in this
case. The argument above may then be simplified by dropping all terms
involving \(\gamma_a\) or \(\tilde{\sigma}\), \(\sigma\). In this case
\(|\psi|^2\leq |\nu|+\zeta r^{-1}\) over \(X\).
\epf

\subsection{Estimates for \(|\beta|^2\)}

Coming up next is an analog of Proposition 3.1 of \cite{Ts}.

\begin{prop}\label{T:prop3.1-}
There exist positive constants \(\sO\geq 8\), \(c\), \(c'\)
\( \zeta_0, \zeta '_0\geq 1\)  that depend only on
the metric, \(\nu \) \(\varsigma_w\), and \(\zeta _\grp\), 
such the following hold: Suppose \(r>1\), \(\delta >0\) are such that
\(\delta\geq\textsc{o}r^{-1/3}\), then 
\begin{equation}\label{ineq:beta0}
\begin{split}
|\beta|^2& \leq 2 c \, \ts^{-3}r^{-1} (
|\nu |-|\alpha|^2)+\zeta _0\, \ts^{-5} r^{-2};\\
|\beta|^2& \leq 2 c' \ts^{-3}r^{-1} (
|\nu |-|\psi|^2)+\zeta '_0\, \ts^{-5} r^{-2}
\end{split}
\end{equation}
on  \(X^{'a}_\delta\), \(0\leq a\leq \frac{25}{3}\). 
\end{prop}


\pf 
Proceed as in \cite{Ts} to get the following analog of (3.19) of
\cite{Ts}. \footnote{There is a sign error in (3.19) of \cite{Ts}.}
\begin{eqnarray}
& 2^{-1} d^*d |\beta|^2 &+|\nabla_A\beta|^2+4^{-1} r|\nu|\,
|\beta|^2+ 4^{-1} r (|\alpha|^2\, |\beta|^2+|\beta|^4)\nonumber\\
&& \leq (\zeta +\zeta_0b^2)|\beta|^2 +\zeta_1
\big(|b|\, |\nabla_A\alpha|+ \zeta _1'
\big)|\beta|, 
\label{ineq:beta}\\
& 2^{-1} d^*d |\alpha|^2 &+|\nabla_A\alpha|^2-4^{-1} r|\nu|\,
|\alpha|^2+ 4^{-1} r (|\alpha|^2\, |\beta|^2+|\alpha|^4)\nonumber\\
&& \leq \zeta'|\alpha|^2 +\zeta_2
\big(|\nabla_A(b\beta)|+\zeta _2'
\big)|\alpha|
\label{ineq:alpha}
\end{eqnarray}
on \(X''\), where \(b\) arises from
\(\nabla J\) and can be bounded by
\(|b|\leq\zeta_0\, \tilde{\sigma}^{-1}\), \(|\nabla b|\leq
\zeta'_0\, \tilde{\sigma}^{-2}\) on \(X''\). In the above, the positive constants,
\(\zeta _i, \zeta _i'\), \(i=0,1,2\) depend only on
the metric, \(\nu \) \(\varsigma_w\), and \(\zeta _\grp\). 

Judicious uses of the
triangle inequality shows that there exist constants \(\epsilon
_1<1\), \(c\geq 1\), \(\zeta _1\), \(\zeta _1'\) depending only on
the metric, \(\nu \) \(\varsigma_w\), and \(\zeta _\grp\), such that 
when \(r^{-1}\delta^{-3}<\epsilon _1\),  the inequalities
(\ref{3.20}), 
(\ref{3.21}) below hold over \(X''_{\delta/3}\):  
\begin{equation}\label{3.20}\begin{split}
2^{-1} d^*d |\beta|^2 &+|\nabla_A\beta|^2+8^{-1} r|\nu|\,
|\beta|^2+ 4^{-1} r (|\alpha|^2\, |\beta|^2+|\beta|^4)\\
& \leq c\, (r^{-1}\ts^{-3} |\nabla_A \alpha |^2+r^{-1}\ts^{-1})
\end{split}\end{equation}
Set \(\varpi:=|\nu|-|\alpha|^2\). Combine (\ref{ineq:alpha}) and
(\ref{eq:omega-ineq}) and use Lemma \ref{lem:4dPsi:l-infty} 
to get 
\begin{equation}\label{3.21}\begin{split}
2^{-1}d^*d(-\varpi) +|\nabla_A\alpha|^2+8^{-1} &
r|\nu|\,(-\varpi)+4^{-1} r\varpi^2\\
&\leq \zeta_1|\nabla_A\beta|^2+\zeta_1' \ts^{-1}  \quad \text{over
  \(X^{'a}_{\delta/3}\), \(0\leq a\leq 9\).}
\end{split}\end{equation}
(To get the last term above, we used Proposition \ref{T:lem3.2} 
to bound
\(\tilde{\sigma}^{-2}|\alpha|^2\leq \zeta_3\, \ts^{-1}\) and invoked
\(r^{-1}\delta^{-3}<\epsilon _1\) to simplify terms.)

A combination of the previous two inequalities then yield that when
\(r^{-1}\delta^{-3}<\epsilon _2:=\min\,  (\epsilon _1, (4c\zeta
_1)^{-1})\), for \(0\leq a\leq 9\)
\begin{equation}\label{de:beta}
2^{-1}|\nabla_A\beta|^2+cr^{-1}\delta ^{-3} |\nabla_A \alpha |^2+2^{-1}d^*du_1+8^{-1} r|\nu|u_1\leq 0,\quad \text{over
  \(X^{'a}_{\delta/3}\)}
\end{equation}
for a suitable constant \(\zeta \geq 1\) and
\begin{equation}\label{def:u_1}
u_1:=|\beta|^2-2cr^{-1}\delta^{-3}\varpi-\zeta r^{-2}\delta^{-5}.
\end{equation}

Fix an \(x\in X^{'a}_\delta\), \(0\leq a\leq \frac{25}{3}\), and use the abbreviation \(\ts_x=\ts (x)\) below. Then \(B(x, \ts_x/3)\subset B(x,
2\ts_x/3)\subset X_{\delta/3}^{'9}\). Define the function
\(\lambda_x(\cdot):=\chi(3\dist (x, \cdot)/\ts_x)\). Then \(\lambda_x
u_1\) is supported on \(B(x, 2\ts_x/3)\) and satisfies
\[
2^{-1}d^*d(\lambda_xu_1)+8^{-1} r|\nu|(\lambda_xu_1)\leq \xi,
\]
where \(\xi\) is supported
on the shell \(A_{\ts_x}:=B(x, 2\ts_x/3)-B(x, \ts_x/3)\), and is
bounded by 
\[
|\xi|\leq \ts_x^{-2}(\zeta '_3 |\psi |^2+\zeta '_4 r^{-1}\delta ^{-3}+\zeta '_5
r^{-2}\delta ^{-5}). 
\] 
By Lemma \ref{lem:4dPsi:l-infty}, this means that \(|\xi|\leq \zeta
'_2\ts_x^{-2}\) when \(r^{-1}\delta ^{-3}<\epsilon _3\) for a certain
\(\epsilon _3\leq \epsilon _2\). 
Let \(\zeta>0\) be a
constant so that \(|\nu|\big|_{X_{\ts_x/3}}\geq 4\zeta\ts_x\), and
  let \(\mu_x\) be the
solution to the equation
\[
2^{-1}d^*d\mu_x+\zeta r\ts_x \mu_x=|\xi|  
\]
on \(B(x, 2\ts_x/3)\) with Dirichlet boundary condition. Then
\(\lambda_x u_1\leq \mu_x\) by the comparison principle. Meanwhile,
using \(G_r\) to denote the integral kernel of the operator
\(2^{-1}d^*d+\zeta r\ts_x \) on \(B(x, 2\ts_x/3)\) with Dirichlet
boundary condition, 
\[
\begin{split}
\mu_x(x) &=\int_{B(x, 2\ts_x/3)}G_r(x, y)|\xi|(y) dy\\ & \leq \zeta'_1\ts_x^{-2}
\int_{A_{\ts_x}} \big( \dist (x, y)^{-2}\exp \big(-\zeta_2 \dist (x, y) (r\ts_x)^{1/2}\big) \Big) dy\\
&\leq  \zeta_4 \ts_x^2\exp \big(-\zeta'_2
(r\ts_x^3)^{1/2}\big) \leq 
\zeta_5 r^{-2}\ts_x^{-4} 
\end{split}
\]
when \(r^{-1}\ts_x^{-3}\leq r^{-1}\delta ^{-3}<\epsilon _4\leq
\epsilon _3\) for certain constant \(\epsilon _4\). 
Combining this with the bound
\(u_1(x)\leq \mu_x(x)\), we have 
\begin{equation}\label{ineq:beta3-}
|\beta|^2 \leq 2 c \delta ^{-3}r^{-1} (
|\nu |-|\alpha|^2)+\zeta _0\, \delta ^{-5} r^{-2}\quad \text{over \(X^{'a}_\delta \)}
\end{equation}
when \(r^{-1}\delta ^{-3}\leq\epsilon _5\), where 
 \(\zeta _0\) and \(\epsilon _5\leq \epsilon _4\) are certain constants
which  are independent of \(r,
\delta \), and \((A, \Psi )\). Given \(r\), \(\delta _0\) satisfying
\(r^{-1}\delta ^{-3}_0<\epsilon _5\), 
 fix \(x\in X^{'a}_{2\delta_0}\) and set
\(\delta=\ts_x/2\) in (\ref{ineq:beta3-}). It follows that 
\[
|\beta|^2 (x) \leq 2 c \ts_x^{-3}r^{-1} (
|\nu |-|\alpha|^2)+\zeta _0 \, \ts^{-5}_x r^{-2}\quad \text{over \(X^{'a}_{2\delta _0}\)}. 
\]
This implies the first inequality in (\ref{ineq:beta0}), with
\(\sO=2\epsilon _5^{-1/3}\). 
The second inequality in (\ref{ineq:beta0}) follows directly from the first. \epf


\subsection{Estimating \(|F_{A}|\)} 

 
We use Lemma \ref{T:lem3.1}  to obtain a preliminary bound on
\(|F_A|\) on \(X^{'a}\). 
Let \(\td{r}\) denote the function \[x\mapsto
  \td{r}_x:=r\ts (x)\] on \(X\). Note that \(\frac{1}{2}\ln r\leq \ln
\td{r}_x\leq \ln r\) 
when \(x\in X_{\delta }\) for \(\delta \leq 1\) and
\(r\geq\delta ^{-2}\).

\begin{prop}\label{T:prop3.2}
{\bf (a)} There exist positive constants \(\zeta_2\),
\(\zeta'_2\), that depend only on \(\varsigma_w\) and \(\zeta
_\grp\), 
such that 
\begin{equation}\label{eq:F+bdd}
|F_A^+| \leq 2^{-3/2}r\, (|\nu|-|\psi|^2)+\zeta _2r
|\beta|^2+\zeta'_2\quad \text{over \(X\).  }
\end{equation}

{\bf (b)} 
There exist positive constants \(r_0>8\), 
\(\sO\geq 8\), \(\zeta \), \(\zeta '\), \(\zeta_1\), \(\zeta _1'\)
that depend only on the metric, \(\nu \), \(\varsigma_w\) and \(\zeta
_\grp\), 
 so that the following holds: Suppose \(r>r_0\), and let \(\delta
 _0:=\textsc{o}r^{-1/3}\).  Then for \(0\leq a\leq \frac{17}{2}\), 
 \begin{equation}\label{bdd:ss}
\begin{split}
\text{(i)} \quad  |F_A^-| &\leq
(2^{-3/2}+\varepsilon_0
) \, r\, (|\nu|-|\psi|^2) +K_0 
\quad \text{over \(X^{'a}_{\delta_0}\)};  \\
\text{(ii)} \quad |F_A^-| &\leq  
\zeta _1 r  \delta +\zeta _1' \delta^{-2}\ln\,  (\delta /\sigma )
\quad \text{over \(Z^{'a}_{\delta}\) for any \(\delta\geq \delta_0\)},\\
\end{split}
\end{equation}
where 
\begin{equation}\label{def:e_0K_0}
  \begin{split}
    \varepsilon_0& :=
    \begin{cases}
    \big(r^{-4/3}\ts^{-1} (\ln  r+\smE)\big)^{2/7} & \text{where \(\sigma
     \geq 3  r^{-1/6} (\ln r+\smE )^{1/6} \)};\\
   r^{-1/6}\, \ts\, (\ln r+\smE )^{-1/2}+\zeta '_5\sO_r^{-3}  &\text{otherwise};
 \end{cases}\\
 K_0 &:= \begin{cases}
      \zeta' (r\ts^{-1})^{5/7}\, (\ln r+\smE )^{2/7} & \text{where \(\sigma
      \geq 3  r^{-1/6} \, (\ln r+\smE )^{1/6} \)};\\
  \zeta '' r^{1/2}\ts^{-2}(\ln r+\smE )^{1/2}   &\text{otherwise}.
   \end{cases}
\end{split}
\end{equation}
{\bf (c)} 
For \(r>r_0\), there exist positive constants \(\zeta _3\), \(\zeta
_3'\) that depend only on the metric, \(\nu \), \(\varsigma_w\) and \(\zeta
_\grp\), 
 so that the following holds over \(X^{'a}_{\delta _0}\): 
 \[\begin{split}
|F_A^-| &\leq\zeta _3 \, r\ts; \\
|F_A^-| &\leq
(2^{-3/2}+\varepsilon_0)\,  r\, (|\nu|-|\psi|^2) +K_1,
\end{split}
\]
where \begin{equation}\label{def:K_1}
  K_1:=\min \, (K_0, \zeta _3  \, r\ts) \leq
    \zeta_3'\, r^{5/6}\, (\ln r +\smE )^{1/6}.
  \end{equation}
 \end{prop}

\pf The estimate for \(|F_A^+|\) is a direct consequence of the the
Seiberg-Witten equation \(\grS_{\mu_r,\hat{\grp }}(A, \Psi)=0\) and (\ref{eq:nonlocalptwise}).

Let \(\textsc{s}:=|F_A^-|\). 
The arguments leading to \cite{T}'s (I.2.19), together with
(\ref{eq:nonlocalptwise}) and (\ref{psi-infty})  give: 
\begin{equation}\label{eq:DE-s}
\begin{split}
& \big(\frac{ d^*d}{2}+\frac{r|\psi|^2}{4}\big) \ss\\
& \quad \leq z_4\, \ss +2^{-3/2}r|\nabla_A\psi|^2
+\zeta _1r\ts^{-2}|\psi|\, |\beta|\\
& \qquad \qquad \quad +\zeta_2r\ts^{-1}(|\nabla_A\psi|\, |\beta |+|\alpha |\,
|\nabla_A\beta |)+\zeta_0\\
& \quad \leq z_4\, \ss +2^{-3/2} r|\nabla_A\psi|^2 
+\zeta_2r\ts^{-1}|\nabla_A\psi|\, |\beta| \\
& \qquad \qquad \quad +\zeta_2\, r|\nabla_A\beta |^2+\zeta '_1r\ts^{-2}|\psi |^2 +\zeta_0\quad \text{over \(X^{'9}\),}\\
\end{split}
\end{equation}
where the constants \(z_4\), \(\zeta _1\), \(\zeta _1'\), \(\zeta _2\)
depend only on the metric and \(\nu \); the constant \(\zeta_0\)
depends on the same parameters, and additionally on \(\varsigma_w\)
and \(\zeta _\grp\). 
(This is a refinement of \cite{Ts}'s
(3.32).) 

It follows from (\ref{ineq:beta}) that 
\[
\begin{split}
&  \big(\frac{ d^*d}{2}+\frac{r|\psi|^2}{4}\big)  |\beta|^2 +|\nabla_A\beta|^2+4^{-1} r|\nu|\,
|\beta|^2\\
& \qquad \quad \leq \zeta \ts^{-2}|\beta|^2
+\zeta'\ts^{-1}|\nabla_A\alpha|\, |\beta|\qquad 
\text{ over \(X''\)},
\end{split}
\]
where \(\zeta \), \(\zeta '\) depend only on the metric, \(\nu \), \(\varsigma_w\)
and \(\zeta _\grp\). So
\begin{equation}\label{ineq:beta2}
\begin{split}
& \big(\frac{ d^*d}{2}+\frac{r|\psi|^2}{4}\big) (\ss +\zeta _2 r|\beta
|^2)+ 4^{-1} \zeta _2r^2|\nu|\,|\beta|^2\\
& \quad \leq z_4\, \ss +2^{-3/2} r|\nabla_A\psi|^2 
+\zeta_3r\ts^{-1}|\nabla_A\psi|\, |\beta| +\zeta '_3r\ts^{-2}|\psi
|^2 +\zeta _0\\
& \quad \leq z_4 \, \ss +c_\varepsilon\, r |\nabla_A\psi|^2 +4^{-1}\zeta _3^{2}\varepsilon^{-1}r \ts^{-2}\,|\beta|^2+\zeta '_3r\ts^{-2}|\psi
|^2  +\zeta _0\quad \text{over \(X^{'9}\),}\\
\end{split}
\end{equation}
where \(\varepsilon \) be an arbitrary positive number small than \(8\), and  
\[
c_\varepsilon := 2^{-3/2}+\varepsilon .
\]  
In the above inequalities as well as for the
rest of this proof, all constants denoted in the form of \(\zeta _*\),
\(z_*\) depend only on the metric, \(\nu \), \(\varsigma_w\)
and \(\zeta _\grp\); in particular, they are independent of \(\varepsilon \) 

Let  \(\delta _0:=\sO r^{-1/3}\), 
where \(\sO\) is as in Proposition
\ref{T:prop3.1-}, and  let \(\sO_r:=r^{1/3}\ts\).  Writing
\(u:=|\psi|^2-|\nu|\) and 
appealing to
Propositions \ref{T:lem3.2} and 
\ref{T:prop3.1-},  Equation (\ref{ineq:beta2}) implies:
{\small
\begin{equation}\label{ineq:beta3}
\begin{split}
& \big(\frac{ d^*d}{2}+\frac{r|\psi|^2}{4}\big) (\ss +\zeta _2 r|\beta
|^2) \leq z_4 \ss +c_\varepsilon\, r |\nabla_A\psi|^2 \\
& \quad + 
(1-\chi(\delta_0^{-1}\sigma))\,\Big(  \zeta _3'r\ts^{-2}|\psi|^2
+\zeta \, \varepsilon ^{-1}\sO_r^{-3} r\ts^{-2}
 (-u) +\zeta '\varepsilon ^{-1}\sO_r^{-6}r \ts^{-1}+\zeta_0\Big)\\ 
& \qquad\qquad 
+ \zeta'_2 \chi (\delta_0^{-1}\sigma)\, \varepsilon^{-1}(r\ts^{-1}+r^{2/3} \ts^{-2}\big)
\quad \text{over \(X^{'9}\).}\\
\end{split}
\end{equation}
}
Meanwhile, by (\ref{u-ineq1}), (\ref{eq:omega-ineq}) we have
\begin{equation}\label{u-ineq2a}
\begin{split}
& \big(\frac{ d^*d}{2}+\frac{r|\psi|^2}{4}\big) \,
u+|\nabla_A\psi|^2\leq \zeta'_0\, \ts^{-1}\qquad
\text{ over \(X''\).}
\end{split}
\end{equation}
Let 
\begin{equation}\label{def:q0}
q_0=q_0^{(\varepsilon )}:=\ss +c_\varepsilon  r\, u+\zeta
_2r|\beta |^2.
\end{equation}
A combination of (\ref{u-ineq2a}) and (\ref{ineq:beta3}) then gives:  
\begin{equation}\begin{split}\label{ineq:q_0}
  \big(\frac{ d^*d}{2}+\frac{r|\psi|^2}{4}\big) q_0  
 & \leq z_4q_0+ \zeta'_1 \chi (\delta_0^{-1}\sigma)\,
 \varepsilon^{-1}(r\ts^{-1}+r^{2/3} \ts^{-2}\big)\\
&  \quad + \zeta ' \, (1-\chi
(\delta_0^{-1}\sigma))\, \Big( r\ts^{-2}(-u)  \big(c_\varepsilon\ts^2+\zeta \, \varepsilon ^{-1}\sO_r^{-3} )\\
& \qquad \quad \quad + \zeta _3'r\ts^{-2}|\psi|^2+r \ts^{-1}( \zeta _0'+\zeta '\varepsilon ^{-1}\sO_r^{-6})+\zeta_0\Big)\\
\end{split}
\end{equation}
on \(X^{'9}\).  In order to to get rid of the term  \((1-\chi
(\delta_0^{-1}\sigma))\, r \ts^{-1}( \zeta _0'+\zeta '\varepsilon
^{-1}\sO_r^{-6})\) above, we will work with  a modified version of
\(q_0\). Note that  for any \(\delta
\geq\sO r^{-1/3}\) and \(k>-1\), 
{\small 
\[\begin{split}
& \big(\frac{ d^*d}{2}+\frac{r|\psi|^2}{4}\big) \big(-(1-\chi (\delta^{-1}\sigma ))\, \ts^{-k}\big)
\leq -\frac{r|\psi|^2}{4} \big(1-\chi (\delta^{-1}\sigma )\big)\, \ts^{-k}\\
&  \qquad \qquad +(1-\chi (2\delta^{-1}\sigma
))\,\big(\zeta_7\, \ts^{-k-2}+\zeta_7'  \delta^{-2} \chi (\delta
^{-1}\sigma /2)\,\ts^{-k}\big)\\
& \quad \leq \big(1-\chi (2\delta^{-1}\sigma )\big)r\ts^{-k}\, \big( -(1/4-\zeta
\sO^{-3}_r)\, |\psi|^2+\zeta
'\sO^{-3}_r(-u)+\zeta_7'  \delta^{-2} \chi (\delta
^{-1}\sigma /2)\big)\\
 & \quad \leq \big(1-\chi (2\delta^{-1}\sigma )\big)\, r \Big(
 -\zeta_1\ts^{-k+1} +\ts^{-k}\big(\zeta _1'(-u) +\zeta_7'  \delta^{-2} \chi (\delta
^{-1}\sigma /2)\big)\Big).
\end{split}
\]
}
Combining (\ref{ineq:q_0}) with the preceding inequality and making
use of Proposition \ref{T:lem3.2},
we see that  there exist constants \(\zeta _5\), \(\zeta _5'\) that
are  independent
of \(r\) and \(\delta \), such that with 
\begin{equation}\label{def:q}
\begin{split}
q=q^{(\varepsilon )}& :=q_0^{(\varepsilon )}- (1-\chi (\delta_0 ^{-1}\sigma ))\, \ts^{-2}(\zeta _5+\zeta_5'
\varepsilon ^{-1} \sO_r^{-6}),\\
\end{split}
\end{equation}
one has: 
 \begin{equation}\label{DE:q}
\begin{split}
  &  \big(\frac{d^*d}{2}+\frac{r|\psi|^2}{4}\big) q\\
  & \qquad \leq z_4
q+ (\zeta _8'+\zeta \, \varepsilon ^{-1}\sO_r^{-3} )(1-\chi (2\delta_0 ^{-1}\sigma ))\, r\ts^{-2}(-u+z_8r^{-1}) \\
& \qquad \qquad +\zeta _8\, \varepsilon ^{-1}\, \chi (\delta_0
^{-1}\sigma /2) \,(r\ts^{-1} +\delta _0^{-2} \ts^{-2}) 
\quad \text{over \(X^{'9}\).}
\end{split}
\end{equation}
for  certain \(r\)-independent constants \(\zeta \), \(\zeta
_8\), \(\zeta '_8\), \(z_8\). 
Thus, writing 
\[
\begin{split}
& \eta_1:=z_8\, (\zeta _8'+\zeta \, \varepsilon ^{-1}\sO_r^{-3} )\, (1-\chi
(2\delta_0 ^{-1}\sigma ))\, \ts^{-2}\\
& \eta_2:=(\zeta _8'+\zeta \, \varepsilon ^{-1}\sO_r^{-3} )\, (1-\chi (2\delta_0 ^{-1}\sigma ))\, r\ts^{-2}(-u);\\
& \xi  := \zeta_8\, \varepsilon ^{-1}\chi (\delta_0 ^{-1}\sigma /2)\,
( r\ts^{-1}+\delta _0^{-2}\ts^{-2}),
\end{split}
\]
one has
\[
d^*dq_+-z_4q_+\leq  \eta_1+\eta_2+   \xi\quad 
\text{over \(X^{'9}\).}
\]
\footnote{
The formula above is regarded as
  inequalities between distributions. See e.g. \cite{Ts} p.187 for justification.
}
Fix \(x\in X^{'a}\), \(0\leq a\leq \frac{17}{2}\). Let \(\rho_*>0\) be such that \(\rho _*\leq 1/4\)
and the operator \(d^*d-z_4\) with Dirichlet
boundary condition on \(B(x, \rho)\) has positive spectrum \(\forall
\rho\leq 2\rho_*\). Suppose that \(\rho _0>0\) is no larger than
\(\frac{1}{2}\rho _*\). 
Multiplying both sides of the preceding differential inequality by \(\chi (\rho_0^{-1}\dist (x,\cdot))\) times \(G(x, \cdot)\), the
Dirichlet Green's function for \(d^*d-z_4\) on a ball of radius \(2\rho_0\), and integrate
over \(B (x,2\rho_0)\). We then have
\begin{equation}\label{eq:q-est}
\begin{split}
 q(x)  & \leq c_{0}\, \rho_0^{-4} \int_{B(x, 2\rho_0)-B(x,\rho_0)}q_++c_1\int_{B(x,2\rho_0)}
\xi \dist (x,\cdot)^{-2}\\
& \qquad \quad 
+c_2\int_{B(x,2\rho_0)}
\eta _1\dist (x,\cdot)^{-2}+c_3\int_{B(x,2\rho_0)}
\eta _2\dist (x,\cdot)^{-2}.
\end{split}
\end{equation}

Note that by Propositions \ref{T:prop3.1-} and \ref{T:lem3.2}, 
\begin{equation}\label{bdd:q-0}
\begin{split}
q& \leq \ss +\zeta \chi (\delta _0^{-1}\sigma )\, r\ts+\zeta '\ts^{-2}
\quad \text{over \(X^{'9}\).}\\
\end{split}
\end{equation}
Thus, the first term on the right hand side of (\ref{eq:q-est}) may be
bounded via  the facts that 
 \[
 \int_{B(x, 2\rho_0)}\ss  \leq \zeta \, r^{1/2}\, (\ln r+\smE )^{1/2}\quad \text{by
   (\ref{eq:L^2_1});}
\]
and 
 \begin{equation}\label{bdd:B_rho0}
\int_{B(x,2\rho_0)}\chi (d_0^{-1}\sigma )\, \ts^{k} \leq \zeta
d_0^{k+3}\rho_0\quad \text{for any \(d_0>0\) and \(k>-3\)}.
\end{equation}
This gives: 
\[
\rho ^{-4}_0\int_{B(x,2\rho_0)-B(x,\rho_0)}q_+\leq \zeta_0'\, r^{1/2}\, (\ln
r +\smE )^{1/2}\rho_0^{-4}.
\]

Now suppose \(x\in X^{'a}_{8\delta }\), \(\delta \geq \delta _0:=\sO
r^{-1/3}\). Choose \(\rho _0>0\) to be such that
\(\rho _0< 2\rho _*\).  
The remaining integrals (\ref{eq:q-est}) are bounded in this case as
follows.  We adopt the shorthand
\(\sigma _y:=\sigma (y)\) and \(\ts_y:=\ts (y)\) in what follows. 
  
When \(\rho _0\) is chosen to be sufficiently small, the second integral on the right hand side of
(\ref{eq:q-est}) can be bounded by the following computation: 
Use \(z\) to parametrize \(\nu ^{-1}(0)\cap B(x,2\rho_0)\), and for \(y\in B(x,2\rho_0)\), let \(z_y:= (4\rho _0^2-(\sigma _y-\sigma
_x)^2)^{1/2}\). 
\[\begin{split}
\int_{B(x,2\rho_0)} \xi \dist (x,
\cdot)^{-2}
 & \leq \zeta \varepsilon ^{-1}\int_0^{4\delta _0}\int^{z_y}_{-z_y}\frac{r  \ts_y+\delta_0^{-2}}{(\sigma _y-\sigma _x)^2+z^2}\, dz\, d\sigma _y, \\
& \leq \zeta ' \varepsilon ^{-1}r\delta_0^{2} \, \sigma _x^{-1}
\quad \text{when \(x\in X^{'a}_{8\delta}\)}. 
\end{split}
\]

To bound the third integral in (\ref{eq:q-est}), take \(\rho
_1=\min \, (\ts(x)/4, \rho _0)\) for \(x\in X^{'a}_{8\delta }\) and
separate \(B(x,2\rho_0)\) into two regions: \(\dist (x, \cdot)\leq
\rho _1\) on the first region, and \(\dist (x, \cdot)>
\rho _1\) on the second. Integrate over the two regions separately,
using the facts that on the first region,
\(\ts\geq\frac{3}{4}\ts_x\), while on the second region, \(\dist (x,
\cdot)\geq \min (\rho _0, \ts_x/4)\) and \(\delta _0/2\leq \ts\leq
2\rho _0\). We get:
\[
\int_{B(x,2\rho_0)} \eta_1\dist (x,\cdot)^{-2}\leq 
\zeta   _1+\zeta_2'\, \varepsilon ^{-1}\sO_r^{-3}+\zeta_3 \, \ts^{-2}\, ( \rho
_0+\varepsilon ^{-1} r^{-1}\delta_0^{-2})\, \rho _0
 \quad 
\]
when \(x\in X^{'a}_{8\delta}\).
(In the above, the first two terms on the RHS comes from the first
region, while the last term comes from the second region.)

To bound the last integral in (\ref{eq:q-est}),
choose small positive numbers \(\rho_1, \delta_2\), such that
\(\rho _1\leq\rho _0 \) 
and \(\sigma (x)/2 >\delta_2\geq\delta_0 \). 
Separate \(B(x,2\rho_0)\) into the
three regions: \(\scrR_1:=B(x, \rho_1)\cap X^{'a}_{\delta _2}\), \(\scrR_2:=X^{'a}_{\delta_2}\cap (B(x,2\rho_0)-B(x,\rho_1)) \),
 \(\scrR_3:=(Z^{'a}_{\delta_2}-Z^{'a}_{\delta _0/2})\cap
B(x,2\rho_0)\), and integrate separately. Using 
the fact that 
\(
\eta _2\leq  \zeta r\ts^{-1}(1+\varepsilon ^{-1} \sO_r^{-3})
\) over \(\scrR_1\)  and \(\scrR_3\), and Lemma \ref{T:lem3.1}
over \(\scrR_2\), we get: 
\begin{equation}\label{int:eta2}
\begin{split}
& \int_{B(x,2\rho_0)}\eta _2\dist (x,\cdot)^{-2} \\
& \qquad \leq \zeta \,
r\ts^{-1}_x\, \rho _1^2 \, (1+\varepsilon ^{-1} \sO_r^{-3})
+ \zeta '\, \rho_1^{-2} \delta _2^{-3}(\ln r)\, (1+\varepsilon ^{-1} r^{-1}\delta_2^{-3})\\
&\qquad \qquad  \quad 
+ \zeta _3'\int_{\delta _0/2}^{\delta
  _2}\int^{z_y}_{-z_y}\frac{r  \ts_y+\varepsilon^{-1} \ts^{-2}_y
   }{(\sigma _y-\sigma _x)^2+z^2}\, dz\, d\sigma _y, \\
& \qquad \leq \zeta \,
r\ts^{-1}_x\, \rho _1^2 \, (1+\varepsilon ^{-1} \sO_r^{-3})
+ \zeta '\, \rho_1^{-2} \delta _2^{-3}(\ln r +\smE )\,
(1+\varepsilon ^{-1} r^{-1}\delta_2^{-3})\\
&\qquad \qquad \quad +\zeta _4'\, 
\big( r\, \delta _2^2+\varepsilon^{-1} \delta_0^{-1}\big) \, \sigma _x^{-1}.\\
\end{split}
\end{equation}
Put together, we have for \(q(x)\leq f_1+\varepsilon ^{-1} f_2\)  for
\(x\in X^{'a}_{8\delta }\), where \(f_1\), \(f_2\) are given as follows. For \(x\in X^{'a}_{8\delta _0}\), 
\[\begin{split}
f_1& =\zeta _1\ts^{-2}_x+\zeta _2 r\ts^{-1}_x\rho _1^2+\zeta _3\, \rho
_1^{-2}\delta _2^{-3} (\ln r +\smE )+\zeta _4\, r \delta _2^2\, \sigma
_x^{-1}\\
\end{split}
\]
Minimizing the RHS by choosing  \(\rho _1=r^{-1/4}\ts^{1/4}_x\delta _2^{-3/4} (\ln
r+\smE )^{1/4}\), one has:
\[
f_1\leq \zeta _1\ts^{-2}_x+\zeta _2'\,  r^{1/2}\ts^{-1/2}_x\delta
_2^{-3/2} (\ln r+\smE )^{1/2}+\zeta _4\, r \delta _2^2\, \sigma
_x^{-1}\]
Next, choose
\[
  \delta _2=\begin{cases}
    r^{-1/7}\ts_x^{1/7} (\ln r+\smE )^{1/7} & \text{when \(\sigma
      _x\geq 3  r^{-1/6} (\ln r+\smE )^{1/6} \)};\\
    \sigma_x/3 &\text{otherwise}.
  \end{cases}
\]
One gets:
\[\begin{split}
f_1&\leq  \begin{cases}
  \zeta _3'\,  (r\ts_x^{-1})^{5/7}(\ln r+\smE )^{2/7}\quad & \text{when \(\sigma
  _x\geq 3  r^{-1/6} (\ln r+\smE )^{1/6} \)};\\
\zeta '_4 r^{1/2}\ts_x^{-2}(\ln r+\smE )^{1/2}
&\text{otherwise}
\end{cases}
\end{split}
\]
Note that such
choice of \(\delta _2\) satisfies the constraint that \(\sigma (x)/2
>\delta_2\geq\delta_0 \) for all sufficiently large \(r\) and \(x\in X^{'a}_{8\delta _0}\). 

Meanwhile, \(f_2\) is of the form
\[
\begin{split}
f_2& =\zeta '_1r^{1/3}\ts^{-1}_x+\zeta '_2 \ts^{-4}_x\rho _1^2+\zeta '_3\, \rho
_1^{-2}r^{-1}\delta _2^{-6} (\ln r+\smE )+\zeta _4\, \delta _0^{-1}\, \sigma _x^{-1}\\
&\leq \zeta _5\,  \ts^{-1}_x\big( r^{1/3}+\zeta _5'r^{-1/2}(\ln r+\smE )^{1/2}\ts_x^{-4}\big), 
\end{split}
\]
with \(\rho _1^2=r^{-1/2}\delta _2^{-3}(\ln r+\smE )^{1/2}\ts^{2}_x\), and \(\delta
_2=\sigma _x/3\). 

For a fixed \(x\in X^{'a}_{8\delta _0}\), set \[\begin{split}
    \varepsilon& =\varepsilon _x\\
    & := \begin{cases}
    \big(r^{-4/3}\ts^{-1}_x (\ln  r+\smE)\big)^{2/7} & \text{when \(\sigma
      _x\geq 3  r^{-1/6} (\ln r+\smE )^{1/6} \)};\\
   r^{-1/6}\, \ts_x\, (\ln r+\smE )^{-1/2}+\zeta '_5\sO_r^{-3}  &\text{otherwise},
  \end{cases}
\end{split}
\]
we then have:
\begin{equation}\label{ineq:q}
  \begin{split}
  2  q^{(\varepsilon _x)}(x) & \leq \begin{cases}
      \zeta' (r\ts^{-1}_x)^{5/7}\, (\ln r+\smE )^{2/7} & \text{when \(\sigma
      _x\geq 3  r^{-1/6} \, (\ln r+\smE )^{1/6} \)};\\
  \zeta '' r^{1/2}\ts_x^{-2}(\ln r+\smE )^{1/2}   &\text{otherwise}
   \end{cases}\\
   &   =: K_0\, \quad \text{over \(X^{'a}_{8\delta }\). }
\end{split}
\end{equation}
Noting that for all sufficiently large \(r\), \(\ts^{-2}<K_0\) on 
\(X^{'a}_\delta \), we have  
\[
\begin{split}
\ss +\zeta_2r|\beta |^2\leq (2^{-3/2}+\varepsilon _0)\, r\, (-u)+K_0 
\quad \text{over  \(X^{'a}_{8\delta _0}\). }
\end{split}
\]
The first inequality in (\ref{bdd:ss}) now follows with \(\sO\)
renamed as \(8\sO\).

To verify the second inequality in (\ref{bdd:ss}), fix \(x\in
Z^{'a}_{\delta }\). Notice that  (\ref{DE:q}) also holds with the
constant \(\delta _0\) therein replaced by \(4\delta \), with the
definition of \(q\) in (\ref{def:q}) correspondingly modified. For the
rest of this proof, let \(q=:q_\delta \), \(\eta_1\), \(\eta_2\), \(\xi\) denote
the correspondingly modified versions, with \(\varepsilon \) set to be
1. In particular, \(\xi\) is now supported on \(Z^{'a}_{8\delta }\), and
\(\eta_1\), \(\eta_2\) are both supported on \(X^{'a}_{2\delta }\).  
Then in this case \(q\, (x)\) is still bounded by (\ref{eq:q-est}), and the first
term on its right hand side is bounded as before. Namely, it is bounded
by a positive multiple of \(r^{1/2}( \ln r+\smE )^{1/2}\). The 
second to the fourth terms on the right hand side are bounded differently as
follows. 

To bound the second integral in  (\ref{eq:q-est}), divide \(B(x, 2\rho
_0)\cap Z^{'a}_{8\delta }\) into 
three regions: the first region  \(\scrR_1:=B(x, \sigma _x/2)\); the
second region \(\scrR_2:=\{y\, |\, \sigma _y\leq \frac{3}{4}\sigma _x,
y\in B(x, 2\rho
_0)-\scrR_1\); the third region \(\scrR_3:=B(x, 2\rho
_0)\cap Z^{'a}_{8\delta }-\scrR_1-\scrR_2\). 
 Make use the following facts: For \(y\in \scrR_1\), \(\ts_y\geq \ts_x/2\)
 and so \(\xi\leq \zeta  r\ts^{-1}_x +\zeta ' \delta
 ^{-2}\ts^{-2}_x\); for \(y\in \scrR_2\), \(|\sigma _y-\sigma
 _x|\geq\sigma _x/4\);   for \(y\in \scrR_3\), \(|\sigma _y-\sigma
 _x|\geq \sigma _y/2\). One has:  
\[\begin{split}
\int_{B(x,2\rho_0)} \xi \dist (x,
\cdot)^{-2}& \leq   \zeta  _2r\delta +\zeta '_2 \, \delta ^{-2}\ln \, (\delta
/\sigma _x),
\end{split}
\]
where \(\zeta _2\), \(\zeta _2'\) are independent of \(r\), \(\delta
\), and \(x\in Z^{'a}_\delta \). 

To bound the remaining two integrals in  (\ref{eq:q-est}), note that
the distance from \(x\) to the support of either \(\eta_1\) or \(\eta_2\) is no
less than \(\delta \geq \sigma _x\).
\[
\begin{split}
\int_{B(x,2\rho_0)} \eta_1\dist (x,
\cdot)^{-2}& \leq  \zeta  \int_{2\delta }^{\sigma _x+2\rho _0}\int _{-2\rho _0}^{2\rho _0}\frac{\sigma _y^{-2} dz}{(\sigma _y-\sigma
  _x)^2+z^2}\, \sigma _y^2 \, d\sigma _y\\
& \leq \zeta '\int_{2\delta }^{\sigma _x+2\rho _0}\frac{d\sigma _y}{\sigma _y-\sigma _x}\\
& \leq \zeta  _3\, \ln \, (\delta ^{-1})
\quad \text{when \(x\in Z^{'a}_{\delta}\)}. 
\end{split}
\]
To bound the last integral in (\ref{eq:q-est}), choose \(\delta _2>0\), and devide 
\(B(x,2\rho_0)\cap X^{'a}_{2\delta }\) into two regions: \(\sigma \geq\delta _2\) on the
first region, and \(\sigma \leq\delta _2\) on the second. (The second
region is empty when \(\delta _2\leq 2\delta \).) 
Over the first region, use Lemma \ref{T:lem3.1} together with the
observation that \(\dist (y, x)\geq  \sigma _y/2\geq\delta _2/2\) for any \(y\) in
this region. 
Over the second region, use the fact that 
by Proposition 
 \ref{T:lem3.2}, 
\(
\eta _2\leq  \zeta   r\ts^{-1} 
\) over \(X^{'a}_{2\delta }\). 
Then, setting \(\rho _2=r^{-1/7}\delta^{1/7}(\ln r+\smE )^{1/7}\), one has: 
\[
\begin{split}
& \int_{B(x,2\rho_0)} \eta_2\dist (x,
\cdot)^{-2} \\
& \qquad \leq  
\zeta '\, \delta _2^{-2} \delta _2^{-3}(\ln r+\smE)+\zeta
_1\int_{2\delta}^{\delta _2}\int _{-2\rho _0}^{2\rho _0}\frac{r\sigma _y
  \, dz \, d\sigma _y}{(\sigma _y-\sigma
  _x)^2+z^2}\,\\
& \qquad \leq \zeta '\, \delta _2^{-5}(\ln r+\smE)+ \zeta '_1\delta^{-1}\int_{2\delta }^{\delta _2}r\sigma _y
  \, d\sigma _y\\
& \qquad \leq \zeta '\, \delta _2^{-5}(\ln r+\smE)+\zeta r \, \delta^{-1}\delta
_2^2
\\
& \qquad \leq \zeta  _4\, r^{5/7}\delta^{2/7}(\ln r+\smE )^{2/7}\leq \zeta_4' r\delta
\quad \text{when \(x\in Z^{'a}_{\delta}\)}. 
\end{split}
\]
Gathering all the termwise bounds obtained, 
we have: 
\[
q_\delta \leq \zeta  \, r \delta +\zeta ' \delta^{-2}\ln\,  (\delta /\sigma ) 
\quad \text{over \(Z^{'a}_{\delta}\)},
\]
where \(\zeta \), \(\zeta '\) are positive constants independent of
\(r\) and \(\delta \). Hence, by Proposition \ref{T:lem3.2}
\[
\ss =q_\delta +c_1r\, u-\zeta _2 r|\beta |^2\leq \zeta _1r \delta +\zeta '_1 \delta^{-2}\ln\,  (\delta /\sigma)
\quad \text{over \(Z^{'a}_{\delta}\), }
\]
as asserted by 
the second inequality in (\ref{bdd:ss}). 

Having established Item (b) of the assertions of the proposition, item
(c) follows directly. Given \(x\in X^{'a}_{\delta _0}\), set the parameter
\(\delta \) in  (\ref{bdd:ss}) (ii) as \(\delta =\sigma _x\), and
observe that over \(X^{'a}_{\delta _0}\), \(r^{2/3}\leq \sO^{-1}r\ts_x\),
one arrives at the first inequality in Item (c). The second inequality
asserted in Item (c) follows from a combination of the first
inequality and  (\ref{bdd:ss}) (ii). 
\epf

We shall repeatedly improve the estimate for \(|F_A^-|\). The first of 
such improvements replaces the occurrence of \(K_1\) in Item c) of
Proposition \ref{T:prop3.2} with a constant multiple of \(\ts^{-2}\).
Recall the notation \(\td{r}=r\ts\). 
\begin{prop}\label{prop:curv-varpi}
Let \(u:=|\psi|^2-|\nu|\), and let \(\varepsilon _0\), \(K_0\),
\(K_1\) be as
in Proposition \ref{T:prop3.2}. 
There exist \(r\)-independent positive constants \(\zeta_O\),
\(r_0>8\), \(\zeta \), \(\zeta'\), \(\zeta_1\) 
that only depend on the metric, \(\nu \), \(\varsigma_w\), \(\zeta
_\grp\) with the following significance:
 Let  \(\delta _0':=\zeta_O \, r^{-1/3}(\ln r)^{2/3}\),  \(\varepsilon _1:=K_1 \,
 \td{r}^{-1}
\), and \(\varepsilon':= \varepsilon
_0+\zeta   \varepsilon_1\). Then for any \(r\geq r_0\) one has: 
\begin{equation}\label{eq:curv-varpi}
\begin{split}
|F_A^-|& \leq (2^{-3/2}+\varepsilon ')\, 
r \, (-u)_+ +\zeta'
\,  \ts^{-2} 
\quad
\text{over \(X^{'a}_{\delta'_0}\), \(0\leq a\leq 7\)}.
\end{split}
\end{equation}
The constants 
\(\sO\), \(r_0\), \(\zeta \)
and \(\zeta '\) above 
depend  on the \(\Spin^c\) structure and the relative homotopy
class of \((A, \Psi)\). 
\end{prop}

Note that \(\varepsilon _1+\varepsilon '\leq \zeta _0'\) over
\(X^{'a}_{\delta _0'}\) for an \(r\)-independent constant \(\zeta _0'\).

The proof of the preceding proposition makes use a comparison
function named \(v_2\), which we now describe. 

 Let \(X^*_{\delta}\subset  X''_{\delta}\) 
denote the subspace consisting of points where
\(|\psi|^2/|\nu|\geq 1/2\), and let \(\sW:=|\nu |/2-|\psi |^2\). 
\begin{lemma}\label{lem:v_2}
Let \(u:=|\psi|^2-|\nu|\). 
There exist positive constants \(\sO\geq 8\), \(\zeta_i\), \(i=1,
\ldots, 5\), \(\zeta'_3, \zeta _4'\), that are independent of
\(r\), \(\delta \), and \((A, \Psi )\), with the following
significance: Suppose \(r>1\), \(\delta >0\) are such that
\(\delta\geq\textsc{o}r^{-1/3}\), then 
given a 
positive constant \(\epsilon<1\), there is a function \(v_2\) on \(X''_{\delta}\) which
satisfies: 
\begin{itemize}
\item Over \(X''_\delta \), 
\[\begin{split}
v_2& \geq \zeta_1\ts^{-\epsilon}(-u)_+\geq
  \zeta_1\ts^{-\epsilon}\, (|\nu |/2+\sW).
\end{split}
\] 
In particular,
  \(v_2\geq \zeta'_1\ts^{1-\epsilon}\) over \(X_\delta ''-X^*_\delta
  \). 
\item \(v_2\geq \zeta _2r^{-1}\ts^{-2-\epsilon}\). 
\item Over \(X''_\delta \), 
\[\begin{split}
(d^*d+r|\psi|^2/2)\, v_2 & \geq
  \ts^{-\epsilon}\big( \zeta _3 \epsilon r
  |\psi|^2\,(-u)_+-\zeta _3'(1-\epsilon )\, r \ts^{-2}\sW\big)\\
&\geq  \ts^{-\epsilon}\big( \zeta _3 \epsilon r
  |\nu |\,(-u)_+/2-\zeta _4'r \ts^{-2}\sW\big).
  \end{split}
\]
In particular, \((d^*d+r|\psi|^2/2)\, v_2\geq\zeta _3 \epsilon 
  \ts^{-\epsilon}r
  |\psi|^2\,(-u)_+\) over  \(X^*_\delta \). 
\item \(v_2\leq\zeta_4r^{\epsilon}\ts^{2\epsilon }(-u+\zeta r^{-1}\ts^{-2})\) on \(X''_\delta\).
\item \(v_2\leq \zeta_5 \ts^{1-\epsilon}\) on \(X''_\delta\).
\end{itemize}
\end{lemma}
\pf 
From (\ref{u-ineq2a})  we have
\begin{equation}\label{u-ineq2}
(2^{-1}d^*d+\frac{r}{4}|\psi|^2)(-u)\geq-\zeta _0\, \ts^{-1}\quad
\text{over \(X''\)}.
\end{equation}
As \(|\psi|^2\geq|\nu|/2\) on \(X^*_{\delta}\), it follows from
(\ref{ts-loBdd}) that 
\begin{equation}\label{ts-loBdd-p}
\begin{split}
& d^*d (r^{-1}\ts^{-k})+\frac{r|\psi|^2}{2}(r^{-1}\ts^{-k})\\
& \quad =
\big(d^*d +\frac{r}{4}|\nu|-\frac{r\sW}{2} \big)(r^{-1}\ts^{-k}) 
\geq\begin{cases}
  \zeta''\ts^{-k+1}-\frac{r\sW}{2} \big)(r^{-1}\ts^{-k}) &  \text{on \(X''_\delta\).}\\
  \zeta''\ts^{-k+1}&  \text{on \(X^*_\delta\).}
\end{cases}
\\
\end{split}
\end{equation}
Adding a suitable multiple of the \(k=2\)'s version of the preceding inequality to
(\ref{u-ineq2}) and combining with  Proposition \ref{T:lem3.2}, one
can find  a positive constant
\(\zeta_2\) such that 
\[
v_1:=-u+\zeta _2r^{-1}\ts^{-2}
\]
satisfies:
\BTitem\label{ineq:v1}
\item \((d^*d+r |\psi|^2/2)\, v_1\geq -\zeta r\ts^{-2} \sW\) on
  \(X''_\delta\). In particular, \((d^*d+r |\psi|^2/2)\, v_1\geq 0\) on \(X_\delta^*\);
\item \(v_1\geq \zeta r^{-1}\ts^{-2}\) on \(X''_\delta\),
\item \(v_1\geq (-u)_+\geq |\nu |/2+\sW\) on \(X''_\delta\),
\item \(v_1<\zeta\ts\) on \(X''_\delta\).
\ETitem
Now take \(v_2:=v_1^{1-\epsilon}\).
\epf



\noindent {\it Proof of Proposition \ref{prop:curv-varpi}.} Let \(\sO\) be the larger of that in Proposition \ref{T:prop3.2}
and that in the preceding lemma. Set \(\delta _0=\sO r^{-1/3}\).  Fix \(\delta \geq \delta _0\). Let \(X^*_{\delta}\subset  X_{\delta}''\)
denote the subspace consisting of points where
\(|\psi|^2/|\nu|\geq 1/2\) as before. Note that over \(X''_\delta -X_\delta ^*\), \(|\nu |\leq -2u\) and \(
\ts^{-1}\leq \sO^{-1} r^{1/3}\). Let \( \varepsilon_{0, \delta }\),
\(K_{0, \delta } \), \(K_{1, \delta } \) be defined by replacing every
occurrence of \(\ts\) and \(\sigma \) in the respective formulas
(\ref{def:e_0K_0}), (\ref{def:K_1}) defining \(\varepsilon_0\),
\(K_0\), \(K_1\), by the number \(\delta \); and let \(K'_{1, \delta } :=\min\, 
(K_{0,\delta }, r\ts)\). 
Let \(q_0=q_0^{(\varepsilon)}\)  be as in (\ref{def:q0}) with 
 \(\varepsilon \) set to be \(\varepsilon_{0, \delta }\). Then by the
 proof of Proposition
 \ref{T:prop3.2} 
 and the fact that \(\delta ^{-2}< K'_{1,\delta }\) over \(X_\delta ''\)
 for \(\delta \geq \delta _0\) and all sufficiently large \(r\), there is a positive constant
 \(\zeta '\) that depends only on the metric, \(\nu \),
 \(\varsigma_w\), and \(\zeta _\grp\), 
 such that 
\begin{equation}\label{bdd:q_0b}
\begin{split}
q_0& \leq \zeta '
K'_{1, \delta } \quad \text{ over \(X^{'a}_{\delta }\)}, \quad 0\leq a\leq \frac{17}{2}
\end{split}
\end{equation}
for all sufficiently large \(r\) and and \(\delta \geq \delta
_0\).   
Combined with (\ref{ineq:q_0}) (replacing \(\delta _0\) therein by \(\delta _0/2\)), this gives: 
\begin{equation}\label{DE:q_0-1}
\begin{split}
  \big(\frac{ d^*d}{2}+\frac{r|\psi|^2}{4}\big) q_0  
 & \leq  (\zeta_4 +\zeta _4' \varepsilon^{-1}_{0,\delta }\,\sO_r^{-3} \ts^{-2}) \, r\, (-u)
 +r \ts^{-1}( \zeta _1+\zeta '_1\varepsilon ^{-1}_{0,\delta }\sO_r^{-6})\\
 & \leq  \zeta''_4 \, r\ts^{-2} (-u)
+\zeta ''_1r \ts^{-1}\quad \text{on \(X^{'a}_{\delta }\).}
\end{split}
\end{equation}
Given \(i\in \grY_m\), let \(q_{0,i}\) denote the \(Y_i\)-end
limit of \(q_0\). This is a function on \(Y_i\), and is bounded via
the 3-dimensional Seiberg-Witten equation \(\grF_{\mu _{i,r}}(B_i,
\Phi _i)=0\) and Lemma \ref{lem:3d-ptws} by 
\begin{equation}\label{bdd:q_0i}
q_{0,i}\leq \zeta _i\,(1- \chi (\sigma _i/\delta _{0,i}))\, \ts_i^{-2}+\zeta _i'\chi (\sigma _i/\delta _{0,i})
\, r^{2/3}
\end{equation}
for certain \(r\)-independent constants 
\(\zeta _i\), \(\zeta _i'\). 
In the above, \((B_i, \Phi _i)\) denotes the \(Y_i\)-end limit of
\((A, \Psi )\), and \(\sigma _i\), \(\ts_i\), \(\delta _{0,i}\) are
respectively what were denoted by \(\sigma \), \(\ts\), \(\delta _0\)
in Lemma \ref{lem:3d-ptws}. According to observations in Section \ref{sec:end-limit},   \(q_0\Big|_{Y_{i:L}}\)
approaches \(q_{0,i}\) as \(L\to \infty\) in \(C^k(Y_i)\) topology.

Combining (\ref{DE:q_0-1}) with (\ref{ts-loBdd}), 
and making use of (\ref{bdd:q_0i}) and
(\ref{bdd:q_0b}), one may find  constants \(\zeta _2, \zeta_2'\) independent
of \(\delta \), and \(r\), such that the function 
\[
q':=q_0-\zeta _2\ts^{-2} 
\]
satisfies: 
\begin{equation}\label{DE:q'}
\begin{split}
& \big(\frac{ d^*d}{2}+\frac{r|\psi|^2}{4}\big) \, q' 
 \leq   \zeta _2' \, r\ts^{-1} =:\xi'\,
\quad \text{on \(X^{'a}_\delta \)};\\  
& \quad q'\leq 0 \quad \text{over \(\hat{Y}_{i, L}\), \(\forall i\in
  \grY_m\); }\\
& \quad q'  \leq \zeta 
K'_{1,\delta } \quad \text{ over \(X^{'a}_{\delta }\)}, \quad 0\leq
a\leq \frac{17}{2}, \quad \delta \geq\delta _0. 
\end{split}
\end{equation}
(The number \(L\) above may depend on \(r\) and \((A, \Psi )\), but is
independent of \(\delta \). This dependence does not affect our
subsequent discussion, and unless otherwise specified, all constants
below are independent of \(L\).) Now let 
\[
q'_\delta :=\gamma _{a, \delta } \, q', \quad 0\leq a\leq \frac{15}{2}
\]
where \(\gamma _{a,\delta }\) is the cutoff function introduced  in the beginning of this
section. The function \(q'_\delta \) satisfies: 
\begin{equation}\label{DE:q'_d}
 \big(\frac{ d^*d}{2}+\frac{r|\psi|^2}{4}\big) \, q'_\delta  
 \leq  \gamma _{a,\delta }\, \xi' +
 \frac{ d^*d \, q'_\delta }{2}-\gamma _{a,\delta }
 \frac{ d^*dq'}{2} =: \xi'_\delta . 
\end{equation}
\((q'_\delta )_+\) is supported within the compact space \(U:=\scrX^{a+1}_{\delta
}-\bigcup_{i\in \grY_m} \hat{Y}_{i, L}\subset X''_\delta \). Let
\[
\scrX^a_{\delta, l}\subset \scrX^a_{\delta}-\bigcup_{i\in \grY_m} \hat{Y}_{i, l}
\]
denote the manifold with boundary obtained by ``rounding the corners''
of \(\scrX^a_{\delta}-\bigcup_{i\in \grY_m} \hat{Y}_{i, l}\). 
More precisely, it satisfies:  \(\partial\scrX^a_{\delta, L}\) is
smooth, and  \(
\scrX^a_{\delta}-\bigcup_{i\in \grY_m}
\hat{Y}_{i,l}-\scrX^a_{\delta, l}\subset \hat{Y}_{i, [l-\epsilon ,
  l]},\) 
where \(\epsilon \ll 2^{-8}\). 
Let \(V:=\scrX^{a+1}_{\delta, L+1}\supset U\). 
Then \(q'_\delta \leq 0\)
on \(\partial V\). 

Let \(q_1\) be a solution to the following
Dirichlet boundary value problem:
\[
\big(\frac{ d^*d}{2}+\frac{r|\nu |}{4}\big) \, q_1
=\xi'_\delta \, \, \, \text{ over \(V\)};\qquad 
q_1|_{\partial V}=0.
\]
The (Dirichlet) Green's function for \(\frac{ d^*d}{2}+\frac{r|\nu |}{4}\), denoted \(G_r\) below, satisfies  
\[
|G_r(x, \cdot)|+\dist(x,\cdot)\, |d G_r(x, \cdot)|\leq \zeta_g'\dist(x,\cdot)^{-2}e^{-\zeta_g
  (r\delta )^{1/2}\dist (x,\cdot)}
\] 
for certain constants \(\zeta _g'\), \(\zeta _g\).  Thus,
{\small \[
\begin{split}
 &  |q_1(x)|\\
  & \quad \leq \zeta _1\int_U \dist(x,\cdot  )^{-2}e^{-\zeta_g
  (r\delta )^{1/2}\dist (x,\cdot )} |\xi'| \\
&\quad  \,\, +\zeta _2\int_{X^{'a+1}_\delta -X^{'a}_\delta }K'_{1,\delta }\dist(x,\cdot  )^{-3}e^{-\zeta_g
  (r\delta )^{1/2}\dist (x,\cdot )} \\
& \quad \, \, +\int_{(Z^{'a}_{2\delta }-Z^{'a}_\delta)\cap U}K'_{1,\delta }\big( \zeta _3\,
\delta^{-2}\dist(x,\cdot  )^{-2}+\zeta '_3\, \delta^{-1}\dist(x,\cdot
)^{-3}\big) e^{-\zeta_g
  (r\delta )^{1/2}\dist (x,\cdot )}\\
& \quad \leq \zeta _1' \delta ^{-1}\ts^{-1}(x)+ \zeta _2'\, K'_{1,\delta }(x) \, (r\delta
)^{-1/2} e^{-\zeta _g(r\delta )^{1/2}\dist (x,\, X^{' a+1}_\delta
  -X^{'a}_\delta  )}\\
&\quad \qquad \qquad +\zeta _3'' K'_{1,\delta }(x)\,  (r\delta ^{3})^{-1/2}
e^{-\zeta _g(r\delta )^{1/2}\dist (x,\, Z^{'a}_{2\delta }-Z^{'a}_\delta) }. \\
& \quad \leq \zeta _1' \delta ^{-1}\ts^{-1}(x)+ \zeta _4  \, K'_{1, \delta }(x) \, (r\delta ^{3})^{-1/2},
\end{split}
\]
}
where the terms in the second and third lines above are obtained via
integration by parts. It follows that there exist constants \(r_0>8\), 
\(\zeta _O\), \(\zeta ''_1\) and , \(\zeta _4'\), which are independent of \(r\), \(\delta \), and
\((A, \Psi )\), such that for all \(r\geq r_0\), and \(\delta \geq
\delta _0':=\zeta _O\, r^{-1/3}\, (\ln r)^{2/3}>\delta _0\), 
\[
\begin{split}
|q_1| & \leq \zeta '_4  \, K'_{1,\delta }(\ln r)^{-1}\quad \text{over \(V\);} \\
|q_1| & \leq \zeta ''_1 \delta ^{-1}\ts^{-1} \quad \text{over \(U':=X^{'a-\frac{1}{2}}_{4\delta }\cap V\).} 
\end{split}
\]

Assume  \(\delta
>\delta _0'\) from now on. Then 
the preceding estimates for \(|q_1|\) and \(q'\),
  the function \(q_2:= q'_\delta -q_1\) satisfies: 
\[\begin{split}
 \big(\frac{ d^*d}{2}+\frac{r|\psi|^2}{4}\big) \, q_2& =-\frac{ru}{4}
q_1\leq \frac{\zeta '_4}{4}  \, K'_{1,\delta }(\ln r)^{-1} r\, |u|; \\
 q_2|_{\partial V} & \leq 0;\\
q_2 & \leq \zeta _5 \, K'_{1,\delta }. 
\end{split}
\]
Choose \(\epsilon =(\ln (r\delta ))^{-1}\), and note that \(\delta
^{-1+\epsilon }K_{1,\delta }\geq\zeta \ts^{-1+\epsilon }K'_{1,\delta
}\).  Appeal to the first and the third bullets in Lemma \ref{lem:v_2} 
to find positive constants \(\zeta _1\), \(\zeta _2\), \(\zeta _3\)
(independent of \(r\), \(\delta \)), such that with 
\[
q_3:=q_2
-\zeta  _1\delta ^{-1+\epsilon } K_{1, \delta } \, v_2, 
\]
one has: 
\[\begin{split}
\big(2^{-1}d^*d+\frac{r}{4}|\psi|^2\big)\, q_3& \leq \zeta   _2\, 
r\ts^{-2} \delta
^{-1} K_{1, \delta } \sW +\frac{\zeta '_4}{4}  \, K_{1,\delta }\, (\ln r)^{-1} r\,
u_+\\
  q_3|_{\partial V}& \leq 0; \\ 
q_3 & \leq  -\zeta   _3 \, \delta ^{-1} K_{1, \delta } \sW. \\
\end{split}
\] 
It follows from the preceding three inequalities that \(q_3\leq 0\)
over \(V\). To see this, note that by the third inequality, \(q_3\leq
0\) over \(V-X^*_\delta \), where \(\sW\geq 0\). Together with the
second inequality, this implies that \(q_3|_{\partial (V\cap
  X^*_\delta )}\leq 0\). Meanwhile, the RHS of first inequality is
nonpositive over \(V\cap X^*_\delta \). 
By the maximum principle, \(q_3\leq 0\) over \(V\cap X^*_\delta \) as well. Now, applying the fourth
bullet in Lemma \ref{lem:v_2}, one has:
\[\begin{split}
q_0& \leq \zeta  _1' \delta ^{-1} \td{r}^\epsilon
K_{1, \delta }\, (-u+\zeta  r^{-1}\ts^{-2}) +\zeta _2' \delta ^{-1}\ts^{-1}
\\
& \leq \zeta  _1''\, \delta ^{-1}K_{1, \delta }\, (-u+\zeta
r^{-1}\ts^{-2}) +\zeta _2' \delta ^{-1}\ts^{-1}
\end{split}
\]
over \(X^{'a}_{4\delta }\), \(0\leq a\leq 7\), for all sufficiently large $r$. All the
constants \(\zeta _*\) above are independent of \(r\), \(\delta \),
and \((A, \Psi )\). Given \(x\in X^{'a}_{8\delta _0}\), set \(\delta
=\ts(x)/2\) in the preceding expression. It then gives 
\[\begin{split}
q_0
& \leq \zeta  _7
\, \varepsilon _{1}r\, (-u)+\zeta _8\, \ts^{-2}. 
\end{split}
\]
This leads directly to the conclusion of the proposition. \epf

\subsection{Estimating \(|\nabla_A\alpha|\) and \(|\nabla_A\beta|\)}

The next proposition is an analog of \cite{Ts}'s Proposition 3.3 and
\cite{T}'s Proposition I.2.8, and the proof is an adaptation of the latter. Let 
 \(\underline{\alpha}:=|\nu|^{-1/2}\alpha\) and
 \(\varpi=|\nu|-|\alpha|^2\). 

\begin{prop}\label{lem:est-1st-der}
There exist positive constants \(r_1\), \(\zeta _O\), \(\zeta', \zeta''\), that are independent of
\(r\) and \((A, \Psi )\), with the following significance: 
Let \(\delta _0'=\zeta _O r^{-1/3}(\ln r)^{2/3}\). For any \(r>r_1\), one has: 
\[\begin{split}
|\nabla_A\underline{\alpha}|^2+r\ts^2|\nabla_A\beta|^2&\leq \zeta'r
\varpi+\zeta''\ts^{-2} \quad \text{over \(X^{'a}_{\delta'_0}\)}, \quad
0\leq a\leq \frac{13}{2}.
\end{split}\]
\end{prop}

\pf Let \(\delta _0'\) be as in Proposition \ref{prop:curv-varpi}. Argue as in p.191  of
\cite{Ts} (which is itself a modification of the arguments in p.20 of 
\cite{T}),\footnote{Caveat: Equations (3.50)--(3.53) in \cite{Ts}
  contain typos and errors.} replacing the use of
Propositions 3.1, 3.2 therein by their counterparts in our setting, 
Propositions \ref{T:prop3.1-} and  \ref{prop:curv-varpi}. 
This leads to the following
variant of (3.53) in \cite{Ts} (which in turn is based on \cite{T}'s
(I.2.38) and (I.2.40)):
\begin{equation}\label{eq:DE-alpha'}\begin{split}
2^{-1}d^*d &|\nabla_A\underline{\alpha}|^2+2^{-1}|\nabla_A\nabla_A\underline{\alpha}|^2+8^{-1}r|\nu||\nabla_A\underline{\alpha}|^2\\
&\leq \zeta \, (\tilde{\sigma}^{-2}+r\varpi_+)
|\nabla_A\underline{\alpha}|^2+
\zeta_2 (r\ts^{-1}|\beta|^2+ r^{-1}\ts^{-6}) \, |\nabla_A\beta|^2\\
& \qquad \quad +
\zeta_3 r^{-1}\tilde{\sigma}^{-4}|\nabla_A\nabla_A\beta|^2
+\zeta_5r^{-1}\tilde{\sigma}^{-7}+\zeta_6r\ts^{-3}\varpi^2;\\
2^{-1}d^*d &
|\nabla_A\beta|^2+|\nabla_A\nabla_A\beta|^2+8^{-1}r|\nu|\, |\nabla_A\beta|^2\\
& \leq \zeta  (\tilde{\sigma}^{-2}+r\varpi_+)
|\nabla_A\beta|^2+
\zeta_2'(r|\beta|^2+ r^{-1}\ts^{-5}) \, |\nabla_A\alpha|^2\\
& \qquad \quad +\zeta_3'r^{-1}\tilde{\sigma}^{-3}|\nabla_A\nabla_A\alpha|^2+\zeta_4'r\tilde{\sigma}^{-1}|\beta|^2+\zeta_5'r^{-1}\ts^{-1}
\end{split}
\end{equation}
on \(X^{'a}_{\delta_0'}\), \(0\leq a\leq 7\). Re-introduce the notation \(\ts_x=\ts(x)\) and
\(\td{r}_x=\td{r}(x)=r\ts_x\), and note that 
\begin{equation}\label{est-ts}
1/2\leq \ts /\tilde{\sigma}_x\leq 3/2\quad
\text{over \(B\, (x, \ts_x/2) \). }
\end{equation}
Fix \(x\in X^{'a}_{2\delta _0'}\), \(0\leq a\leq\frac{13}{2}\). Use (\ref{eq:DE-alpha'}), (\ref{3.20}),
(\ref{3.21}), Proposition \ref{T:prop3.1-} and (\ref{est-ts}) to find constants \(\zeta \), \(\zeta
'\), \(\zeta ''\), \(\zeta _1\) that are independent of \(r\), \(x\),
and \((A, \Psi )\),  such that 
\begin{equation}\label{def:y}
\begin{split}
& \quad d^*d y_x+4^{-1}r|\nu|y_x\leq 0 \quad \text{over \(B\, (x,
  \ts_x/2) \),  with }\\ 
&y_x:=\max \, \big(|\nabla_A\underline{\alpha}|^2+\zeta
r\ts^2_x|\nabla_A\beta|^2-\zeta' r\varpi+\zeta'' r^2\ts^3_x|\beta|^2-\zeta_1\tilde{\sigma}^{-2}, 0\big).
\end{split}
\end{equation}
Note that by Propositions \ref{T:lem3.2} and \ref{T:prop3.1-} and (\ref{est-ts}), 
\begin{equation}\label{bdd:y0}
y_x\leq \zeta _2\, \ts^{-1}|\nabla_A\alpha|^2+\zeta_2'
r\ts^2|\nabla_A\beta|^2+\zeta_3' \ts^{-2}+\zeta _1'r \varpi . 
\end{equation}

Multiply both sides of the differential inequality (\ref{def:y}) by
the function 
\(\chi  (4\dist (x, \cdot)/\ts_x)\, G_r\), 
where \(G_{r}\) denotes the 
Green's function of \(d^*d+4^{-1}r|\nu|\) on \(B\, (x, \ts_x/2)\) with Dirichlet
boundary condition, then integrate over \(B\, (x, \ts_x/2)\). This
Green's function satisfies a bound of the form: 
\[
\Big|G_r(x, x')\Big|+\dist (x, x')\, \Big|\, dG_r(x,
x')\Big|\leq  \zeta_0\dist(x,x')^{-2}\, e^{-2\zeta'_0
  \td{r}_x^{1/2}\dist (x,x')}. 
\]
It follows from integration by parts that 
\begin{equation}\label{eq:y-bdd0}
\begin{split}
y_x(x)& \leq\zeta _3 \ts_x^{-4}e^{-\zeta'_0 \td{r}_x^{1/2}\ts_x/2}\int_{A_x}y_x, \\
\quad 
\end{split}
\end{equation}
where \(A_x= B(x, \ts_x/2)-B(x, \ts_x/4)\). 

By (\ref{bdd:y0}), (\ref{est-ts}), Propositions \ref{prop:SW-L2-bdd},
\ref{T:prop3.1-}, 
and Lemma  \ref{T:lem3.1}, 
\begin{equation}\label{eq:y-bdd1}
\int_{A_x}y_x\leq\zeta  '\, \ts_x^{-1} \, (\ln r+\smE)
+\zeta r\ts^2_x\int_{A_x}|\nabla_A\beta|^2. 
\end{equation}

To estimate the last term above, 
multiply both sides of (\ref{3.20}) by \(\chi\big(2\dist (x, \cdot)/\ts_x\big)\), then
integrate over \(B(x,\ts_x)\). 
This gives 
\[\begin{split}
\int_{A_x}|\nabla_A\beta|^2& \leq  \zeta  _4\ts_x^{-2} \int_{B\,
  (x, \ts_x)}|\beta|^2 +\zeta  _4'r^{-1}\ts_x^{-3}\int_{B\,
  (x, \ts_x)}|\nabla_A\alpha |^2+ \zeta  _4''r^{-1}\ts_x^3\\
& \leq \zeta _5 \, r^{-1}\ts_x^{-3}(\ln r+\smE). 
\end{split}
\]
In the above,  (\ref{est-ts}), Propositions  \ref{T:prop3.1-},
\ref{prop:SW-L2-bdd},
and Lemma  \ref{T:lem3.1} are used to derive the second line from the
first line. Inserting the preceding inequality into
(\ref{eq:y-bdd1}) and combining with (\ref{eq:y-bdd0}),
we have 
\[
y_x(x) \leq\zeta _6\, \ts_x^{-5}e^{-\zeta_6'(r\ts_x^3)^{1/2}} (\ln
r+\smE ), 
\]
where \(\zeta _6\), \(\zeta _6'\) are independent of \(r\), \(x\), and \((A, \Psi )\).  
This implies the existence of constants \(r_1\), \(\zeta _7\), \(\zeta _7'\),
also independent of \(r\), \(x\), and \((A, \Psi )\), 
such that \(y_x(x)\leq \zeta _7' \ts^{-2}_x\) when \(r\geq r_1\) and \(\sigma_x\leq
\zeta _7 \, r^{-1/3}(\ln r)^{2/3}\).  Rename \(\zeta _O\)
to be the larger of this \(\zeta _7\) and twice of the version from Proposition
\ref{prop:curv-varpi}, and re-set \(\delta _0'=\zeta _O\, r^{-1/3}(\ln
r)^{2/3}\). \epf


\begin{rem}
It is basically equivalent to estimate either of
 \(|\nabla_A\underline{\alpha}|\) and \(|\nabla_A\alpha|\), as they
 are related  by 
\[
\Big| |\nabla_A\alpha|^2-|\nu|\,
|\nabla_A\underline{\alpha}|^2\Big|\leq \zeta \ts^{-1}\quad \text{on
  \(X_{\delta_0}^{'a}\supset X_{\delta_0'}^{'a}\)}.
\]
In fact, a slight modification of the argument in the preceding proof
would yield a direct, similar pointwise bound for
\(|\nabla_A\alpha|\). However, the bound for \(|\nabla_A\ud{\alpha}|\)
given in the previous proposition is slightly better, and it is more
amenable to our applications later.
\end{rem}

\section{Monotonicity and its consequences}\label{sec:mono}

This section contains variants and counterparts of  Steps (3)  and (4) of
Taubes' proof, as outlined in 
Section \ref{sec:synopsis} above. Intermediate steps of these
arguments are used to improve various integral bounds and pointwise
estimates in the previous two sections. 

Let \((A, \Psi)=(A_r, \Psi_r)\) be as in the statement of  Lemma
\ref{lem:Etop-bdd1}, 
and let \(B\subset X\) be an open set. Following Taubes, we define the
(local) {\em energy} of \((A, \Psi)\in \Conn (\bbS^+)\times \Gamma(\bbS^+)\) on \(B\) to
be:
\[
\mathcal{W}_B(A, \Psi):=4^{-1}r\int_B |\nu| \, \Big||\nu|-|\psi|^2\Big|.
\]
Taubes' proof of the monotonicity formula hinges on a bound for
\(\mathcal{W}_X\). In contrast to \cite{Ts}, \(\mathcal{W}_X\) can be infinite in our
situation. However, Lemma \ref{T:lem3.1} provides us with bounds on
\(\mathcal{W}_{X_\bullet}\) for any compact \(X_\bullet\subset
X'\). 


\subsection{Bounding \(\mathcal{W}_{Z_\delta }\)}\label{sec:W_z}
Recall the notation
\(\td{r}=r\ts\). 

Let \(\delta _1:= 2r^{-1/6} (\ln r+\smE)^{3/4}\). Over \(X^{'a}_{\delta _1/2}\), \(\ts^{-5}\leq
\td{r}\, 
(\ln r+\smE)^{-9/2}\). This ensures that there is an \(r\)-independent
constant \(\zeta _e\) such that  the functions \(\varepsilon _1\) and \(\varepsilon
'\leq \zeta \varepsilon _1\) from Proposition
\ref{prop:curv-varpi} have the property that 
\begin{equation}\label{bdd:eps_1}
(\ln r+\smE ) \, (\varepsilon _1+\varepsilon ') \leq \zeta _e \quad \text{over \(X^{'a}_{\delta
    _1/2}\)}. 
\end{equation}

Fix an  \(x\in \nu ^{-1}(0)\cap X^{'a}\). Choose local coordinates \((t,
x_1, x_1, x_3)\) centered at \(x\) so that
\(\nu^{-1}(0)\) is identified with the \(t\)-axis, and that
\(\omega=2\nu ^+=d\eta\), where 
\begin{equation}\label{eq:eta0}
\begin{split}
\eta& =-\big(\frac{x_1^2+x_2^2}{2}\big) \,
dt-x_3(x_1dx_2-x_2dx_1)+O\, (\sigma ^3), \\
|\eta|& \leq \sigma \, |\nu |/2+O\, (\sigma ^3). 
\end{split}
\end{equation}
Let \(Z_x(\delta , l)\subset
X''\) denote the cylinder consisting of points with distance no greater
than \(\delta \) from \(\nu ^{-1}(0)\), and whose \(t\)-coordinate
lies in the interval \([-l, l]\). The admissibility of \(\nu \)
implies that there exist \(x\)-independent constants \(1>\delta _c>0\),
\(1\geq l_c>0\) such that for any given \(x\in \nu ^{-1}(0)\cap X^{'8}\), local coordinates of the type described above are
defined over \(Z_x(\delta , l)\) for all \(\delta \leq \delta _c\),
\(l\leq 2l_c\). 

\begin{lemma}\label{lem:monotone0}
Let \(r_0\), and \(\grt(r)\) be as in Lemmas
\ref{lem:Etop-bdd1} and \ref{lem:Etop-bddf}. 
There exist positive constants \(r_1\geq r_0\), 
\(\zeta\), \(\zeta '\) such that \(\zeta , \zeta '\) depend only on
\(\nu \) and the parameters listed in (\ref{parameters}), 
and \(r_1\) depends only on all the
above as well as on \(\smE\) (in particular, they are independent of
\(r\), \(\delta \), \(x\), and \((A, \Psi )\)), with the following
significance: 
For any given \(r\geq r_1\),  \(x\in \nu ^{-1}(0)\cap X^{'a}\),
\(0\leq a\leq \frac{47}{8}\), \(\delta
\) and \(l\) such that \(\delta _c\geq\delta \geq \delta _1/2\), 
\(l_c\geq l\geq l_c/4\), and  \((A, \Psi)=(A_r, \Psi_r)\) satisfying
the assumption of Lemma
\ref{lem:Etop-bdd1}, 
one has: 
\[
\mathcal{W}_{Z_x(\delta , l)}(A, \Psi)\leq 
\begin{cases}\zeta \delta ^2(\ln r+\smE ) & \\
\zeta '\delta ^2(\grt(r)+\smE ) &\text{when Lemma \ref{lem:Etop-bddf} applies.}
\end{cases}
\]
\end{lemma}
\pf Let \(l\) be such that \(l_c\geq l\geq l_c/4\), and let \(\lambda _l(t)\) be a smooth cutoff function which equals 1 over
\([-l, l]\), and vanishes outside \([-l-l_c/8, l+l_c/8]\). Let
\(\eta_l:=\lambda _l\eta\) and 
\[
\scrW_\delta :=4^{-1}\int_{Z_x(\delta , l+l_c/8)}\lambda _l\, r |\nu| \,
\Big||\nu |-|\psi|^2\Big|.
\]
Note that 
\(\mathcal{W}_{Z_x(\delta , l+l_c/8)}(A, \Psi )\geq\scrW_\delta \geq \mathcal{W}_{Z_x(\delta , l)}(A, \Psi )\).

Suppose \(r\) is sufficiently large so that both \(r\geq r_0\) and
\(\delta _1=\delta _1(r)\leq \delta _c\). 
Fix \(\delta \geq \delta _1/2\). 
By the Seiberg-Witten equation \(\grS_{\mu,
  \hat{\grp}}(A, \Psi)=0\), (\ref{eq:eta0}), Propositions
\ref{T:lem3.2}, \ref{prop:curv-varpi}, \ref{T:prop3.2}, and 
\ref{T:prop3.1-}, for \(x\in \nu ^{-1}(0)\cap X^{'\frac{47}{8}}\), 
\begin{equation}\label{eq:W-ubdd-0}
\begin{split}
\frac{i}{2}\int_{Z_x(\delta , l+l_c/8)}F_A\wedge d \eta_l
&\geq \frac{i}{2}\int_{Z_x(\delta , l+l_c/8)}\lambda _l\, F_A\wedge
\omega-\zeta _1\delta \int_{Z'_l}
|F_A|\, |\nu |
\\
& \geq  \scrW_\delta -\zeta  _2\delta ^2-\zeta _1'\delta \, \mathcal{W}_{Z'_l}(A,
\Psi ),
\end{split}
\end{equation}
where \(Z'_l\subset Z_x(\delta ,l +l_c/8)-Z_x(\delta , l)\) is the support of
\(\lambda '_l\). 
Meanwhile, by the Stokes' theorem,  (\ref{bdd:eps_1}), Propositions
\ref{prop:curv-varpi}, \ref{T:prop3.2}, \ref{T:lem3.2} and
\ref{T:prop3.1-}, one has: {\small
\begin{equation}\label{eq:DE-mono-0}
\begin{split}
& \frac{i}{2}\int_{Z_x(\delta , l+l_c/8)}F_A\wedge d
\eta_l =\frac{i}{2}\int_{\partial Z_x(\delta , l+l_c/8)}F_A\wedge  \eta_l\\
& \, \, \leq \frac{1}{4}\int_{\partial Z_x(\delta , l+l_c/8)}\left( \lambda _l |\nu|\Big(r\Big||\nu |-|\psi|^2\Big|
  \big(1+\frac{\zeta _3}{\ln r+\smE}\big)+\zeta_3'\,
  \tilde{\sigma}^{-2}\Big) \sigma (1+\zeta' \sigma  )\right) \\
&\, \, \leq \big(1+\zeta _3 (\ln r+\smE)^{-1}\big)\frac{\delta +\zeta '\delta ^2}{2}\, 
\frac{d}{d\delta }\scrW_\delta +\zeta _4\delta ^2\quad \text{for \(\delta \geq\delta _1/2\).}
\end{split}
\end{equation}}
Combining (\ref
{eq:DE-mono-0}) and (\ref{eq:W-ubdd-0}) and appealing to Lemma
\ref{T:lem3.1}, one obtains: 
\begin{equation}\label{Dineq:W_0}
\begin{split}
\frac{d}{d\delta }\left(e^{f_0}\scrW_\delta \right) & \geq- \big(\zeta _1''\, \mathcal{W}_{Z'_l}(A,
\Psi )+ \zeta _4'\delta \big)\, e^{f_0}, \\
& \geq- \zeta _5 \, (\ln r+\smE )\, e^{f_0}, \qquad\text{when \(\delta \geq \delta _1/2\), where}\\
f_0(\delta ) & :=-2\, (1+\zeta _3(\ln r+\smE)^{-1})^{-1}\ln \, (\zeta '\delta /(1+\zeta '\delta )).
\end{split}
\end{equation}
Integrating the preceding differential inequality, one finds that there exist positive constants \(\zeta
_w\), \(\zeta ''\) independent of \(r\), \(R\), \(s\), \(x\), and \((A,
\Psi )\),  such that for \(s>\delta \geq \delta _1\), 
\[
\zeta _w\scrW_s/s^2\geq \scrW_\delta/\delta ^2-\zeta''(\delta
^{-1}-s^{-1})\, (\ln r+\smE ).
\]
Setting \(s=\delta _c\) and appealing to Lemma 
\ref{T:lem3.1}, 
this gives: 
\[
\mathcal{W}_{Z_x(\delta , l)}(A, \Psi)\leq \zeta _0 \, \delta (\ln r+\smE)
\quad \text{for \(\delta _c\geq \delta \geq\delta _1/2\), \(l\leq l_c\)}, 
\]
where \(\zeta _0\) is a positive constant \(\zeta _0\)  independent of \(r\),
\(\delta \), \(x\), \(l\), and \((A, \Psi)\). In particular, this implies
that 
\[
 \mathcal{W}_{Z'_l}(A,
\Psi )\leq \mathcal{W}_{Z_x(\delta , l+l_c/8)}(A, \Psi)\leq \zeta _0 \, \delta (\ln r+\smE)
\quad \text{for \(\delta \geq\delta _1/2\), \(l\leq l_c\)}. 
\]
Reinsert this into the first line of (\ref{Dineq:W_0}) to get
\[
\frac{d}{d\delta }\left(e^{f_0}\scrW_\delta \right) \geq- \zeta _5'\, 
\delta ^{-1+\zeta _6(\ln r+\smE )^{-1}} (\ln r+\smE ). 
\]
Integrating as before, one gets 
\[
\mathcal{W}_{Z_x(\delta , l)}(A, \Psi)\leq \zeta \, \delta ^2(\ln r+\smE)
\quad \text{for \(\delta _c\geq \delta \geq\delta _1/2\), \(l\leq l_c\)} 
\]
as asserted. The same argument shows that when
Lemma \ref{lem:Etop-bddf} applies, the factor of \(\ln r\) in the
preceding inequality may be replaced by \(\grt(r)\). 
\epf

The preceding lemma has the following corollary: 
\begin{lemma}\label{lem:u^a-int}
Suppose \(X_\bullet\subset X^{'a}\), \(0\leq a\leq \frac{47}{8}\), has
length \(|X_\bullet|=1\). Let \(\delta _c\) be as defined in the
paragraph preceding Lemma \ref{lem:monotone0}. Then
for \(-1<b\leq 1\), 
\[
\int_{X_\bullet\cap X^{'a}_{\delta }}r |\nu|^b\, \Big||\nu
|-|\psi|^2\Big|\leq \zeta  \delta _c^{-1+b} (\ln r +\smE )\quad \text{\(\forall r\geq r_1\) and
  \(\delta \geq \delta _1/2\).} 
\]
In the above, \(r_1\), \(\delta _1\) are as in the previous lemma, and
\(\zeta \) depend only on
\(\nu \) and the parameters listed in (\ref{parameters}). 
When Lemma \ref{lem:Etop-bddf}
  applies, the factors \(\ln r\) in the preceding inequality can be replaced
  by \(\grt(r)\). 
\end{lemma}

\pf Fix \(r\geq r_1\) and \(\delta \), 
  \(\delta _c\geq\delta \geq \delta _1/2\). Decompose \(X^{'a}_\delta\) as the union \(X^{'a}_{\delta _c}\cup
\bigcup_{n=1}^{N}(Z^{'a}_{\delta_n}-Z^{'a}_{\delta_{n-1}})\), where
\(\delta_n:=2^{n}\delta \), and \(N\) is the
smallest \(n\in \bbZ^+\) such that \(\delta_n\geq\delta _c\). Then
according to (\ref{eq:v-sig}), \(\zzz_v \, \delta _n\geq|\nu|\geq \zzz_v  ^{-1}\delta _{n-1}\) over \(Z^{'a}_{\delta_n}-Z^{'a}_{\delta_{n-1}}\). 
Make use of Lemma \ref{lem:Etop-bddf} over \(X^{'a}_{\delta _c}\cap
X_\bullet\), and Lemma \ref{lem:monotone0} over each
\((Z^{'a}_{\delta_n}-Z^{'a}_{\delta_{n-1}})\cap X_\bullet\), one has: 
\begin{equation}\label{int:u^a}
\begin{split}
& \int_{X^{'a}_\delta\cap X_\bullet}
 r |\nu|^b\, \Big| |\nu|-|\psi|^2\Big|\\
&\quad \leq \zzz_v  \sum_{n=1}^{N}\delta
_{n-1}^{-1+b}\, \mathcal{W}_{Z^{'a}_{\delta _n}\cap X_\bullet}(A, \Psi
)+\zzz_v  \, \delta _c^{-1+b}\, \mathcal{W}_{X^{'a}_{\delta _c}\cap  X_\bullet}(A, \Psi )\\
&\quad \leq\zeta _1(\ln r+\smE )\sum_{n=1}^{N}\delta
_{n-1}^{-1+b}\delta _n^2+\zeta'_1\delta _c^{-1+b}(\ln r+\smE)\\
&\quad \leq 4\zeta _1(\ln r+\smE )\, \delta ^{1+b}\sum_{n=0}^{N-1}(2^{1+b})^n
+\zeta'_1\delta _c^{-1+b}(\ln r+\smE)\\
 & \quad \leq \zeta\delta _c^{-1+b}(\ln r+\smE).
\end{split}
\end{equation}
The same argument shows that when Lemma \ref{lem:Etop-bddf}
  applies, the factors \(\ln r\) in (\ref{int:u^a}) can be replaced
  by \(\grt(r)\). 
\epf

\subsection{Integral bounds redux}\label{sec:improved}

We are now in a position to eliminate the undesirable factor of \(\ln
r\), which made its first appearance in Lemma  \ref{lem:Etop-bdd1} and
has propagated throughout our discussion so far. 
A closer look at the proof of Lemma  \ref{lem:Etop-bdd1} reveals that
this factor  originates from the preliminary lower bound on the last
term of (\ref{bdd:CSD-lower}), namely
\[
-\frac{ir}{4}\int_{\hat{Y}_{[l,L]}}ds\,(
*_Y\xi_\nu)\wedge(F_B-F_{B_0}). 
\]
The bound in Section \ref{sec:4.2} made use of \(L^2\)-bounds of (components of)
\(F_A-F_{A_0}\). In what follows we argue differently, using instead a
better \(L^1\)-bound on \(F_A-F_{A_0}\). The latter bound is in turn
made possible by the pointwise estimates from Section \ref{sec:pt-est} and
the preliminary local energy bounds from the previous
subsection.

\begin{lemma}\label{lem:M_bdd+}
Let \(r_0\), \(\grt(r)\) be as in Lemmas
\ref{lem:Etop-bdd1} and \ref{lem:Etop-bddf} and
let \(X_\bullet\subset X^{'a}\), \(0\leq a\leq \frac{47}{8}\), be compact and  \(|X_\bullet|=1\). 
Then there exist \(r\)-independent constants \(\zeta\), \(\zeta '\),
\(r_1\geq r_0\), such that for all \(r\geq r_1\), 
\[
\|F_A-F_{A_0}\|_{L^1(X_\bullet)}\leq 
\begin{cases}
\zeta (\ln r+\smE )&\\
 \zeta '(\grt(r)+\smE ) &\text{when Lemma \ref{lem:Etop-bddf}
  applies}.
\end{cases}
\]
In the above,  \(\zeta\), \(\zeta '\) depend only on
\(\nu \) and the parameters listed in (\ref{parameters}), 
and \(r_1\) depends only on all the
above as well as on \(\smE\). 
\end{lemma}
\pf Let \(\delta _0'\), \(\delta _1\)
be respectively as in Proposition \ref{prop:curv-varpi} and Lemma
\ref{lem:monotone0}.  Set \(\delta_0''=r^{-1/4}\) and assume that \(r\) is
sufficiently large so that \(\delta _0''>\delta _0'\).  Decompose
\(X_\bullet\) as the union of three regions, \(\scrR_i\), \(i=1, 2,
3\): \(\sigma \leq \delta _0''\) on \(\scrR_1\);  \(\delta _0''\leq
\sigma \leq \delta _1\) over \(\scrR_2\); \(\sigma \geq\delta _1\)
over \(\scrR_3\). Then according to Propositions \ref{prop:curv-varpi}
and \ref{T:prop3.2}, \ref{T:lem3.2} and
\ref{T:prop3.1-}, 
\[
|F_A|\leq \begin{cases} \zeta  _1r\ts +\zeta _1'\ts^{-2}& \text{over
    \(\scrR_1\)};\\
\zeta _2r \, (-u)_++\zeta _2' \ts^{-2} & \text{over
    \(\scrR_2\cup \scrR_3\)}, 
\end{cases}
\]
where \(-u=|\nu |-|\psi |^2\) as before. 
Consequently, 
\begin{equation}\label{eq:M-bdd+}
\begin{split}
& \|F_A-F_{A_0}\|_{L^1(X_\bullet)}\\
& \quad \leq r\int_{\scrR_1}(\zeta _1\ts+\zeta _1'r^{-1}\ts^{-2})+
\zeta_3\, (\delta _0'')^{-1}\mathcal{W}_{\scrR_2}(A, \Psi )+\zeta _3'
\int_{\scrR_3} r\, | -u|\\
& \qquad \qquad +\zeta _2'\int_{\scrR_2\cup \scrR_3}\ts^{-2}\\
&\quad \leq\zeta _1'' r\, (\delta _0'')^4+\zeta_4\, (\delta
_0'')^{-1}\delta _1^2\, (\ln r +\smE) +\zeta _4' \, (\ln r+\smE ) +\zeta _2''.\\
 & \quad \leq \zeta_5+\zeta _5'(\ln r+\smE ).
\end{split}
\end{equation}
In the above, Lemma \ref{lem:monotone0} was invoked to bound the
second term in the second line in terms of the the second term in the
third line; and Lemma \ref{lem:u^a-int} was used (with \(b=0\)) to
bound the third term in the second line in terms of the the third term in the
third line.  When Lemma \ref{lem:Etop-bddf}
  applies, the same argument shows that the factors of \(\ln r\) in
  (\ref{eq:M-bdd+}) can be replaced by \(\grt(r)\). \epf


The preceding lemma is now used to amend the bounds in Section
\ref{sec:4}. 
\begin{prop}\label{prop:integral-est2}
The conclusions of Lemma \ref{lem:Etop-bdd1}, Proposition \ref{prop:SW-L2-bdd},
Proposition \ref{est-L^2_1},
Lemma \ref{co:E-omega-bdd},
Lemma \ref{co:E-omega-bdd3}, and 
Lemma \ref{T:lem3.1} 
 remain valid with all factors of \(\ln r\) replaced by a positive
 constant \(z_0>1\). This constant depends only on \(\nu \), the
 parameters listed in (\ref{parameters}), and \(\smE\). 
\end{prop}
\pf 
Recall  from (\ref{bdd:CSD-lower}) that when \(\hat{Y}_{i, [l, L]}\subset
X^{'a}\) 
\[
\begin{split} 
&  -2\op{CSD}_{\mu_+, \hat{\grp}}^{\partial \hat{Y}_{i,
     [l,L]}}(B,\Phi)
\geq -\frac{ir}{4}\int_{\hat{Y}_{i, [l,L]}}ds\,(
*_Y\xi_\nu)\wedge(F_B-F_{B_0}). 
\end{split}
\]
Combining the preceding inequality with 
(\ref{diff:CSD-E}), (\ref{eq:xi-exp}), and  Lemma \ref{lem:M_bdd+}, we
see that if \(\grt(r)\co [r_0, \infty)\to (1, \infty)\) is such that
(\ref{assume:t(r)}) holds, 
then for any \(\hat{Y}_{i, [l,L]}\subset X^{'a}\)
\begin{equation}\label{bdd:E-lo}
 \scrE^{'\mu_r}_{top}(\hat{Y}_{i, [l,L]})(A,\Psi )\geq -\zeta  \, r
 e^{-\kappa _il} (\smE +\grt (r))
\end{equation}
for some positive constant \(\zeta \) independent of \(r\) and
\(\hat{Y}_{i, [l,L]}\subset X^{'a}\). 
But this means that if (\ref{assume:t(r)}) holds for some 
\(\grt(r)=\grt_k(r)\), then it holds also for
\(\grt(r)=\grt_{k+1}(r)\), where 
\begin{equation}\label{def:t_k}
  \grt_{k+1}(r)=
  \ln \, (1+\grt_k(r)) +\ln\,  (\zeta /\zeta'_5).
\end{equation}
We already know that (\ref{assume:t(r)}) holds for
\(\grt(r)=\grt_1(r)\), 
\[
\grt_1(r):=\ln r. 
\] 
So in this way we get a sequence of functions \(\{\grt_1, \grt_2,
\cdots\}\) and (\ref{assume:t(r)}) holds for \(\grt(r)=\grt_k(r)\)
\(\forall k\in \bbN\). Now, there is a number \(z_0>1\) such that 
\[
\ln\,  (\smE +z)+\ln\,  (\zeta /\zeta
'_5)<z/2 \quad\forall z\geq z_0. 
\]
Then according to (\ref{def:t_k}), 
\[
\grt_{k+1}(r)<\frac{1}{2}\grt_k(r) \quad \text{if \(\grt_k(r)\geq z_0\)}.
\]
Thus, for any fixed \(r\geq r_0>1\), there exists an \(n\in \bbN\)
(depending on \(r\)) such that \(\grt_{n+1}(r)\leq
z_0\). Consequently, 
\begin{equation}\label{assume:t(r)k}
\scrE^{'\mu_r}_{top} (X_\bullet)(A, \Psi ) \geq  -\zeta '_5 \, r\quad \forall X_\bullet\subset X^{'a}-X^{'a}_{z_0\, 
  \pmb{\hatl}}, r\geq r_0. 
\end{equation}
That is to say, the assumption (\ref{assume:t(r)}) in Lemma
\ref{lem:Etop-bddf} holds when \(\grt(r)\) therein is taken to be
the constant function \(z_0\). So the conclusion of Lemma
\ref{lem:Etop-bddf}, 
(\ref{eq:E_top-Mf}), also  holds with
\(\grt(r)\) replaced by \(z_0\). This is the same as saying that Lemma
\ref{lem:Etop-bdd1} holds with all factors of \(\ln r\) replaced by
\(z_0\). It follows that all appearances of \(\ln r\) or \(\grt(r)\)
in Proposition \ref{prop:SW-L2-bdd},
Proposition \ref{est-L^2_1},
Lemma \ref{co:E-omega-bdd},
Lemma \ref{co:E-omega-bdd3}, and 
Lemma \ref{T:lem3.1} can be likewise be replaced by \(z_0\). \epf

\begin{lemma}\label{lem:amend-int}
The conclusions of  Lemmas \ref{lem:M_bdd+}, 
\ref{lem:monotone0} and  \ref{lem:u^a-int} 
 remain valid with all factors of \(\ln r\) replaced by a positive
 constant \(z_0>1\). This constant depends only on  \(\nu \), the
 parameters listed in (\ref{parameters}), and \(\smE\). 
In particular, there exist constants \(r_1\), \(\zeta \), \(\zeta '\),
\(\zeta_1\), \(\zeta _1'\), such that \(\forall r\geq r_1\) and for
any \(X_\bullet\subset X^{'a}\), \(0\leq a\leq \frac{47}{8}\), 
\begin{equation}\label{bdd:L^14d}
\|F_A\|_{L^1(X_\bullet)}\leq \zeta  |X_\bullet|+\zeta _1\smE; \quad
\|A'-A_0\|_{L^1(X_\bullet)}\leq (\zeta '  |X_\bullet|+\zeta
'_1\smE)\, |X_\bullet |,
\end{equation}
where \((A', \Psi ')=u\cdot (A, \Psi )\) is in the normalized
Coulomb-Neumann gauge over \(X_\bullet\), and the constants \(\zeta
\), \(\zeta '\), \(\zeta
_1\), \(\zeta _1'\) depend only on \(\nu \) and the
 parameters listed in (\ref{parameters}), and \(r_1\) depends on
 \(\smE\) in addition to all the above. 
\end{lemma}
\pf Simply replace all usage of  Lemma \ref{T:lem3.1} in the proofs of
Lemmas \ref{lem:M_bdd+}, 
\ref{lem:monotone0} and  \ref{lem:u^a-int} by the amended
version of  Lemma \ref{T:lem3.1} in Proposition
\ref{prop:integral-est2}. The \(L^1\)-bound on \(A'-A_0\) follows from
that for \(F_A\) via a standard elliptic estimate, noting that the lowest absolute
value of non-zero eigenvalues of \(*d+d*\) on 1-forms on \(X_\bullet\)
with Neumann boundary condition has a
lower bound propotional to \(|X_\bullet|^{-1}\)). 
\epf


The next two lemmas will be useful in Section \ref{sec:7}.
\begin{lemma}\label{bdd:A3d}
Adopt the assumptions and notation of Lemma \ref{lem:Etop-bdd1}. Let
\(\hat{Y}_{[l, l+1]}\subset X^{'a}\), \(0\leq a\leq \frac{47}{8}\), and suppose \(\eta\) is an exact
2-form on \(\hat{Y}_{[l, l+1]}\). Then there exist positive constants
\(r_0\), \(\zeta _2\), \(\zeta _2'\) such that \(\forall r\geq r_0\),
\[
  \Big|\int_{Y:l}i(A-A_0)\wedge
\eta\Big|+\Big|\int_{Y:l+1}i\, (A-A_0)\wedge \eta\Big|\leq (\zeta
_2+\zeta _2'\smE)\, \|\eta \|_{L^\infty(\hat{Y}_{[l, l+1]})}.\]
(Note that both integrals in the preceding formula are gauge-invariant
due to the exactness of \(\eta\).) In the above, the constants \(\zeta _2\), \(\zeta _2'\) depend only on \(\nu \) and the
 parameters listed in (\ref{parameters}), and \(r_0\) depends on
 \(\smE\) in addition to all the above. 
\end{lemma}
\pf Let \(\chi _l(s)\co [ l, l+1]\to
[0, 1]\) be a non-increasing function that equals 1 on \([l, l+1/4]\),
and vanishes over \([l+3/4, l+1]\). Use the same notation to denote its
pull back over \(\hat{Y}_{i,
  [l, l+1]}\) under the projection \(\hat{Y}_{[l,l+1]}\simeq [l, l+1]\times Y\to [l, l+1]\). 
\begin{equation}
\begin{split}
& i\int_{Y_{i:l}} (A-A_0)\wedge\eta =i\int_{\hat{Y}_{i,
  [l, l+1]}}d\big( \chi_l (s)\, (A-A_0)\wedge\eta\big)\\
&\quad = i\int_{\hat{Y}_{i,
  [l, l+1]}}  \chi_l (s)\, (F_{A}-F_{A_0})\wedge\eta +i\int_{\hat{Y}_{i,
  [l, l+1]}} \chi_l' (s) \, ds\, (A'-A_0)\wedge\eta ,
\end{split}
\end{equation}
where \((A', \Psi ')=u\cdot (A, \Psi )\) is in the normalized
Coulomb-Neumann gauge over \(\hat{Y}_{i,
  [l, l+1]}\).  
The asserted bound for \(\Big|\int_{Y:l}i(A-A_0)\wedge
\eta\Big|\) then follows from Equation (\ref{bdd:L^14d}).  The same
argument can be used to bound \(\Big|\int_{Y:l+1}i(A-A_0)\wedge
\eta\Big|\), by simply replacing \(\chi _l\) in the preceding formula
by \(1-\chi _l\). 
\epf 

\begin{lemma}\label{lem:F_v}
  Adopt the assumptions and notation of Lemma \ref{lem:Etop-bdd1}. 
Suppose \(X_\bullet\) is such that \(\partial X_\bullet\subset
X^{'a}\), \(0\leq a\leq \frac{47}{8}\). 
Then there exist positive constants \(r_0\), \(\zeta _h, \zeta '_h\)
such that 
\begin{equation}\label{eq:F-nu0}
\Big| \int_{X_{\bullet}}i F_{A_r}\wedge \nu\Big| \leq  \zeta _h\smE +\zeta_h'.
 \end{equation}
In the above, the positive constants \(\zeta _h\), \( \zeta '_h\)
depend  only on \(\nu \) and the
 parameters listed in (\ref{parameters}), and \(r_0\) depends on
 \(\smE\) in addition to all the above.  
\end{lemma}
\pf Write 
\begin{equation}\label{eq:F-nu}
\begin{split}
& i\int_{X_{\bullet}} F_{A}\wedge\nu\\
& \quad =4r^{-1}\scrE_{top}^{\mu_r}(X_{\bullet})(A,
 \Psi)-4r^{-1}\scrE_{top}^{w_r}(X_{\bullet})(A, \Psi)\\
& \quad =4r^{-1}\scrE_{top}^{'\mu _r}(X_{\bullet})(A,
\Psi)-4r^{-1}\scrE_{top}^{w_r}(X_{\bullet})(A, \Psi)\\
& \qquad \qquad -ir^{-1}\int_{\partial
   X_\bullet} (A-A_0) \wedge (*_4\mu _r)\\
& \quad =4r^{-1}\scrE_{top}^{'\mu _r}(X_{\bullet})(A,
 \Psi)-4r^{-1}\scrE_{top}^{w_r}(X_{\bullet})(A, \Psi)-i\int_{\partial
   X_\bullet} (B-B_0) \wedge \grv
\\
& \quad =4r^{-1}\scrE_{top}^{'\mu _r}(X_{\bullet})(A,
 \Psi)-4r^{-1}\scrE_{top}^{w_r}(X_{\bullet})(A, \Psi) -i\int_{\partial
   X_\bullet} (B'-B_0) \wedge \grv, 
\end{split}
\end{equation}
where \((B, \Psi ):=(A,\Psi)|_{\partial
   X_\bullet}\), and \((B', \Psi ')\) is the representative of
 \([(B, \Psi )]\) in the normalized Coulomb gauge. We used the fact
 \(\grv|_{\partial X_\bullet}\) is exact to derive the last line above
 from the penultimate line. Using Lemma \ref{lem:Etop-bdd1} (or Lemma
 \ref{lem:Etop-bddf}), Lemma 
 \ref{co:E-omega-bdd}, 
and their amendments in Proposition
 \ref{prop:integral-est2} to estimate the first two terms in the last
 line of the previous formula, and using the preceding lemma to bound
 the last term, we arrive at the asserted bound (\ref{eq:F-nu0}). 
\epf

\subsection{The monotonicity formula}

We begin by improving the pointwise estimates of
\(|F_A^-|\) given in Propositions \ref{prop:curv-varpi}
and \ref{T:prop3.2}. The estimates in the next lemma is only better
than that in Proposition \ref{prop:curv-varpi} where \(3  r^{-1/6} \,
(\ln r+\smE )^{1/6}\geq \sigma \geq r^{-2/9}(\ln r)^{1/4}\). However,
this improvement will be useful later. 
Let \(u:=|\psi|^2-|\nu|\) as before. 

\begin{lemma}\label{lem:curv-varpi-3}
   Adopt the assumptions and notation of Lemma \ref{lem:Etop-bdd1}. 
   Let \(\delta _1':= r^{-2/9}(\ln r)^{1/4}\). 
There exist positive constants 
\(r_0>8\), \(\zeta \), \(\zeta'\), with the following significance:
For any \(r\geq r_0\) one has: 
\begin{equation}\label{eq:curv-varpi-3}
\begin{split}
|F_A^-|& \leq (2^{-3/2}+\zeta (\ln r)^{-1} )\, 
r \, (-u)_+ +\zeta'
\,  \ts^{-2} 
\quad
\text{over \(X^{'a}_{\delta'_1}\), \(0\leq a\leq \frac{47}{8}\).}
\end{split}
\end{equation}
In the above, the constants  \(\zeta \)
and \(\zeta '\) depend  only on \(\nu \) and the
 parameters listed in (\ref{parameters}), and \(r_0\) depends on
 \(\smE\) in addition to all the above.  
More generally, given any positive \(\epsilon <1/2\),
there exists constants \(r_\epsilon \), \(\zeta _\epsilon \), \(\zeta'_\epsilon \),
such that \(\forall r\geq r_\epsilon \), 
\begin{equation}\label{eq:curv-varpi-4}
\begin{split}
|F_A^-|& \leq (2^{-3/2}+\zeta _\epsilon (\ln r)^{-1} )\, 
r \, (-u)_+ +\zeta'_\epsilon 
\,  \ts^{-2} 
\quad
\text{over \(X^{'a}_{\delta _{1,\epsilon }}\),}
\end{split}
\end{equation}
where \(\delta _{1, \epsilon }:= r^{-(1-\epsilon )/4}\). 
\end{lemma}
\pf Let \(\delta _1\), \(\delta _0'\) \(\delta _0\) be respectively
as in the beginning of Section \ref{sec:W_z}, Proposition
\ref{prop:curv-varpi}, 
and Proposition \ref{T:prop3.2}. Assume below that \(r\) is sufficiently large so that \(\delta
_1>\delta '_1>\delta _0'>8\delta _0\).  
Return to the proof of Proposition \ref{T:prop3.2}. Re-examine the
bounds for each term on the right hand side of (\ref{eq:q-est}),
taking \(\varepsilon =(8\ln r)^{-1}\). The same argument, with
(\ref{eq:L^2_1}) replaced by its amended version from Proposition
 \ref{prop:integral-est2}, show that
the first term on the RHS of (\ref{eq:q-est}) is bounded by a constant
multiple of \(r^{1/2}\). Meanwhile, the second and the third terms
are respectively bounded by 
constant multiples of \(r^{1/3}(\ln r)\,\ts^{-1}\)  and \(\ts^{-2}\) when \(x\in
X^{'a}_{\delta _0'}\). Together, we have
\begin{equation}\label{est:qa}
  \begin{split}
    & q(x) \leq \zeta _1 r^{1/3}(\ln r)\, \ts^{-1}+\zeta _2r^{1/2}+\zeta _3\ts^{-2}+\text{the last term of (\ref{eq:q-est})}
\quad 
\end{split}
\end{equation}
for \(x\in X^{'a}_{\delta _0'}\).

A little more work is required to improve the bound
on the last term. 
Note that by Proposition \ref{prop:curv-varpi} and (\ref{bdd:eps_1}), 
(\ref{eq:curv-varpi-3}) already holds where \(\sigma \geq \delta
_1/2\). Thus, it suffices to consider \(x\in Z_{\delta _1/2}^{'a}\cap X_{\delta_0'}^{'a}\). Fix an
\(x\in Z_{\delta _1/2}^{'a}\cap X_{\delta_0'}^{'a}\) and reintroduce the notation \(\sigma
_x:=\sigma (x)\); \(\ts_x:=\ts (x)\). This time, we divide \(B(x, 2\rho
_0)\cap X^{'a}_{\delta _0/2}\) into four regions defined from three parameters \(\rho _1,
\delta _2\), and \(\delta _3\). These parameters are chosen such that: 
\[
\rho _1\leq \sigma _x/4; \quad \delta _2\geq \max \, (2\sigma _x,
\delta _1/2); \quad \delta _0/2\leq \delta _3<\sigma _x/2. 
\]
Since \(\sigma _x\leq \delta _1/2\) by assumption, we take \(\delta
_2=\delta _1\). 
Then \(B(x, 2\rho_0)=\bigcup_{1\leq i\leq 4}\scrR'_i\), \(\scrR'_1:=B (x, \rho _1)\subset
Z^{'a}_{\delta _2}-Z^{'a}_{\delta _3}\); \(\scrR'_2:=X^{'a}_{\delta_2}\cap B(x,2\rho_0)\), \(\scrR'_3:=(Z^{'a}_{\delta _2}-Z^{'a}_{\delta
  _3})\cap B(x,2\rho_0)-\scrR'_1\);
\(\scrR'_4:=(Z^{'a}_{\delta_3}-Z^{'a}_{\delta _0/2})\cap B(x,2\rho_0)\). Note 
that with the present choice of \(\varepsilon \) and the assumption
that \(\delta >\delta _0'\), we have simpler pointwise bounds for \(\eta_2\): 
\[
\eta _2\leq  \zeta _2r\ts^{-2}(-u)_+\leq \zeta _2' r\ts^{-1}\quad
\text{on \(X^{'a}_{\delta _0'}\)}. 
\]
Integrate over each  of the four regions separately to replace
(\ref{int:eta2}) by: 
\begin{equation}\label{int:eta2-a}
\begin{split}
& \int_{B(x,2\rho_0)} \eta_2\dist (x,
\cdot)^{-2}\\
&\quad  \leq   \zeta '_1\,
r\ts^{-1}_x\, \rho _1^2  +
\zeta _2\, \delta _2^{-2} \delta _2^{-3/2}
+\zeta _3 \delta _3^{-3}
\delta _2^2\rho _1^{-2}+\zeta_4 r\sigma _x^{-1} \delta _3^2
\, \, \text{when \(x\in Z^{'a}_{\delta_1/2}\)}. 
\end{split}
\end{equation}
In the above, the first and the last terms on the right hand side
respectively bound the integrals over the first and the last regions,
following the same computations that respectively yield the first and
the last term on the RHS of (\ref{int:eta2}). The second term on the
RHS of (\ref{int:eta2-a}) bounds the integral over \(\scrR'_2\), which follows from the
amended version of Lemma
\ref{lem:u^a-int} (with \(b=-1/2\)) in Lemma
\ref{lem:amend-int}, together with the observation that 
 \(\dist (x, \cdot)\geq \delta _2/2\) over \(\scrR'_2\).  The third term on the
RHS of (\ref{int:eta2-a}) bounds the integral over \(\scrR'_3\), which follows from the
amended version of Lemma \ref{lem:monotone0} in Lemma
\ref{lem:amend-int}. Now we take \(\rho _1=\delta _3/2\), and \(\delta
_3=(r^{-1}\ts_x \delta _2^2)^{1/7}\). This choice of \(\delta _3\)
meets the requirement that \(\delta _3<\sigma _x/2\) when \(r\) is
sufficiently large. With such choice, 
\begin{equation}\label{est:qb}\begin{split}
\int_{B(x,2\rho_0)} \eta_2\dist (x,
\cdot)^{-2} 
& \leq   \zeta  \, (r\ts_x^{-1})^{5/7}\delta _2^{4/7}+\zeta
_2\, \delta _2^{-7/2}\\
& \leq \zeta  \, (r\ts_x^{-1})^{5/7}\delta _1^{4/7}+\zeta
' r^{7/12}(\ln r)^{-21/8}. 
\end{split}
\end{equation}
Combining (\ref{est:qa}), (\ref{est:qb}) and make use of Propositions 
\ref{T:lem3.2} and \ref{T:prop3.1-} as in the proof of Proposition \ref{T:prop3.2}, we have: 
\[
|F_A^-|\leq \big(2^{-3/2}+(8\ln r)^{-1} \big)\, 
r \, (-u)_+ + K_1' \quad \text{over \(X^{'a}_{\delta '_0}\),}
\]
where \(K_1':=\zeta  _2\, (r\ts_x^{-1})^{5/7}\delta _1^{4/7}+\zeta_2
' r^{7/12}(\ln r)^{-21/8}+\zeta_1\ts^{-2}+\zeta _1' r^{1/3}(\ln r)\, \ts^{-1}.\)
(As observed previously, by Proposition
\ref{prop:curv-varpi} this inequality already hold over \(X_{\delta
  _1/2}\).) 
Next, re-run the proof of  Proposition
\ref{prop:curv-varpi},  replacing the use of Proposition
\ref{T:prop3.2} by the preceding inequality. In particular, the role
of \(K_1\) therein is now played by \(K_1'\). The
arguments in the proof then yield: 
\[
|F_A^-|\leq \big(2^{-3/2}+(8\ln r)^{-1} +\zeta _5 \varepsilon _1'\big)\, 
r \, (-u)_+ +\zeta_5'\ts^{-2} \quad \text{over \(X^{'a}_{\delta '_0}\)}, 
\]
where 
\[
\varepsilon '_1:=K_1'/\td{r}\leq \zeta  \, \td{r}^{-2/7}\ts^{-10/7}\delta
_1^{4/7}+\zeta _6(\ln r)^{-1}\quad \text{over \(X^{'a}_{\delta _0\ln r}\)}.
\] 
This leads to the asserted inequality
(\ref{eq:curv-varpi-3}). The second assertion of the lemma is obtained
by iterating the preceding arguments. Let \(\delta _{1:0}:=\delta _1\)
and \(\delta _{1:1}:=\delta _1'\). Suppose an inequality of the form
(\ref{eq:curv-varpi-3}) holds over \(X^{'a}_{\delta _{1:k}}\) and
\(\delta _{1:k}\geq r^{-1/4}\), then
the preceding arguments imply that it also holds over \(X^{'a}_{\delta
  _{1:k+1}}\) (but with different constants \(\zeta \), \(\zeta '\)),
where \(\delta_{1:k+1}=r^{-1/6}\delta_{1:k}^{1/3}(\ln
r)^{7/12}\). This implies the assertion (\ref{eq:curv-varpi-4}) by iteration. 
\epf 

An immediate consequence of the preceding lemma is: 
\begin{lemma}\label{lem:monotone00}
 Lemmas  \ref{lem:monotone0}, \ref{lem:u^a-int}, and their amended
 versions in Lemma
\ref{lem:amend-int} also hold with all appearances of \(\delta _1\)
therein replaced by \(\delta _1'\) or more
generally, \(\delta _{1,\epsilon }\) for any given
\(0<\epsilon <1/2\) (but with the various constant coefficients
\(\zeta _*\) changed). 
\end{lemma}
\pf Simply replace the use of Proposition  \ref{prop:curv-varpi} in the proof of Lemma
\ref{lem:monotone0} by the preceding lemma. \epf 

We also have the following variant of Lemma  \ref{lem:monotone0}. This
may be viewed as a Seiberg-Witten version of Taubes' monotonicity formula for 
pseudo-holomorphic curves near \(\omega^{-1}(0)\) (cf.  Section 2 of
\cite{T99}). 

\begin{lemma}\label{lem:SW-monotone}
Let \(\delta '_1\) \(\delta _{1, \epsilon }\) be as in Lemma 
\ref{lem:curv-varpi-3} and let \(\delta '_*:=\delta '_1\) or
\(\delta _{1, \epsilon }\).
There exist positive constants \(r_1\),  \(\zeta\), \(\zeta '\),
\(\zeta ''\)  that are independent of
\(r\), \(R\),  \(x\), and \((A, \Psi )\) with the following
significance: 
Given any \(r\geq r_1\),  \(x\in \nu ^{-1}(0)\cap X^{'a}\), \(0\leq
a\leq \frac{43}{8}\), 
and \(R\) such that \(1/4\geq R\geq 2\delta _*'\), one has: 
\[
\mathcal{W}_{B_x(R)}(A, \Psi)\leq
\zeta\,  R^3+\zeta '( \delta _*')^2. 
\]
In particular, when \(R\geq (\delta '_*)^{2/3}\), then \(\mathcal{W}_{B_x(R)}\, (A, \Psi)\leq
\zeta''\,  R^3\). 
\end{lemma}
\pf 
Fix \(x\in \nu ^{-1}(0)\cap X^{'a}\) and \(R\) with \(1/8\geq R\geq 2\delta '_*\).
Let \(\rho \) denote the function \(\dist (x, \cdot)\) on \(X\). 
 Write \(\omega=dq\) on small neighborhood of \(x\), where
\(q\) is of the form 
\begin{eqnarray}\label{eq:loc-omega}
q&=&  q_0+O(\sigma^2\rho),\quad \text{with}\nonumber\\
&q_0&= \frac{t}{2}d \, (x_1^2+x_2^2-2x_3^2)-(x_1^2+x_2^2-2x_3^2) \,
      dt-3x_3\, 
      (x_1 dx_2-x_2 dx_1)
\end{eqnarray}
in a certain local coordinate chart at \(x\), in terms of which
\(\nu^{-1}(0)\) is identified with the \(t\)-axis: \(\{(t, x_1, x_2, x_3)|\,
x_1=x_2=x_3=0\}\), and \(|q_0| =2^{-1/2}\rho|\omega|/3\). Let \(\rho
_0\) be the maximal value of \(R\) such that \(B_x(R)\) is contained in this
local coordinate chart, and suppose \(R\leq \rho _0\) from now on. (If
\(R\geq \rho _0\), the asserted inequality already hold
for certain \(\zeta \).) 

Recall again the notation \(-u=|\nu |-|\psi |^2\), and let
\[
  \scrW_{R}:=4^{-1}\int_{B_x(2R)}\chi (\rho /R)\, r |\nu| \,|-u|\geq
  \mathcal{W}_{B_x(R)}(A, \Psi ).
  \]
Use the Stokes' theorem, the Seiberg-Witten equation \(\grS_{\mu,
  \hat{\grp}}(A, \Psi)=0\), and Proposition \ref{T:lem3.2}  to get: 
 \begin{equation}\label{eq:W-ubdd-z}
\begin{split}
\scrW_{R}-\zeta  _1R^3& \leq  \frac{i}{2}\int_{B_x(2R)}\chi (\rho /R)\, F_A\wedge \omega   \leq
-\frac{1}{2R}\int_{B_x(2R)}\chi '(\rho /R)\, |F_A|\, |q|\\
& \leq  -\frac{1}{2R}\int_{B_x(2R)}\chi '(\rho /R)\, |F_A|\, |\nu |\, (\rho/3+\zeta_2\rho\sigma).
\end{split}
\end{equation}
The last line in the preceding expression is in turn bounded via Lemma
\ref{lem:curv-varpi-3}, Propositions \ref{prop:curv-varpi},
\ref{T:prop3.2}, \ref{T:lem3.2} and  \ref{T:prop3.1-} as: 
\begin{equation}\label{eq:DE-mono-z}
\begin{split}
& \frac{1}{2R}\int_{B_x(2R)}\chi '(\rho /R)\, |F_A|\, |\nu |\,
(\rho/3+\zeta_2\rho\sigma)\\
&\quad \leq \big(1+\zeta _8(\ln r)^{-1}\big) \frac{R+\zeta _4\, R^2}{3}
\frac{d\scrW_{R}}{dR}+\zeta _4'\, R^3\\
& \qquad \quad +\zeta _2\int_{B_x(2R)-B_x(R)}\chi (\sigma /\delta
_0')\, r\,\ts^2 \\
&\qquad \quad +\int_{B_x(2R)-B_x(R)}\chi (\sigma /\delta
_*') \big(1-\chi (\sigma /\delta
_0')\big) \big( \zeta _3 \ts^{-1}+\zeta _3' r\, |\nu |\, |-u|\big).
\end{split}
\end{equation}
Combining (\ref{eq:DE-mono-z}) and (\ref{eq:W-ubdd-z}), and applying the amended
version of Lemma \ref{lem:monotone0} in Lemma \ref{lem:monotone00} to
estimate the last integral above, one obtains: 
\[\begin{split}
\scrW_{R}& \leq  \big(1+\zeta _8(\ln r)^{-1}\big) \frac{R+\zeta _4\, R^2}{3}
\frac{d\scrW_{R}}{dR}  \\
& \qquad +\zeta  '_1R^3+\zeta _5 r(\delta _0')^5 R+\zeta
_5' (\delta _*')^2 R+\zeta _6 (\delta _*')^2.
\end{split}
\]
This in turn leads to: 
\[
\begin{split}
\frac{d}{dR}\left(e^{f_z}\scrW_{R}\right) & \geq
- \big(\zeta_7\,R^2+\zeta_7' \delta''^2/R\big)\, e^{f_z}, \qquad\text{when \(R\geq 2\delta '_*\), where}\\
f_z(R) & :=-3(1+\zeta _8(\ln r)^{-1})^{-1}\ln \, (\zeta _4R/(1+\zeta _4R)).
\end{split}
\]
Integrating, one finds that there exist positive constants \(\zeta
_w\), \(\zeta ''\) independent of \(r\), \(R\), \(s\), \(x\), and \((A,
\Psi )\),  such that for \(s>R\geq 2\delta '_*\), 
\[
\zeta _w\scrW_{s}/ s^{3/(1+z_r)}\geq \scrW_{
  R}/R^{3/(1+z_r)}-R^{3z_r/(1+z_r)}\big(\zeta_9 +\zeta_9' \delta''^2R^{-3}\big),
\]
where \(z_r:=\zeta _8(\ln r)^{-1}\).  
Now take \(s=\rho _0\), and appeal to the amended
version of Lemma \ref{lem:monotone0} in Lemma \ref{lem:monotone00}
again  to obtain the asserted inequality. 
\epf


We are now ready to state and prove the counterpart of Proposition 4.1 of \cite{Ts} and
Proposition I.3.1 of \cite{T} in our context. 

Let \(\ud{\alpha }=|\nu |^{-1/2}\alpha \), and recall the notations 
\(\td{r}=r\ts\); \(-u=|\nu |-|\psi |^2\). 

\begin{prop}\label{rem:monotone} {\rm (the monotonicity formula) }
 Adopt the assumptions and notation of Lemma \ref{lem:Etop-bdd1}. 
Let \(\delta '_1\) \(\delta _{1, \epsilon }\) be as in Lemma
\ref{lem:curv-varpi-3}, and let \(\delta _*:=(\delta _*')^{2/3}=(\delta '_1)^{2/3}\) or
\((\delta _{1, \epsilon })^{2/3}\). There exist positive constants
\(r_0\),  \(\zeta\) 
with the following
significance: Given any \(r\geq r_0\),  \(x\in X^{'a}_{\delta_*}\),
\(0\leq a\leq \frac{39}{8}\), and an 
\(R\in \bbR^+\) satisfying \((r\tilde{\sigma}(x))^{-1/2}\leq
R<\tilde{\sigma}(x)/4\), one has:  
\begin{equation}\label{W-u-bdd+}
\mathcal{W}_{B_x (R)}\, (A, \Psi)\leq \zeta
\tilde{\sigma}(x)\, R^2. 
\end{equation}
Fix a number \(v\),  \(0\leq v<1\). If in addition,
\(|\underline{\alpha}(x)|\leq v\), then there exist constants \(\zeta
_v\), \(\kappa _v<1/4\)  that are independent of
\(r\), \(x\), and \((A, \Psi )\) (but possibly depend on \(v\)), such that
for and \(R\) satisfying \((r\tilde{\sigma}(x))^{-1/2}\leq R<\kappa _v\, \tilde{\sigma}(x)\),
one has: 
\begin{equation}\label{W-l-bdd+}
\mathcal{W}_{B_x(R)}\, (A, \Psi) \geq \, \zeta_v  
\, \tilde{\sigma}(x)\, R^2. 
\end{equation}
In the above, the constants  \(\zeta \), \(\kappa _v\)
and \(\zeta _v\) (for each fixed \(v\)) depend  only on \(\nu \) and the
 parameters listed in (\ref{parameters}), and \(r_0\) depends on
 \(\smE\) in addition to all the above.  
\end{prop}
\pf 
Let  \(\delta '_*\) is as in the preceding lemma, and choose \(r\) to
be sufficiently large so that \(\delta _*>2\delta _*'\). Fix \(x\in
X_{\delta_*'}\). 
Let \(\rho\, (y) : =\dist ( x, y)\) as before. For simplicity, we frequently write
\(\mathcal{W}_{x, R}=\mathcal{W}_{B_x(R)}(A, \Psi)\) and \(\tilde{\sigma}_x=\tilde{\sigma}(x)\) below.
Let \(R<\tilde{\sigma}(x)/4\), so that \(B_x(2R)\subset
X^{'a+1/2}_{\delta_*/2}\subset X^{'a+1/2}_{\delta _*'}\), \(0\leq a\leq
\frac{39}{8}\), 
and \(\tilde{\sigma}(x)/2<\tilde{\sigma}(y)<3\tilde{\sigma}(x)/2\)
\(\forall y\in B_x(2R)\). 
As usual, use the Seiberg-Witten equation \(\grS_{\mu,
  \hat{\grp}}(A, \Psi)=0\) and Proposition \ref{T:lem3.2} to get:  
\begin{equation}\label{eq:W-ubdd}\begin{split}
\frac{i}{2}\int_{B_x(R)}F_A\wedge \omega   +\zeta _1'\ts^{-1} R^4&\geq 4^{-1}\int_{B_x(R)}r |\nu| \,
\Big||\nu |-|\psi|^2\Big|\\
&= \mathcal{W}_{x, R}.
\end{split}
\end{equation}
This time, write \(\omega=2\nu ^+=d\theta\) over \(B_x(2R)\),
where \(\theta =\frac{1}{2} (y_1dy_2-y_2dy_1+y_3dy_4-y_4dy_3)+O(\rho
^2)\) in a symplectic local coordinate system \((y_1, y_2, y_3, y_4)\)
where \(\omega=dy_1\wedge d y_2+dy_3\wedge dy_4 +O(\rho )\). 
In particular, it satisfies  
\(|\theta|\leq 2^{-3/2}|\omega|(\rho+O(\rho^2))\). Then the usual 
argument via  the Stokes'
theorem,  
Lemma \ref{lem:curv-varpi-3}, 
Propositions \ref{T:prop3.2} (a), \ref{T:lem3.2},  and 
\ref{T:prop3.1-} 
gives: 
\begin{equation}\label{eq:DE-mono}\begin{split}
& \frac{i}{2}\int_{B_x(R)}F_A\wedge \omega   
=\frac{i}{2}\int_{\partial B_x(R)}F_A\wedge \theta\\
&\quad \leq \frac{1}{2}\int_{\partial B_x(R)}
(|F_A^+|^2+|F_A^-|^2)^{1/2}|\nu|(\rho+\zeta\rho^2) \\
&\quad  \leq \frac{1}{4}\int_{\partial B_x(R)}\left(r \, (1+\zeta
  _1(\ln r)^{-1})\Big||\nu |-|\psi|^2\Big|
 +\zeta'\, \tilde{\sigma}^{-2}\right) |\nu|\, (\rho+\zeta\rho^2)/2 \\
& \qquad \qquad +\zeta_2\tilde{\sigma}_x^{-1}R^4\\
&\quad \leq 2^{-1}(1+\zeta _1(\ln r)^{-1})\, (R+\zeta R^2)
\frac{d\mathcal{W}_{x, R}}{dR}+\zeta_2\tilde{\sigma}_x^{-1}R^4\quad
\text{for \(x\in X^{'a}_{\delta _*}\).} 
\end{split}
\end{equation}
Combined with (\ref{eq:W-ubdd}), this leads to: 
\[
\begin{split}
\frac{d}{dR}\left(e^f\mathcal{W}_{x, R}\right)\geq
-\zeta'_2\tilde{\sigma}^{-1}_xR^3e^f, \qquad\text{for \(x\in X^{'a}_{\delta_*}\), where}\\
f(R):=-2\big(1+\zeta _1(\ln r)^{-1}\big)^{-1}\ln (R/(1+\zeta R)).
\end{split}
\]
Fix a positive constant \(c\leq 1\). Integrating the preceding
differential inequality, we have for  \(x\in X^{'a}_{\delta
    _*}\) and \(c\, (r\ts_x)^{-1/2}\leq R\leq s\leq\tilde{\sigma}_x/4\), 
\begin{equation}\label{eq:W-comparison}
\zeta _c\, \mathcal{W}_{x, s}/ s^2\geq \mathcal{W}_{x,
  R}/R^2-\zeta''\tilde{\sigma}_x^{-1}(s^2-R^2)\geq 0,  
\end{equation}
 where \(\zeta '', \zeta
_c>0\) are certain \(r\)-independent constants. (\(\zeta _c\) may depend on \(c\)). 
Now 
take \(s=\tilde{\sigma}_x/4\) in
the preceding inequality. The assertion (\ref{W-u-bdd+}) then
follows from  Lemma \ref{lem:SW-monotone}. 

To prove the second inequality claimed, i.e. (\ref{W-l-bdd+}), suppose that \(x\) satisfies \(|\underline{\alpha}(x)|<v\). Take
\(R=c\, (r\tilde{\sigma}_x)^{-1/2}\) in (\ref{eq:W-comparison}), one has:
\begin{equation}\label{eq:W-lower}
\zeta_c\mathcal{W}_{x, s}/s^2\geq c^{-2}r\ts_x\mathcal{W}_{x, c\, (r\ts_x)^{-1/2}}-\zeta''\tilde{\sigma}_x^{-1}(s^2-c^2(r\ts_x)^{-1}) .
\end{equation}
Repeat the argument in Step 3 of the proof of Proposition I.3.1 in \cite{T} with the following modifications:
\begin{itemize}
\item[(1)]  Replace the use of Proposition I.2.8 in \cite{T} by its
analog in our setting, Proposition \ref{lem:est-1st-der};
\item[(2)] Replace \(r\) by \(r\ts_x\),
\item[(3)] Replace \(\alpha\) by \(\underline{\alpha}\).
\item[(4)] Replace the condition \(\ud{\alpha }<1/2\) therein by
  \(\ud{\alpha }<v\). 
\item[(5)]  
Proposition \ref{lem:est-1st-der} 
implies that there is a constant \(c<1\) independent of \(r\) and \(v\),
so that \(\ud{\alpha
}(y)<1-(1-v)/2\) when \(\dist (y, x)<c \, (r\ts_x)^{-1/2}\). 
\item[(6)] Consequently, the counterpart of \cite{T}'s (I.3.11) in our
  context says: 
\[
\mathcal{W}_{x, c\, (r\ts_x)^{-1/2}}\geq \zeta _6
\big(1-(1+v)^2/4\big)\, r\tilde{\sigma}^2_x c^4 (r\ts_x)^{-2}. 
\]
\end{itemize}
Inserting the preceding inequality into (\ref{eq:W-comparison}), setting the constant \(c\) therein to be the constant
\(c\) in Items (5) and (6) above. We get: 
the first term on the right hand side of (\ref{eq:W-lower})
has a lower bound by a constant \(\zeta'\); therefore
\[
\mathcal{W}_{x, s}\geq \zeta _c^{-1}s^2\Big(\zeta_6\big(1-(1+v)^2/4\big)c^2\ts_x-\zeta''\tilde{\sigma}_x^{-1}s^2\Big).
\] 
This implies that the second inequality asserted, (\ref{W-l-bdd+}),
with \[
\zeta _v=\zeta _c^{-1}\zeta_6\, c^2\big(1-(1+v)^2/4\big)/2,\quad 
\kappa_v=\min \, \big(1/4, (\zeta _v/\zeta '')^{1/2}\big).\] 
\epf

The  monotonicity formula is used primarily by way of the next corollary. 

Given a set \(\Lambda ^*=\{B_k\}_k\) consisting of balls \(B_k\subset
X\), and a subset \(U\subset X\), we use \(\Lambda ^*_U\subset \Lambda
^*\) to denote the subset \(\Lambda ^*_U=\{B_k\, |\, B_k\cap U\neq
\emptyset\}_k\). Given \(v\), \(0<v<1\), let  \(X_{\delta, v}^{'a}:=\{x\,
|\, |\ud{\alpha}|(x)\leq v, \, x\in
  X_\delta^{'a}\}\). 
\begin{cor}\label{cor:ball-covering}
 Adopt the assumptions and notation of Lemma \ref{lem:Etop-bdd1}, 
 and let \(\delta_*\), \(\kappa _v\) be as in
Proposition \ref{rem:monotone}. 
There exist
constants \(r_0\), \(\zeta _0\), \(\zeta \), \(\zeta '\)  that are independent of \(r\),
\(\rho \),  \(\delta \),  \(X_\bullet\), \(\grB_R\), 
and \((A, \Psi )\),  such that the following hold for any given \(r\geq
r_0\), \(1>\delta \geq \delta _*\), 
\begin{itemize}\itemsep -2pt
\item Fix \(\rho\in [r^{-1/2},
\kappa _0\, \delta^{3/2}]\). 
Let \(\Lambda=\{B(x_k, \rho _k)\}_k\) be a set of mutually disjoint
balls in \(X_\delta ^{'a}\), \(0\leq a\leq \frac{39}{8}\), with 
\(\rho _k:=\tilde{\sigma}(x_k)^{-1/2}\rho\) and  \(x_k\in
\alpha^{-1}(0)\). Then given a compact \(X_\bullet \subset X^{'a}\), the
subset \(\Lambda_{X_\bullet}\subset \Lambda \) has no more than
\(\zeta _0\, \rho^{-2}(1+|X_\bullet|)\) elements.
\item Let \(\rho \) be as in the previous bullet. The set \(\alpha^{-1}(0)\, \cap X^{'a}_\delta\)
 can be covered by a set of balls, \(\Lambda^\rho=\{B \, (x_k,\rho
 _k)\}_k \),  where \(\rho _k=\tilde{\sigma}(x_k)^{-1/2}\rho\) and 
 \(x_k\in \alpha^{-1}(0)\cap X^{'a}_{\delta/2}\), such that the concentric
 balls \(B (x_k, \rho_k/2)\) are mutually disjoint. Any such set
 \(\Lambda ^\rho \)  has the following
 properties: 
\begin{itemize}
\item[(i)] Given any compact \(X_\bullet \subset X^{'a}\), the set \(\Lambda^\rho_{X_\bullet}\)  consists of no
  more than \(\zeta  \, \rho^{-2}(1+|X_\bullet|)\) elements.
\item[(ii)] Given any ball \(\grB_R=B\, (x, R)\subset X^{'a}_\delta \) of radius
  \(R\), with \\ \((r\tilde{\sigma}(x))^{-1/2}\leq
R<\tilde{\sigma}(x)/4\),  the set \(\Lambda_{\grB_R}^{\rho}\)
consists of no more than \(\zeta '\ts (x)\, \rho^{-2}R^2\) elements.
\end{itemize}
\item  Fix \(v\), \(0<v<1\) and  \(\rho\in [r^{-1/2},
\kappa _v\delta^{3/2}]\). Then there exists positive constants
  \(c_v\), \(c_v'\) that are independent of \(r\), \(\delta
  \), \(\rho \), \((A, \Psi )\), \(X_\bullet\) and \(\grB_R\) (but possibly on \(v\)), such that
  the following holds: 
  The set \(X^{'a}_{\delta, v}\) has a cover by a set of balls,
  \(\Lambda^{v, \rho}=\{B (x_k, \rho _k)\}_k\), where \(\rho
  _k=\tilde{\sigma}(x_k)^{-1/2}\rho\), and the concentric balls
  \(B(x_k, \rho _k/2)\) are mutually disjoint.
  Any such set \(\Lambda^{v, \rho}\) has the following properties:
\begin{itemize}
\item[(i)] Given any compact \(X_\bullet \subset X^{'a}\), the set \(\Lambda_{X_\bullet}^{v, \rho}\) consists of no
  more than \(c_v\, \rho^{-2}(1+|X_\bullet|)\) elements.
\item[(ii)] Given any ball \(\grB_R=B\, (x, R)\subset X^{'a}_\delta \) of radius
  \(R\), with \\ \((r\tilde{\sigma}(x))^{-1/2}\leq
R<\tilde{\sigma}(x)/4\),  the set \(\Lambda_{\grB_R}^{v, \rho}\)
consists of no more than \(c_v'\, \ts (x)\, \rho^{-2}R^2\) elements.
\end{itemize}
\end{itemize}
In the above, the constants  \(\zeta \), \(\zeta '\)
 \(\zeta _0\), and \(c_v\), \(c'_v\) (for each fixed \(v\)) depend  only on \(\nu \) and the
 parameters listed in (\ref{parameters}), and \(r_0\) depends on
 \(\smE\) in addition to all the above.  
\end{cor}
\pf This is an 
analog of Lemmas I.3.6 and I.3.8 in \cite{T} and follows from the same
argument with the use of Proposition I.3.1 by its counterpart in our
context, Proposition 
 \ref{rem:monotone} above. 
\epf

\subsection{(Further) refined curvature estimate}

The next proposition is the counterpart of Proposition I.3.4 
in \cite{T} (and also \cite{Ts}'s Proposition 4.2) in our context. 
It improves the estimate for \(|F_A^-|\)
given in Proposition \ref{prop:curv-varpi} (and Lemma
\ref{lem:curv-varpi-3}) over \(X^{'a}_{\delta _*}\) by eliminating the coefficient \(\varepsilon '\) therein, thereby making the first
terms in the RHS of (\ref{eq:curv-varpi}) and (\ref{eq:F+bdd})
coincide.  

The aforementioned improvement is made possible by the following
consequence of the monotonicity formula. This is in turn the
counterpart of \cite{T}'s Lemma I.3.5 in our context.  
\begin{lemma}\label{T:lem4.2}
Adopt the assumptions and notation of Lemma \ref{lem:Etop-bdd1}, 
and let \(\delta_*\) be as in Proposition \ref{rem:monotone}. 
Then there exist positive constants \(r_0\), \(\zeta \), \(\zeta '\)
with the following significance: for \(r\geq r_0\), there is a smooth function \(\gru\) on \(X^{'a}_{\delta _*}\),
\(0\leq a\leq \frac{39}{8}\), 
satisfying the following properties:
\begin{itemize}
\item \(d^*d\gru\geq \zeta r\ts \) where \(|\ud{\alpha}|\leq
  2^{-1}\);
\item \(\Big|d^*d\gru\Big|\leq \zeta 'r\ts \);
\item\(|\gru|\leq \zeta'\). 
\end{itemize}
In the above, the constants  \(\zeta \), \(\zeta '\) depend  only on \(\nu \) and the
 parameters listed in (\ref{parameters}), and \(r_0\) depends on
 \(\smE\) in addition to all the above.  
\end{lemma}
\pf Modify 
the proof of Lemma I.3.5 in  \cite{T} as follows.
\begin{itemize}
\item Replace the parameter \(r\)
in \cite{T} by \(\td{r}=r\ts \), and replace balls of the form \(B\, 
(x,\rho_*)\) in \cite{T} by their counterparts, 
\(B\, (x, \ts^{-1/2}(x) \, \rho_*)\), in our context. 
\item Replace all appearances of \(\alpha \) in \cite{T}'s proof of
  its Lemma I.3.5  by \(\ud{\alpha }\). Replace the use of Lemmas I.3.6, I.3.8 in \cite{T} by Corollary
  \ref{cor:ball-covering} above, taking \(\delta =\delta _*\). 
\item In Step 2 of the proof of Lemma I.3.5 in  \cite{T}, replace the
set of balls \(\{B_i\}_i\) therein by the set of balls 
\(\Lambda ^{v, \rho }\) in the third
bullet of our Corollary
  \ref{cor:ball-covering}; with \(v=1/2\), and \(\rho =r^{-1/2}\) as in
  \cite{T}. 
\item  In Step 3 of the proof of Lemma I.3.5 in  \cite{T}, a function
  \(u_i\) is introduced for each element \(B_i\) in \(\{B_i\}_i\). As
  explained in the previous bullet, The
  counterpart of \cite{T}'s \(B_i\) is a typical element, \(B\, (x_k,
  \rho _k)\), where \(\rho _k=r^{-1/2}\ts (x_k)^{-1/2}\),  in
  \(\Lambda ^{1/2, r^{-1/2}}\).  To this we associate a function
  \(\gru_k\) on \(X\), which plays the part of \cite{T}'s \(u_i\): In
  parallel to the defining equation of \cite{T}'s \(u_i\) in (I.3.32)
  of \cite{T}, \(\gru_k\) is defined to be the solution to: 
\[\begin{cases}
 \, \, \frac{1}{2} d^*d \gru_k=
r \ts(x_k)\, (s_k-\kappa _k\, \grs_k); & \\
\,\, \text{\(\forall i\in \grY\), \(\gru_k\) has \(Y_i\)-end limit \(0\).} &
\end{cases}
\]
where \(s_i\), \(\grs_i\), \(\kappa _i\) are defined as in \cite{T},
modified according to the rules set forth by the first two bullets
above. The argument in \cite{T} leading to its (I.3.36) shows that
\(\gru_k\) obeys: 
\begin{equation}\label{bdd:rmu_k}
|\gru_k|+|d\gru_k|\leq \zeta _4 \, (r\ts (x_k))^{-3/2} \dist (\cdot,
x_k)^{-3}
\end{equation}
where \(\dist (\cdot x_k)\geq \zeta ''\rho _k\). 
\item Take \[\gru=\sum_{k\in \Lambda ^{1/2, r^{-1/2}}}\chi \big(8\dist
  (\cdot, x_k)/\ts (x_k)\big)\, \gru_k; \qquad \gru':=\sum_{k\in \Lambda
    ^{1/2, r^{-1/2}}}\, \gru_k. 
\]  
The arguments in Step 5 of \cite{T}'s proof of its Lemma I.3.5 (again
  modified according to the previously described rules) then confirm
  the bound on \(\gru\) asserted in the last bullet of the lemma, and
  also that: 
\begin{equation}\label{bdd:ddu}
\begin{split}
& \text{\(d^*d\, \gru'\geq \zeta  r\ts \) where \(|\ud{\alpha}|\leq
  2^{-1}\);} \qquad 
  \Big|d^*d\, \gru'\Big|\leq \zeta ' r\ts.
\end{split}
\end{equation}
Note that given \(x\) and \(k\), the value of \(\Big|d^*d\Big(\chi \big(8\dist
  (\cdot, x_k)/\ts (x_k)\big)\, \gru_k-\gru_k\Big)\Big|\) at \(x\)
  vanishes unless \(k\) belongs to the subset
  \[
    \Lambda _x:=\{k\, |\,
  x_k\in B_x(x, 3\ts(x))-B(x, 2\ts(x)/5),\, k\in \Lambda ^{1/2,
    r^{-1/2}}\}.
\]
When \(k\in \Lambda _x\), (\ref{bdd:rmu_k}) may be used to bound this value by \\ \(\zeta _4' \,
  (r\ts (x))^{-3/2} \ts (x)^{-5}\). Meanwhile, by Corollary
  \ref{cor:ball-covering} \(\Lambda _x\) has at most \(\zeta _5\,r
  \ts(x)^3\) elements. Together, these imply that  
   \(\big|\, d^*d\, \gru-d^*d\, \gru'\big|\, (x)\leq \zeta _6 \, (r\delta
   _*^3)^{-3/2} r\ts (x)\). Combined with (\ref{bdd:ddu}), this
   confirms the first two bullets of the lemma. \epf
\end{itemize}

Here is the promised improved curvature estimate. 
\begin{prop}\label{prop:F-refined}
Re-introduce the notations \(\varpi =|\nu |-|\alpha |^2\) and
\(-u=|\nu |-|\psi |^2\). 
Adopt the assumptions and notation of Lemma \ref{lem:Etop-bdd1}, 
and let  \(\delta_*\)  be as in
Proposition \ref{rem:monotone}.  Then there exist \(r\)-independent constants \(r_0\),
\(\zeta \), \(\zeta '\) such that for all \(r\geq r_0\), 
\begin{equation}\label{eq:F-refined}
\begin{split}
|F_A^-| & \leq 2^{-3/2}r\varpi+\zeta \ts^{-2};\\
|F_A^-| & \leq 2^{-3/2}r\, (-u)+\zeta '\ts^{-2}\quad \text{over
  \(X^{'a}_{8\delta_*}\), \(0\leq a\leq \frac{15}{4}\)}. 
\end{split}
\end{equation}
In the above, the constants  \(\zeta \), \(\zeta '\) depend  only on \(\nu \) and the
 parameters listed in (\ref{parameters}), and \(r_0\) depends on
 \(\smE\) in addition to all the above.  
\end{prop}
\pf Note that the two asserted inequalities are equivalent over
\(X^{'a}_{\delta _0}\) by way of Proposition \ref{T:prop3.1-}; so it
suffices to establish either of them. 

To begin, recall (\ref{ineq:beta2}), and this time take \(\varepsilon
=4\zeta _3^{-2}\ts^{-2} r^{-1}|\nu |^{-1}\leq \zeta ''r^{-1}\ts^{-3}\). Fix
\(\delta \geq\delta _0'\), where  \(\delta _0'\) is as in Proposition
\ref{prop:curv-varpi}, and let
\(\varepsilon _\delta :=\zeta '' r^{-1}\delta ^{-3}\). Introduce
the following variant of the function \(q_0\) in (\ref{def:q0}): 
\[
q_0=q_{0, \delta }:=\ss +(2^{-3/2}+\varepsilon _\delta ) \, r\, u+\zeta
_2r|\beta |^2. 
\]
By Propositions
\ref{prop:curv-varpi}, \ref{T:lem3.2}, \ref{T:prop3.1-}, it satisfies the following variant of
(\ref{DE:q_0-1}) and (\ref{bdd:q_0b}): 
\[
\begin{split}
  \big(\frac{ d^*d}{2}+\frac{r|\psi|^2}{4}\big) q_{0, \delta }  
 & \leq \zeta  _0r\, \ts^{-1}; \\
q_{0,\delta } & \leq \zeta  (\ln r)^{-1} r\, (-u)_++\zeta ' \ts^{-2}\quad \text{over \(X^{'a}_\delta \)}. 
\end{split}
\]
Combine this and Lemma
\ref{T:lem4.2} to find a constant \(\zeta_1
\) that is independent of \(r\), \(\delta \), and \((A, \Psi )\), such
that with \(\grq_{0,\delta }:=q_{0,\delta
 }-\zeta _1\delta ^{-2}\gru\), 
one has for \(\delta \geq \delta _*\): 
\[
\begin{split}
 & \big(\frac{ d^*d}{2}+\frac{r|\psi|^2}{4}\big) \, \grq_{0,\delta}\leq \zeta  _0'
 r\, \delta ^{-2}|\psi |^{2}; \\
& \qquad \grq_{0,\delta }\leq \zeta  \, (\ln
r)^{-1} r\, (-u)_++\zeta '' \delta ^{-2}\quad \text{over \(X^{'a}_\delta
  \), \(0\leq a\leq \frac{39}{8}\). }\\
 \end{split}
\]
Together with (\ref{bdd:q_0i}) and the last bullet of Lemma
\ref{T:lem4.2}, 
this in turn leads to the existence of a constant \(\zeta _1'\) (again 
independent of \(r\), \(\delta \), and \((A, \Psi )\)), such that
the function \(\grq_\delta :=\grq_{0, \delta }-\zeta _1'\delta ^{-2}\)
on \(X^{'a}_\delta \) satisfies: 
\[
\begin{split}
 & \big(\frac{ d^*d}{2}+\frac{r|\psi|^2}{4}\big) \, \grq_{\delta}\leq 0; \\
&  \grq_{\delta } \leq \zeta  \, (\ln
r)^{-1} r\, (-u)_+\quad \text{over \(X^{'a}_\delta
  \), \(0\leq a\leq \frac{39}{8}\). }\\
& \text{\((\grq_\delta )_+\) is supported on \(U\subset V\),}
 \end{split}
\]
where \(U\), \(V\) are compact spaces defined as in the proof of
\ref{prop:curv-varpi}. 
Now repeat the remaining argument in the proof of Proposition
\ref{prop:curv-varpi}, with \(\grq_\delta \) playing the role of \(q'\)
therein. To give more details: We now have
\[
\gamma _{a, \delta } \grq_\delta =q_1+q_2, \quad 0\leq a\leq \frac{31}{8},
\]
where \(\gamma _{a, \delta } \) is as in  (\ref{DE:q'_d}), and
\(q_1\), \(q_2\) respectively satisfy:
\[
\big(\frac{ d^*d}{2}+\frac{r|\nu |}{4}\big) \, q_1
=\xi'_\delta \, \, \, \text{ over \(V\)};\qquad 
q_1|_{\partial V}=0,
\]
where \[
 \xi'_\delta =
\frac{ d^*d}{2} (\gamma _{a, \delta } \grq_\delta )-\gamma _{a,\delta }
 \frac{ d^*d\grq_\delta }{2}; 
\]
and 
\begin{equation}\label{def:q_2new}
\begin{split}
 & \big(\frac{ d^*d}{2}+\frac{r|\psi|^2}{4}\big) \, q_2=-\frac{ru}{4}
q_1\\
& \qquad q_2|_{\partial V} \leq 0. \\
 \end{split}
\end{equation}
The present version of \(q_1\) is bounded by 
\[
\begin{split}
|q_1|& \leq  \zeta _2'(\ln r)^{-1} r\, (-u)_+ (r\delta
)^{-1/2} e^{-\zeta _g(r\delta )^{1/2}\dist (\cdot,\, X^{' a+1}_\delta
  -X^{'a}_\delta  )} \\
 &\qquad \qquad 
+\zeta _2 (\ln
r)^{-1} r\, (-u)_+ (r\delta ^{3})^{-1/2}
e^{-\zeta _g(r\delta )^{1/2}\dist (\cdot,\, Z^{'a}_{2\delta }-Z^{'a}_\delta) }. \\
\end{split}
\]
Plugging this back in (\ref{def:q_2new}), we have: 
\[
\begin{split}
 & \frac{ d^*d}{2} \, (q_2)_+\leq  \zeta _3'(\ln r)^{-1} (r\ts)^2 (r\delta
)^{-1/2} e^{-\zeta _g(r\delta )^{1/2}\dist (\cdot,\, X^{' a+1}_\delta
  -X^{'a}_\delta  )} \\
 &\qquad \qquad +\zeta _3' (\ln
r)^{-1} (r\ts)^2 (r\delta ^{3})^{-1/2}
e^{-\zeta _g(r\delta )^{1/2}\dist (\cdot,\, Z^{'a}_{2\delta }-Z^{'a}_\delta) }\\
\\
& \qquad (q_2)_+|_{\partial V} =0, \qquad 0\leq a\leq \frac{31}{8}. \\
 \end{split}
\]
Noting again that the Green's function \(G\) for the differential
operator \( \frac{ d^*d}{2}\) with
Dirichlet boundary condition on \(V\) satisfies: 
\[
|G(x, \cdot)|+\dist(x,\cdot)\, |d G (x, \cdot)|\leq
\zeta\dist(x,\cdot)^{-2}, 
\] 
this leads to the following bound for \(q_2\) over \(X^{'a}_{4\delta
}\subset V\), \(0\leq a\leq \frac{15}{4}\): 
\[
q_2\leq \zeta _4\, \ts^{-2} (\ln
r)^{-1} r^{-1}\delta ^{-4}. 
\]
Combing this with the previously obtained bound for \(q_1\) and the last bullet of Lemma
\ref{T:lem4.2}, we have: 
\[
|F_A| \leq 2^{-3/2}\, r\, (-u)_++\zeta
_4'\delta ^{-2}\quad \text{\(X^{'a}_{4\delta
}\), \(0\leq a\leq \frac{15}{4}\)},
\]
when \(\delta \geq\delta _*\), where \(\zeta _4'\) is independent of
\(r\), \(\delta \), and \((A, \Psi )\). Given a fixed \(x\in X^{'a}_{8\delta
}\), we take \(\delta =\ts(x)/2\) in the preceding formula. This leads
to the second line of (\ref{eq:F-refined}), which in turn implies the
first line. \epf

\begin{rem}
The preceding improvement of Proposition \ref{prop:curv-varpi} and and Lemma
\ref{lem:curv-varpi-3}) is not strictly necessary for the
proof of our main theorems. (Proposition \ref{prop:curv-varpi} and/or Lemma
\ref{lem:curv-varpi-3} suffice for our
purposes). Similarly, Proposition I.3.4 in \cite{T} improves the
curvature estimate in Proposition I.2.4 therein, but the proof of the
main Theorem of part I of \cite{T} uses only the
latter. Nevertheless, the more precise estimate in 
\cite{T}'s Proposition  I.3.4 plays a role in part IV of \cite{T}. As
the present article is the counterpart of part I of \cite{T},
Proposition \ref{prop:F-refined} should be useful for further studies
on the relationship between the Seiberg-Witten-Floer theory and
ECH-type Floer theories.
 \end{rem}

\subsection{Approximation by local models}\label{sec:local_m}

Proposition I.4.2 of \cite{T} (and correspondingly, Proposition 5.2 in
\cite{Ts}) has the following counterpart in our context.

Let \(\delta _*\) be as in be as in Proposition \ref{rem:monotone},
and let \(x\in X^{'a}_{2\delta _*}\). Re-introduce the notation \(\ts_x=\ts
(x)\) and \(\td{r}=r\ts\); \(\td{r}_x=r\ts_x\).  Given an \( (A,
\Psi)=(A_r, \Psi _r)=(A, (\alpha , \beta ))\), let \(\underline{A}_x,
(\underline{\alpha}_x, \underline{\beta}_x))\) be the rescaled version
of \((A, (\alpha , \beta ))\) defined in \cite{T}'s (I.4.7) (cf. also
Equation (5.1) of 
\cite{Ts}), with the parameter \(\lambda \) therein chosen to be
\((\td{r}_x)^{1/2}\) instead. (\(\underline{A}_x,
(\underline{\alpha}_x, \underline{\beta}_x))\) was called \((\underline{A},
 (\underline{\alpha}, \underline{\beta}))\) in \cite{T}. 
\begin{lemma}\label{lem:loc-mod}
Adopt the assumptions and notation of Lemma \ref{lem:Etop-bdd1}. There is an \(r_0>8\) depending on the same parameters listed in
Proposition I.4.2 of \cite{T} as well as \(\smE\), such that for all \(r\geq r_0\), \(x\in
X^{'a}_{2\delta _*}\), \(0\leq a\leq \frac{39}{8}\), and \( (A, (\alpha , \beta ))=(A, \Psi)=(A_r, \Psi _r)\), 
 there is a Seiberg-Witten solution on \(\bbR^4\),
\((A_0, (\alpha_0,0))\), as described in 
Proposition I.4.1 of \cite{T}, which approximates \(\underline{A}_x,
(\underline{\alpha}_x, \underline{\beta}_x))\) in the sense that it satisfies Items (1)--(4) of
Proposition I.4.2 in \cite{T}. 
\end{lemma}

\pf The proof of Proposition I.4.2 in \cite{T} can be copied, with the
following modifications:
\begin{itemize}  
\item As usual, the factors of \(r\) in \cite{T} are replaced by
  \(\td{r}_x\). 
\item Equation (I.4.8) of \cite{T} has the following counterpart in
  our setting:
\[\begin{split}
& |\ud{\alpha }_x| +|F_{\ud{A}_x}|\leq \zeta _1\\
& |\nabla_{\ud{A}_x}\ud{\alpha }_x|^2\leq \zeta  \, (1- |\ud{\alpha }_x|^2
)_++\zeta ' \ts_x^{-2}\td{r}^{-1}\\
& |\ud{\beta }_x|+|\nabla_{\ud{A}_x}\ud{\beta }_x|\leq \zeta _1'
\ts_x^{-1}\td{r}^{-1/2}, 
  \end{split}
\]
where \(\zeta \), \(\zeta '\), \(\zeta _1\), \(\zeta _1'\) are
independent of \(r\), \(x\), and \((A, \Psi )\). Note that
\(\ts_x^{-1}\td{r}^{-1/2}\leq \zeta _2r^{-5/36}\to 0\) as \(r\to
\infty\). These inequalities follow from Propositions \ref{T:lem3.2},
\ref{T:prop3.1-}, \ref{prop:curv-varpi}, and  \ref{lem:est-1st-der}. 
\item Proposition \ref{rem:monotone} is used in place of \cite{T}'s
  Proposition I.3.1. Proposition \ref{prop:curv-varpi} (or its
  improved version, Proposition \ref{prop:F-refined}) is used in place  of \cite{T}'s
  Proposition I.3.4.
\item The necessary integral bound on \(r\, (|\nu |-|\alpha |^2)\) (which, by Proposition \ref{T:lem3.2},
\ref{T:prop3.1-}, is equivalent to an integral bound on \(r\, (|\nu
|-|\psi |^2)\))  is supplied by Lemma  \ref{lem:u^a-int}
and its amended versions in Lemmas
\ref{lem:amend-int} and \ref{lem:monotone00}. The necessary integral bound on
\(|F_A^-|^2-|F_A^+|^2\) is supplied by 
Proposition \ref{prop:SW-L2-bdd} and its amendment in    Proposition
\ref{prop:integral-est2}. 
\end{itemize}
\epf

The preceding lemma leads to the following variant of
Lemma I.4.5 in \cite{T}: 
\begin{lemma}\label{lem:T1.4.5}
Adopt the assumptions and notation of Lemma \ref{lem:Etop-bdd1}, and
let \(\delta _*\) be as in be as in Proposition \ref{rem:monotone}. 
Then  the statement of Lemma I.4.5 in \cite{T} holds for \((A, \Psi)\) over \(X^{'a}_{2\delta
  _*}\), \(0\leq a\leq \frac{39}{8}\), with \(\alpha \) and \(r\) therein respectively replaced by
\(\ud{\alpha }\) and  \(\td{r}\). In other words, given \(0<v<1\), there
is a constant \(z_v\) independent of \(r\), \(x\), and \((A, \Psi
)\) (but depends on \(v\)), such that 
\[
\dist \, (\cdot, \alpha ^{-1}(0))\leq z_v\, \td{r}^{-1/2} \quad
\text{over \(X_{2\delta _*, v}\)},
\]
where \(X_{\delta  , v}\) is as in Corollary
\ref{cor:ball-covering}. Namely, \(X_{\delta, v}^{'a}:=\{x\,
|\, |\ud{\alpha}|(x)\leq v, \, x\in
  X_\delta^{'a}\}\). 
\end{lemma}
\pf The proof identical to the proof of Lemma I.4.5 in \cite{T}, with
the role of \cite{T}'s Proposition I.4.2 therein replaced by the
preceding lemma. 
\epf

\subsection{Exponential decay}\label{sec:exp}

The purpose of this subsection is to establish the counterpart of
Lemma 6.1 in \cite{Ts} (which is a variant of 
Proposition I.4.4 in \cite{T}) in our context, the upcoming Proposition
\ref{prop:exp-decay}. The proof of this proposition makes use
of a corollary of the monotonicity formula, Lemma \ref{lem:exp-decay}
below, 
which is in turn the counterpart of  Lemma I.4.6 in \cite{T}.

Before embarking on the aforementioned tasks, we introduce a notion of
``\(Y_i\)-end limits'' for functions over (subdomains of)
\(X^{'a}_\delta \) (or more generally, sections of bundles over
\(X^{'a}_\delta \)). The definition is parallel to the notion of
\(Y_i\)-end limits for functions/sections over (subdomains of) \(X\),
previously introduced in Section \ref{sec:convention} (6).

Given \(\delta >0\) and \(i\in \grY_m\), let \(\sigma _i:=\dist
(\cdot,\nu _i^{-1}(0))\), and \(Y_{i,
  \delta }:=\{y\in Y_i\, |\, \sigma _i\geq\delta \}\). In what
follows, \(Y_{i,
  \delta }\) is often identified with \(Y_{i: \infty, \delta
}:=Y_{i:\infty}\cap X_\delta \). 
Observe that the admissibility of \(\nu \) implies the
following: Let \(g\) denote the metric on \(X\) and recall \(\gri\co
\coprod_{i\in \grY}[0, \infty)\times Y_i\to X-X_c^\circ\), the
isometry from Definition \ref{def:mce}. Let \(\bfY_{i, \delta ,
  L}:=\gri \big([L, \infty)\times Y_{i, \delta }\big)\) when
\(i\in \grY_m\). 
Given \(i\in \grY_m\) and a small
\(\delta<1/8\),  there is an \(L=L_\delta >8\) and a diffeomorphism 
\(
\gri _\delta \co [L, \infty)\times Y_{i, \delta }\to \hat{Y}_{i, L, \delta }:=X_\delta \, \cap \hat{Y}_{i,
  L}\), such that: 
\[\begin{split}
& \qquad \bfY_{i, 2\delta , L}\subset 
\hat{Y}_{i, L, \delta } \subset \bfY_{i,  \delta /2, L}; \\
& \qquad \gri_\delta \big|_{\bfY_{i,
    2\delta, L}}=\gri\big|_{\bfY_{i, 2\delta, L}};\quad \mathfrc{s}_i=\pi
_\bbR\circ\gri_\delta , \\
& \|\gri_\delta ^*\, g-\gri^*g\|_{C^2([l, \infty)\times Y_{i, \delta}
  )}+\|\gri_\delta ^*\, \nu -\pi _Y^*\, \nu _\infty\|_{C^2([l, \infty)\times
  Y_{i,\delta })}\leq \zeta  e^{-\kappa  _i\, l/2} \, \,  \forall l\geq
L, 
\end{split}
\]
where \(\kappa _i\) is the constants  from 
(\ref{eq:xi-exp}),  \(\pi_\bbR, \pi _Y\) respectively denote
the projection to the first factor and the second factor of the
product \(\bbR\times Y\), and  \(\zeta \) is a constant depending only on the metric and
\(\nu \). 
Let  \(\rmV\)  be a bundle over
\(\ov{X_\delta ''}\) (or its subdomains such as \(\ov{X^{'a}_\delta }\)
or \(\hat{Y}_{i, L}\cap  \ov{X}_\delta \)), and let 
\(\rmV_i:=\rmV\big|_{Y_{i: \infty, \delta }}\). 
A section \(q\in \Gamma
(\rmV\big|_{X''_\delta })\) is said to have \(q_i\in \Gamma (\rmV_{i})\) as a {\em
  \(Y_i\)-end limit} if  \(q\) extends to be defined over
\(\ov{X_\delta ''}\) with (the extended) \(q\big|_{Y_{i: \infty, \delta
  }}=q_i\); and 
\[
\|\gri_\delta ^*\, q-\pi _Y^*\, q_{i}\|_{C^2([l, \infty)\times
  Y_{i,\delta })}+\|\gri_\delta ^*\, q-\pi _Y^*q_{i}\|_{L^2_1([l, \infty)\times
  Y_{i,\delta })}\to 0 \quad\text{as \(l\to \infty\).}
\]

Let \(\rmd_\alpha \, (\cdot):=\dist (\cdot, \alpha^{-1}(0))\). When
\(\alpha ^{-1}(0)=\emptyset\), we set \(\rmd_\alpha \, (\cdot)=\dist
(\cdot, \alpha^{-1}(0))\equiv \infty\).  Recall the
notation \(\td{r}=r\ts\). Let \(\td{r}_x:=\td{r}(x)\) and \(\ts_x=\ts
(x)\). 

\begin{lemma}\label{lem:exp-decay}
Adopt the assumptions and notation of Lemma \ref{lem:Etop-bdd1}, and
let \(\delta _*\) be as in be as in Proposition \ref{rem:monotone}. 
There exist positive constants \(r_0>8\), \(c\), \(c'\geq 2^5\, c\),
\(\zeta _0\), \(\zeta , \zeta '\), 
that are independent of
\(r\) and \((A, \Psi )\), such that for any \(r\geq r_0\), there exists a function
\(h_r\) over \(X^{'a}_{\delta _*}\), \(0\leq a\leq \frac{39}{8}\), satisfying  
\begin{itemize}
\item[(1)] \(d^*d \, h_r+ \frac{r}{16} |\nu | h_r\geq -\zeta _0 \exp \big(-(r\ts^3)^{1/2}/c')\) where \(\rmd_\alpha \geq  c \td{r}^{-1/2}\).
\item[(2)] \(h_r\geq\zeta \rmd_\alpha ^{-2}\exp \big(-\td{r}^{1/2}\rmd_\alpha /c \big)\) 
\item[(3)] \(h_r \leq \zeta ' \td{r}\exp \big(-\td{r}^{1/2}\rmd _\alpha /c \big)\)
  where \(\rmd_\alpha \geq c \td{r}^{-1/2}\).
\item[(4)] \(h_r\) has a \(Y_i\)-end limit, denoted \(h_{i,r}\), for
  every Morse end \(i\in \grY_m\). 
 This is a function over \(Y_{i, \delta
    _*}\), which  satisfies  3-dimensional versions of (1)--(3) above. 
\end{itemize}
In the above, the constants  \(c\), \(c'\),
\(\zeta _0\), \(\zeta , \zeta '\) depend  only on \(\nu \) and the
 parameters listed in (\ref{parameters}), and \(r_0\) depends on
 \(\smE\) in addition to all the above.  
\end{lemma}
\pf  To construct \(h_r\) and to verify that it satisfies Items
(2) and (3) above,  repeat the arguments in the proof of (and the
paragraph preceding) Lemma I.4.6 in \cite{T}. 
 (\(h_r\) is called \(h\) in \cite{T}.) Modify the arguments of
 \cite{T}  in like manner to what was done in the proof of Lemma \ref{T:lem4.2}: 
\begin{itemize}
\item Replace the parameter \(r\)
in \cite{T} by \(\td{r}=r\ts \), and replace balls of the form \(B\, 
(x,\rho_*)\) in \cite{T} by their counterparts, 
\(B\, (x, \ts^{-1/2}_x\, \rho_*)\), in our context. 
\item Use respectively 
Corollary \ref{cor:ball-covering} and Proposition \ref{rem:monotone}
 in place of Lemma I.3.6 and Proposition I.3.1 in
\cite{T}. 
\item Associate to each \(y\in X''_{\delta _*/2}\) a function,
  \(H_y\), on \(X\). This function \(H_y\) is given by the formula in
  \cite{T}'s (I.4.17), but with \(r\) replaced by \(\td{r}_y\), and
  with the constant \(c\) therein chosen so that \(H_y\) satisfies the
  following analog of \cite{T}'s (I.4.18): 
\begin{equation}\label{DE:H_y}
d^*d \, H_y+ \frac{r}{32} |\nu  (y)|\,  H_y\geq 0\qquad \text{where \(\rmd_\alpha \geq  c \td{r}_y^{-1/2}\).}
\end{equation}
(Recall the relation between \(|\nu |\) and \(\ts\) from
(\ref{eq:v-sig}).) Choose \(\rmm\geq 8\) such that for any given
\(y\in X''_{\delta _*}\), 
\begin{equation}\label{ratio-v}
\text{
\(\Big||\nu |/|\nu (y)|-1\Big|<1/4\) over \(B(y,
2\rmm^{-1}\ts_y)\).}
\end{equation}
Note that given \(y\in X''_{\delta _*}\), the function \(\Big|d^*d\Big(\chi \big(\rmm\dist
  (\cdot, y)/\ts_y\big)\, H_y-H_y\Big)\Big|\) is supported on \(B(y,
  \rmm^{-1}\ts_y)-B(y, (2\rmm)^{-1}\ts_y)\), and by the defining
  formula for \(H_y\),  there is constant \(\zeta _\rmm\) depending on
  \(\rmm\), but is independent of \(r\), \(y\), \((A, \Psi )\), such
  that for all sufficiently large \(r\) and \(y\in X''_{\delta _*}\), 
\begin{equation}\label{dchi-H}
\Big|d^*d\Big(\chi \big(\rmm\dist
  (\cdot, y)/\ts_y\big)\, H_y-H_y\Big)\Big|\leq \zeta _\rmm\, \td{r}_y^{1/2}
  \ts_y^{-3}\exp \big(-\td{r}_y^{1/2}(2\rmm c)^{-1}\ts_y \big). 
\end{equation}
\item Let \(\Lambda _r:=\Lambda ^{-r^{1/2}}=\{ B\, (x_k,
  \rho _k)\}_k\) 
  denote the version of  \(\Lambda ^\rho \) from Corollary
  \ref{cor:ball-covering} with \(\rho
  =r^{-1/2}\) and \(\delta =\delta _*\). Set 
\[
h_r:=\sum_{k\in \Lambda _r}\chi \big(\rmm\dist
  (\cdot, x_k)/\ts (x_k)\big) H_{x_k}.
\] 
Items (2) and (3) asserted by 
  the lemma follow from the arguments in the proof of \cite{T}'s
  Lemma I.4.6, modified as instructed in the first two bullets above. 
\end{itemize}

To verify that \(h_r\) satisfies Item (1) of the lemma, argue as in
the last bullet in the proof of  Lemma \ref{T:lem4.2}:  By
(\ref{DE:H_y}), (\ref{ratio-v}), (\ref{dchi-H}), and Corollary
  \ref{cor:ball-covering}, 
\[
\begin{split}
  & \big(d^*d \, h_r+ \frac{r}{16} |\nu | h_r\big) (x)\,\\
  &\qquad  \geq -\zeta'
_\rmm\,r\ts_x^3\sum_{k\in \Lambda _r\cap B(x, 3\ts_x)} \td{r}_{x_k}^{1/2}
  \ts_{x_k}^{-3}\exp \big(-\td{r}_{x_k}^{1/2}(2\rmm c)^{-1}\ts_{x_k}
  \big), \\
& \qquad \geq -\zeta_\rmm''\,\exp \big(-(r\ts_x^3)^{1/2}(4\rmm c)^{-1} \big) \quad \text{where \(\rmd_\alpha \geq  c \td{r}_x^{-1/2}\),}
\end{split}
\]
when \(x\in X''_{\delta _*}\) and \(r\) is sufficiently large. 

Lastly, with Items (1)-(3) now confirmed, Item (4) follows from the admissibility assumption on \((A,
\Psi )\) and  the standard compactness/properties results in
Seiberg-Witten theory (cf. \cite{KM}'s Section 10.7). 
\epf

\begin{prop}\label{prop:exp-decay}
Adopt the assumptions and notation of Lemma \ref{lem:Etop-bdd1}, and let \(\delta _*\) be as in Proposition \ref{rem:monotone}.
There exist positive constants \(r_0>8\), 
\(\zeta_i\), \(\zeta _i'\), \(i=1,2, 3\), and \(c\) that are independent of
\(r\), and \((A, \Psi )\), with the following
significance: Suppose \(r\geq r_0\), Then over \(X^{'a}_{2^5\delta
  _*}\), \(0\leq a\leq \frac{3}{2}\), 
\begin{eqnarray}
&&r|\varpi|+|F_A| \leq \zeta_2 r\tilde{\sigma}\exp
  \big(-( r\tilde{\sigma})^{1/2} \rmd_\alpha /c\big)+\zeta_2'\ts^{-2};\label{F-decay}\\
\label{eq:T-6.4}
&&|\nabla_A\underline{\alpha}|^2+r\ts^2 |\nabla_A\beta|^2\leq
\zeta_3\, r \tilde{\sigma}\exp \big(-(
  r\tilde{\sigma})^{1/2} \rmd_\alpha /c\big)+\zeta_3'\, 
r^{-1}\ts^{-5};\\
\label{beta-decay}&&
|\beta|^2\leq \zeta _1r^{-1}\ts^{-2}\exp\, (-(r\tilde{\sigma})^{1/2}\rmd_\alpha /c)+\zeta_1'r^{-2}\ts^{-5}. 
\end{eqnarray}
In the above, the constants  \(c\),
\(\zeta _i, \zeta '_i\) depend  only on \(\nu \) and the
 parameters listed in (\ref{parameters}), and \(r_0\) depends on
 \(\smE\) in addition to all the above.  
\end{prop}
\pf 
Return to the beginning of the proof of Proposition \ref{lem:est-1st-der}. 

Fix \(\delta \) with \(1\geq\delta \geq 2\delta _0'\).  
Use (\ref{eq:DE-alpha'}) 
to find constants \(\zeta _0\), 
\(\zeta _1'\), \(\zeta _1\) that are independent of \(r\), \(\delta \)
and \((A, \Psi )\),  such that the function 
\[
 \mu=\mu _\delta :=\delta ^{-2} \,
 |\nabla_A\underline{\alpha}|^2+\zeta_0\, 
 r\, |\nabla_A\beta|^2
\] 
satisfies 
\begin{equation}\label{ineq:mu_delta}
  \begin{split}
     \Big(d^*d+\frac{r}{8}|\nu|\Big)\, \mu_\delta & \leq  \zeta_1
    r\varpi_+\mu _\delta   \, +\zeta_2\delta^{-1}r^{-1}\ts^{-7}\\
    & \quad +\zeta_2'\ts^{-6}+\zeta_3\delta^{-2}r\ts^{-3}\varpi^2+\zeta_3'r\ts^{-4}\varpi_+
\end{split}
\end{equation}
over \(X^{'a}_\delta\), \(0\leq a\leq \frac{39}{8}\). In particular, this implies that for all sufficiently large \(r\) and
\(\delta \geq 2\delta'_0\), 
\begin{equation}\label{mu-ineq}
\big(d^*d+\frac{r}{16} |\nu | \big)\, \mu_\delta  \leq  \zeta' r\delta
^{-3} 
\quad \quad \text{over
  \(
X_{\delta }^{'a}-X_{\delta , v}^{'a}\), \(0\leq a\leq \frac{39}{8}\),}
\end{equation}
where \(v=1-(9\zeta _1)^{-1}\) and  \( X^{'a}_{\delta, v}:=\{y\,
|\, |\ud{\alpha}|(y)\leq v, \, y\in
  X^{'a}_\delta\}\)  as in Corollary \ref{cor:ball-covering}. 

The following observation will be useful.
Given  \(z>0\), let 
  \[
U^{'a}_{\delta , z}:=\{y\, |\,  \, \rmd_\alpha \, ( y) \geq
 z \, \td{r}^{-1/2}, y\in X^{'a}_\delta \}.\] 
According to Lemma \ref{lem:T1.4.5},
  \(X^{'a}_{\delta }-X^{'a}_{\delta , v}\supset U^{'a}_{\delta , z_v}\), where \(z_v\) is the constant from 
  Lemma \ref{lem:T1.4.5}. 
These noted, 
(\ref{mu-ineq}) implies: 
\begin{equation}\label{mu-ineq2}
\big(d^*d+\frac{r}{16} |\nu | \big)\, \mu_\delta  
< \zeta ' r\delta ^{-3}\qquad  \text{over \(U^{'a}_{\delta , z_v}\). }
\end{equation}

Note that by Propositions \ref{lem:est-1st-der} and \ref{T:prop3.1-}, there is a positive
constant \(\zeta _4\) independent of \(r\), \(\delta \), \((A, \Psi
)\), such that 
\begin{equation}\label{bdd:mu-in}
\begin{split}
\mu _\delta & \leq \zeta _4 \, r\delta ^{-1} \quad \text{over
  \(X_\delta;\quad 
  \) in particular, }\\
\mu _\delta & \leq \zeta _4\,e^{ z/c} \,  r\delta ^{-1}
\exp \big(-  \td{r}^{1/2} \rmd_\alpha /c\big) \quad \text{over
  \(X^{'a}_\delta -U^{'a}_{\delta , z}\) for arbitrary \(c , z>0\).}
\end{split}
\end{equation}
Now take \(c\) to be the constant from Lemma
\ref{lem:exp-decay}, and set \(z=\max \, ( c, z_v)\). 
Then (\ref{mu-ineq2}), (\ref{bdd:mu-in}), Lemma
\ref{lem:exp-decay} 
can be used to find constants \(\zeta _5\), \(\zeta
_5'\), which are independent of \(r\), \(\delta \), and \((A, \Psi
)\), so that: 
\begin{equation}\label{DE:h_rd}
\begin{split}
& \big(d^*d + \frac{r}{16} |\nu |\big) \, ( \mu _\delta -\zeta _5\,\delta ^{-2}
h_{r}-\zeta _5'\delta ^{-4}
) \leq 0 \quad \text{over \(U^{'a}_{\delta , z}\);}\\
& \quad \mu _\delta -\zeta _5\,\delta ^{-2} h_{r}
\leq  0 \quad \text{over \(X_\delta ^{'a}-U^{'a}_{\delta , z}\)}. 
\end{split}
\end{equation}

To be able to apply the maximum principle to (\ref{DE:h_rd}), we 
still need another comparison function, 
 \(q_1\) below,  in order  to offset the positive values of \(\mu _\delta \big|_{\partial
  U_{\delta , z}}\) along the \((\partial X_\delta )\, \cap 
  U_{\delta , z}\)  part of \(\partial \, U_{\delta , z}\). 

Recall the definition of  \(\scrX^a_\delta \subset X_\delta ^{'a}\)
and \(\gamma _{a,\delta }\) from the beginning of Section \ref{sec:pt-est}. 
Given any \(i\in \grY_m\), let \( \mu _{\delta , i}\) denote the
\(Y_i\)-end limit of \(\mu _\delta \) , and let  \(q_{1,i}\) be the function on  \(Y_{i,
  \delta }\) solving the following Dirichlet boundary value problem:
\begin{equation}\label{def:q_1i}
\begin{split}
\big(d^*d + \frac{r}{16} |\nu _i|\big) \,q_{1,i} & = d^*d\big(( 1-\chi (\sigma_i
/\delta ))\, \mu _{\delta , i}\big)-( 1-\chi (\sigma_i
/\delta ))\, d^*d \, \mu _{\delta , i};\\
q_{1,i}\big|_{\partial Y_{i, \delta }}& =0. 
\end{split}
\end{equation}
Then \(q_1\) is defined to be  the function on \(\scrX^{a+1}_\delta
\), \(0\leq a\leq \frac{31}{8}\), solving:
\begin{equation}\label{def:q_1exp}
\begin{split}
& \big(d^*d + \frac{r}{16} |\nu |\big) \,q_1= d^*d\big(\gamma
_{a,\delta }\, \mu _{\delta}\big)-\gamma _{a,\delta }\, d^*d \, \mu _{\delta };\\
& \qquad \qquad q_{1}\big|_{\partial \scrX^{a+1}_\delta }=0;\\
& \text{\(\forall i\in \grY_m\), \(\, \,  q_{1,i}\) is the \(Y_i\)-end limit
  of \(q_1\). }
\end{split}
\end{equation}

Let 
\[\begin{split}
\rmu_\delta & :=\gamma _{a, \delta }\,  \mu
_\delta -q_1-\zeta _5\,\delta ^{-2}h_{r}-\zeta _5'\, 
\delta ^{-4}; \\
\rmu_{i, \delta } & :=( 1-\chi (\sigma_i
/\delta ))\, \mu _{\delta , i}-q_{1,i}-\zeta _5\,\delta ^{-2}h_{i,r}-\zeta _5'\, 
\delta ^{-4}, \quad i\in \grY_m.  
\end{split}
\]
Then (\ref{def:q_1exp}) 
and (\ref{DE:h_rd}) together imply:
\[
\begin{split}
 & \big(d^*d + \frac{r}{16} |\nu |\big) \, \rmu_\delta   \leq 0 \quad  \text{over
  \(\scrX_{\delta }^{a+1}\cap U^{a+1}_{\delta , z} \) 
;}\\
& \qquad \rmu_\delta <0   \qquad \text{over \(\, \partial\,
  (\scrX_{\delta }^{a+1} \cap  U^{a+1}_{\delta , z})\)};\\
 &  \text{\(\forall i\in \grY_m\), \(\rmu_{i, \delta}\) is the
   \(Y_i\)-end limit of \(\rmu_\delta \),}
\end{split}
\] 
while \(\rmu_{i, \delta}\) satisfies: 
\[
\begin{split}
& \big(d^*d + \frac{r}{16} |\nu _i|\big) \, \rmu_{i,\delta } \leq 0
\quad \text{over \(U_{i, \delta , z}\);}\\
& \quad \rmu_{i,\delta } \big|_{\partial (U_{i,\delta ,
    z})}<0\quad \forall i\in \grY_m, 
\end{split}
\]
where \(U_{i, \delta , z}\subset Y_{i,\delta }\) is the 3-dimensional
analog of \(U^{a+1}_{\delta , z}\):\\ \(U_{i, \delta , z}:=\{y\, |\,  \, \rmd_{\alpha _i}\, ( y) \geq
 z \, (r\ts_i(y))^{-1/2}, y\in Y_{i, \delta }\}\), and \(\ov{U^{a+1}_{\delta,z}
 }\cap Y_{i: \infty, \delta }=U_{i, \delta , z}\). 
Thus, according to the maximum principle, 
\[\begin{split}
  \sup_{U_{i, \delta ,z}} \rmu_{i,\delta   } &\leq 0;  
 \\
\sup_{\scrX_{\delta }^{a+1}\cap U^{a+1}_{\delta , z}} \rmu_\delta    &
\leq \max_{i\in \grY_m}\,   ( \sup_{U_{i, \delta ,
     z}}  \rmu_{i,\delta
})\leq 0,  
\end{split}
\]
Together with  the second line of (\ref{DE:h_rd}), this implies in
particular that 
\begin{equation}\label{bdd:mu-d-q}
\mu _\delta \leq |q_1|+\zeta _5\,\delta ^{-2}h_{r}+\zeta _5'\, 
\delta ^{-4}  \quad  \text{over \(X^{'a}_{2\delta }\subset \scrX_\delta ^{a+1}\).} 
\end{equation}

We next give a pointwise bound for \(|q_1|\), beginning by 
estimating its \(Y_i\)-end limits, \(|q_{1,i}|\). 
First, note that the (Dirichlet) Green's function for the elliptic
system (\ref{def:q_1i}),  
denoted \(G_{i,r}\), satisfies a bound of the form: 
\[\begin{split}
    |G_{i, r}(x, \cdot)|& +\dist(x,\cdot)\, |d G_{i,r}(x, \cdot)|\\
    & \leq
\rmc_i\dist(x,\cdot)^{-1}\exp\,  \big(-
  (r\delta )^{1/2}\dist (x,\cdot)/\rmc_i\big),
\end{split}
\] 
where \(\rmc_i\) is a constant depending only on the metric and
\(\nu _i\). Multiply both sides of (\ref{def:q_1i}) by \(G_{i,r}\) and integrate over \(Y_{i,
 \delta}\). The preceding bounds for \(G_{i,r}\) and \(dG_{i,r}\),
together with the (\(Y_i\)-end limit of) the first line in (\ref{bdd:mu-in}) then leads to: 
\begin{equation}\label{bdd:q_1i0}
\begin{split}
|q_{1,i} | \leq 
\zeta '_i r^{-1/2}\delta ^{-5/2} \quad  \text{over
  \(Y_{i,\delta }\)}, \\
\end{split}
\end{equation}
where \(\zeta '_i\)  is a constant independent of \(r\), \(\delta \),
and \((A, \Psi )\). Moreover, using  (the \(Y_i\)-end limit of the bound from)
Proposition \ref{lem:est-1st-der} and (\ref{eq:3d-E_t}), 
one has for all sufficiently large \(r\) and \(\delta
\geq\delta _*\): 
\begin{equation}\label{bdd:q_1i}
\begin{split}
|q_{1,i} |& \leq \zeta \delta ^{-2}\sigma _i^{-1}\exp\, \big(-(2\rmc _i)^{-1}
  (r\delta )^{1/2}\sigma _i\big)\int_{A_{i,\delta }}\mu _{\delta ,i}\\
& \leq \zeta '\delta ^{-2}\sigma _i^{-1}\exp\, \big(-(2\rmc _i)^{-1}
  (r\delta )^{1/2}\sigma _i\big)\, \delta ^{-3} \\
& \leq \zeta ''\, \exp\, \big(-
  (r\delta ^3)^{1/2}/\rmc_i\big)\quad \text{over \(Y_{i,
    4\delta }\)}, 
\end{split}
\end{equation}
where \(A_{i,\delta }\) denotes the support of \(d\chi (\sigma_i
/\delta ))\), and \(\zeta ''\)  is a constant independent of \(r\), \(\delta \),
and \((A, \Psi )\).

These \(q_{1,i} \) being the \(Y_i\)-end limits of
\(q_1\), one may find an \(L>8\), such that 
\begin{equation}\label{diff:q_1-i}
\begin{split}
& |q_{1}-(\gri_\delta ^{-1}\circ \pi _Y)^*q_{1,i} | \leq r^{-1}\quad \text{over
  \(\hat{Y}_{i,L,\delta }\) \(\forall \, i\in \grY_m\), and hence by
  (\ref{bdd:q_1i}), }\\
& |q_1|\leq \zeta _7 r^{-1} \quad \text{over
  \(\hat{Y}_{i,L,4\delta }\) \(\forall \, i\in \grY_m\).} 
\end{split}
\end{equation}
In the above, \(\zeta _7\) is a constant independent of   \(r\),
\(\delta \), and \((A, \Psi )\) (though \(L\) may be
dependent). Let \(\scrX^{a+1}_{\delta , L}\subset \scrX^{a+1}_{\delta
}\) be a manifold with boundary with \(\scrX^{a+1}_{\delta
}-\bigcup_{i\in \grY_m}\hat{Y}_{i,L+2, \delta
}\subset\scrX^{a+1}_{\delta , L}\subset \scrX^{a+1}_{\delta
}-\bigcup_{i\in \grY_m}\hat{Y}_{i,L+3, \delta
}\), and let \(\lambda _L\) be a smooth cutoff function supported on
\(\scrX^{a+1}_{\delta , L}\) which equals 1 over  \(\scrX^{a+1}_{\delta
}-\bigcup_{i\in \grY_m}\hat{Y}_{i,L+1, \delta
}\). By (\ref{def:q_1exp}), the function \(\lambda _Lq_1\) now solves the Dirichlet
boundary valued problem: 
\[
\begin{split}
& \big(d^*d + \frac{r}{16} |\nu |\big) \,(\lambda _Lq_1)=
d^*d\big(\gamma _{a,\delta }\, \mu _{\delta }\big)-\gamma _{a,\delta
}\, d^*d \, \mu _{\delta }+d^*d\, (\lambda _Lq_1)-\lambda _Ld^*dq_1;\\
& \qquad \qquad (\lambda _Lq_{1})\big|_{\partial \scrX^{a+1}_{\delta , L}}=0.\\
\end{split}
\]
The Green's function for the preceding elliptic system satisfies a
bound of the form:   
\[\begin{split}
    |G_{r}(x, \cdot)|& +\dist(x,\cdot)\, |d G_{r}(x, \cdot)|\\
    & \leq
c_0\dist(x,\cdot)^{-2}\exp\, \big(-
  (r\delta )^{1/2}\dist (x,\cdot)/c_0\big), 
\end{split}
\]
where \(c_0\) is a constant depending only on the metric and \(\nu
\). As usual, multiply both sides of the first line in (\ref{def:q_1exp}) 
by \(G_{r}\) and integrate over \(\scrX_\delta^{a+1}\).
Integration by parts then gives: 
\[\begin{split}
|q_{1} |& \leq \zeta _1\delta ^{-2}\sigma ^{-2}\exp\, \big(-(2c_0)^{-1}
  (r\delta )^{1/2}\sigma \big)\int_{X^{'a+1}_\delta-X^{'a+1}_{2\delta}}\mu
_{\delta }\\
& \qquad +\zeta _1'\fs_*^{-2}\exp\, \big(-(r\delta
)^{1/2}\fs_*/c_0\big)\int_{X^{'a+1}_\delta-X^{'a}_\delta}\mu_{\delta
}\\
& \qquad +\zeta _2'\fs_L^{-2}\exp\, \big(-(r\delta
)^{1/2}\fs_L/c_0\big)\int_{\bigcup_{i\in \grY_m}(\hat{Y}_{i,L+1, \delta
}-\hat{Y}_{i,L+2, \delta
})}|q_1|\quad 
\end{split}
\]
over \(X^{'a-\frac{1}{8}}_{4\delta }-\bigcup_{i\in \grY_m}\hat{Y}_{i,L, \delta
}\) for all sufficiently large \(r\) and \(\delta
\geq\delta _*\), where $\fs_*\geq 1$, $\fs_L\geq 1$ respectively
denote distances to $ X^{'a+1}_\delta-X^{'a}_\delta$ and $\bigcup_{i\in \grY_m}(\hat{Y}_{i,L+1, \delta
}-\hat{Y}_{i,L+2, \delta
})$; and $\zeta _1$, $\zeta _1'$, $\zeta _2'$ are constants
independent of $r$, $\delta $, $L$, and $(A, \Psi )$. Use Propositions 
\ref{T:prop3.1-}, 
\ref{lem:est-1st-der}, Lemma
\ref{lem:u^a-int} and its amendments in Lemmas \ref{lem:monotone0},
\ref{lem:monotone00} to estimate the first two  integrals in the
preceding expression, and use (\ref{bdd:q_1i0}), (\ref{diff:q_1-i}) to
bound the last integral. One has: 
\[\begin{split}
|q_{1} |& \leq \zeta _7'\exp\, \big(-
(r\delta ^3)^{1/2}/\rmc_0\big)\\
& \leq \zeta _8 r^{-1}\quad \text{over \(X^{'a-\frac{1}{8}}_{4\delta }-\bigcup_{i\in \grY_m}\hat{Y}_{i,L, 4\delta
}\), \(0\leq a\leq \frac{31}{8},\) }
\end{split}
\]
for all sufficiently large \(r\) and \(\delta
\geq\delta _*\), where $\zeta _7'$, $\zeta _8$  are constants 
independent of $r$, $\delta $, $L$, and $(A, \Psi )$. 
Combine the preceding bound with the second line of (\ref{diff:q_1-i}), 
(\ref{bdd:mu-d-q}), the second line of (\ref{DE:h_rd}), Item (3) of
Lemma \ref{lem:exp-decay}, and the second line of (\ref{bdd:mu-in}) to see that 
\[\begin{split}
& 
 \delta ^{-2}|\nabla_A\underline{\alpha}|^2+\zeta_0\, 
 r\, |\nabla_A\beta|^2
 =\mu _\delta
 \\
& \quad \leq \zeta _5\,\delta ^{-2} h_r+\zeta _6'\, 
\delta ^{-4}  \\
& \quad \leq \zeta _6\,r\delta ^{-1}\exp \big(-(r\delta )^{1/2}\rmd _\alpha /c \big)+\zeta _6'\, 
\delta ^{-4}  \quad  \text{over \(X^{'a}_{4\delta }\), \(0\leq a\leq \frac{15}{4}\).}
\end{split}
\]
Fix $x\in X^{'a}_{8\delta _*}$ and take $\delta =\ts (x)$ in the preceding
  expression. It follows from Item (3) of Lemma \ref{lem:exp-decay} that 
\begin{equation}\label{eq:T-6.4+}\begin{split}
& |\nabla_A\underline{\alpha}|^2+\zeta_0\, r\ts^2
  |\nabla_A\beta|^2/4
 \\
& \qquad \leq \zeta _5\,h_r+\zeta _6'\, 
\ts^{-2}\quad \\
& \qquad \leq \zeta _6\,\td{r}\exp \big(-\td{r}^{1/2}\rmd _\alpha  /c\big)+\zeta _6'\, 
\ts^{-2} \quad  \text{over \(X^{'a}_{8\delta _*}\), \(0\leq a\leq \frac{15}{4}\)} 
\end{split}
\end{equation}
where $\zeta _6$, $\zeta _6'$ are constants 
independent of $r$, $\delta $, $L$, and $(A, \Psi )$. Hold this for
the moment before we return to show
how this leads to (\ref{eq:T-6.4}). 

To prove the remaining assertions of the proposition, (\ref{F-decay})
and (\ref{beta-decay}), it suffices to verify the bound for
\(|\varpi|\); in fact, we shall show: 
\begin{equation} \label{w-decay}
|\varpi|\leq \zeta \tilde{\sigma}\exp
  \big(-( r\tilde{\sigma})^{1/2} \rmd_\alpha
  /c\big)+\zeta'r^{-1}\ts^{-2}\quad
\text{over \(X^{'a}_{16\, \delta
    _*}\), \(0\leq a\leq \frac{21}{8}\).}
\end{equation}
Once (\ref{w-decay}) is established, the asserted bounds for \(|F_A|\) and \(|\beta |^2\) would follow from
 Propositions \ref{prop:F-refined} and \ref{T:prop3.1-}. 

To verify (\ref{w-decay}), note that the argument for
(\ref{ineq:alpha}) actually shows that the {\em absolute value} of the
left hand side of (\ref{ineq:alpha}) is no greater than its right hand
side. Together with (\ref{eq:omega-ineq}), Proposition
\ref{T:lem3.2}, and (\ref{eq:T-6.4+}),  this implies: 
\begin{equation}\label{DE:varpi|}
\begin{split}
2^{-1}d^*d\, |\varpi| +4^{-1}r|\alpha |^2 |\varpi |
& \leq |\nu |\,
|\nabla_A\ud{\alpha}|^2+\zeta_1'|\nabla_A\beta|^2+\zeta_2'
\ts^{-2}|\beta |^2+\zeta _2 \ts^{-1} \\
& \leq \zeta  \ts h_r +\zeta ' \ts^{-1} \quad  \text{over
  \(X^{'a}_{8\delta _*}\), \(0\leq a\leq \frac{15}{4}\).}
\end{split}
\end{equation}
Meanwhile, by Proposition
\ref{T:lem3.2}, 
\begin{equation}
\label{bdd:w-in}
\begin{split}
|\varpi |& \leq |\nu |\quad \text{over
  \(X^{'a}_{\delta_*}\);\quad  in particular, }\\
|\varpi | & \leq \zeta _4\,e^{ z/c} \,  \ts
\exp \big(-  \td{r}^{1/2} \rmd_\alpha /c\big) \quad \text{over
  \(X^{'a}_{\delta _*}-U^{'a}_{\delta _*, z}\) for arbitrary \(c , z>0\).}
\end{split}
\end{equation}
Recall that   \(X^{'a}_{\delta }-X^{'a}_{\delta , v}\supset U^{'a}_{\delta , z_v}\), where \(z_v\) is the constant from 
  Lemma \ref{lem:T1.4.5}. 
Let \(z_*=\max\, (z, z_{1/4})\), where \(z\) is as in the paragraph
after (\ref{bdd:mu-in}), and 
combine (\ref{DE:varpi|}),  (\ref{bdd:w-in}) with Lemma \ref{lem:exp-decay} to 
 find constants \(\zeta _3\), \(\zeta
_3'\) that are independent of \(r\), and \((A, \Psi
)\), so that: 
\begin{equation}\label{DE:u_rd}
\begin{split}
& \big(d^*d + \frac{r}{8} |\nu |\big) \, \big( |\varpi |-\zeta
_3\,r^{-1}h_r-\zeta _3' r^{-1}\ts^{-2}\big) \leq 0 \,\,  \\
& \qquad \qquad \qquad \text{over
  \(X^{'a}_{8\delta _*}-X^{'a}_{8\delta _*, 1/4}\supset X^{'a}_{8\delta
    _*}\cap U^{'a}_{8\delta _*, z_*}\);}\\
& \quad |\varpi |-\zeta
_3\,r^{-1}h_r\leq  0 \quad  \text{over \(X^{'a}_{8\delta
    _*}-U^{'a}_{8\delta _*, z_*}\), \(0\leq a\leq \frac{15}{4}\)}. 
\end{split}
\end{equation}
Noting that (\ref{DE:u_rd}), (\ref{bdd:w-in}) are respectively
parallel to (\ref{DE:h_rd}) and (\ref{bdd:mu-in}), a 
repetition of the arguments following  (\ref{DE:h_rd})  leads to
(\ref{w-decay}).

We now return to the verification of (\ref{eq:T-6.4}). Given
(\ref{eq:T-6.4+}), it suffices to verify (\ref{eq:T-6.4}) over the
region where \(\td{r}\exp \big(-\td{r}^{1/2}\rmd _\alpha  /c\big)\leq
\ts^{-2}\). That is, where  \(\rmd_\alpha \geq\frac{c}{\td{r}}\ln
(\td{r}\ts^2)\). Denote this region in \(X^{'a}_{16\delta _*}\),
\(0\leq a\leq \frac{21}{8}\) by
\(V_r\). Reconsider the inequality 
(\ref{ineq:mu_delta}), and the arguments that follow. Armed with (\ref{DE:varpi|}), (\ref{mu-ineq2})
may be improved as:
\[\begin{split}
 & \big(d^*d+\frac{r}{16} |\nu | \big)\, \mu_\delta  
 < \zeta '_1\ts^{-6} 
 +\zeta_4'r\ts^{-1}\delta^{-2}\exp
  \big(-2( r\tilde{\sigma})^{1/2} \rmd_\alpha
  /c\big)\\
  & \qquad +\zeta_5'r\ts^{-3}\exp
  \big(-( r\tilde{\sigma})^{1/2} \rmd_\alpha
  /c\big)\quad  \text{over \(U^{'a}_{\delta , z_v}\). }
\end{split}
\]
In particular,
\[
 \big(d^*d+\frac{r}{16} |\nu | \big)\, \mu_\delta  
 < \zeta '_2\ts^{-6}\quad  \text{over \(V_r\). }
\]
Noting that over \(\partial V_r\), \(\mu _\delta \leq \zeta \delta
^{-2}\ts^{-2}\) and appealing to Lemma
\ref{lem:exp-decay} again, we have the following analog of (\ref{DE:h_rd}): 
\[\begin{split}
& \big(d^*d + \frac{r}{16} |\nu |\big) \, ( \mu _\delta -\zeta _6\,\delta ^{-2}
h_{r}-\zeta _6'r^{-1}\delta ^{-7}
) \leq 0 \quad \text{over \(V_r\);}\\
& \quad \mu _\delta -\zeta _5\,\delta ^{-2} h_{r}
\leq  0 \quad \text{over \(\partial V_r\)}. 
\end{split}\]
Repeating the arguments following (\ref{DE:h_rd}) with straightforward
modifications then leads to (\ref{eq:T-6.4}). 
\epf

\section{Proofs of the main theorems}\label{sec:7}

The section brings together the intermediate results obtained in
Sections 4-6 above to prove the theorems claimed in Section 1.  In
this section, the parameter \(a\) in \(X^{'a}\) is set to 
satisfy \(0\leq a\leq 3/2\).

\subsection{Proving Theorem \ref{thm:l-conv}}\label{sec:l-conv}
Let \(\mathfrc{c}_r:=(A_r, \Psi_r)\), \(r\in \{r_n\}_n=:\Gamma\) be as in the statement of
Theorem \ref{thm:l-conv}. As observed before, in this case by Lemma \ref{lem:E_topX}, both inequalities
in (\ref{assume:EtopX-ubdd}) in Lemma \ref{lem:Etop-bdd1} hold with
the constant \(\smE\) determined by (\ref{eq:CSD-est}) via
(\ref{eq:bbE}). In this section, we shall repeated invoke this fact in
combination with other results in previous sections, sometimes
implicitly. 

\subsubsection*{\it Proof of Theorem \ref{thm:l-conv}. }

{\bf (a):} To find the desired t-curve \(\mathbf{C}\) on \(X^{'a}\), associate
to each \((A_r, \Psi_r)\) a 2-current \(\mathcal{F}_r\) on \(X^{'a}\)
as follows.  Given a compactly supported smooth
2-form \(\mu\) on \(X^{'a}\), let
\[
\mathcal{F}_r(\mu):=\frac{1}{2}\int_{X^{'a}}\frac{i}{2\pi}
(F_{A_r}-F_{A_K})\wedge \mu,
\]
where \(F_{A_K}\) is the curvature form of the canonical connection on
\(K^{-1}\). (As observed in \cite{Ts}, the 2-form \(iF_{A_K}\) is
singular at \(Z:=\nu ^{-1}(0)\), but for any compact \(X_\bullet\subset X''\), the
integral \(\int_{X_\bullet}iF_{A_K}\wedge \mu\) is still finite because
\(|F_{A_K}|\leq \zeta _{X_\bullet}\dist (\cdot, Z)^{-2}\) near \(Z\cap
X_\bullet\)). 
In fact, by our definition of MCE, 
\begin{equation}\label{F_K}
\|F_{A_K}\|_{L^1(X_\bullet)}\leq
\zeta_K|X_\bullet|.
\end{equation}
In the above, \(\zeta _K\) is a positive
constant that depends only on the metric on \(X\) and \(\nu \). Together with
(\ref{bdd:L^14d}) and Lemma \ref{lem:E_topX}, this implies that there is a positive constant
\(\zeta _E\) that depends only on
\BTitem\label{dep-list}
 \item the metric and \(\Spin^c\)-structure on
\(X\),  
\item the relative homology class \(\rmh(\mathfrc{h})\in \scrH
\big((X,\nu),\{\tilde{\gamma}_i\}_i\big)\) (which in turn depends
implicitly on \(\{\tilde{\gamma}_i\}_i\), 
\item \(\nu \), \(\varsigma_w\), \(\zzz_\grp\), 
\ETitem
such that 
\begin{equation}\label{def:zeta-E}
\|F_{A_r^E}\|_{L^1(X_\bullet)}\leq
\zeta_E\, (|X_\bullet|+1).
\end{equation}
Thus, by
Alaoglu's theorem, there is a subsequence of \(\{\mathcal{F}_r\}_r\)
which converges on \(X_\bullet\). 
A diagonalization
argument then shows that  a subsequence, \(\Gamma _0\),  of the former converges to a
current \(\mathcal{F}\) on \(X^{'a}\). Abusing notation, this subsequence
is also denoted by \(\{\mathcal{F}_r\}_r\). With this done, the arguments in
Sections 7(b)-(d) and 7(f) in \cite{Ts} carry through directly, with
Lemma 4.1 therein substituted by our Corollary
\ref{cor:ball-covering}. This gives the t-curve \(\mathbf{C}\)
and proves the first claim of the theorem regarding t-convergence.

{\bf (b):} 
By definition, 
\[
\scrF_\omega(\mathbf{C})\leq 
\scrF_\nu(\mathbf{C}|_{X'})+\sup_{X_*\in \grX_1}\scrF_\omega (\mathbf{C}|_{X_*}),
\]
where \(\grX_1=\{ X_*\, |\, X_*\subset X^{'a},
|X_*|\leq 1\}\). Consequently, it suffices to bound
\(\scrF_\nu(\mathbf{C}|_{X^{'a}})\) and \(\scrF_\omega
(\mathbf{C}|_{X_*})\). 

We begin by bounding \(\scrF_\nu(\mathbf{C}|_{X^{'a}})\).
Fix a compact \(X_\bullet\subset X^{'a}\).
Since \(\{\mathcal{F}_r\}_r\) converges to
\(\tilde{C}\) as currents, 
\begin{equation}\label{scrF1}
\scrF_\nu(\mathbf{C}|_{X_\bullet})=\lim_{r\to \infty}\frac{1}{2}\int_{X_\bullet}\frac{i}{2\pi}
(F_{A_r}-F_{A_K})\wedge \nu
\end{equation}
Note that by (\ref{F_K}) and (\ref{eq:xi-exp}), 
\begin{equation}\label{int:F_K}
\Big|\int_{X_\bullet}\frac{i}{2\pi}
F_{A_K}\wedge \nu\Big|\leq \zeta'_K
\end{equation}
for a positive constant \(\zeta'_K\) independent of
\(X_\bullet\). Combine this with (\ref{scrF1}) and Lemmas
\ref{lem:F_v},  \ref{lem:E_topX} 
to get 
\[
\scrF_\nu(\mathbf{C}|_{X_\bullet})\leq \zeta_v, 
\]
where \(\zeta_v>0\) is independent of \(X_\bullet \) and \(r\), and is determined
by the items listed in (\ref{dep-list}). This implies that
\[\scrF_\nu(\mathbf{C}|_{X^{'a}})\leq \zeta_v\] as well.

Fix \(X_*\in \grX_1\). We now proceed to bound \(\scrF_\omega
(\mathbf{C}|_{X_*})\).
By the convergence of \(\{\mathcal{F}_r\}_r\) again, 
\begin{equation}\label{scrF2}
\scrF_\omega(\mathbf{C}|_{X_*})=\lim_{r\to \infty}\frac{1}{2}\int_{X_*}\frac{i}{2\pi}
(F_{A_r}-F_{A_K})\wedge \omega
\end{equation}
Combine this with (\ref{F_K}), Lemma \ref{lem:E_topX}, Lemma \ref{co:E-omega-bdd3} and its
amendment in Proposition
\ref{prop:integral-est2}, Lemma \ref{co:E-omega-bdd3}, we have 
\[
\scrF_\omega(\mathbf{C}|_{X_*})\leq \zeta_v' 
\]
for any \(X_*\in \grX_1\). Like \(\zeta_v\),  \(\zeta_v'\) only depends 
on the items listed in (\ref{dep-list}). Together with the  bound 
\(\scrF_\nu(\mathbf{C}|_{X^{'a}})\leq \zeta _v\), this implies Assertion (b) of 
Theorem \ref{thm:l-conv}. \epf

\subsection{SW and Gr versions of Floer-theoretic energies}\label{sec:energy-comp} 


Compactness
results in Floer theory typically hinge on bounds on certain
``topological energy'', whose value only depends on the relative
homotopy class of the solution. In Morse theory, this is just the
change of the value of the Morse function along the gradient flow. 
In the Seiberg-Witten setting, the
relevant notion of topological energy is \(\scrE_{top}\)
in \cite{KM} (cf. e.g. Equation (5.1) of \cite{KM}). For
\(Gr\)-variants of Floer theories (such as ECH/PFH, SFT) that do
not really have a Morse-theoretic model, local convergence and global
convergence are respectively controlled by \(\scrF_\omega\) and
\(\scrF_\nu\). In Theorem \ref{thm:l-conv}, we obtained a bound on
\(\scrF_\omega\) of the limiting t-curve \(\mathbf{C}\). For
convenience of later discussions, rename this \(\mathbf{C}\) as
\(\mathbf{C}_0\).
Before moving on to the proof of Theorem
\ref{thm:g-conv}, we need to better understand the relation between
\(\scrF_\nu(\mathbf{C}_0)\) and its counterpart on the SW side,
\(\scrE_{top}\). The two are related via the quantity: 
\[\scrV_r(X_\bullet):=\int_{X_\bullet}\frac{i}{2\pi}
F_{A^E_r}\wedge \nu=\int_{X_\bullet}\frac{i}{4\pi}
F_{A_r}\wedge \nu-\int_{X_\bullet}\frac{i}{4\pi}
F_{A_K}\wedge \nu .\] 
Let 
\[
\grE_r(X_\bullet):=
\frac{1}{\pi r}\scrE_{top}^{\mu_r, \hat{\grp}}(X_\bullet)(A_r, \Psi_r)-\int_{X_\bullet}
\frac{iF_{A_K}}{4\pi }\wedge \nu . 
\]

It is easier to compare averaged versions of the three notions
of energy, \(\grE_r\), \(\scrV_r\), \(\scrF_\nu \), defined as
follows. Suppose \(X_\bullet\) is of the following type: 
\begin{equation}\label{ave-X}
X_\bullet=\begin{cases}
\hat{Y}_{i, [l, L]}, & \text{ with \(|\hat{Y}_{i, [l,
  L]}|\geq 2\), \(l\geq 1/2\) or}\\
X_\bfl\subset X^{'a}& \text{ \(\bfl(i)\geq 1/2\) for every \(i\in
  \grY_m\); \(\bfl(i)=\ul_i\) for \(i\in \grY_v\)}.
\end{cases}
\end{equation}
Let
\[\begin{split}
\bar{\grE}_r(\hat{Y}_{i,[l, L]}) & :=\int_{-1/2}^{1/2}
\grE_r(\hat{Y}_{i,[l+t, L+t]}) \, dt.\\
\bar{\grE}_r(X_\bfl) & :=\int_{-1/2}^{1/2}\grG _r(X_{\bfl+\bft}) \,
dt, 
\end{split}\]
where \(\bft\) is the function given by \(\grt(i)=t\) \(\forall i\in
\grY_m\); \(\bft(i)=0\) \(\forall i\in \grY_v\).
Let \(\bar{\scrV}_r(X_\bullet)\), \(\bar{\scrF}_\nu
(\bfC|_{X_\bullet})\) be similarly defined. 
Thus defined, \(\bar{\grE}_r\) (and also \(\bar{\scrV}_r\), \(\bar{\scrF}_\nu\)) are
additive: Given \(l<L<L'\) with \(L-l\), \(L'-L\) no less than 2, for
\(i\in \grY_m\)
\[
 \bar{\grE}_r(\hat{Y}_{i,[l, L]})+\bar{\grE}_r(\hat{Y}_{i,[L,
  L']})=\bar{\grE}_r(\hat{Y}_{i,[l, L']});
\]
given functions  \(\bfl\co \grY\to [0,\infty]\),  \(\bfL\co
\grY\to [0,\infty]\) such that \(\bfL(i)-\bfl(i)\geq 2\) \(\forall i \in
\grY_m\) and that \(X_\bfL\), \(X_\bfl\) are of the type
(\ref{ave-X}), 

\[
\bar{\grE}_r(X_\bfL)-\bar{\grE}_r(X_\bfl)=\sum_{i\in
  \grY_m}\bar{\grE}_r(\hat{Y}_{i, [\bfl(i), \bfL(i)]}). 
\]

From Theorem \ref{thm:l-conv} we already know
that 
\[\begin{split}
& \lim_{r\to \infty}\scrV
_r(X_\bullet)=\scrF_\nu(\mathbf{C}_0|_{X_\bullet}) \quad \text{for all
  compact
\(X_\bullet\subset X^{'a}\). }\\
& \lim_{r\to \infty}\bar{\scrV}
_r(X_\bullet)=\bar{\scrF}_\nu(\mathbf{C}_0|_{X_\bullet}) \quad \text{for all
  compact
\(X_\bullet\subset X^{'a}\) of the type (\ref{ave-X}). }
\end{split}
\]
Meanwhile, by (\ref{eq:F-nu}), Lemma \ref{lem:E_topX}, Lemma
\ref{co:E-omega-bdd} and its amendment in Proposition \ref{prop:integral-est2}, we
know that for all \(X_\bullet\), the difference between \(\grE_r(X_\bullet)-\scrV_r(X_\bullet)\) is bounded by a constant
independent of \(r\) and \(X_\bullet\). Similarly for
\(\bar{\grE}_r(X_\bullet)-\bar{\scrV}_r(X_\bullet)\). We now show that
the latter goes to 0 as \(r\to \infty\).

\begin{prop}\label{cor:E-cp}
Let \(\{(A_r, \Psi_r)\}_{r\in \Gamma }\) be as in the statement of Theorem
\ref{thm:l-conv}. Let \(X_\bullet\subset X^{'a}\) be of the type
(\ref{ave-X}). 
Then there exists positive constants  \(R_S>1\), 
\(\zeta _S\) depending only on (\ref{dep-list}), 
such that
for all \(r\geq R_S\), 
\[
\Big|\bar{\grE}_r(X_\bullet)-\bar{\scrV}_r(X_\bullet)\Big|\leq \zeta
_S\, r^{-1/4}.
\]
\end{prop}
\pf We verify the preceding inequality only for the case
\(X_\bullet=X_\bfl\). The case when \(X_\bullet=\hY_{i, [l, L]}\) is
similar. Write \(\partial\ov{X_\bfl}=\bigcup_{i\in \grY} Y_{i:
  l_i}\), \(l_i=\ul_i\) when \(i\in \grY_v\). Then 
\[
\pi \big(\bar{\grE}_r(X_\bullet)-\bar{\scrV}_r(X_\bullet) \big)=r^{-1}\int_{-1/2}^{1/2}
\scrE_{top}^{w_r}(X_\bullet)\, dt. 
\]
The right hand side of the preceding identity can be expressed as a
sum of three terms: 
\begin{equation}\label{E-V:sum}
\begin{split}
& \pi  \big(\bar{\grE}_r(X_\bullet)-\bar{\scrV}_r(X_\bullet)\big)\\
& \quad =\frac{1}{4r}\int_{-1/2}^{1/2}\int_{X_{\bfl+\bft}}F_{A_0}\wedge
(F_{A_0}+iw_r)\, dt \\
& \qquad -\sum_{i\in \grY}\int_{\hY_{i,[
    l_i-1/2, l_i+1/2]}} \big(2 \langle \psi , \slp_B\psi
\rangle +\partial _s |\nu |\big)\\
& \qquad + \frac{1}{4r}\sum_{i\in \grY}\int_{\hY_{i,[
    l_i-1/2, l_i+1/2]}}  (B-B_0)\wedge
(F_B+F_{B_0}+iw_r) \wedge ds.\\
\end{split}
\end{equation}
These three terms are estimated as follows. 

(i) For the first term on the RHS of (\ref{E-V:sum}), 
\[
\Big|\frac{1}{4r}\int_{-1/2}^{1/2}\int_{X_{\bfl+\bft}}F_{A_0}\wedge
(F_{A_0}+iw_r)\, dt\Big| \leq \zeta  _1r^{-1}, 
\]
where \(\zeta _1\) depends only on the choice of \(A_0\) and
\(w_r\). (In fact, it equals 0 by the assumptions (\ref{eq:A_0}) and
Assumption \ref{assume} (3).)

(ii) To estimate the second term on the RHS of (\ref{E-V:sum}), use
the Seiberg-Witten equation to rewrite: 
\begin{equation}\label{ii-}
\begin{split}
& -\int_{\hY_{i,[
    l_i-1/2, l_i+1/2]}} \big( 2 \langle \psi , \slp_B\psi
\rangle +\partial_s |\nu |\big) \\
& \quad\qquad =\int_{\hY_{i,[
    l_i-1/2, l_i+1/2]}} \big( 2\langle \psi , \partial _s\psi
\rangle -\partial_s |\nu |\big)\\
\\
& \quad \qquad =\int_{\hY_{i,[
    l_i-1/2, l_i+1/2]}} \partial _s (|\psi |^2-|\nu |).
\end{split}
\end{equation}
Let \(\delta \geq\sO r^{-1/3}\), where \(\sO\) is as in Proposition
\ref{T:prop3.1-}. 
Let \(v_r:=(\zzz_v\delta )^{-1}|\nu |\), where \(\zzz_v\)
is the positive
constant  from (\ref{eq:v-sig}). 
In particular, \(\sigma \geq \delta \) where \(v_r\geq 1\). Write 
\[
|\psi |^2-|\nu |=(1-\chi
(v_r))\, (|\nu |\, (|\ud{\alpha}|^2 -1)+|\beta|^2)+\chi (v_r)\, (|\psi |^2 -|\nu |).
\]
We have:
\begin{equation}\label{ii-c}\begin{split}
& \int_{\hY_{i,[l_i-1/2, l_i+1/2]}} \partial _s (|\psi |^2-|\nu |) \\
& \quad =\int_{\hY_{i,[
    l_i-1/2, l_i+1/2]}} \partial _s \big((1-\chi (v_r)) \, (|\nu |\, (|\ud{\alpha}|^2 -1)+|\beta|^2)\big)\\
& \qquad \qquad +\int_{\partial\hY_{i,[
    l_i-1/2, l_i+1/2]}} \big(\chi (v_r)\, (|\psi |^2 -|\nu |)\big)\\
&\quad =2\int_{\hY_{i,[
    l_i-1/2, l_i+1/2]}} \big((1-\chi (v_r))\, \big(|\nu |\,
\langle\ud{\alpha}, \partial_s\ud{\alpha}\rangle+\langle \beta, \partial_s\beta\rangle\big)\big) \\
&\qquad \qquad + \int_{\hY_{i,[
    l_i-1/2, l_i+1/2]}} \big((1-\chi (v_r)) \partial_s|\nu |
(|\ud{\alpha}|^2 -1)\big)\\
&\qquad \qquad -(\zzz_v\delta)^{-1} \int_{\hY_{i,[
    l_i-1/2, l_i+1/2]}} \big( \chi
'(v_r)\, \partial_s|\nu |  (|\psi |^2 -|\nu |)\big)\\
&\qquad \qquad +\int_{\partial\hY_{i,[
    l_i-1/2, l_i+1/2]}} \big(\chi (v_r)(|\psi |^2 -|\nu |)\big). 
\end{split}
\end{equation}
By Proposition \ref{T:lem3.2} and the
assumptions on \(\nu \), the last two terms above are bounded together
by 
\begin{equation}\label{ii-b}
\begin{split}
& (\zzz_v\delta)^{-1} \Big|\int_{\hY_{i,[
    l_i-1/2, l_i+1/2]}} \big( \chi
'(v_r)\, \partial_s|\nu |  (|\psi |^2 -|\nu |)\big)\Big|\\
& \quad \quad \qquad +\Big|\int_{\partial\hY_{i,[
    l_i-1/2, l_i+1/2]}} \big(\chi (v_r)(|\psi |^2 -|\nu |)\big)\Big|\\
&\quad \leq 
\zeta  _0\, \delta ^{3}, 
\end{split}
\end{equation}
where \(\zeta _0\) depends only on the metric, \(\nu \),
\(\varsigma_w\), \(\zeta _\grp\). 
The other two terms on the RHS of the second equality of
(\ref{ii-c}) are together bounded by 
\[
\int_{\bar{Y}_\delta } \Big(2|\nu |\, |\ud{\alpha
}|\, | \partial_s\ud{\alpha}|+ 2|\beta|\, |\partial_s\beta|+\Big|\partial_s|\nu |\Big|\, 
\Big||\ud{\alpha}|^2 -1\Big|\Big),
\]
where \(\bar{Y}:=\hY_{i,[
    l_i-1/2, l_i+1/2]}\). Now take  \(\delta =2^4 \delta_*\geq\sO r^{-1/3}
\), \(\delta _*\) being as in Propositions \ref{rem:monotone} and 
  \ref{prop:exp-decay}.  The preceding
  integral may be bounded using Proposition
\ref{prop:exp-decay} 
by {\small
\begin{equation}\label{ii-a}
\begin{split}
& \int_{\bar{Y}_\delta } \left( \zeta_1\sqrt{r\ts^3}\exp \big(-(
  r\tilde{\sigma})^{1/2} 
  \dist (\cdot, \alpha^{-1}(0))/(2c)\big)
   +\zeta' 
  r^{-1/2}\ts^{-3/2} +\zeta ''
r^{-1}\ts^{-3}\right)\\
& \quad \quad \leq  \zeta_1\int_{\bar{Y}_\delta } \Big( \sqrt{r\ts^3}\exp \big(-(
  r\tilde{\sigma})^{1/2} \dist (\cdot, \alpha^{-1}(0))/(2c)\big)\Big)
  +\zeta'_1 
 r^{-1/2},
\end{split}
\end{equation}
}
where \(c\) is the constant from Proposition
\ref{prop:exp-decay}. To bound the first term on the RHS of
(\ref{ii-a}), first note that it vanishes when \(\alpha
^{-1}(0)=\emptyset\). Assume consequently that \(\alpha
^{-1}(0)\neq\emptyset\). 
Invoke the second bullet of Corollary \ref{cor:ball-covering} with
\(\rho =r^{-1/2}\) to
cover \(\alpha^{-1}(0)\, \cap \, X_\delta \) 
by a set of balls, \(\Lambda ^{\rho }=\{B_k\}_k\), where \(B_k\) is centered at
\(x_k\in\alpha^{-1}(0)\, \cap \, X_\delta \), and is of radius
\((\ts (x_k))^{-1/2}\rho \). Let \(B'_k\supset B_k\) denote the concentric
ball with twice the radius, and let \(U=\bigcup_kB'_k\). Define the
function \(d_k(\cdot):=\dist (\cdot, x_{k})\) on \(X\). Given \(x\in X_\delta 
\), there exists \(k_x\in \Lambda ^\rho \) such that
\(B_{k_x}\) contains a point whose distance to \(x\) equals \( \dist (x,
\alpha^{-1}(0))\). Then 
\[\begin{split}
\dist (x, \alpha^{-1}(0))\geq d_{k_x}(x)-(\ts (x_{k_x}))^{-1/2}\rho & \geq
d_{k_x}(x)/2 \\
& \geq (\ts (x_k))^{-1/2}\rho \quad \forall x\in X_\delta -U.
\end{split}
\]
Let \(\bar{\Lambda }^\rho :=\Lambda^\rho _{\hY_{i,[l-1, l+1],
    \delta}}\). By bullet 2 of by Corollary \ref{cor:ball-covering},
\(\bar{\Lambda }^\rho \) has at most \(\zeta   r\) elements. If \(\bar{\Lambda }^\rho =\emptyset\), then \(\dist (\cdot,
\alpha^{-1}(0))\geq 1/2\) over \(\bar{Y}_\delta\), and in this case 
\[\begin{split}
& \int_{\bar{Y}_\delta } \Big(\sqrt{r\ts^3}\exp \big(-(
  r\tilde{\sigma})^{1/2} \dist (\cdot, \alpha^{-1}(0))/(2c)\big)\Big)
  \leq \zeta  _2' (r\delta )^{-1/2}. 
\end{split}
\]
Otherwise, \(\bar{\Lambda }^\rho \neq\emptyset\), and \(k_x\in \bar{\Lambda }^\rho \) 
for any 
\(x\in \bar{Y}_\delta\). In this case
\[
\begin{split}
& \int_{\bar{Y}_\delta } \Big(\sqrt{r\ts^3}\exp \big(-(
  r\tilde{\sigma})^{1/2} \dist (\cdot, \alpha^{-1}(0))/(2c)\big)\Big)\\
& \qquad \leq \int_{U\cap \bar{Y}_\delta } \sqrt{r\ts^3}+\int_{\bar{Y}_\delta -U} \Big(\sqrt{r\ts^3}\exp \big(-(
  r\tilde{\sigma})^{1/2} \dist (\cdot, \alpha^{-1}(0))/(2c)\big) \\
& \qquad \leq 
  \sum_{k\in  \bar{\Lambda }^\rho } \int_{B_k'} \sqrt{r\ts^3}+\sum_{k\in  \bar{\Lambda }^\rho } \int_{\bar{Y}_\delta -B_k'} \sqrt{r\ts^3}\exp \big(-(
  r\ts)^{1/2}d_k/(4c)\big)\\
& \qquad \leq  \zeta _2\, (r\delta )^{-1/2}.
\end{split}
\]

Insert 
the preceding bounds  for the first integral on the right hand side of
(\ref{ii-a}), and combine with 
Equations 
(\ref{ii-b}), (\ref{ii-c}) and (\ref{ii-}), 
setting \(\delta =r^{-1/7}\), which is no less than \(2^4\delta _*\)
when \(r\) is sufficiently large and when \(\epsilon \leq 1/7\) in the
definition of \(\delta _{1,\epsilon }\).  
We then have: 
\[
\text{the second term on the RHS of (\ref{E-V:sum})} \leq \zeta _6\,
r^{-3/7}. 
\]

(iii) To estimate the last term in (\ref{E-V:sum}), note that \(F_B+F_{B_0}+iw_r\) is exact by
Assumption \ref{assume} (2), and therefore
\begin{equation}\label{eq:bF}\begin{split}
& \int_{\hY_{i,[
    l_i-1/2, l_i+1/2]}}  (B-B_0)\wedge
(F_B+F_{B_0}+iw_r) \wedge ds\\& \qquad =\int_{\hY_{i,[
    l_i-1/2, l_i+1/2]}}  (B'-B_{i,0})\wedge
(F_B+F_{B_0}+iw_r) \wedge ds
\end{split}
\end{equation}
where \(B'(s)=A'|_{Y_{i:s}}\), and \(A'=u\cdot A\) is in the normalized
Coulomb-Neumann gauge on \(\hY_{i,[
  l_i-1/2, l_i+1/2]}\).

Next, arguing as in the proof of Lemma
  \ref{lem:CSD-est}, one has: 
\begin{equation}\label{bdd:A-inf}
\|A'-A_0\|_{L^\infty(\hY_{i,[
    l_i-1/2, l_i+1/2]})}\leq \zeta _2 r^{3/4},
\end{equation}
where \(\zeta _3\)
depends only on the parameters in (\ref{parameters}), \(\nu \), and
\(\smE\).

To see this, let \(\Omega'\) denote the space of 1-forms \(a\) on \(\hY_{i,[
    l_i-1/2, l_i+1/2]}\) satisfying the conditions \(\iota_{\partial s}a=0\) over \(\partial(\hY_{i,[
    l_i-1/2, l_i+1/2]})\) and \(c_k(a):=\int_{\hY_{i,[
    l_i-1/2, l_i+1/2]}} q_k\wedge a=0\), where \(q_k=d(\chi _+(s) *_3 h_k)\), \(k=1, \ldots,
b^1(Y_i)\), are as chosen in Definition \ref{rem:CNgauge}.  Note that the operator \(d+d^*\co \Omega'\to
B^2\oplus\Omega ^0\) has a Green's function \(G(x,y)\) that satisfies
the inequality
\(|G(x, y)|\leq \zeta _g \dist (x, y)^{-3}\), where \(B^2, \Omega ^0\)
are respectively the spaces of exact 2-forms and smooth functions on \(\hY_{i,[
    l_i-1/2, l_i+1/2]}\), and \(\zeta _g\) is a positive constant that
  only depends on the metric on \(Y_i\) and the choice of \(q_k\). Let
  \(q_k^*:=\pi _2^*h_k\), where \(h_k\) are as in Definition
  \ref{rem:CNgauge} and (\ref{eq:deltab1}). Write \(A'-A_0=\hata+ \sum_kc_k(A'-A_0) \, q_k^*\), where
\(\hata\in i \Omega'\). 
Then for \(x\in \hY_{i,[
    l_i-1/2, l_i+1/2]}\), 
\begin{equation}\label{bdd:hata}\begin{split}
 |\hata (x)| & \leq \zeta \int_{\hY_{i,[
    l_i-1/2, l_i+1/2]}} \dist (x, \cdot)^{-3} \big(
|F_{A}|+|F_{A_0}|\big)\\
&\leq \zeta   \|F_{A}\|_{L^\infty(\hY_{i,[
    l_i-1/2, l_i+1/2]})}\int_{B_x(\rho )\cap\hY_{i,[
    l_i-1/2, l_i+1/2]}} \dist (x, \cdot)^{-3} \\
&\qquad +\zeta  \rho ^{-3}\int_{\hY_{i,[
    l_i-1/2, l_i+1/2]}-B_x(\rho )} |F_A| +\zeta ', 
\end{split}
\end{equation}
where \(\zeta \) depends only on the metric on \(Y_i\) and the choice
of \(q_k\), and \(\zeta '\) depends on the above as well as the choice
of \(A_0\). Use Propositions \ref{T:prop3.2} and \ref{T:lem3.2} to
bound \(\|F_{A}\|_{L^\infty(X'')}\leq \zeta _1 r\), where \(\zeta _1\) depends only on the metric, \(\nu \), \(\varsigma_w\) and \(\zeta
_\grp\), and apply this bound to bound the first term on the RHS of
(\ref{bdd:hata}). Use the \(L^1\)-estimate for \(F_A\) in
(\ref{bdd:L^14d}) to bound the second term on the RHS of
(\ref{bdd:hata}). We have, taking \(\rho =r^{-1/4}\), \(\|\hata\|_{L^\infty(\hY_{i,[
    l_i-1/2, l_i+1/2]})}\leq \zeta _3 r^{3/4}\), where \(\zeta _3\)
depends only on the parameters in (\ref{parameters}), \(\nu \), and \(\smE\). 
Combining this in turn with the fact that since \(A'\) is in the {\rm
  normalized} Coulomb-Neumann gauge, \(|c_k(A'-A_0)|<2\pi \) \(\forall
k\), the asserted bound (\ref{bdd:A-inf}) follows.

To proceed, combining  (\ref{bdd:A-inf}) with the \(L^1\)-estimate for
\(F_A\) in  (\ref{bdd:L^14d}) with (\ref{eq:bF}), 
one has:
\[
 \frac{1}{4r}\Big|\int_{\hY_{i,[
    l_i-1/2, l_i+1/2]}}  (B-B_0)\wedge
(F_B+F_{B_0}+iw_r) \, ds\Big|\leq  \zeta _3 r^{-1/4}. 
\]
Finally, the claimed inequality in the proposition follows by inserting the bounds from (i), (ii), (iii) above into
(\ref{E-V:sum}).
\epf

In the upcoming proof of Theorem \ref{thm:g-conv}, we make use of a
lower bound on \(\grE_r\), which follows from  an improved lower bound on \(\scrE_{top}^{\mu _r, \hat{\grp}}\). 
\begin{lemma}\label{lem:E-lbdd}
There are positive constants \(r_0\), \(\zeta _i'\), \(\zeta _i''\)
depending only on (\ref{dep-list}) with the following significance: For all \(i\in \grY\) and \(l, L\geq 0\), 
\[\begin{split}
\scrE_{top}^{\mu_r, , \hat{\grp}}(\hat{Y}_{i,[l,L]})(A_r, \Psi_r)
& \geq -r \zeta '_i e^{-\kappa _il}, \\
\grE_r(\hat{Y}_{i,[l,L]})(A_r, \Psi_r) & \geq
- \zeta _i'' e^{-\kappa _il} \quad \text{if \(\hat{Y}_{i,[l,L]}\subset
  X^{'a}\).}
\end{split}
\]
\end{lemma}
\pf By (\ref{bdd:CSD-lower}), (\ref{diff:CSD-E}), and (\ref{bdd:L^14d}), 
\[
\scrE^{'\mu_r,  \hat{\grp}}(
\hat{Y}_{i,[l,L]}) (A_r, \Psi _r)\geq -r \zeta _i e^{-\kappa _il}.
\]
By Lemma \ref{bdd:A3d} and (\ref{eq:xi-exp}), 
\[\begin{split}
\scrE_{top}^{\mu_r, \hat{\grp}}(\hat{Y}_{i,[l,L]})(A_r,
\Psi_r)& =\scrE^{' \mu_r, \hat{\grp}}(
\hat{Y}_{i,[l,L]}) (A, \Psi )-\frac{ir}{4}\int_{\partial
  \hat{Y}_{i,[l,L]}}(B_r-B_0)\wedge \grv\\ & \geq -r \zeta _i' e^{-\kappa
  _il}.
\end{split}
\]
So 
\[
\grE_r(\hat{Y}_{i,[l,L]})=
\frac{1}{\pi r}\scrE_{top}^{\mu_r, \hat{\grp}}(\hat{Y}_{i,[l,L]})(A_r, \Psi_r)-\int_{\hat{Y}_{i,[l,L]}}
\frac{iF_{A_K}}{4\pi }\wedge \nu\geq - \zeta _i'' e^{-\kappa _il}.
\]
\epf

\begin{lemma}\label{lem:grE-X}
 Let \((A_r, \Psi _r)\) be as in Theorem \ref{thm:l-conv} and 
let \(\smE\), \(\bbE\) be respectively the constants from
(\ref{eq:CSD-est}) and (\ref{eq:bbE}) determined (via Lemma \ref{lem:E_topX}) by (\ref{dep:rel-htpy}). 
There exists a positive constant \(\zeta \) depending only on
(\ref{dep-list} 
such that 
\[
\grE_r(X^{'a})\leq \frac{\smE}{\pi } +\zeta ; \quad \bar{\grE}_r(X^{'a})\leq
\frac{\smE}{\pi } +\zeta  . 
\]
If \(X\) has no vanishing ends, \(\grE_r(X)=\bar{\grE}_r(X^{'a})\) and 
\[
\lim_{r\to\infty}\grE_r(X)=\frac{\bbE}{\pi }-\int_{X}
\frac{iF_{A_K}}{4\pi }\wedge \nu. 
\]
\end{lemma}
\pf Write 
\[\begin{split}
\grE_r(X^{'a})=& 
\frac{1}{\pi r}\scrE_{top}^{\mu_r, \hat{\grp}}(X)(A_r,
\Psi_r)-\frac{1}{\pi r}\scrE_{top}^{\mu_r, ,
  \hat{\grp}}(X-X^{'a})(A_r, \Psi_r)\\
&\qquad -\int_{X^{'a}}
\frac{iF_{A_K}}{4\pi }\wedge \nu. 
\end{split}
\]
The assertions of the lemma then follow directly from Lemma
\ref{lem:E-lbdd}, (\ref{eq:CSD-est}), and (\ref{eq:bbE}). 
\epf

\subsection{Proving Theorem \ref{thm:g-conv} and Proposition
  \ref{cor:F-positive}, the non-cylindrical case}\label{sec:g-conv:a}

We introduce some notations for convenience before embarking on the proofs. 

Given two t-curves \(\mathbf{C}=[C, \tilde{C}]\), \(\mathbf{C}'=[C',
\tilde{C}']\) in \((X, \nu)\), and a compact \(X_\bullet\subset X^{'a}\), let 
\[
\mathfrc{d}_{X_\bullet}(\mathbf{C}', \mathbf{C}) :=\dist_{X_\bullet}
(C', C)+
\|\tilde{C}'-\tilde{C}\|_{op, X_\bullet},
\]
where \(\|\cdot\|_{op, X_\bullet}\) denotes the norm as currents on
\(X_\bullet\). 
The same expression is also used when either of the 
the entries \(\mathbf{C}\) or \(\mathbf{C}'\) on the left  is replaced by a
general pseudo-holomorphic subvariety \(C\) or \(C'\). In this case
the \(\tilde{C}\) or \(\tilde{C}'\) on the right simply denotes the
associated current of the pseudo-holomorphic subvariety.

Given a t-orbit \(\pmb{\gamma}\) or a general multi-orbit
\(\gamma \) in \((Y, \theta)\), let \(\mathfrc{d}_{Y}(\pmb{\gamma},
\cdot)\), \(\mathfrc{d}_{Y}(\gamma,
\cdot)\) be similarly defined. 

Let \(\mathfrc{c}_r:=(A_r, \Psi_r)\) be an admissible solution to
\(\grS_{\mu _r, \hat{\grp}}(A_r, \Psi _r)=0\)  on the \(\Spin^c\) admissible pair 
\((X, \nu)\) as
before, and let \(\mathbf{C}=[C, \tilde{C}]\) be a t-curve. Given a
compact \(X_\bullet\subset X^{'a}\), denote: 
\[\begin{split}
\mathfrc{d}_{X_\bullet}(\mathfrc{c}_r, \mathbf{C})& :=\dist_{X_\bullet}
(\alpha_r^{-1}(0), C)+
\Big\|\frac{i}{2\pi}F_{A^E_r}-\tilde{C}\Big\|_{op, X_\bullet}.
\end{split}
\]

Recall from Section \ref{sec:gr-spin}(a) the definition of  the set of
t-orbits, \(\bbP(Y_i, \grs_i)\). 
Let 
\[
\mathfrc{d}_i:=\min_{\pmb{\gamma}, \pmb{\gamma}'\in \bbP(\grs_i), \pmb{\gamma}\neq\pmb{\gamma}'}\big(\dist_{Y_i}
(\gamma', \gamma)+\|\tilde{\gamma}'-\tilde{\gamma}\|_{op, Y_i}\big)>0.
\]

\subsubsection*{\it Proof of Theorem \ref{thm:g-conv}: the non-cylindrical
  case.} Assume that \((X, \nu )\)
is non-cylindrical. The proof for the cylindrical case requires only little
modification and will be postponed to the next subsection. 

\subsubsection*{(a):} Let \(\{\fc_r:=(A_r, \Psi _r)\}_{r\in \Gamma _0}\) be the sequence from the
conclusion of Theorem \ref{thm:l-conv} (a). Rename the t-curve
\(\bfC\) in the statement of Theorem \ref{thm:l-conv} by
\(\bfC_0\). That is to say, \(\{(A_r, \Psi _r)\}_{r\in \Gamma _0}\)
t-converges to \(\bfC_0\) over \(X^{'a}\). Given \(i\in \grY_m\), let \(\pmb{\gamma }_{0,i}\) denote the
\(Y_i\)-end limit of \(\bfC_0\). Write the sequence \(\Gamma _0\) as \(\{r_n\}_n\). Given \(l\in [0, \infty]\),
let  \(\vec{l}\co \grY_X\to [0, \infty]\) denote the
function with \(\vec{l}_i=l\) when \(i\in \grY_m\),
\(\vec{l}_i=\ul_i\) when \(i\in \grY_v\).  

Recalling the construction of
\(\bfC_0\) in the proof of Theorem \ref{thm:l-conv} (a) via diagonalization, we choose the  
diagonalization process such that 
the subsequence \(\Gamma _0=\{r_n\}_n\subset \Gamma \) so obtained (\(\Gamma , \Gamma _0\) being
as in the statement of Theorem \ref{thm:l-conv})  satisfies the
following: Let \(\varepsilon _n:=e^{-n}\), \(L_n:=e^n\), \(r_n\) increases with
\(n\), \(r_n\to \infty\) as \(n\to \infty\), and \(\fd_{X_\bullet}(\fc_{r_n}, \bfC_0)\leq\varepsilon _n\) \(\forall
X_\bullet\subset X_{\vec{L}_n}\). We find it convenient to restate  the
preceding statement in a re-parametrized fashion: Let \(n(r)\co \Gamma _0\to \bbN\)
denote the inverse function of the map \(n\mapsto r_n\), and write
\(\varepsilon _T(r):=\varepsilon _{n(r)}\),
\(L_T(r):=L_{n(r)}=\varepsilon _T(r)^{-1}\). Then
\begin{equation}\label{def:L_T}
\fd_{X_\bullet}(\fc_{r}, \bfC_0)\leq\varepsilon _T(r) \quad \forall
X_\bullet\subset X_{\vec{L}_{T}(r)}, r\in \Gamma _0.
\end{equation}

Given \(\varepsilon > 0\), let \(\ell_i(\varepsilon)\geq 0\) be the minimal
\(L\) such that 
\begin{equation}\label{def:ell}
\text{
\(\fd_{\hat{Y}_{i, [l,l+3]}}(\bfC_0, \bbR\times
\pmb{\gamma}_{i,0})\leq \varepsilon \) \(\forall l\geq L\) and \(\int_{\tilde{C}_0|_{\hY_{i, L}}}\nu 
\leq  \varepsilon /3\). }
\end{equation}
Such \(\ell_i\) exists because of Proposition
\ref{prop:t-curve-asymp}, and because of (\ref{eq:xi-exp}).

Let \(R(\varepsilon )>1\) denote the minimal \(r\) such that all 
of the following hold: 
\begin{equation}\label{def:R_e}
\begin{split}
& \text{(i) \( L_T(r)\geq\ell_i(\varepsilon )+3\) \(\forall i\in \grY_m\);}\\
& \text{(ii)  \(\varepsilon _T(r)=L_T(r)^{-1}\leq
(3\|\nu \|_{\infty})^{-1}\varepsilon \);}\\
& \text{(iii) \(\zeta _S r^{-1/4}\leq \varepsilon /3\) and \(r\geq R_S\), where \(R_S\),
  \(\zeta _S\) are as in Proposition \ref{cor:E-cp}.}
\end{split}
\end{equation}
The preceding conditions for \(R(\varepsilon )\), (\ref{def:L_T}), the
condition (\ref{def:ell}) for \(\ell(\varepsilon )\), and Proposition
\ref{cor:E-cp} 
together ensure that 
\begin{equation}\label{borderE}
\text{ \(\forall r\geq R(\varepsilon )\), \(\bar{\grE}_r(\hat{Y}_{i, [L_T(r)-5/2, L_T(r)-1/2]})
\leq \varepsilon \).}
\end{equation}

To proceed, we need the next two preliminary lemmas. 
\begin{lemma}\label{lem:exhausting intervals}
Fix \(i\in \grY_m\). Suppose that \(I_1\subset I_2\subset \cdots \) is a sequence of
intervals in \(\bbR\) so that \(\bigcup_n I_n=:I_\infty\). Assign to
each \(n\in \bbZ^+\) a triple \((r=r_n, \mathfrc{c}_r, m_r)\), so
that:
\begin{itemize}
\item \(\fc_r\)  is a solution to the \(r\)-th version of
  Seiberg-Witten equation\\ \(\grS_{\mu _r, \hat{\grp}}(\fc_r)=0\);
\item \(m_r\in \bbR^{\geq 0}\) is unbounded, and 
\item \(I_n\) is isomorphic by some translation
\(\tilde{\tau}_n\) to an interval
\(I'_n\subset [m_r, \infty)\).
\end{itemize}
Let \(\mathfrc{b}_n:=\tilde{\tau}_{n}^*(\mathfrc{c}_{r_n}|_{\hat{Y}_{i,
    I'_n}})\). Then there is a t-curve \(\mathbf{C}\subset
I_\infty\times Y_i\) for the cylindrical admissible pair
\((I_\infty\times Y_i, \pi_2^*\nu_i)\), and a subsequence of \(\{\mathfrc{b}_n\}_n\),
denoted by the same notation below, which (locally) t-converges to
\(\mathbf{C}\) over \(I_\infty\times Y_i\). Moreover, there is an
upper bound on \(\scrF_\omega(\mathbf{C})\) depending only on
\(\mathfrc{h}\), \(\grs\), and other data mentioned in the end of
Theorem \ref{thm:l-conv}. 
\end{lemma}
\pf 
Passing to a subsequence if necessary, we may assume that \(m_r\) increases with \(r\) and goes to infinity
as \(r\to \infty\). Meanwhile, \(\mathfrc{b}_n=:(\tilde{A}_n, \tilde{\Psi}_n)\) is a solution to the Seiberg-Witten
equation \(\grS_{\tilde{\tau}_{n}^*\mu_r, 0}(\mathfrc{b}_n)=0\) over 
\((\tilde{\tau}_{n})^{-1}\hat{Y}_{i,   I'_n}\simeq I_n\times
Y_i\). Over this, 
\begin{equation}\label{eq:F_KDecay}
\|\tilde{\tau}_{n}^*(\nu|_{\hat{Y}_{i,
    I'_n}})-\pi_2^*\nu_i\|_{C^k}\leq \zeta_i e^{-\kappa_i m_r}\to 0
\end{equation}
as \(n\to \infty\), while \(\|\tilde{\tau}_{n}^*(w_r|_{\hat{Y}_{i,
    I'_n}})-\pi_2^*w_{i,r}\|_{C^k}=0\) by assumption. 

Because of (\ref{def:zeta-E}) 
(and a diagonalization argument), there is a subsequence of
\(\{\mathfrc{b}_n\}_n\) (again denoted by the same notation), and a
current \(\tilde{C}\) on \(I_\infty\times
Y_i\), so that over any compact subset \(V\subset I_k\times
Y_i\) the corresponding sequence of currents
\(\{\frac{i}{2\pi}F_{\td{A}^E_n}|_V\}_{n\geq k}\) converge to
\(\tilde{C}|_V\). Following Steps 5 and 6 of Section  \ref{sec:synopsis}'s
summary of Taubes' arguments, the support of \(\tilde{C}\) is
a \(J_i\)-holomorphic subvariety \(C\), where \(J_i\) is the almost complex
structure on \(I_\infty\times
(Y_i-\nu_i^{-1}(0)\) defined by \((\pi_2^*\nu_i)^+\). Moreover, 
\(\mathbf{C}=(C, \tilde{C})\) is a t-curve for the cylindrical
admissible pair 
\((I_\infty\times Y_i, \pi_2^*\nu_i)\) and \(\mathfrc{b}_n\)
t-converges to \(\mathbf{C}\). 
The claim about the \(\omega\)-energy of \(\mathbf{C}\) follows from
the same argument as in the proof of Theorem \ref{thm:l-conv}. \epf

\begin{lemma}\label{lem:t-close} 
Fix \(i\in \grY_m\) and a 
positive number \(\epsilon\), 
\(\epsilon\leq\mathfrc{d}_i/8\). 
Let \(\{\fc_r\}_{r\in\Gamma _0}\) be the sequence from the conclusion
of Theorem \ref{thm:l-conv}, chosen in the aforementioned manner, and
let \(L_T(r)\) be as in (\ref{def:L_T}).                                 
Then there
exist positive constants \(h_i=h_i(\epsilon)\leq\mathfrc{d}_i/8\), \(\scrR_i=\scrR_i(\epsilon)\) such that the
following holds: For any \(r\geq \scrR_i\), \(r\in \Gamma _0\), and \(l\geq L_T(r)\), if 
\[
\bar{\grE}_r(\hat{Y}_{i, [l+1/2, l+3/2]})
\leq h_i,\]
then \(\mathfrc{d}_{\hat{Y}_{i, [l,
    l+3]}}(\mathfrc{c}_r, \bbR\times \pmb{\gamma})<\epsilon\) for some
t-orbit \(\pmb{\gamma}\) on \((Y_i, \nu_i)\).
\end{lemma}
\pf Suppose the contrary. Then for any pair of positive numbers \(\delta>0\) and \(\rmr\geq 1\), there exist
\(r\geq \rmr\), \(r\in\Gamma _0\), and \(l\geq L_T(r)\) such that 
\begin{equation}\label{eq:con1}
\bar{\grE}_r(\hat{Y}_{i, [l+1/2, l+3/2]})\leq \delta
\end{equation}
and \begin{equation}\label{eq:con2}\mathfrc{d}_{\hat{Y}_{i, [l,
    l+3]}}(\mathfrc{c}_r, \bbR\times \pmb{\gamma})\geq\epsilon\quad \text{for all \(\pmb{\gamma}\in \bbP(\grs_i)\).}
\end{equation}
Choose a sequence of
\((\delta, \rmr)\), denoted \(\{(\delta_n, \rmr_n)\}_n\), with
\(\lim_{n\to \infty}\delta_n\to 0\) and \(\rmr_n\to \infty\) as \(n\to
\infty\), and denote the corresponding  \(r, l\) by
\(r'_n\), \(l_n\) respectively. Note that \(l_n\to \infty\)
as \(n\to \infty\), since \(l_n\geq L_T(r'_n)\) and \(L_T(r'_n)\to
\infty\) as \(n\to \infty\). Apply Lemma
\ref{lem:exhausting intervals} with \(I_1=I_2\cdots =I_\infty=[0,3]\),
and with \(I_n'=[l_n, l_n+3]\). This gives a t-curve \(\mathbf{C}=[C, \td{C}]\) on
the cylindrical admissible
pair \((\hat{Y}_{i, [0,
    3]}, \pi_2^*\nu_i)\) and a subsequence of 
\(\{(\delta_n, \rmr_n)\}_n\) (denoted by the same
notation), such that with \(\fc_n':=\fc_{r_n}\), 
\begin{equation}\label{lim_C}
\mathfrc{d}_{\hat{Y}_{i, [0,
    3]}}(\tau^*_{l_n}(\mathfrc{c}'_n|_{\hat{Y}_{i, [l_n,
    l_n+3]}}), \mathbf{C})=\mathfrc{d}_{\hat{Y}_{i, [l_n,
    l_n+3]}}(\mathfrc{c}'_n, \tau_{-l_n}\mathbf{C})\to 0\quad \text{as \(n\to \infty\). }
\end{equation}
Combining with Proposition \ref{cor:E-cp}, 
 this implies that 
\[\begin{split}
0\leq  & \int_{-1/2}^{1/2}\scrF_\nu (\tilde{C}|_{\hY_{i, [1/2+t, 3/2+t]}})\,
dt\\
 &\quad  =\lim_{n\to \infty}\bar{\scrV}_{r_n} (\hY_{i, [l_n+1/2, l_n+3/2]})\\
& \quad =\lim_{n\to
  \infty}\bar{\grE}_{r_n} (\hY_{i, [l_n+1/2, l_n+3/2]})=0. 
\end{split}
\]
The \(v\)-energy \(\scrF_\nu \) for t-curves being nonnegative, this implies that\\
\(\scrF_\nu (\tilde{C}|_{\hY_{i, [1, 2]}})=0\). But by Lemma
\ref{lem:F-curve}, this in turn implies that 
\(\mathbf{C}|_{\hY_{i, [1, 2]}}=[1,2]\times \pmb{\gamma}\subset [1,2]\times Y_i\simeq\hY_{i, [1, 2]}\) for some \(\pmb{\gamma}\in
\bbP(\grs_i)\), and hence \(\mathbf{C}=[0,3]\times \pmb{\gamma}\subset
\hY_{i, [0, 3]}\). This, together with
(\ref{lim_C}), contradict the
assumption (\ref{eq:con2}).
\epf

To continue with the proof of Theorem \ref{thm:g-conv}, fix \(\epsilon \leq\mathfrc{d}_i/8\), and let
\(\varepsilon _0:=h_i(\epsilon )\), where \(h_i(\epsilon )\) is as in Lemma \ref{lem:t-close}. 
Divide \([L_T(r)-5/2, \infty)\) into
intervals \(I^*\) of length 2. For each \(i\in \grY_m\), let \(\invbreve{\scrI}_{i,r}\subset
[L_T(r)-5/2, \infty)\subset \bbR^+\) be the union of those \(I^*\)
satisfying 
\begin{equation}
\bar{\grE}_r(\hat{Y}_{i, I^*})\geq\varepsilon _0.
\end{equation}
Let \(\breve{\scrI}_{i,r}\) be the closure of \([L_T(r)-5/2,
\infty)-\invbreve{\scrI}_{i,r}\). 
By Lemma
\ref{lem:t-close}, this implies that when \(r\geq \scrR_i(\epsilon
)\), for each of the length-2
interval,  \(I^*=[l-5/2, l-1/2]
\subset\breve{\scrI}_{i,r}\), there is a t-orbit \(\pmb{\gamma }_*\in
\bbP(\grs_i)\) so that 
\begin{equation}\label{eq:7.20a}
  \fd_{\hY_{i, [l-3, l]}}
  (\fc_r, \bbR\times \pmb{\gamma }_*)\leq \epsilon .
\end{equation}
Since
\(\bbP(\grs_i)\) is discrete, this t-orbit \(\pmb{\gamma
}_*\) only depends on the connected component of
\(\breve{\scrI}_{i,r}\) that \(I^*\) lies in. 
Let \(\Lambda _{i,r}+1\) be the number of connected components of
\(\breve{\scrI}_{i,r}\). The set of 
connected components of \(\breve{\scrI}_{i,r}\) is  naturally ordered,
and we denote it by \(\{\breve{\scrI}_{i,r}^{(k)}\}_{k=0}^{\Lambda _{i,r}}\):  for any pair \(k<k'\) and 
any \(l\in \breve{\scrI}_{i,r}^{(k)}\) and \(l'\in
  \breve{\scrI}_{i,r}^{(k')}\),   \(l<l'\). Denote the t-orbit
associated to the component \(\breve{\scrI}^{(k)}_{i,r}\) by
\(\pmb{\gamma}_{i,r}^{(k)}\). By (\ref{borderE}), the interval 
  \([L_T(r)-5/2, L_T(r)-1/2]\subset\breve{\scrI}_{i,r}^{(0)}\), and it follows from the
  definition of \(R(\varepsilon )\) that 
  \(\pmb{\gamma}_{i,r}^{(0)}=\pmb{\gamma }_{0, i}\) \(\forall r\geq
  R(\varepsilon _0)\). Let \(r_0\) denote the maximum among
  \(R(\varepsilon _0)\) and \(\scrR_i(\epsilon )\) \(\forall i\in
  \grY_m\), and assume that \(r\geq r_0\) from now
  on. 

Correspondingly, label the connected components of \(\invbreve{\scrI}_{i,r}\) by
  \(\invbreve{\scrI}_{i,r}^{(k)}\), \(k\in \bbZ^+\), so that \(\invbreve{\scrI}_{i,r}^{(k)}\) lies between
\(\breve{\scrI}_{i,r}^{(k-1)}\) and
\(\breve{\scrI}_{i,r}^{(k)}\). We claim that each
\(\invbreve{\scrI}_{i,r}^{(k)}\) has finite length, and the number of connected components of
\(\invbreve{\scrI}_{i,r}\), 
namely \(\Lambda _{i,r}\), is bounded independently of \(r\). 
To see this, write:
\begin{equation}\label{eq:grE}\begin{split}
\bar{\grE}_r(X^{'a})(\mathfrc{c}_r)&
=\bar{\grE}_r(X^{'a}_{\vec{L}_T(r)-\vec{\frac{5}{2}}})(\mathfrc{c}_r)\\
& \quad +
\sum_{i\in \grY_m}\Big(\bar{\grE}_r(\hat{Y}_{i,
  \breve{\scrI}_{i,r}})(\mathfrc{c}_r)+\sum_{k=1}^{\Lambda _{i,r}}\bar{\grE}_r(\hat{Y}_{i,
  \invbreve{\scrI}_{i,r}^{(k)}})(\mathfrc{c}_r)\Big)
\end{split}
\end{equation}
By Proposition \ref{cor:E-cp} and (\ref{def:L_T}),  \(\forall r\geq r_0\) 
\[\begin{split}
\bar{\grE}_r(X^{'a}_{\vec{L}_T(r)-\vec{\frac{5}{2}}})(\mathfrc{c}_r)& \geq
\bar{\scrV}_r(X^{'a}_{\vec{L}_T(r)-\vec{\frac{5}{2}}})(\mathfrc{c}_r)-\zeta  _S\, r^{-1/4}\\
& \geq \int_{-1/2}^{1/2}\scrF_\nu
(\tilde{C}_0|_{X^{'a}_{\vec{L}_T(r)-\vec{\frac{5}{2}}+\vec{t}}})\, dt-\|\nu
\|_\infty \varepsilon _T(r)-\zeta  _S\, r^{-1/4}\\
& \geq -\|\nu
  \|_\infty \, \varepsilon _T(r)-\zeta  _S\, r^{-1/4}.
\end{split}
\]
By definition, for each pair of \(i, k\), \(\bar{\grE}_r(\hat{Y}_{i,
  \invbreve{\scrI}_{i,r}^{(k)}})(\mathfrc{c}_r)\geq |\invbreve{\scrI}_{i,r}^{(k)} |\varepsilon _0/2\). By Lemma
\ref{lem:E-lbdd}, for each \(i\in
\grY_m\)
\[
\bar{\grE}_r(\hat{Y}_{i,
  \breve{\scrI}_{i,r}})(\mathfrc{c}_r)\geq-\zeta _i'' e^{-\kappa _iL_T(r)}
\]
Combine the above bounds with (\ref{eq:grE}) to get:
\begin{equation}\label{bdd:Lambda_i}
\begin{split}
& \sum_{i\in \grY_m}\sum_{k=1}^{\Lambda _{i,r}}\, |\invbreve{\scrI}_{i,r}^{(k)}
|\, \varepsilon _0/2 \\
&\qquad  \leq  
\sum_{i\in \grY_m}\sum_{k=1}^{\Lambda _{i,r}}\bar{\grE}_r(\hat{Y}_{i,
  \invbreve{\scrI}_{i,r}^{(k)}})(\mathfrc{c}_r)\\
&\qquad \leq  \bar{\grE}_r(X^{'a})(\mathfrc{c}_r)+\|\nu
\|_\infty \varepsilon _T(r)+\zeta  _Sr^{-1/4}
+\sum_{i\in
    \grY_m}\zeta _i'' e^{-\kappa _iL_T(r)}. 
\end{split}
\end{equation}
By Lemma \ref{lem:grE-X}, the right hand side of the preceding
inequality has an \(r\)-independent upper bound. This in turn gives an \(r\)-independent upper
bound on the left hand side for all \(r\geq r_0\). Recalling that every
\(|\invbreve{\scrI}_{i,r}^{(k)} |\) is a positive even integer, this
gives an upper bound for each \(\Lambda _{i,r}\) and \(|\invbreve{\scrI}_{i,r}^{(k)} |\).

The preceding uniform bounds on \(\Lambda _{i,r}\) enables us to pass
to a subsequence \(\Gamma _1\) of \(\Gamma _0\), 
such that \(\Lambda_{i,r}=\Lambda_i\) \(\forall r\in
\Gamma_1\) for some \(\Lambda_i\in  \bbZ^{\geq 0}\). Given \(k, k'\), \(1\leq k<k'\leq \Lambda _i\), Let
\(\fm_{i,r}(k, k')\) denote the distance between the centers of
intervals 
\(\invbreve{\scrI}_{i,r}^{(k)}\) and
\(\invbreve{\scrI}_{i,r}^{(k')}\). Let
\(\fm_{i,r}(0,k)\) denote the distance between the center of
\(\invbreve{\scrI}_{i,r}^{(k)}\) and \(L_T(r)\). 
We choose the subsequence \(\Gamma
_1\) so that:
\BTitem\label{def:Gamma1}
\item For every \(i\) and \(k\in \Lambda _i\), the t-orbit \(\pmb{\gamma
  }_{i,r}^{(k)}\) is independent of \(r\). This is possible because
  \(\bbP(\grs_i)\) consists of finitely many elements. 
\item Given \(i\in \grY_m\), and any pair \(k, k'\), \(0\leq k<k'\leq
  \Lambda _i\), the distance 
\(\fm_{i,r}(k, k')\) is non-decreasing with respect to \(r\). 
\ETitem
Fix \(i\in \grY_m\). Let \(k_1\) be the smallest \(k\) such that \(\fm_{i,r}(0, k)\) is 
unbounded with respect to \(r\). Let \(\grI_1\) denote the finite interval  \(\grI_1\supset
\invbreve{\scrI}_{i,r}^{(1)}\cup  \invbreve{\scrI}_{i,r}^{(k_1)}\) with \(\partial\grI_1\subset
\partial\invbreve{\scrI}_{i,r}^{(1)}\cup
\partial\invbreve{\scrI}_{i,r}^{(k_1)}\) if \(k_1\) exists; otherwise
set \(\grJ_1=\emptyset\). Given
any pair \(k\), \(k'\) such that 
\(k<k'\leq \Lambda _i\) and \(\fm_{i,r}(k, k')\) is 
bounded with respect to \(r\), we form a finite interval \(\grI_*\supset
\invbreve{\scrI}_{i,r}^{(k)}\cup  \invbreve{\scrI}_{i,r}^{(k')}\) with \(\partial\grI_*\subset
\partial\invbreve{\scrI}_{i,r}^{(k)}\cup  \partial\invbreve{\scrI}_{i,r}^{(k')}\). Given
\(i\) and \(r\), let \(\grJ_{i,r}\subset [L_T(r)-1/2, \infty)\) denote
\(\bigcup_*\grJ_*\cup\grJ_1\cup\bigcup_{k}\invbreve{\scrI}_{i,r}^{(k)}\),
where the index \(k\) runs over all \(k\) such that 
 \(
 \invbreve{\scrI}_{i,r}^{(k)}\subset[L_T(r)-1/2,\infty)-\bigcup_*\grJ_*-\grJ_1\). 
Let \(\grK_i\) denote the number of connected components of
\(\grJ_{i,r}\), and when \(\grK_i>0\), let
\(\{\grI_{i,r}^k\}_{k=1}^{\grK_i}\) denote the set of connected components of
\(\grJ_{i,r}\). These connected components are ordered such that 
\(\grI_{i,r}^{k'}\) lies to the right of \(\grI_{i,r}^k\subset
[L_T(r)-1/2, \infty)\) for any \(k, k'\) such that \(1\leq k<k'\leq \grK_i\). For any
\(1\leq k<\grK_i\), the intervals \(\grI_{i,r}^k\)
and \(\grI_{i,r}^{k+1}\) are separated by one of the connected
components of \(\breve{\scrI}_{i,r}\), say \(\breve{\scrI}_{i,r}^{(k')}\). 
We rename the
t-orbit associated to this connected component,  \(\pmb{\gamma
}_{i,r}^{(k')}\), by \(\pmb{\gamma
}_{i}^k\). (Recall the second bullet of (\ref{def:Gamma1}).) Meanwhile, \(\grJ_{i,r}^1\) is adjacent to the first
connected component of \(\breve{\scrI}_{i,r}\), i.e.
\(\breve{\scrI}_{i,r}^{(0)}\),  on the left; and \(\grJ_{i,r}^{\grK_i}\)
is adjacent to the last connected component of
\(\breve{\scrI}_{i,r}\), i.e. \(\breve{\scrI}_{i,r}^{(\Lambda _i)}\), which is a half-infinite interval of the form
\([L_+, \infty)\). The t-orbits associated to
\(\breve{\scrI}_{i,r}^{(0)}\) and \(\breve{\scrI}_{i,r}^{(\Lambda
  _i)}\) are respectively \(\pmb{\gamma }_{0,i}\) and \(\pmb{\gamma
}_i\). (\(\pmb{\gamma }_i\) is as in Condition (1) in the
statement of Theorem \ref{thm:l-conv}.) Let \(\pmb{\gamma
}_{i}^0:=\pmb{\gamma }_{0,i}\) and \(\pmb{\gamma
}_{i}^{\grK_i}:=\pmb{\gamma }_{i}\). Note that when \(\grK_i=0\),
\(\grJ_{i,r}=\emptyset\) and \(\pmb{\gamma }_{0,i}=\pmb{\gamma }_{i}\).

Recall also that by (\ref{bdd:Lambda_i}), the length 
\(|\invbreve{\scrI}_{i,r}^{(k)}|\) has an upper bound independent of
\(r\), \(i\), and \(k\). This means that there is a finite number
\(\fm \geq 1\) so that \(|\grI_{i,r}^{k}|\leq 2\fm\) for all \(i, k\)
and \(r\). 
Let \(L^i_{k,r}\in [L_T(r)-1/2, \infty)\) denote the mid point of the
interval \(\grI_{i,r}^{k}\), and let \(\grY'_m=\{i\, |\, i\in \grY_m, \grK_i>0\}\). By construction, we have: 
\BTitem\label{def:L_kr}
\item Given \(i\in \grY'_m\), for every \(k\) with \(1\leq k\leq \grK_i-1\), 
 \(\Delta ^i_{k,r}:=L^i_{k+1,r}-L^i_{k,r}\)  is non-decreasing with respect to
 \(r\). Also, \(\Delta ^i_{0,r}:=L^i_{1,r}-L_T(r)\)  is non-decreasing with respect to
 \(r\). (cf. the second bullet in (\ref{def:Gamma1}).)
 \item Given \(i\in \grY'_m\), for every \(k\) with \(0\leq k\leq \grK_i-1\), 
 \(\lim_{r\to \infty}\Delta ^i_{k,r}\to\infty\) as \(r\to \infty\). 
\ETitem
In particular, there exists an \(R_1\geq R_0\) such that 
\(\Delta ^i_{k,r}>2\fm +4\) for all \(i\in \grY'_m\), \(k\in \{0,
\ldots, \grK_i-1\}\), \(\forall
r\geq R_1\). Assume that \(r\geq R_1\) from now on. 

Fix \(i\in \grY'_m\).  
Given \(r\) and \(k\), \(1\leq k\leq \grK_i\), let
\(I^{'i}_{k,r}\supset \grJ^k_{i,r}\) be defined by 
\[
I^{'i}_{k,r}:=\begin{cases}
[L^i_{k-1, r}+\fm, L^i_{k+1,
  r}-\fm] & \text{when \(2\leq k\leq \grK_i-1\)};\\
[L_T(r)-3, L^i_{2,
  r}-\fm] & \text{when \(k=1\)};\\
[L^i_{\grK_i-1, r}+\fm, \infty)  & \text{when \(k=\grK_i\)}.
\end{cases}
\]
Let \(I_{k,r}^i:=\tau
_{-L^i_{k,r}}I^{'i}_{k,r}\). That is, 
\[
I^i_{k,r}:=\begin{cases}
[-\Delta ^i_{k-1, r}+\fm, \Delta ^i_{k,r}-\fm] & \text{when \(2\leq k\leq \grK_i-1\)};\\
[-\Delta ^i_{0,r}-3, \Delta ^i_{1,r}-\fm] & \text{when \(k=1\)};\\
[-\Delta {}^i_{\grK_i-1, r}+\fm, \infty)  & \text{when \(k=\grK_i\)}.
\end{cases}
\]
Then by (\ref{def:L_kr}),
for every fixed \(1\leq k\leq \grK_i\), \(I^i_{k,r}\subset I^i_{k,r'}\) when
\(r'>r\), and together they 
form a  nested sequence of intervals exhausting \(\bbR=I^i_{k,
  \infty}\). Let \[\begin{split}
    & \breve{\grI}^i_{k,r}:=\\
    & \quad \begin{cases}
I^{'i}_{k,r}\cap I^{'i}_{k+1,r}=[L_{k, r}^i+\fm, L_{k+1, r}^i-\fm]\subset
\breve{\scrI}_{i,r}& \text{when \(1\leq k\leq \grK_i-1\);} \\
[L_T(r)-5/2, L_{1, r}^i-\fm]\subset
\breve{\scrI}_{i,r}^{(0)}&
\text{when \(k=0\);} \\
[L_{\grK_i, r}^i+\fm,
\infty) \subset
\breve{\scrI}_{i,r}^{(\Lambda_i)}& \text{when \(k=\grK_i\).} \\
\end{cases}
\end{split}
\]
Then by construction (recall (\ref{eq:7.20a})), for \(r\geq R_1\), 
\begin{equation}\label{chain-ends}
\fd_{\hY_{i, I^*}}
(\fc_r, \bbR\times \pmb{\gamma }_i^k)\leq \epsilon \quad \forall\, 
I^*\subset\breve{\grI}_{k,r}^i, \,  |I^*|=2. 
\end{equation}
Consider now the sequence \(\{\mathfrc{b}_r:=\tau
_{-L^i_{k,r}}\fc_r|_{\hY_{i, I_{k,r}'}}\}_{r\in \Gamma _1}\). \(\mathfrc{b}_r\) is defined
on \(I_{k,r}^i\times Y_i\). Apply Lemma
\ref{lem:exhausting intervals} to \(\{\mathfrc{b}_r\}_{r\in \Gamma
  _1}\) to get a t-curve \(\bfC_{i,k}\) on \((\bbR\times Y_i, \pi
_2^*\nu _i)\) and a subsequence \(\Gamma ^i_k\) of \(\Gamma _1\) so
that \(\{\mathfrc{b}_r\}_{r\in \Gamma
  _k^i}\) t-converges to \(\bfC_{i,k}\). More specifically, the
subsequence \(\Gamma ^i_k\) is chosen via diagonalization so that 
\begin{equation}\label{appr:b-C}
\fd_{I^*\times Y_i} (\mathfrc{b}_{r_n}, \bfC_{i,k})\leq  e^{-n}
\quad \forall 
 I^*\subset I^i_{k,r_n},  \forall n, 
\end{equation}
where \(r_n\) now denotes the \(n\)-th element in the sequence \(\Gamma
^i_k=\{r_n\}_{n\in \bbZ^+}\). Comparing the preceding inequality with
(\ref{chain-ends}), we see that the \(-\infty\)-limit of
\(\bfC_{i,k}\) is \(\pmb{\gamma }_i^{k-1}\), and its \(+\infty\)-limit
is \(\pmb{\gamma }_i^{k}\). Thus, when \(\grK_i>0\), \(\{\bfC_{i,k}\}_{k=1}^{\grK_i}\)
forms a chain of t-curves on \((\bbR\times Y_i, \pi _2^*\nu _i)\), \(\grC_i\),  with \(-\infty\)-limit \(\pmb{\gamma
}_i^0=\pmb{\gamma }_{i,0}\), and with \(+\infty\)-limit \(\pmb{\gamma
}_i^{\grK_i}=\pmb{\gamma }_{i}\). When \(\grK_i=0\),
\(\breve{\scrI}_{i,r}=[L_T(r)-5/2, \infty)\) is connected; by
(\ref{eq:7.20a}), 
\(\pmb{\gamma }_{0,i}=\pmb{\gamma }_i,\) and 
\[
\fd_{\hY_{i, I^*}}
(\fc_r, \bbR\times \pmb{\gamma }_i)\leq \epsilon \quad \forall\, 
I^*\subset[L_T(r)-5/2, \infty), \,  |I^*|=2. 
\]
In this case, let \(\grC_i\) be the 0-component chain of t-curves with
a single rest orbit \(\pmb{\gamma }_{0,i}=\pmb{\gamma }_i\). Now 
\(\grC:=\{\bfC_0, \{\grC_i\}_{i\in
  \grY_m}\}\) forms a chain of t-curves on \((X^{'a}, \nu )\),
with \(Y_i\)-end limit \(\pmb{\gamma }_i\) for each \(i\in \grY_m\). Finally, choose an order of the finite set
\[\{(i, k)\, | \, \text{\(i\in \grY_m\) with \(\grK_i>0\)} , \, 1\leq k\leq \grK_i\},\] and choose
consecutive subsequences \(\Gamma ^i_k\subset \Gamma _1\) in this order:
\(\Gamma ^i_k\subset \Gamma ^{i'}_{k'}\) whenever \((i, k)>(i',
k')\) to get a subsequence \(\Gamma '=\bigcap_{\{(i, k)\}}\Gamma
^i_k\subset \Gamma _1\). Then \(\Gamma
'\) and \(\grC\) together satisfies Item (a) asserted by Theorem
\ref{thm:g-conv}.


\subsubsection*{(b) and (d):} Recall that
\(\scrF_\omega(\grC)=\scrF_\omega(\bfC_0)+\sum_{i\in
  \grY'_m, \grK_i>0}\sum_{k=1}^{\grK_i}\scrF_\omega(\bfC_{i,k}) \). We already
provided a bound on \(\scrF_\omega(\bfC_0)\) in 
the proof of Theorem \ref{thm:l-conv} (b). Thus, to verify Item (b) in
the statement of the theorem, it suffices to bound
each \(\scrF_\omega(\bfC_{i,k})\). As in the proof of Theorem
\ref{thm:l-conv} (b), this reduces to the following two tasks: 
(i) obtaining a \(l\)-independent bound on
\(\scrF_\omega(\bfC_{i,k}|_{[l, l+1]})\); (ii) bounding \(\scrF_\nu
(\bfC_{i, k})\). 

(i): Fix \(i\) with \(\grK_i>0\) and \(k\). Given \(l\in \bbR\), there exists \(\grr_l\geq R_1\)
such that \([l, l+1]\subset I^i_{k, r}\) \(\forall r\geq \grr_l\). Let
\(r_n\) denote the \(n\)-th element of \(\Gamma _k^i\), and consider an \(r_n\geq\grr_l\). Then
according to (\ref{appr:b-C}) and (\ref{bdd:L^14d}), 
\[\begin{split}
& \scrF_\omega(\bfC_{i,k}|_{[l, l+1]})-\frac{1}{2}c_1(\grs_i)\cdot
[\theta ]+\zeta _\theta  \\
& \qquad \leq \int_{\hY_{i,
    [L^i_{k,r_n}+l, L^i_{k,r_n}+l+1]}}\frac{iF_{A_{r_n}^E}}{2\pi
}\wedge \nu _i+e^{-n}\|\nu _i\|_\infty\\
& \qquad \leq \zeta '+e^{-n}\|\nu _i\|_\infty,
\end{split}
\]
where \(\theta =*\nu _i\) and \(\zeta _\theta \) is as in Proposition
\ref{prop:t-conv3d}, and \(\zeta '>0\) is independent of \(n\)
(equivalently, \(r=r_n\)), \(k\), and
\(l\). This gives us an \(l\)-independent bound
\begin{equation}\label{bdd:F_wik}
\scrF_\omega(\bfC_{i,k}|_{[l, l+1]})\leq \frac{1}{2}c_1(\grs_i)\cdot
[\theta ]+\zeta '-\zeta _\theta  .
\end{equation}

(ii):  
When \(\grK_i>1\), let \(l^i_{k,r}:=(L^i_{k,
  r}+L^i_{k+1, r})/2\) denote the mid point of the interval
\(\breve{\grI}^i_{k,r}\) for  \(1\leq k\leq \grK_i-1\). Define \(J^{'i}_{k,r}\),
\(I^{'i}_{k,r}\supset J^{'i}_{k,r}\supset\grI^k_{i,r}\) by
\[
J^{'i}_{k,r}:=\begin{cases} [l^i_{k-1,r}, l^i_{k,r}]& \text{when
    \(\grK_i>1\) and \(2\leq k\leq
  \grK_i-1\);}\\
[L_T(r)-\frac{5}{2}, l^i_{1,r}]& \text{when   \(\grK_i>1\) and \(k=1\);}\\
[l^i_{\grK_i-1,r}, \infty\big) & \text{when   \(\grK_i>1\) and \(k=
  \grK_i\);}\\
[L_T(r)-\frac{5}{2},\infty)& \text{when   \(\grK_i=1\) and \(k=1\).}
\end{cases}
\]

Let \(J^{i}_{k,r}=\tau
_{-L^i_{k,r}}J^{'i}_{k,r}\). In other words, 
\[
J^i_{k,r}:=\begin{cases} \big[-\frac{\Delta ^i_{k-1,r}}{2},\frac{
\Delta ^i_{k,r}}{2}\big] & \text{when \(\grK_i>1\) and \(2\leq k\leq
  \grK_i-1\);}\\
\big[-\Delta ^i_{0,r}-\frac{5}{2},
\frac{\Delta ^i_{1,r}}{2}\big] & \text{when when \(\grK_i>1\) and \(k=1\);}\\
\big[-\frac{\Delta ^i_{\grK_i-1,r}}{2}, \infty\big) & \text{when \(\grK_i>1\) and \(k=
  \grK_i\).}\\
\big[-\Delta ^i_{0,r}-\frac{5}{2},
\infty) & \text{when when \(\grK_i=1\) and \(k=1\);}\\
\end{cases}
\]
Then \(J^i_{k,r}\subset I^i_{k,r}\), yet
according to (\ref{def:L_kr}),
\(\{J^i_{k,r}\}_{r\in \Gamma '}\) still forms a nested sequence of
intervals exhausting \(\bbR\). Meanwhile, for every fixed \(i\in \grY'_m\) and
\(r\),  \(\bigcup_{k=1}^{\grK_i}J^{'i}_{k,r}=[L_T(r)-5/2, \infty)\),
and the interiors of \(J^{'i}_{k,r}\) are mutually disjoint. Now
\begin{equation}\label{eq:grEI}
\begin{split}
  & \bar{\grE}_r(X^{'a}) \\
  &=\bar{\grE}_r(X^{'a}_{\vec{L}_T(r)-\vec{\frac{5}{2}}})+
\sum_{i\in
  \grY'_m}\sum_{k=1}^{\grK_i}\bar{\grE}_r(\hat{Y}_{i,J^{'i}_{k,r}})+\sum_{i\in
  \grY_m-\grY'_m}\bar{\grE}_r(\hat{Y}_{i, L_T(r)-\frac{5}{2}})
\\
& =\bar{\scrV}_r(X^{'a}_{\vec{L}_T(r)-\vec{\frac{5}{2}}})+
\sum_{i\in
  \grY'_m}\sum_{k=1}^{\grK_i}\bar{\scrV}_r(\hat{Y}_{i,J^{'i}_{k,r}})+\sum_{i\in
  \grY_m-\grY'_m}\bar{\scrV}_r(\hat{Y}_{i, L_T(r)-\frac{5}{2}})\\
& \qquad +\scrO(r^{-1/4})\\
& =\bar{\scrF}_\nu
(\bfC_0|_{X^{'a}_{\vec{L}_T(r)-\vec{\frac{5}{2}}}})+\sum_{i\in
  \grY'_m}\sum_{k=1}^{\grK_i}\bar{\scrF}_\nu
(\bfC_{i, k}|_{J^i_{k,r}\times Y_i})+\scrO(e^{-n_r})+\scrO(r^{-1/4})\\
& \qquad +
\sum_{i\in
  \grY'_m}\sum_{k=1}^{\grK_i}\int_{-1/2}^{1/2}\int_{\tau _t \hY_{i,
  J^{'i}_{k,r}}}\frac{i}{2\pi }F_{A^E}\wedge(\nu -\pi _2^*\nu _i)\,
dt\\
& \qquad +
\sum_{i\in\grY_m-
  \grY'_m}\int_{-1/2}^{1/2}\int_{\hY_{i,
  L_T(r)-\frac{5}{2}+t}}\frac{i}{2\pi }F_{A^E}\wedge(\nu -\pi _2^*\nu _i)\,
dt,\\
\end{split}
\end{equation}
where \(r\in \Gamma '\) and \(n_r\in\bbZ^+\) are related as follows:
\(r\) is the \(n_r\)-th element in \(\Gamma '\). In the above, the
second equality follows from 
Proposition \ref{cor:E-cp}, and the third equality is a consequence of
the construction of \(\Gamma '\) and \(\bfC_{i,k}\); see in particular 
(\ref{def:L_T}) and  (\ref{appr:b-C}). 
Using (\ref{bdd:L^14d}) and (\ref{eq:xi-exp}) to estimate the last two
terms above, we have:
\begin{equation}\label{eq:grEJ}
\begin{split}
\bar{\grE}_r(X^{'a}) & =\bar{\scrF}_\nu
(\bfC_0|_{X^{'a}_{\vec{L}_T(r)-\vec{\frac{5}{2}}}})+\sum_{i\in
  \grY_m'}\sum_{k=1}^{\grK_i}\bar{\scrF}_\nu
(\bfC_{i, k}|_{J^i_{k,r}\times Y_i})\\
& \qquad \quad +\scrO(e^{-n_r})+\scrO(r^{-1/4})+\sum_{i\in
  \grY_m} \scrO(e^{-\kappa _i L_T(r)})
\end{split}
\end{equation}
Together with Lemma \ref{lem:grE-X}, the preceding identity  then
implies that 
\[\bar{\scrF}_\nu
(\bfC_0|_{X^{'a}_{\vec{L}_T(r)-\vec{\frac{5}{2}}}})+\sum_{i\in
  \grY'_m}\sum_{k=1}^{\grK_i}\bar{\scrF}_\nu
(\bfC_{i, k}|_{J^i_{k,r}\times Y_i})\leq \frac{\smE}{\pi }+\zeta 
\]
for all sufficiently large \(r\). 
Note that
as 
functions of \(r\), \(\bar{\scrF}_\nu
(\bfC_0|_{X^{'a}_{\vec{L}_T(r)-\vec{\frac{5}{2}}}})\), \(\bar{\scrF}_\nu
(\bfC_{i, k}|_{J^i_{k,r}\times Y_i})\) are all  nondecreasing, nonnegative functions, and recall
that \(\{J^i_{k,r}\}_r\) exhausts \(\bbR\) when \(i\in \grY_m'\). Consequently,\\ \(\lim_{r\to\infty}\bar{\scrF}_\nu
(\bfC_0|_{X^{'a}_{\vec{L}_T(r)-\vec{\frac{5}{2}}}})=\bar{\scrF}_\nu(\bfC_0)\), 
and for all \(i\in \grY'_m, k\), the limit
\(\lim_{r\to \infty}\bar{\scrF}_\nu
(\bfC_{i, k}|_{J^i_{k,r}\times Y_i})\) exists and equals \(\scrF_\nu
(\bfC_{i, k})\). 
Take the \(r\to \infty\) limit on both sides of the 
identity (\ref{eq:grEJ}) and apply Lemma \ref{lem:grE-X}. We get: 
\begin{equation}\label{bdd:F_vgrC}
\bar{\scrF}_\nu(\grC):=\bar{\scrF}_\nu
(\mathbf{C}_0)+\sum_{i\in
  \grY'_m}\sum_{k=1}^{\grK_i}\scrF_\nu
(\bfC_{i, k})\leq \smE/\pi +\zeta , 
\end{equation}
where \(\zeta \) is as in Lemma \ref{lem:grE-X}. Also,  when \(X\) has
no vanishing ends, by (\ref{eq:bbE}) 
\begin{equation}\label{eq:F_vgrC}
\begin{split}
\scrF_\nu(\grC)& =\bar{\scrF}_\nu(\grC)\\
& =\lim_{r\to \infty}\big( (\pi  r)^{-1}\scrE_{top}^{\mu _r}(A_r, \Psi _r)\big)-\int_{X}
\frac{iF_{A_K}}{4\pi }\wedge \nu\\
& = \frac{i}{2\pi }\int_{X}F_{A_0^E}\wedge\nu 
+\frac{1}{\pi }\sum_{i\in \grY_m}[\nu _i]\cdot \jmath_h (\pmb{\gamma }_i)+ \frac{1}{2\pi }i^*[\nu ]\cdot h_{A_0}(\fh),
\end{split}
\end{equation}
as claimed in Item (d) in the statement of the theorem. Meanwhile,
Item (b) in the statement of the theorem follows from a combination of
(\ref{bdd:F_vgrC}), (\ref{bdd:F_wik}), and Theorem \ref{thm:l-conv}
(b). 

\subsubsection*{(c):} Let 
\(l^i_{0,r}:=L_T(r)-3/2\). Let
\(u_{r}\in \scrG\) be such that
\((A_r', \Psi _r'):=u_{r}\cdot ( A_r,\Psi _r)\) is in a temporal gauge on \(\hY_{i,
  \breve{\grI}^i_{k,r}}\)  \(\forall i\in \grY_m\) and \(k\in \{0, 1, \ldots,
\grK_i\}\). 
Given \(i\in \grY'_m\)  and \(k\in \{0, 1, \ldots,
\grK_i-1\}\), let \(B^i_{k,r}\in \Conn (Y_i)\), \(\Phi ^i_{k,r}\in \Gamma (\bbS_i)\) be
\[
B^i_{k,r}:=B_{0}+\int_{-1/2}^{1/2}(B_r(l^i_{k,r}+t)-B_0)\, dt ,\quad \Phi ^i_{k,r}:=\int_{-1/2}^{1/2}\Phi _r(l^i_{k,r}+t)\, dt,
\] 
where
\((B_r(s), \Phi _r(s))=(A_r', \Psi '_r)|_{Y_{i: s}}\).  
When \(i\in \grY_v\), set \(\grK_i=0\), and for all \(i\in \grY\), let
\(\grc^i_{\grK_i,r}:=\grc_{i,r}\). When \(i\in \grY'_m\), 
let \(\grc^i_{k,r}\in \scrB_{Y_i}\) denote the gauge-equivalence class
of \((B^i_{k,r}, \Phi ^i_{k,r})\) when \(k\in \{0, 1, \ldots,
\grK_i-1\}\). 
Recall that given \(\grc=[(B, \Phi )]\in \scrB_{Y_i}\), \(\td{\grc}\)
denotes the 1-current 
\[
\td{\grc}=
\frac{iF_{B^{E}}}{2\pi}=\frac{iF_{B}}{4\pi}-\frac{iF_{B^{K_i}}}{4\pi}, 
\]
where \(K_i^{-1}\) denotes the anti-canonical bundle over \(Y_i\) defined
by \(\nu _i\), and \(B^{K_i}\) denotes the connection on \(K_i^{-1}\) determined by
the Levi-Civita connection. 
Let \([x]_a\subset \bbR\) denote the length \(a\)  interval centered
at \(x\). Given a 1-form \(\mu \) on
\(Y_i\), the evaluation of the exact 1-current
\(\td{\grc}_{k,r}^i-\td{\gamma }_i^k\)
 on \(\mu \) is equal to the evaluation of the 2-current 
\(\frac{iF_{A_r}}{4\pi }-\frac{iF_{B^{K_i}}}{4\pi }-\pi _2^*\td{\gamma
}_i^k\) on the 2-form \(\pi _2^*\mu \wedge ds\) over \(\hY_{i,
  [l^i_{k,r}]_1}\). The absolute value of this is bounded by 
\[
\|\mu \|_\infty\Big(\fd_{\hY_{i,
  [l^i_{k,r}]_1}} (\fc_r, \bbR\times\pmb{\gamma }^k_i)+\Big\|\frac{iF_{A_K}}{4\pi }-\frac{iF_{B^{K_i}}}{4\pi }\Big\|_{L^1(\hY_{i,
  [l^i_{k,r}]_1})}\Big),
\]
and therefore 
\[\begin{split}
 & \|\td{\grc}_{k,r}^i-\td{\gamma }_i^k\|_{op}\leq  \fd_{\hY_{i,
   [l^i_{k,r}]_1}} (\fc_r, \bbR\times \pmb{\gamma }^k_i)
+\Big\|\frac{iF_{A^K}}{4\pi }-\frac{iF_{B^{K_i}}}{4\pi }\Big\|_{L^1(\hY_{i,
  [l^i_{k,r}]_1})}\\
& \qquad \leq  \fd_{\hY_{i,
  [-s^i_{k,r}]_1}} (\tau _{-L^i_{k+1, r}}\fc_r, \bfC_{i,k+1})+\fd_{\hY_{i,
  [-s^i_{k,r}]_1}} (\bfC_{i,k+1},  \bbR\times \pmb{\gamma }^k_i)\\
& \qquad \qquad +\zeta _i
e^{-\kappa _il^i_{k,r}}\\
& \qquad \to  0\quad \text{as \(r\to \infty\)},
\end{split}
\]
where \(s^i_{k,r}=\Delta ^i_{k,r}/2\) when \(k\in \{1, \ldots,
\grK_i-1\}\), and \(s^i_{0,r}:=\Delta ^i_{0,r}+3/2\). In the above, we
used (\ref{eq:xi-exp}), (\ref{appr:b-C}) and the facts that \(s^i_{k,r}\), \(s^i_{k,r}\)
both go to \(\infty\) as \(r\to \infty\), and that the
\(-\infty\)-limit of \(\bfC_{i,k+1}\) is \(\pmb{\gamma
}^k_i\). Given that \(\td{\grc}_{k,r}^i=\frac{i}{2\pi }F_{B^i_{k,r}}\)
converges to \(\td{\gamma }_i^k\), The arguments in the proof of  Theorem \ref{thm:strong-t} in
Section \ref{sec:strong-t} may
be combined with (\ref{appr:b-C}) to show that \([B^i_{k,r}]\in \Conn
(Y_i)/\scrG\subset \op{C}(Y_i)\) converges in the current topology. We
denote this limit by \([B^i_{k,\infty}]\). 
Thus, 
the arguments in Lemma \ref{rem:rel_class} shows that there
exists an \(R_2\geq R_1\) such that for all  \(r, r'\) satisfying
\(r>r'\geq R_2\), we have canonical isomorphisms of affine spaces
\[\begin{split}
& \pi _{Y_i}(\grc^i_{k-1, r'}, \grc^i_{k, r'})  \stackrel{\sim}{ \to}
\pi _{Y_i}(\grc^i_{k-1, r}, \grc^i_{k, r})\quad \forall i\in \grY'_m,
\, k\in
\{1, \ldots, \grK_i\};\\
& \pi _0(\scrB_X(\{\grc^i_{0,r'}\}_{i\in \grY})  
\stackrel{\sim}{\to } 
\pi _0(\scrB_X(\{\grc^i_{0,r}\}_{i\in \grY})
\end{split}
\]
via concatenations with the distinguished elements
\(o_{Y_i}(\grc^i_{k, r'}, \grc^i_{k, r})\)\\ \(\in \pi _{Y_i}(\grc^i_{k, r'},
\grc^i_{k, r})\) for \(0\leq k\leq \grK_i\), \(i\in \grY_m\). (When
\(i\in \grY_v\), \(\grc^i_{0,r}=\grc^i_{0,r'}:=[(B_i, \Phi _i)]\) is
the \(Y_i\)-end limit of \([(A_r, \Psi _r)]\), and
\(o_{Y_i}(\grc^i_{0, r'}, \grc^i_{0, r})\) is the identity in this
case.) 
Recall from Lemma \ref{lem:htpy} the definition of \(\grh'\); and that 
\(\grh=i_\infty\circ\grh'\). 
Each of the preceding maps fits into its respective version of the diagram
(\ref{CD:grh}):
{\small
\[
  \xymatrix@=18pt{
   \pi_0(\scrB_{X}(\{\grc^i_{0,r'}\}_{i\in\grY})\ar@{->}[r]^{\op{c}_{\{o_{Y_i}(\grc^i_{0,r'},\grc^i_{0,r})\}_i}}
   \ar@{->}[d]^{\grh'}
     &\pi_0(\scrB_{X}(\{\grc^i_{0,r}\}_{i\in \grY})
  \ar@{->}[d]^{\grh'}
   \ar@{->}[rd]^{\grh}\\
\scrH ^\bbR\big((X^{'a}, \nu),\{\tilde{\grc}^i_{0,r'}\}_{i\in \grY_m}\big)
\ar@{->}[r]^{i_{r',r}} & \scrH ^\bbR\big((X^{'a},
\nu),\{\tilde{\grc}^i_{0,r}\}_{i\in \grY_m}\big) \ar@{->}[r]^{i_\infty} & \scrH ^\bbR\big((X^{'a}, \nu),\{\tilde{\gamma }^0_i\}_i\big);
 }
\]}
and for \(i\in \grY'_m\),
\( k\in\{1, \ldots, \grK_i\}\), 
{\small\[
\xymatrix@=15pt{
 \pi _{Y_i}(\grc^i_{k-1, r'}, \grc^i_{k, r'})   \ar@{->}[r]^{*}
   \ar@{->}[d]^{\grh'}
     &\pi _{Y_i}(\grc^i_{k-1, r}, \grc^i_{k, r}) 
  \ar@{->}[d]^{\grh'}
   \ar@{->}[rd]^{\grh}\\
\scrH^\bbR(Y_i, \nu_i,\grs_i; \td{\grc}^i_{k-1, r'}, \td{\grc}^i_{k, r'}) 
\ar@{->}[r]^{i_{r',r}} &  \scrH^\bbR(Y_i, \nu_i,\grs_i; \td{\grc}^i_{k-1, r}, \td{\grc}^i_{k, r}) \ar@{->}[r]^{i_\infty} & \scrH^\bbR(Y_i, \nu_i,\grs_i; \tilde{\gamma }_i^{k-1}, \tilde{\gamma }_i^{k}),
 }
\]}
where the horizonal map in the top row labeled \(*\) is given by
\(c\mapsto  o_{Y_i}(\grc^i_{k-1,r},\grc^i_{k-1,r'})*c*
o_{Y_i}(\grc^i_{k,r'},\grc^i_{k,r})\). Recalling that in the cylindrical
case, \(\grh'\) factors as \(\ud{\grh}\circ\Pi _*\) and combing this
with the fact that \(\grh'\), \(\grh\), \(\ud{\grh}\) preserve
concatenations, we see that the horizonal maps \(i_{r',r}\),
\(i_\infty\) are given by {\small 
\begin{equation}\begin{split}
    &  i_{r',r}=\\
    & \begin{cases}
    \op{c}_{\{\ud{\grh}\, (\td{o}_{Y_i}([B^i_{0,r'}],[B^i_{0,r}]))\}_i}
    &\text{\footnotesize in 1st CD;}\\
    \text{the map \(h\mapsto \ud{\grh}\, (\td{o}_{Y_i}([B^i_{k-1,r}],[B^i_{k-1,r'}])*h*
\ud{\grh}\, (\td{o}_{Y_i}([B^i_{k,r'}],[B^i_{k,r}]))\)}
&\text{\footnotesize 
in 2nd CD;}
  \end{cases}
\end{split}
\end{equation}
\begin{equation}\label{eq:i_infty}\begin{split}
    & i_{\infty}=\\
    &\begin{cases}
    \op{c}_{\{\ud{\grh}\, (\td{o}_{Y_i}([B^i_{0,r}],[B^i_{0,\infty}]))\}_i}
    &\text{\footnotesize in 1st CD;}\\
    \text{the map \(h\mapsto \ud{\grh}\, (\td{o}_{Y_i}([B^i_{k-1,\infty}],[B^i_{k-1,r}])*h*
\ud{\grh}\, (\td{o}_{Y_i}([B^i_{k,r}],[B^i_{k,\infty}]))\)}
&\text{\footnotesize in 2nd CD.}
  \end{cases}
\end{split}
\end{equation}
}
(In the above, ``CD'' stands for
``commutative diagram''.)

Recall that the relative homology class of \(\grC\) is defined to be
the composition
\[
[\grC]= \op{c}_{\{[\grC_i]\}_{i}\in \grY_m} ([\bfC_0]),
\]
and 
\[
[\grC_i]:=\begin{cases}
[\tilde{C}_{i, \grK_i}]* \cdots *[\tilde{C}_{i,2}]*[\tilde{C}_{i,1}]
&\text{when \(i\in \grY'_m\);}\\
0\in H_2(Y_i;\bbZ) &\text{when \(\grK_i=0\).} 
\end{cases}
\]
Denote by \(\mathfrc{h}_r\) the relative homotopy class of \((A_r,
\Psi _r)\) and recall that by assumption, versions of
\(\mathfrc{h}_r\) corresponding to different, sufficiently large \(r\)
are identified via the canonical isomorphisms in
Lemma \ref{rem:rel_class} and collectively
denoted as  \(\mathfrc{h}\). To compare \(\grh (\fh)\) with \([\grC]\), first fix \(r>R_2\) and
identify \(\fh=\fh_r\) with the relative homotopy class of \((A_r, \Psi
_r)\) in \(\pi_0\scrB_X(\{\grc_{i,r}\}_{i\in \grY})\). 
We shall  express
\(\fh_r\) as a composition, 
\begin{equation}\label{def:fh_i}
  \fh_r=\op{c}_{\{\fh_{i,r}\}_i} (\fh_{0,r}), \quad
\end{equation}
where \(\fh_{0,r}\in \pi_0\scrB_{X}(\{\grc^i_{0,r}\}_{i\in \grY})\),
and
\begin{equation}\label{def:fh_i-}\begin{split}
    & \fh_{i,r}= \\
    & \begin{cases}
    \fh^r_{i,\grK_i}, * \cdots*\fh^r_{i,2}*\fh^r_{i,1},\quad \fh^r_{i,
  k}\in \pi _{Y_i}(\grc^i_{k-1, r}, \grc^i_{k, r}) &\text{when 
      \(i\in \grY'_m\)};\\
 0\in H_2(Y;\bbZ)\simeq \pi  _{Y_i}(\grc_{i,r}, \grc_{i,r})&
 \text{when \(i\in \grY-\grY'_m\)},\\
\end{cases}
\end{split}
\end{equation}
such that
\begin{equation}\label{eq:fh-C}
\text{
\(\grh(\fh_{0,r})=[\bfC_0]\) and \(\grh(\fh^r_{i,
  k})=[\bfC_{i,k}]\) for all \(i\in \grY'_m\) and \(k\in \{1, \ldots,
\grK_i\}\). }
\end{equation}
Once (\ref{def:fh_i}), (\ref{def:fh_i-}), (\ref{eq:fh-C}) are
established, modulo Lemma
\ref{lem:htpy}(b), the proof of Item (c) in the statement of the
theorem follows directly from the naturality of \(\grh\) under concatenation/composition maps. 

To define \(\fh_{0,r}\), \(\fh^r_{i,k}\), we begin with a modification of
\([(A_r, \Psi _r)]\in \scrB_X(\{\grc_{i,r}\}_{i\in \grY})\). 
Introduce the cutoff function \(\lambda (s):=\chi (2s)\chi (-2s)\) and let \((\hat{B}^i_{k,r},
\hat{\Phi }^i_{k,r})\) denotes the pull back of \((B^i_{k,r}, \Phi ^i_{k,r})\) to
\(\hY_i\). Set 
\[\begin{split}
\td{A}_r:& =A'_r+\sum_{i\in \grY'_m}\sum_{k=0}^{\grK_i-1} \big(\lambda
(s-l^i_{k,r})(\hat{B}^i_{k,r}-A'_r)\big),\\
\td{\Psi }_r:& =\Psi '_r+\sum_{i\in \grY'_m}\sum_{k=0}^{\grK_i-1} \big(\lambda
(s-l^i_{k,r})(\hat{\Phi }^i_{k,r}-\Psi '_r)\big).
\end{split}
\]
Then \([(\td{A}_r, \td{\Psi }_r)]\) represents the same relative homotopy class
as \([(A_r, \Psi _r)]\), namely \(\fh_r\). Now, we define \(\fh_{0,r}\) to
be the relative homotopy class of \([(A_{0,r}, \Psi _{0,r})]\), which
is in turned given as follows: 
Define \((A_{0,r}, \Psi _{0,r})\in \Conn (X)\times \Gamma (\bbS_X^+)\)
by setting \((A_{0,r}, \Psi _{0,r})=(\td{A}_r, \td{\Psi }_r)\) over
\(X^+_{\vec{L}_T(r)-\vec{3/2}}:=X^{'a}_{\vec{L}_T(r)-\vec{3/2}}\cup
(X-X^{'a})\), and then extending over \(X\) by the \(s\)-independent
configurations \((\hat{B}^i_{0,r},
\hat{\Phi }^i_{0,r})\), \(i\in \grY_m'\),  over the complement of 
\(X^+_{\vec{L}_T(r)-\vec{3/2}}\). 


Similarly, when \(i\in \grY'_m\), \(\fh^r_{i,k}\) is defined to be the relative homotopy
class of \([(A_{k,r}^i, \Psi _{k,r}^i)]\), where \((A_{k,r}^i, \Psi
_{k,r}^i)\in \in \Conn (\bbR\times Y_i)\times \Gamma (\bbS_{\bbR\times
  Y_i}^+)\) is given by setting 
\((A_{k,r}^i, \Psi _{k,r}^i)=\tau _{-L^i_{k,r}}(\td{A}_r, \td{\Psi
}_r)\) over  \(J^{i}_{k,r}\times Y_i \), and then extending over
\(\bbR\times Y_i\) by the \(s\)-independent
configurations over the complement, namely \(\hat{B}^i_{k-1,r}\) and \(\hat{B}^i_{k,r}\) respectively over the
 \(-\infty\)-end and the \(+\infty\)-end.  With \(\fh^r_{i,k}\),
 \(\fh_{0,r}\) so defined, the identities in 
(\ref{def:fh_i}), (\ref{def:fh_i-}) hold by construction. 

Let \(\{\mu _j\}_j\) be a set of closed 2-forms on \(X\) that
represents a basis of \(H^2(X;\bbR)\). We choose
\(\mu _j\) so that for every \(j\), 
\begin{itemize}
 \item \(\mu _j|_{\hY_i}=\pi _2^*\mu _{i,j}\) for a closed
2-forn \(\mu _{i,j}\) on \(Y_i\) for all \(i\in \grY\). 
\item When \(i\in \grY_v\), \(\mu _{i,j}=0\). This is possible because
  of the condition (\ref{b_1=0}). 
\item When \(i\in \grY_m\), \(\mu _{i,j}\) is supported on the
  complement of 
  \(\bigcup_{\pmb{\gamma }\in \bbP(\grs_i)} \gamma \subset
    Y_i\). This is possible because \(\bbP(\grs_i)\) consists of
    finitely many elements, and the closure of each \(\pmb{\gamma }\) is a compact
    1-manifold with (possibly empty) boundary. 
\end{itemize}
To verify the first identity in (\ref{eq:fh-C}), it suffices to verify
that 
for every \(j\), the pairing 
\begin{equation}\label{pairing=0}
\big\langle[\mu _j], \grh'
(\fh_{0,r})-(i_\infty)^{-1}[\bfC_0]\big\rangle =0, \quad r> R_2.
\end{equation}
Given (\ref{eq:i_infty}), the left hand side of the preceding expression equals {\small
\[
\begin{split}
&
\int_{-1/2}^{1/2}\Big(\int_{X^{'a}_{\vec{L}_T(r)-\vec{3/2}+\vec{t}}}\frac{iF_{A_{0,r}^E}}{2\pi}\wedge\mu_j-\int_{\td{C}_0|_{X^{'a}_{\vec{L}_T(r)-\vec{3/2}+\vec{t}}}}\mu_j\\
& \,\,  -\sum_{i\in
    \grY'_m}\Big(\int_{\td{C}_0|_{\hY_{i,L_T(r)-3/2+t}}}\pi _2^*\mu_{i,j}-\lim_{r'\to
    \infty}\int_{Y_i}\big(\frac{i}{2\pi}(B_{0,r'}^E-B_{0,r}^E)\big)\wedge\mu_{i,j}\Big)\Big)
  dt\\
  &
=\int_{-1/2}^{1/2}\Big(\int_{X^{'a}_{\vec{L}_T(r)-\vec{3/2}+\vec{t}}}\frac{iF_{A_{r}^E}}{2\pi}\wedge\mu_j-\int_{\td{C}_0|_{X^{'a}_{\vec{L}_T(r)-\vec{3/2}+\vec{t}}}}\mu_j\\
& \, \,   -\sum_{i\in
    \grY'_m}\Big(\int_{\td{C}_0|_{\hY_{i,L_T(r)-3/2+t}}}\pi _2^*\mu_{i,j}-\lim_{r'\to
    \infty}\int_{Y_i}\big(\frac{i}{2\pi}(B_{0,r'}^E-B_{0,r}^E)\big)\wedge\mu_{i,j}\Big)\Big)
dt,
\end{split}
\]
}
and is independent of \(r\). However, according to (\ref{def:L_T}),
the fact that \(\bfC_0\) is asymptotic to a t-orbit, and the
convergence of \([B_{0,r}]\) in current topology, the preceding
expression converges to 0 as \(r\to\infty\). We have thus confirmed
the first identity in (\ref{eq:fh-C}). 


The second identity in (\ref{eq:fh-C}) for each given \(i, k\) is verified similarly, with the
roles of \(A_{0,r}\), \(\bfC_0\), \(X^{'a}_{\vec{L}_T(r)-\vec{3/2}}\)
above respectively replaced by \(A_{k,r}^i\), \(\bfC_{i,k}\), and
\(J^{i}_{k,r}\times Y_i \), and with the role of  (\ref{def:L_T}) in
the preceding argument now played by 
(\ref{appr:b-C}). \epf

Finally, we prove 
Lemma \ref{lem:htpy}(b).

\subsubsection*{\it Proof of Lemma \ref{lem:htpy}(b).} Recall that the
currents \(\{\tilde{\gamma }_i\}_{i\in \grY_m}\) are integral, and
\(\scrH\big((X^{'a}, \nu),\{\tilde{\gamma }_i\}_i\big)\) embeds in \(\scrH
^\bbR\big((X^{'a}, \nu),\{\tilde{\gamma }_i\}_i\big)\) as an orbit of
the \(\scrH_X\)-action. The map 
 \(\grh=i_\infty\circ \grh'\) maps
\(\pi_0\scrB_X(\{\grc_{i,r}\}_{i\in \grY})\simeq \pi_0\scrB_{X^{'a}}(\{\grc_{i,r}\}_{i\in \grY_m})\) to another orbit of the
\(\scrH_X\)-action in \(\scrH ^\bbR\big((X^{'a}, \nu),\{\tilde{\gamma
}_i\}_i\big)\). We claim that these two orbits are identical. To see
this, note that by Theorem \ref{thm:g-conv} (a), there is a
subsequence, also denoted \(\{(A_r, \Psi
_r)\}_r\), which weakly converges to a chain of
t-curves, \(\grC\).  The associated current \(\td{\grC}\) is integral,
and has relative homology class in \(\scrH\big((X^{'a},
\nu),\{\tilde{\gamma }_i\}_i\big)\subset\scrH^\bbR\big((X^{'a},
\nu),\{\tilde{\gamma }_i\}_i\big)\).  Meanwhile, \([\grC]=\grh (\fh)\) is in
the image of \(\grh \). This implies that the \(\scrH_X\)-orbits
\(\scrH\big((X^{'a}, \nu),\{\tilde{\gamma }_i\}_i\big)\) and \(\grh \,
(\pi_0\scrB_X(\{\grc_{i,r}\}_{i\in \grY}))\) in \(\scrH
^\bbR\big((X^{'a}, \nu),\{\tilde{\gamma }_i\}_i\big)\) are identical. 
\epf

\subsubsection*{\it Proof of  Proposition \ref{cor:F-positive}, the
  non-cylindrical case.} Adopt the
notations from the proof of Theorem \ref{thm:g-conv} above. 
By the assumption (\ref{def:P}), \(C_0\cap P\), \(C_{i,k}\cap(\bbR\times
\pp_i)\) \(\forall i, k\) all consist of finitely many points. Thus, there exists an \(R_3\geq R_2\) such
that for all \(r\geq R_3\):
\begin{itemize}
\item \(C_0\cap P\) lies in the interior of
\(X^{'a}_{\vec{L}_T(r)-\vec{3}}\);
\item  \(C_{i,k}\cap(\bbR\times
\pp_i)\) lies in the interior of \(J^{i}_{k,r}\times Y_i\) for all
\(i\), \(k\). 
\item For all \(i\in \grY_m\), \(\fd_{\hY_{i, I_*}}(P, \bbR\times \pp_i)< \epsilon /4\) for all
  \(I_*\subset [L_T(r)-5/2, \infty)\) with length \(|I_*|=2\). 
\end{itemize}
Let \(\td{P}\subset X^{'a}\) be a subvariety such that:
\begin{itemize}
\item  \(\td{P}\) agrees
with \(P\) over \(X^{'a}_{\vec{L}_T(r)-\vec{3}}\);
\item  it agrees with
\(\bbR\times \pp_i\) over \(\hY_{i, L_T(r)-5/2}\) \(\forall i\in
\grY_m\), and 
\item \(\fd_{\hY_{Y_{i, [L_T(r)-3, L_T(r)-5/2]}}}(P,
  \td{P}) +\fd_{\hY_{Y_{i, [L_T(r)-3, L_T(r)-5/2]}}}(\td{P}, \bbR\times
    \pp_i)  \leq \epsilon \). 
\item \(P-\td{P}\)  is a boundary. 
  \end{itemize}
  Write 
\[\begin{split}
\int_P\frac{i}{2\pi}F_{A_r^E} 
& =
\int_{\td{P} \cap X^{'a}_{\vec{L}_T(r)-\vec{\frac{5}{2}}}}\frac{i}{2\pi}F_{A_r^E}+
\sum_{i\in  \grY'_m}\sum_{k=1}^{\grK_i}\int_{\td{P}\cap\hat{Y}_{i,J^{'i}_{k,r}}}\frac{i}{2\pi}F_{A_r^E}\\
\end{split}
\]
Taking the \(r\to \infty\) limit, by (\ref{def:L_T}) and 
(\ref{appr:b-C}) we then have
\[
\lim_{r\to \infty}\int_P\frac{i}{2\pi}F_{A_r^E} =\# \, (P\cap
\bfC_0)+\sum_{i\in  \grY'_m}\sum_{k=1}^{\grK_i}\# \, \big((\bbR\times \pp_i)\cap \bfC_{i,k}\big).
\]
Now observe that each term on the right hand side is a non-negative
integer, since \(P\), \(\bbR\times \pp_i\), \(C_0\), \(C_{i,k}\) are
all pseudo-holomorphic. 
\epf

\subsection{Proving Theorem \ref{thm:g-conv} and Proposition
  \ref{cor:F-positive}, the cylindrical case}\label{sec:g-conv:b}

The proofs require only minor modifications of the proofs in the
non-cylindrical case.

\subsubsection*{\it Proof of  Theorem \ref{thm:g-conv}, the cylindrical
  case.} Begin as with the proof of the noncylindrical case in the previous
subsection: 
Rename the t-curve
\(\bfC\) in the statement of Theorem \ref{thm:l-conv} as 
\(\bfC_0\). That is, \(\{(A_r, \Psi _r)\}_{r\in \Gamma _0}\)
t-converges to \(\bfC_0\) over \(X=\bbR\times Y\).
As noted in Remark \ref{rem:cyl-end}, in this case we may regard \(X\) as an MCE with two ends: \(\grY=\grY_m=\{\pm\}\) consists two elements,
with the elements \(+/-\) labelling respectively  the
\(+\infty\)/\(-\infty\)-ends; their corresponding ending
3-manifolds are \(Y_\pm=\pm Y\). Rerunning the arguments in the
previous subsection, we get a subsequence \(\Gamma '\) and two chains of
t-curves \(\grC_\pm \) on \(\bbR\times (\pm Y)\), where \(\grC_\pm\) has
\(-\infty\)-limit \(\pm\pmb{\gamma }_{0,\pm}\) and \(+\infty \)-limit
\(\pm \pmb{\gamma }_\pm\). Let \(\iota\co \bbR\times (-Y)\to
\bbR\times Y\) denote the orientation preserving map \((s, x)\mapsto
(-s,x)\). Write {\footnotesize
\[\begin{split}
    & \grC=\\
    & \begin{cases}
\{\iota(\bfC_{-,\grK_-}), \iota(\bfC_{-,\grK_--1}), \ldots,\iota (\bfC_{-,1}), \bfC_0,
\bfC_{+,1}, \ldots, \bfC_{+, \grK_+}\} & \text{\footnotesize when \(\grK_->0\);
  \(\grK_+>0\)};\\
\{\iota(\bfC_{-,\grK_-}), \iota(\bfC_{-,\grK_--1}), \ldots,
\iota(\bfC_{-,1}), \bfC_0\} & \text{\footnotesize when \(\grK_->0\);
  \(\grK_+=0\)}; \\
\{\bfC_0,\bfC_{+,1}, \ldots, \bfC_{+, \grK_+}\} & \text{\footnotesize when \(\grK_-=0\);
  \(\grK_+>0\)}; \\
\{\bfC_0\} & \text{\footnotesize when \(\grK_-=0\); \(\grK_+=0\)}
  \end{cases}
\end{split}
\]}
if  \(\bfC_0\) is non-constant; and when  \(\bfC_0\) is constant, let {\footnotesize
\[\begin{split}
    &  \grC=\\
    & \begin{cases}
\{\iota(\bfC_{-,\grK_-}), \iota(\bfC_{-,\grK_--1}), \ldots, \iota(\bfC_{-,1})
\bfC_{+,1}, \ldots, \bfC_{+, \grK_+}\} & \text{when \(\grK_->0\);
  \(\grK_+>0\)};\\
\{\iota(\bfC_{-,\grK_-}), \iota(\bfC_{-,\grK_--1}), \ldots, \iota(\bfC_{-,1})\} & \text{when \(\grK_->0\);
  \(\grK_+=0\)}; \\
\{\bfC_{+,1}, \ldots, \bfC_{+, \grK_+}\} & \text{when \(\grK_-=0\);
  \(\grK_+>0\)}. \\
  \end{cases}
\end{split}
\]}
Lastly, when \(\bfC_0\) is constant and \(\grK_\pm=0\), let \(\grC\)
be the 0-component chain of t-curves with a single
rest point \(\pmb{\gamma }_+\). The arguments from the previous
subsections then shows that \(\{(A_r, \Psi
_r)\}_{r\in \Gamma '}\) weakly t-converges to 
 \(\grC\), which have the properties described in Items (b)-(d) in the
 statement of the theorem. (Note that when \((Y, \nu )\) is
 cylindrical, the conditions in Items (c) are (d) are both met.)
The case when \(\bfC_0\) is constant and \(\grK_\pm=0\) occurs only when \(\{(A_r, \Psi _r)\}_{r\in \Gamma _0}\) is such that
there exists \(r_0>1\), such that \(\forall r\geq r_0\), \((A_r, \Psi
_r)\) is constant. This means that \((A_r, \Psi
_r)=(\hat{B}_r, \hat{\Phi }_r)\), where \((B_r, \Phi _r)\) strongly
t-converges to \(\pmb{\gamma }_-=\pmb{\gamma }_+\) by assumption. In
this case, \(\grC\) is the 0-component chain of t-curves with a single
rest point \(\pmb{\gamma }_-=\pmb{\gamma }_+\). \epf

\subsubsection*{\it Proof of  Proposition
  \ref{cor:F-positive}, the cylindrical case.} Let \(\bfC_0\),
\(\grC\) be as in the proof of the cylindrical case of Theorem
\ref{thm:g-conv} above. 
Then simply repeat
the arguments in the proof for the noncylindrical case of Proposition
  \ref{cor:F-positive} in the previous subsection. 
  \epf

\end{document}